\pdfoutput=1
\RequirePackage{ifpdf}
\ifpdf 
\documentclass[pdftex]{sigma}
\else
\documentclass{sigma}
\fi

\numberwithin{equation}{section}

\newtheorem{Theorem}{Theorem}[section]
\newtheorem{Corollary}[Theorem]{Corollary}
\newtheorem{Lemma}[Theorem]{Lemma}
\newtheorem{Proposition}[Theorem]{Proposition}
\newtheorem{Condition}[Theorem]{Condition}
{ \theoremstyle{definition}
\newtheorem{Definition}[Theorem]{Definition}

\newtheorem{Remark}[Theorem]{Remark}}

\newcommand{\nbigc}{\mathcal{C}}
\newcommand{\nbigd}{\mathcal{D}}
\newcommand{\nbige}{\mathcal{E}}
\newcommand{\nbigf}{\mathcal{F}}
\newcommand{\nbigg}{\mathcal{G}}
\newcommand{\nbigh}{\mathcal{H}}
\newcommand{\nbigi}{\mathcal{I}}

\newcommand{\nbigk}{\mathcal{K}}

\newcommand{\nbigm}{\mathcal{M}}

\newcommand{\nbigo}{\mathcal{O}}
\newcommand{\nbigp}{\mathcal{P}}

\newcommand{\nbigr}{\mathcal{R}}
\newcommand{\nbigs}{\mathcal{S}}
\newcommand{\nbigt}{\mathcal{T}}
\newcommand{\nbigu}{\mathcal{U}}
\newcommand{\nbigv}{\mathcal{V}}
\newcommand{\nbigw}{\mathcal{W}}
\newcommand{\nbigx}{\mathcal{X}}
\newcommand{\nbigy}{\mathcal{Y}}

\newcommand{\proj}{\mathbb{P}}
\newcommand{\seisuu}{{\mathbb Z}}

\newcommand{\cnum}{{\mathbb C}}
\newcommand{\real}{{\mathbb R}}
\newcommand{\hyperh}{\mathbb{H}}

\newcommand{\DD}{\mathbb{D}}
\newcommand{\EE}{\mathbb{E}}

\newcommand{\GG}{\mathbb{G}}

\newcommand{\gbigh}{\mathfrak H}

\newcommand{\gbigu}{\mathfrak U}
\newcommand{\gbigv}{\mathfrak V}

\newcommand{\gbigx}{\mathfrak X}

\newcommand{\gbigz}{\mathfrak Z}

\newcommand{\gminia}{\mathfrak a}
\newcommand{\gminib}{\mathfrak b}

\newcommand{\gminis}{\mathfrak s}

\newcommand{\gminiv}{\mathfrak v}

\newcommand{\vecv}{{\boldsymbol v}}

\newcommand{\veczero}{{\boldsymbol 0}}
\newcommand{\vecalpha}{{\boldsymbol \alpha}}
\newcommand{\veca}{{\boldsymbol a}}
\newcommand{\vecb}{{\boldsymbol b}}

\newcommand{\vecdelta}{{\boldsymbol \delta}}
\newcommand{\vecs}{{\boldsymbol s}}

\newcommand{\vecc}{{\boldsymbol c}}
\newcommand{\vecd}{{\boldsymbol d}}

\newcommand{\vecm}{{\boldsymbol m}}

\newcommand{\vecepsilon}{{\boldsymbol \epsilon}}

\newcommand{\vecn}{{\boldsymbol n}}

\newcommand{\lrarr}{\longrightarrow}

\def\Hom{\mathop{\rm Hom}\nolimits}

\def\End{\mathop{\rm End}\nolimits}

\def\Cok{\mathop{\rm Cok}\nolimits}

\def\Image{\mathop{\rm Im}\nolimits}

\def\Gr{\mathop{\rm Gr}\nolimits}

\def\rank{\mathop{\rm rank}\nolimits}

\def\Ker{\mathop{\rm Ker}\nolimits}

\def\Gr{\mathop{\rm Gr}\nolimits}
\def\Sym{\mathop{\rm Sym}\nolimits}

\def\Res{\mathop{\rm Res}\nolimits}

\def\ord{\mathop{\rm ord}\nolimits}
\def\degpar{\mathop{\rm par\textrm{-}deg}\nolimits}

\def\ch{\mathop{ch}\nolimits}

\def\Tr{\mathop{\rm Tr}\nolimits}

\def\dvol{\mathop{\rm dvol}\nolimits}

\def\id{\mathop{\rm id}\nolimits}

\def\gap{\mathop{\rm gap}\nolimits}

\def\ch{\mathop{\rm ch}\nolimits}

\def\Irr{\mathop{\rm Irr}\nolimits}

\newcommand{\del}{\partial}
\newcommand{\delbar}{\overline{\del}}

\newcommand{\pardeg}{\degpar}

\newcommand{\nhom}{{\mathcal H{\rm om}}}

\newcommand{\nend}{{\mathcal End}}

\newcommand{\sankaku}{\triangle}

\newcommand{\barz}{\overline{z}}
\newcommand{\zbar}{\barz}
\newcommand{\zetabar}{\overline{\zeta}}

\newcommand{\barlambda}{\overline{\lambda}}
\newcommand{\lambdabar}{\barlambda}

\newcommand{\xbar}{\overline{x}}

\newcommand{\DDlambda}{\DD^{\lambda}}

\newcommand{\Par}{{\mathcal P{\rm ar}}}

\newcommand{\lefttop}[1]{{}^{#1}\!}

\def\reg{\mathop{\rm reg}\nolimits}

\newcommand{\openclosed}[2]{]#1,#2]}
\newcommand{\openopen}[2]{]#1,#2[}

\newcommand{\nablatilde}{\widetilde{\nabla}}

\newcommand{\wbar}{\overline{w}}

\newcommand{\DDlambdatilde}{\widetilde{\DD}^{\lambda}}
\newcommand{\htilde}{\widetilde{h}}

\newcommand{\Phitilde}{\widetilde{\Phi}}

\newcommand{\atilde}{\widetilde{a}}
\newcommand{\stilde}{\widetilde{s}}

\newcommand{\DDtilde}{\widetilde{\DD}}

\newcommand{\nbigvtilde}{\widetilde{\nbigv}}
\newcommand{\nbigvhat}{\widehat{\nbigv}}

\newcommand{\Utilde}{\widetilde{U}}
\newcommand{\Dtilde}{\widetilde{D}}
\newcommand{\Xtilde}{\widetilde{X}}

\newcommand{\nbigrtilde}{\widetilde{\nbigr}}

\def\ord{\mathop{\rm ord}\nolimits}

\def\Gal{\mathop{\rm Gal}\nolimits}

\def\HE{\mathop{\rm HE}\nolimits}

\def\inn{\mathop{\rm in}\nolimits}
\def\TW{\mathop{\rm TW}\nolimits}

\newcommand{\Ptilde}{\widetilde{P}}

\newcommand{\Phat}{\widehat{P}}

\newcommand{\varphitilde}{\widetilde{\varphi}}
\newcommand{\Ctilde}{\widetilde{C}}

\newcommand{\Stilde}{\widetilde{S}}

\newcommand{\Ytilde}{\widetilde{Y}}

\newcommand{\Lambdabar}{\overline{\Lambda}}

\newcommand{\Htilde}{\widetilde{H}}

\newcommand{\gminiatilde}{\widetilde{\gminia}}

\newcommand{\ttP}{\texttt{P}}

\newcommand{\ttPtilde}{\widetilde{\texttt{P}}}

\newcommand{\btilde}{\widetilde{b}}
\newcommand{\Qhat}{\widehat{Q}}

\newcommand{\Partilde}{\widetilde{\Par}}
\newcommand{\gaptilde}{\widetilde{\gap}}

\begin{document}

\newcommand{\arXivNumber}{1902.08298}

\renewcommand{\thefootnote}{}

\renewcommand{\PaperNumber}{068}

\FirstPageHeading

\ShortArticleName{Good Wild Harmonic Bundles and Good Filtered Higgs Bundles}

\ArticleName{Good Wild Harmonic Bundles\\ and Good Filtered Higgs Bundles\footnote{This paper is a~contribution to the Special Issue on Primitive Forms and Related Topics in honor of~Kyoji Saito for his 77th birthday. The full collection is available at \href{https://www.emis.de/journals/SIGMA/Saito.html}{https://www.emis.de/journals/SIGMA/Saito.html}}}

\Author{Takuro MOCHIZUKI}

\AuthorNameForHeading{T.~Mochizuki}

\Address{Research Institute for Mathematical Sciences, Kyoto University, Kyoto 606-8502, Japan}
\Email{\href{mailto:takuro@kurims.kyoto-u.ac.jp}{takuro@kurims.kyoto-u.ac.jp}}

\ArticleDates{Received July 16, 2020, in final form June 28, 2021; Published online July 17, 2021}

\Abstract{We prove the Kobayashi--Hitchin correspondence between good wild harmonic bundles and polystable good filtered $\lambda$-flat bundles satisfying a~vanishing condition. We~also study the correspondence for good wild harmonic bundles with the homogeneity with res\-pect to a group action, which is expected to provide another way to construct Frobenius manifolds.}


\Classification{53C07; 58E15; 14D21; 81T13}

\begin{flushright}
\begin{minipage}{60mm}
\it Dedicated to Professor Kyoji~Saito \\ on the occasion of his~77th~birthday
\end{minipage}
\end{flushright}

\renewcommand{\thefootnote}{\arabic{footnote}}
\setcounter{footnote}{0}

\section{Introduction}

Let $X$ be a smooth projective variety
with a simple normal crossing hypersurface $H$.
Let $L$ be an ample line bundle on $X$. We~shall prove the following theorem,
that is the Kobayashi--Hitchin correspondence
for good wild harmonic bundles
and good filtered $\lambda$-flat bundles.

\begin{Theorem}[Corollary
\ref{cor;19.2.18.3}]
\label{thm;19.2.2.20}
The following objects are equivalent:
\begin{itemize}\itemsep=0pt
\item
 Good wild harmonic bundles on $(X,H)$.
\item
 $\mu_L$-Polystable filtered $\lambda$-flat bundles
 $\big(\nbigp_{\ast}\nbigv,\DDlambda\big)$
 on $(X,H)$
 satisfying
\begin{gather*}
\int c_1(\nbigp_{\ast}\nbigv)c_1(L)^{\dim X-1}=0,\qquad
\int\ch_2(\nbigp_{\ast}\nbigv)c_1(L)^{\dim X-2}=0.
\end{gather*}
\end{itemize}
\end{Theorem}
We shall recall the precise definitions
of the objects in Section~\ref{section;19.2.18.1}.

In~\cite{Mochizuki-wild},
we have already proved that
good wild harmonic bundles on $(X,H)$
induce $\mu_L$-polystable good filtered $\lambda$-flat bundles
satisfying the vanishing condition.
Note that $0$-flat bundles are equivalent to Higgs bundles,
and $1$-flat bundles are flat bundles in the ordinary sense.
Moreover,
we studied an analogue of Theorem~\ref{thm;19.2.2.20}
in the case $\lambda=1$,
i.e., the correspondence between
good wild harmonic bundles
and $\mu_L$-polystable good filtered flat bundles
satisfying a~similar vanishing condition~\cite[Theorem~16.1.1]{Mochizuki-wild}.
It~was applied to the study of the correspondence
between semisimple algebraic holonomic $D$-modules
and pure twistor $D$-modules.

In this paper, as a complement,
we shall explain the proof for all $\lambda$.
There is no new essential difficulty
to prove Theorem~\ref{thm;19.2.2.20}
after our studies~\cite{mochi4, mochi2-I,mochi2-II, mochi5, Mochizuki-wild}
on the basis of~\cite{s1, s2}.
Moreover, in some parts of the proof,
the arguments can be simplified in the Higgs case.
However,
because the Higgs case is also particularly important,
it would be useful to explain a rather detailed proof. We~shall also explain the correspondences in homogeneous cases
which would be useful in a generalized Hodge theory.

\subsection{Kobayashi--Hitchin correspondences}

\subsubsection{Kobayashi--Hitchin correspondence for vector bundles}

We briefly recall a part of the history of this type of correspondences.
(See also~\cite{Hitchin-note-vanishing, koba, Lubke-Teleman-universal}.)
For a~holomorphic vector bundle $E$ on a compact Riemann surface $C$,
we set $\mu(E):=\deg(E)/\rank(E)$,
which is called the slope of $E$.
A holomorphic bundle $E$ is called stable (resp.~semistable)
if~$\mu(E')<\mu(E)$ (resp.~$\mu(E')\leq\mu(E)$)
holds for any holomorphic subbundle
$E'\subset E$ such that $0<\rank(E')<\rank(E)$.
It~is called polystable if
it is a direct sum of~stable subbundles with the same slope.
This stability, semistability and polystability conditions
were introduced by~Mum\-ford~\cite{Mumford-ICM1962}
for the construction of the moduli spaces
of vector bundles with reasonable properties.
Narasimhan and Seshadri~\cite{narasimhan-seshadri}
established the equivalence between unitary flat bundles
and polystable bundles of degree $0$ on compact Riemann surfaces.

Let $(X,\omega)$ be a compact connected K\"ahler manifold.
For any torsion-free $\nbigo_X$-module $\nbigf$,
the~slope of $\nbigf$ with respect to $\omega$
is defined as
\begin{gather*}
 \mu_{\omega}(\nbigf):=
\frac{\int_Xc_1(\nbigf)\omega^{\dim X-1}}{\rank \nbigf}.
\end{gather*}
If the cohomology class of $\omega$ is
the first Chern class of an ample line bundle $L$,
then $\mu_{\omega}(\nbigf)$ is also denoted by $\mu_L(\nbigf)$.
Then, a torsion-free $\nbigo_X$-module $\nbigf$ is called
$\mu_{\omega}$-stable
if
$\mu_{\omega}(\nbigf')<\mu_{\omega}(\nbigf)$ holds
for any saturated coherent subsheaf $\nbigf'\subset\nbigf$
such that $0<\rank(\nbigf')<\rank(\nbigf)$.
This condition was first studied by
Takemoto~\cite{Takemoto1, Takemoto2}.
It~is also called
$\mu_{\omega}$-stability, or slope stability.
Slope semistability and slope polystability
are naturally defined.

Bogomolov~\cite{Bogomolov}
introduced the $T$-stability condition for torsion-free sheaves
on connected projective surfaces,
and he proved the inequality of the Chern classes
$c_2(E)-(r-1)c_1(E)^2/2r\geq 0$
for any $T$-semistable bundle $E$ of rank $r$. We~do not recall the precise definition of
$T$-stability condition here,
but we note that if a holomorphic vector bundle
on a complex projective manifold is slope semistable,
then it is $T$-semistable.
(See~\cite[Section~7]{Bogomolov} for more details.)
Gieseker~\cite{Gieseker} gave a~different proof
of the inequality for slope semistable bundles.
The inequality is called Bogomolov--Gieseker inequality
or Bogomolov inequality.

Inspired by these works,
Kobayashi~\cite{Kobayashi-Nagoya} introduced
the concept of Hermitian--Einstein condition
for metrics of holomorphic vector bundles.
Let $\big(E,\delbar_E\big)$ be a holomorphic vector bundle on
a K\"ahler manifold $(X,\omega)$.
Let $h$ be a Hermitian metric of $E$.
Let $R(h)$ denote the curvature of the Chern connection
$\nabla_h=\delbar_E+\del_{E,h}$,
associated with $h$ and $\delbar_E$.
Then, $h$ is called Hermitian--Einstein
if $\Lambda R(h)^{\bot}=0$,
where $R(h)^{\bot}$ denotes the trace-free part of $R(h)$. In~par\-ti\-cu\-lar,
he proved in~\cite{Kobayashi-Nagoya} that
if a holomorphic vector bundle on a compact K\"ahler manifold
has a Hermitian--Einstein metric,
then it is $T$-semistable.
Kobayashi~\cite{Kobayashi-Academy, Kobayashi-seminar}
and L\"ubke~\cite{Lubke2}
proved that
a holomorphic vector bundle
on a compact connected K\"ahler manifold
satisfies the slope polystability condition
if it has a Hermitian--Einstein metric.
Moreover,
L\"ubke~\cite{Lubke1} established
the so called Kobayashi--L\"ubke inequality
for the first and the second Chern forms
associated with Hermitian--Einstein metrics,
which is reduced to the inequality
$\Tr\bigl(\big(R(h)^{\bot}\big)^2\bigr)\omega^{\dim X-2}\geq 0$
in the form level. In~particular, it implies
the Bogomolov--Gieseker inequality
for holomorphic vector bundles $\big(E,\delbar_E\big)$
with a Hermitian--Einstein metric $h$ on compact K\"ahler manifolds
$(X,\omega)$.
Moreover,
if $c_1(E)=0$ and $\int_X\ch_2(E)\omega^{\dim X-2}=0$
are satisfied
for such $\big(E,\delbar_E,h\big)$,
and if we impose $\det(h)$ is flat,
then the Kobayashi--L\"ubke inequality implies
that $R(h)=0$,
i.e., $\nabla_h$ is~flat.

Independently, in~\cite{Kotake-Ochiai},
Hitchin proposed a problem
to ask an equivalence of the stability condition
and the existence of a metric $h$ such that $\Lambda R(h)=0$,
under the vanishing of the first Chern class of the bundle.
(See~\cite{Hitchin-note-vanishing} for more precise explanation.)
It~clearly contains the most important essence.
He also suggested possible applications
of the vanishings.
His problem stimulated Donaldson
whose work on this topic
brought several breakthroughs to
whole geometry.

In~\cite{Donaldson-Narasimhan-Seshadri},
Donaldson introduced the method of global analysis
to reprove the theorem of~Nara\-si\-mhan--Seshadri.
In~\cite{don3},
by using the method of the heat flow associated with
the Hermitian--Einstein condition,
he established the equivalence of
the slope polystability condition
and the existence of a Hermitian--Einstein metric
for holomorphic vector bundles on any complex projective surface.
The important concept of Donaldson functional was
also introduced in~\cite{don3}.

Eventually, Donaldson~\cite{don} and
Uhlenbeck--Yau~\cite{uy}
established the equivalence
on any dimensional complex projective manifolds.
Note that Uhlenbeck--Yau proved it for any compact K\"ahler manifolds,
more generally.
The correspondence is called
with various names;
Kobayashi--Hitchin correspondence,
Hitchin--Kobayashi correspondence,
Donaldson--Hitchin--Uhlenbeck--Yau correspondence, etc. In~this paper, we call it the Kobayashi--Hitchin correspondence.

As a consequence of the Kobayashi--Hitchin correspondence
and the Kobayashi--L\"ubke ine\-qu\-a\-lity,
we also obtain an equivalence between
unitary flat bundles
and slope polystable holomorphic vector bundles $E$
satisfying $\mu_{\omega}(E)=0$ and $\int_X\ch_2(E)\omega^{\dim X-2}=0$.
Note that Mehta and Ramanathan~\cite{mehta-ramanathan1, mehta-ramanathan2}
deduced the equivalence on complex projective manifolds
directly from the equivalence in the surface case due to Donaldson
\cite{don3}.

\subsubsection[Higgs bundles and lambda-flat bundles]
{Higgs bundles and $\boldsymbol\lambda$-flat bundles}

Such correspondences have been also studied for
vector bundles equipped with some additional structure,
which are also called Kobayashi--Hitchin correspondences in this paper.
One of the most rich and influential is the case of Higgs bundles,
pioneered by Hitchin and Simpson.

Let $\big(E,\delbar_E\big)$ be a holomorphic vector bundle
on a compact Riemann surface $C$.
A Higgs field of $\big(E,\delbar_E\big)$ is a holomorphic section
$\theta$ of $\End(E)\otimes\Omega^1_C$.
Let $h$ be a Hermitian metric of $E$. We~obtain the Chern connection $\delbar_E+\del_{E,h}$
and its curvature $R(h)$.
Let $\theta^{\dagger}_h$
denote the adjoint of $\theta$.
In~\cite{hitchin},
Hitchin introduced the following equation,
called the Hitchin equation,
\begin{gather}
\label{eq;19.2.16.1}
 R(h)+\big[\theta,\theta_h^{\dagger}\big]=0.
\end{gather}
Such $\big(E,\delbar_E,\theta,h\big)$ is called a harmonic bundle. In~particular,
he studied the case $\rank E=2$.
Among many deep results in~\cite{hitchin},
he proved that
a Higgs bundle $\big(E,\delbar_E,\theta\big)$ has
a Hermitian metric $h$ satisfying (\ref{eq;19.2.16.1})
if and only if
it is polystable of degree $0$.
Here, a Higgs bundle $\big(E,\delbar_E,\theta\big)$ is called stable
(resp.~semistable)
if $\mu(E')<\mu(E)$
(resp.~$\mu(E')\leq\mu'(E)$) holds for any
holomorphic subbundle $E'\subset E$
such that $\theta(E')\subset E'\otimes\Omega^1_C$
and that $0<\rank(E')<\rank(E)$,
and a~Higgs bundle is called polystable
if it is a direct sum of stable Higgs subbundles
with the same slope. By~this equivalence and another equivalence due to Donaldson~\cite{don2}
between irreducible flat bundles
and twisted harmonic maps,
Hitchin obtained that the moduli space
of polystable Higgs bundles of degree $0$
and the moduli space of semisimple flat bundles
are isomorphic.
His~work revealed that the moduli spaces
of Higgs bundles and flat bundles
have extremely rich structures.

The higher dimensional case was studied by Simpson~\cite{s1}.
Note that Simpson started his study independently
motivated by a new way to construct variations of Hodge structure,
which we shall mention later in Section~\ref{subsection;19.2.18.10}.
For a holomorphic vector bundle $\big(E,\delbar_E\big)$
on a complex manifold $X$ with arbitrary dimension,
a Higgs field $\theta$ is defined to be
a holomorphic section of~$\End(E)\otimes\Omega^1_X$
satisfying the additional condition $\theta\wedge\theta=0$.
Suppose that $X$ has a K\"ahler form.
Let $h$ be a Hermitian metric of $E$.
Let $F(h)$ denote the curvature of the connection
$\nabla_h+\theta+\theta_h^{\dagger}$.
A Hermitian metric $h$ of a Higgs bundle $\big(E,\delbar_E,\theta\big)$
is called Hermitian--Einstein
if~$\Lambda F(h)^{\bot}=0$.
When $X$ is compact,
the slope stability, semistability and polystability
conditions for Higgs bundles
are naturally defined
in terms of the slopes of Higgs subsheaves.
Simpson established that
a Higgs bundle $\big(E,\delbar_E,\theta\big)$ on a compact K\"ahler manifold
$(X,\omega)$ has a Hermitian--Einstein metric
if and only if it is slope polystable.
Moreover, he generalized the Kobayashi--L\"ubke inequality
for the Chern forms to the context of Higgs bundles,
which is reduced to the inequality
$\Tr\bigl(\big(F(h)^{\bot}\big)^2\bigr)\omega^{\dim X-2}\geq 0$
in the form level
for any Hermitian--Einstein metric $h$ of $\big(E,\delbar_E,\theta\big)$.
Here, the condition $\theta\wedge\theta=0$ is essential. In~particular, it implies that
if $\big(E,\delbar_E,\theta\big)$ on a compact K\"ahler manifold $(X,\omega)$
satisfies $\mu_{\omega}(E)=0$ and $\int_X\ch_2(E)\omega^{\dim X-2}=0$,
then a Hermitian--Einstein metric $h$ of $\big(E,\delbar_E,\theta\big)$
is {\em a pluri-harmonic metric},
i.e., the connection $\nabla_h+\theta+\theta_h^{\dagger}$ is flat.
It~is equivalent to the following:
\begin{gather*}
 \del_{E,h}\theta=0,\qquad
 \delbar\theta_h^{\dagger}=0,\qquad
 R(h)+\big[\theta,\theta_h^{\dagger}\big]=0.
\end{gather*}
A Higgs bundle $\big(E,\delbar_E,\theta\big)$
with a pluri-harmonic metric is called a harmonic bundle.
This equi\-va\-lence
and another important equivalence due to Corlette~\cite{corlette}
induce an equivalence bet\-ween
semisimple flat bundles and polystable Higgs bundles
$\big(E,\delbar_E,\theta\big)$
satisfying $\mu_{\omega}(E)=0$ and~$\int_X\ch_2(E)\omega^{\dim X-2}=0$
on any connected compact K\"ahler manifold.

After the work of Corlette, Donaldson, Hitchin and Simpson,
it turned out that the moduli space $\nbigm(X)$
of flat bundles on a complex projective manifold $X$
has a hyper-K\"ahler metric. In~particular, it induces the twistor space of the moduli space
$\TW(\nbigm(X))$,
which is a complex analytic space with a fibration
$\TW(\nbigm(X))\lrarr \proj^1$,
such that the fiber over $1$ is the moduli space of
flat bundles,
and that the fiber over $0$ is the moduli space of
Higgs bundles with vanishing rational Chern classes.
The notion of $\lambda$-connections
was introduced and developed by
Deligne and Simpson~\cite{s4,s3}
for a more complex analytic construction
of the twistor space $\TW(\nbigm(X))$.
They obtain the family
of the moduli spaces $\nbigm^{\lambda}(X)$
of $\lambda$-flat bundles on $X$,
and the family of~the moduli spaces
$\nbigm^{\mu}\big(X^{\dagger}\big)$
of $\mu$-flat bundles on
the conjugate $X^{\dagger}$ of $X$.
They proved that the twistor space $\TW(\nbigm(X))$
can be obtained as the gluing of
the two families
$\coprod_{\lambda}\nbigm^{\lambda}(X)$
and~$\coprod_{\mu}\nbigm^{\mu}\big(X^{\dagger}\big)$
by the natural identification of
$\nbigm^{\lambda}(X)
=\nbigm^{\mu}\big(X^{\dagger}\big)$
for $\lambda\mu=1$.

These correspondences are not only really interesting,
but also provide a starting point of~the further investigations.
Simpson pursued the comparison of flat bundles,
Higgs bundles
and more generally $\lambda$-flat bundles
in deeper levels~\cite{s5, s3},
and developed the non-abelian Hodge theory~\cite{s4}. In~particular, he explained that
the Kobayashi--Hitchin correspondences
for $\lambda$-flat bundles can be
studied in a unified way~\cite{s5}.
For more recent study on the moduli spaces of~$\lambda$-connections,
see~\cite{Collier-Wentworth, Hu-Huang, Huang, Simpson-iterated-destability},
for example.

\subsubsection{Filtered case}

It is interesting to generalize such correspondences
for objects on complex quasi-projective mani\~folds. We~need to impose a kind of boundary condition,
that is parabolic structure.

Mehta and Seshadri~\cite{mehta-seshadri}
introduced the concept of parabolic structure of vector bundles
on~compact Riemann surfaces.
Let $C$ be a compact Riemann surface
with a finite subset $D\subset C$.
Let $E$ be a holomorphic vector bundle on~$C$.
A parabolic structure of $E$ at $D$
is a tuple of filtrations
$F_{\bullet}(E_{|P})$ $(P\in D)$
indexed by $\openclosed{{-}1}{0}$
satisfying
$F_{a}(E_{|P})=\bigcap_{b>a} F_{b}(E_{|P})$.
Set
$\Gr^F_a(E_{|P}):=F_{a}(E)/F_{<a}(E)$,
and
\begin{gather*}
\deg(E,F):=\deg(E)
-\sum_{P\in D}
 \sum_{-1<a\leq 0}
 a\dim\Gr^F_a(E_{|P}).
\end{gather*}
We set $\mu(E,F):=\deg(E,F)/\rank(E)$.
For any subbundle $E'\subset E$,
filtrations $F(E'_{|P})$ on $E'_{|P}$
are induced as
$F_a(E'_{|P}):=F_a(E_{|P})\cap E'_{|P}$.
Then, $(E,F)$ is called stable
if $\mu(E',F)<\mu(E,F)$
for any subbundle $E'\subset E$
with $0<\rank(E')<\rank(E)$.
Semistability and polystability conditions
are also defined naturally.
Then, Mehta and Seshadri proved an equivalence of
irreducible unitary flat bundles on $C\setminus D$
and stable parabolic vector bundles $(E,F)$ with $\mu(E,F)=0$
on $(C,D)$.

For some purposes,
it is more convenient to replace
parabolic bundles with filtered bundles
introduced by Simpson~\cite{s1,s2}.
Let $\nbigo_C(\ast D)$ denote the sheaf of
meromorphic functions on~$C$ which may have poles along $D$.
Let $\nbigv$ be a locally free $\nbigo_C(\ast D)$-module.
A filtered bundle $\nbigp_{\ast}\nbigv$ over $\nbigv$
is a tuple of lattices
$\nbigp_{\veca}\nbigv$ $(\veca=(a_P)_{P\in D}\in\real^D)$
such that
$(i)$ $\nbigp_{\veca}\nbigv(\ast D)=\nbigv$,
$(ii)$ the restriction of $\nbigp_{\veca}\nbigv$
to a neighbourhood of $P\in D$ depends only on $a_P$,
$(iii)$ $\nbigp_{\veca+\vecn}\nbigv= \nbigp_{\veca}\nbigv\big(\sum n_PP\big)$
for any $\veca\in\real^D$ and $\vecn\in\seisuu^D$,
$(iv)$ for any $\veca\in\real^D$, there exists
 $\vecepsilon\in\real_{>0}^D$ such that
 $\nbigp_{\veca}\nbigv=\nbigp_{\veca+\vecepsilon}\nbigv$.
Let $\veczero$ denote $(0,\ldots,0)\in\real^D$.
Then, $\nbigp_{\veczero}\nbigv$ is equipped with
the parabolic structure $F$
induced by the images
of $\nbigp_{\veca}\nbigv_{|P}\lrarr
 \nbigp_{\veczero}\nbigv_{|P}$ $(P\in D)$.
It~is easy to observe that
filtered bundles are equivalent to parabolic bundles. We~set
$\mu(\nbigp_{\ast}\nbigv):=
 \mu(\nbigp_{\veczero}\nbigv,F)$
for the filtered bundle $\nbigp_{\ast}\nbigv$.

Simpson~\cite{s1,s2} generalized the theorem of
Mehta-Seshadri to the correspondences of
tame harmonic bundles
and regular filtered $\lambda$-flat bundles
on compact Riemann surfaces.
A harmonic bundle $\big(E,\delbar_E,\theta,h\big)$ on $C\setminus D$
is called tame on $(C,D)$
if the closure of the spectral curve of $\theta$
in~$T^{\ast}C(\log D)$ is proper over $C$.
A regular filtered $\lambda$-flat bundle consists of
a filtered bundle $\nbigp_{\ast}\nbigv$
equipped with a flat $\lambda$-connection
$\DDlambda\colon \nbigv\lrarr\nbigv\otimes\Omega^1_C$
such that
$\DDlambda\cdot\nbigp_{\veca}\nbigv\subset
 \nbigp_{\veca}\nbigv\otimes\Omega^1_C(\log D)$
for~any $\veca\in\real^D$.
Stability, semistable and polystable conditions
are naturally defined in terms of~the slope.
Then, Simpson established the equivalence of
tame harmonic bundles on $(C,D)$
and~polystable regular filtered $\lambda$-flat bundles
$\big(\nbigp_{\ast}\nbigv,\DDlambda\big)$ satisfying
$\mu(\nbigp_{\ast}\nbigv)=0$.
Note that filtered bundles express
the growth order of the norms of holomorphic sections
with respect to the metrics. We~should mention that
the study of the asymptotic behaviour of
tame harmonic bundles is much harder than
that of the asymptotic behaviour of unitary flat bundles.
Hence,
it is already hard to prove that
tame harmonic bundles induce
regular filtered $\lambda$-flat bundles.

There are several directions to generalize.
One is a generalization in the context of
tame harmonic bundles on higher dimensional varieties.
Let $X$ be a smooth connected projective variety
with a simple normal crossing hypersurface $H$
and an ample line bundle $L$.
Then, there should be equivalences of
tame harmonic bundles on $(X,H)$
and $\mu_L$-polystable regular filtered $\lambda$-flat bundles
$(\nbigp_{\ast}\nbigv,\theta)$
on $(X,H)$ satisfying
$\int_X c_1(\nbigp_{\ast}\nbigv)c_1(L)^{\dim X-1}=0$
and $\int_X\ch_2(\nbigp_{\ast}\nbigv)c_1(L)^{\dim X-2}=0$
for each $\lambda\in\cnum$.
In~\cite{b}, Biquard studied the case where $H$ is smooth.
In~\cite{li2, Li-Narasimhan, steer-wren},
Li, Narasimhan, Steer and Wren studied
the correspondence for parabolic bundles
without flat $\lambda$-connections.
In~\cite{JZ2}, Jost and Zuo studied
the correspondence between semisimple flat bundles
and tame harmonic bundles.
In~\cite{mochi4,mochi2-I,mochi2-II, mochi5},
the author obtained the satisfactory equivalences
for tame harmonic bundles.
Note that Donagi and Pantev recently proposed
an attractive application of the Kobayashi--Hitchin correspondence
for tame harmonic bundles to the study of geometric Langlands theory
\cite{Donagi-Pantev-2009}.

In another natural direction of generalization,
we should consider more singular objects
than regular filtered Higgs or flat bundles.
A harmonic bundle $\big(E,\delbar_E,\theta,h\big)$
on $X\setminus H$ is called
wild if the closure of the spectral variety of
$\theta$
in the projective completion of $T^{\ast}X$ is complex analytic.
For the analysis,
we should impose that
the spectral variety of the harmonic bundle
satisfies some non-degeneracy condition along $H$.
(See Section~\ref{subsection;19.2.21.1}.)
This is not essential because the condition is always satisfied
once we replace $X$ by its appropriate blow up.
The notion of regular filtered $\lambda$-flat bundle
is appropriately generalized to the notion of
good filtered $\lambda$-flat bundle.
The results of Simpson should be generalized to
equivalences of good wild harmonic bundles
and $\mu_L$-polystable good filtered $\lambda$-flat bundles
$\big(\nbigp_{\ast}\nbigv,\DDlambda\big)$
satisfying
$\int_Xc_1(\nbigp_{\ast}\nbigv)c_1(L)^{\dim X-1}=0$
and $\int_X\ch_2(\nbigp_{\ast}\nbigv)c_2(L)^{\dim X-2}=0$.
Sabbah~\cite{sabbah3} studied
the correspondence between
semisimple meromorphic flat bundles
and wild harmonic bundles
in the one dimensional case.
Biquard and Boalch~\cite{biquard-boalch}
obtained generalization for wild harmonic bundles
in the one dimensional case.
Boalch informed the author that wild generalization
in the context of the Higgs case was not expected in those days.

As mentioned, the author studied
the wild harmonic bundles
on any dimensional varieties
in~\cite{Mochizuki-wild}. We~obtained that good wild harmonic bundles
induce $\mu_L$-polystable good filtered $\lambda$-flat bundles
satisfying the vanishing conditions.
Moreover,
we proved that the construction induces an equivalence
of good wild harmonic bundles
and slope polystable good filtered flat bundles
satisfying the vanishing condition.
Such an equivalence for meromorphic flat bundles
is particularly interesting
because we may apply it to prove a conjecture of Kashiwara~\cite{k5}
on~semisimple algebraic holonomic $\nbigd$-modules.
See a survey paper
\cite{Mochizuki-ICM} for more details on this application.

In~\cite{Mochizuki-wild},
we did not give a proof of the equivalence
for wild harmonic bundles
in the case $\lambda\neq 1$
because it is rather obvious that a similar argument can work
after~\cite{mochi4,mochi2-I,mochi2-II,mochi5,Mochizuki-wild}
on~the basis of~\cite{s1,s2}.
But, because the Higgs case is also important,
it would be better to have a~reference
in which a rather detailed proof is explained.
It~is one reason why the author writes this manuscript.
As another reason, in the next subsection,
we shall explain an application
to~the correspondence for good wild harmonic bundles
with homogeneity,
which is expected to be useful in the generalized Hodge theory.

\subsection{Homogeneity with respect to group actions}

\subsubsection{Variation of Hodge structure}
\label{subsection;19.2.18.10}

As mentioned, Simpson~\cite{s1} was motivated by
the construction of polarized variation of Hodge structure.
Let us recall the definition of polarized complex variation of
Hodge structure given in~\cite{s1},
instead of the original definition of
polarized variation of Hodge structure due to Griffiths.
A complex variation of Hodge structure of weight $w$
is a graded $C^{\infty}$-vector bundle
$V=\bigoplus_{p+q=w} V^{p,q}$
equipped with a flat connection $\nabla$
satisfying the Griffiths transversality condition,
i.e.,
$\nabla^{0,1}(V^{p,q})
\subset
 \Omega^{0,1}\otimes
 \bigl(
 V^{p+1,q-1}\oplus V^{p,q}\bigr)$
and
$\nabla^{1,0}(V^{p,q})
\subset
 \Omega^{1,0}\otimes
 \bigl(
 V^{p-1,q+1}\oplus V^{p,q}
 \bigr)$,
where $\nabla^{p,q}$ denote the $(p,q)$-part
of $\nabla$.
A polarization of a complex variation of Hodge structure
is a flat Hermitian pairing $\langle\cdot,\cdot\rangle$
satisfying the following conditions:
$(i)$ the decomposition $V=\bigoplus V^{p,q}$ is orthogonal
 with respect to $\langle\cdot,\cdot\rangle$,
$(ii)$ $\big(\sqrt{-1}\big)^{p-q}\langle\cdot,\cdot\rangle$
 is positive definite on $V^{p,q}$.

A polarization of pure Hodge structure
typically appears
when we consider the Gauss--Manin connection
associated with
a smooth projective morphism $f\colon \nbigx\lrarr \nbigy$.
Namely,
the family of vector spaces
$H^{w}\big(f^{-1}(y)\big)$ $(y\in \nbigy)$
naturally induces a flat bundle on $\nbigy$.
With the Hodge decomposition,
it is a variation of Hodge structure of weight $w$.
A relatively ample line bundle
induces a polarization on the variation of Hodge structure.

Simpson discovered a completely different way to construct
a polarized variation of Hodge structure.
Let $(V=\bigoplus V^{p,q},\nabla)$
be a complex variation of Hodge structure.
Note that $\nabla^{0,1}$ induces
holomorphic structures
$\delbar_{V^{p,q}}\colon V^{p,q}\lrarr V^{p,q}\otimes\Omega^{0,1}$
of $V^{p,q}$. We~set
$\delbar_V:=\bigoplus \delbar_{V^{p,q}}$.
Then,
$\big(V=\bigoplus V^{p,q},\delbar_V\big)$ is
a graded holomorphic vector bundle. We~also note that $\nabla^{1,0}$ induces
linear maps $V^{p,q}\lrarr V^{p-1,q+1}\otimes\Omega^{1,0}$,
and hence
$\theta\colon V\lrarr V\otimes\Omega^{1,0}$.
It~is easy to check that
$\theta$ is a~Higgs field of $\big(V,\delbar_V\big)$.
Such a graded holomorphic bundle
$V=\bigoplus_{p+q=w} V^{p,q}$
with a Higgs field $\theta$ such that
$\theta(V^{p,q})\subset V^{p-1,q+1}\otimes\Omega^{1,0}$
is called a Hodge bundle of weight $w$. In~general, we cannot construct
a complex variation of Hodge structure
from a Hodge bundle.
However, Simpson discovered that
if a Hodge bundle $(V=\bigoplus V^{p,q},\theta)$
on a compact K\"ahler manifold
satisfies the stability condition and the vanishing condition,
then there exists a flat connection~$\nabla$
and a flat Hermitian pairing $\langle\cdot,\cdot\rangle$
such that
$(i)$ $(V=\bigoplus V^{p,q},\nabla)$ is a complex variation
of Hodge structure which induces the Hodge bundle,
$(ii)$ $\langle\cdot,\cdot\rangle$ is a polarization
of $(V=\bigoplus V^{p,q},\nabla)$.
Indeed, according to the equivalence of Simpson
between Higgs bundles and harmonic bundles,
there exists a pluri-harmonic metric $h$ of $(V,\theta)$.
It~turns out that
the flat connection
$\nabla_h+\theta+\theta^{\dagger}_h$
satisfies the Griffiths transversality.
Moreover,
the decomposition $V=\bigoplus V^{p,q}$ is orthogonal
with respect to $h$,
and flat Hermitian paring $\langle\cdot,\cdot\rangle$
is constructed by the relation
$\big(\sqrt{-1}\big)^{p-q}
 \langle\cdot,\cdot\rangle_{V^{p,q}}
=h_{|V^{p,q}}$.

Note that a Hodge bundle
is regarded as a Higgs bundle $\big(V,\delbar_V,\theta\big)$
with an $S^1$-homogeneity,
i.e., $\big(V,\delbar_V\big)$ is equipped with an $S^1$-action
such that $t\circ \theta\circ t^{-1}=t\cdot\theta$
for any $t\in S^1$.
It~roug\-hly means that
semistable Hodge bundles correspond to
the fixed points in the moduli space of semistable Higgs bundles
with respect to the natural $S^1$-action
induced by $t\big(E,\delbar_E,\theta\big)=\big(E,\delbar_E,t\theta\big)$.

By the deformation
$\big(E,\delbar_E,\alpha\theta\big)$ $(\alpha\in\cnum^{\ast})$,
any semistable Higgs bundles is deformed to
an~$S^1$-fixed point in the moduli space,
i.e., a semistable Hodge bundle as $\alpha\to 0$.
Note that the Higgs field of the limit is not necessarily $0$.
Hence, by the equivalence between Higgs bundles and flat bundles,
it turns out that
any flat bundle is deformed to a flat bundle underlying
a polarized variation of Hodge structure.

In particular, Simpson~\cite{s1} applied these ideas
to construct uniformizations
of some types of~pro\-jec\=tive manifolds.
He also applied it to prove that
some type of discrete groups cannot be
the fundamental group of any projective manifolds
in~\cite{s5}.

\subsubsection{TE-structure}

We recall that a complex variation of Hodge structure
on $X$ induces
a TE-structure in the sense of Hertling~\cite{Hertling},
i.e.,
a holomorphic vector bundle $\nbigv$
on $\nbigx:=\cnum_{\lambda}\times X$
with a meromorphic flat connection
\begin{gather*}
 \nablatilde\colon\ \nbigv\lrarr\nbigv\otimes\nbigo_{\nbigx}\big(\nbigx^0\big)\otimes
 \Omega_{\nbigx}^1\big(\log\nbigx^0\big),
\end{gather*}
where $\nbigx^0:=\{0\}\times X$.
Indeed,
for a complex variation of Hodge structure
$(V=\bigoplus V^{p,q},\nabla)$,
$F^p(V):=\bigoplus_{p_1\geq p} V^{p_1,q_1}$
are holomorphic subbundles
with respect to $\nabla^{0,1}$.
Thus, we obtain a~decreasing filtration
of holomorphic subbundles
$F^p(V)$ $(p\in\seisuu)$
satisfying the Griffiths transversality
$\nabla^{1,0}F^p(V)\subset F^{p-1}(V)\otimes\Omega^{1,0}$.
Let $p\colon \cnum^{\ast}_{\lambda}\times X\lrarr X$
denote the projection. We~obtain the induced flat bundle $(p^{\ast}V,p^{\ast}\nabla)$. By~the Rees construction,
$p^{\ast}V$ extends to a locally free
$\nbigo_{\nbigx}$-module $\nbigv$,
on which $\nablatilde:=p^{\ast}\nabla$ is a meromorphic flat connection
satisfying the condition
$\nablatilde\nbigv\subset
\nbigv\otimes\nbigo_{\nbigx}\big(\nbigx^0\big)
 \otimes\Omega^1_{\nbigx}\big(\log\nbigx^0\big)$.

It~is recognized that
a TE-structure appears
as a fundamental piece of interesting structures
in~various fields of mathematics.
For instance,
TE-structure is an ingredient of Frobenius manifold,
which is important
in the theory of primitive forms and flat structures due to K.~Saito~\cite{Saito-Takahashi},
the topological field theory of Dubrovin~\cite{Dubrovin93},
the $tt^{\ast}$-geometry of Cecotti--Vafa~\cite{Cecotti-Vafa91, Cecotti-Vafa93},
the Gromov--Witten theory,
the theory of Landau--Ginzburg models,
etc.
For the construction of~Frobenius manifolds,
it is an important step to obtain TE-structures.
Abstractly,
TE-structure is also an important ingredient
of semi-infinite variation of Hodge structure~\cite{Barannikov1, Coates-Iritani-Tseng, Iritani},
TERP structure~\cite{Hertling, hertling-sevenheck, hertling-sevenheck2},
integrable variation of twistor structure~\cite{sabbah2},
etc. (See also~\cite{mochi8,Mochizuki-TodaII}.)

\subsubsection{Homogeneous harmonic bundles}
\label{subsection;19.2.16.10}

As Simpson applied his Kobayashi--Hitchin correspondence
to construct complex variations of~Hodge structure,
we may apply Theorem~\ref{thm;19.2.2.20}
to construct TE-structures
with some additional structure.
It~is done through harmonic bundles with homogeneity
as in the Hodge case.

Let $X$ be a complex manifold equipped with an $S^1$-action.
Let $\big(E,\delbar_E\big)$ be an $S^1$-equivariant holomorphic vector bundle.
Let $\theta$ be a Higgs field of $\big(E,\delbar_E\big)$,
which is homogeneous with respect to the $S^1$-action,
i.e.,
$t^{\ast}\theta=t^m\theta$ for some $m\neq 0$.
Let $h$ be an $S^1$-invariant
pluri-harmonic metric of $\big(E,\delbar_E,\theta\big)$.
Then, as studied in~\cite[Section~3]{Mochizuki-TodaII},
we naturally obtain a TE-structure.
More strongly, it is equipped with
a grading in the sense of~\cite{Coates-Iritani-Tseng, Iritani},
and it also underlies a polarized integrable variation of pure twistor structure
of weight $0$~\cite{sabbah2}.
Moreover,
if there exists an $S^1$-equivariant isomorphism
between
$\big(E,\delbar_E,\theta,h\big)$ and its dual,
the TE-structure is enhanced to a~semi-infinite variation of
Hodge structure with a grading
\cite{Barannikov1, Coates-Iritani-Tseng, Iritani}.
If the $S^1$-action on $X$ is trivial,
this is the same as the construction of
a variation of Hodge structure
from a Hodge bundle with a pluri-harmonic metric
for which the Hodge decomposition is orthogonal.

Let $H$ be a simple normal crossing hypersurface of $X$.
From an $S^1$-homogeneous
 good wild harmonic bundle $\big(E,\delbar_E,\theta,h\big)$
on $(X,H)$,
as mentioned above, we obtain a TE-structure with a~grading
on $X\setminus H$.
Moreover, it extends to a meromorphic TE-structure on $(X,H)$
as studied in~\cite[Section~3]{Mochizuki-TodaII}. We~obtain the mixed Hodge structure
as the limit objects at the boundary,
which is useful for the study of more detailed properties
of the TE-structure.

\subsubsection{An equivalence}

Let $X$ be a complex projective manifold
with a simple normal crossing hypersurface $H$
and an~ample line bundle $L$,
equipped with a $\cnum^{\ast}$-action.
A good filtered Higgs bundle
$(\nbigp_{\ast}\nbigv,\theta)$
is called $\cnum^{\ast}$-homogeneous
if $\nbigp_{\ast}\nbigv$ is $\cnum^{\ast}$-equivariant
and $t^{\ast}\theta=t^m\cdot\theta$ for some $m\neq 0$.
Then,
we~obtain the following theorem by using Theorem~\ref{thm;19.2.2.20}.
(See Section~\ref{subsection;19.2.18.21}
for the precise definition of the stability condition
in this context.)

\begin{Theorem}[Corollary~\ref{cor;19.2.18.20}]\label{thm;19.2.16.20}
There exists an equivalence between the following objects:
\begin{itemize}\itemsep=0pt
\item
 $\mu_L$-polystable
 $\cnum^{\ast}$-homogeneous
 good filtered Higgs bundles $(\nbigp_{\ast}\nbigv,\theta)$
 on $(X,H)$
 satisfying
\begin{gather*}
 \int_X c_1(\nbigp_{\ast}\nbigv)c_1(L)^{\dim X-1}=
 \int_X\ch_2(\nbigp_{\ast}\nbigv)c_1(L)^{\dim X-2}=0.
\end{gather*}
\item
 $S^1$-homogeneous
 good wild harmonic bundles on $(X,H)$.
\end{itemize}
\end{Theorem}

As mentioned in Section~\ref{subsection;19.2.16.10},
Theorem~\ref{thm;19.2.16.20} allows us
to obtain a meromorphic TE-structure on $(X,H)$ with a grading
from a $\mu_L$-polystable
$\cnum^{\ast}$-equivariant good filtered Higgs bundle
satisfying the vanishing condition. We~already applied it to a classification
of solutions of the Toda equations
on $\cnum^{\ast}$~\cite{Mochizuki-TodaI}.
It~seems natural to expect that this construction would be
another way to obtain Frobenius manifolds.

Although we explained the homogeneity with respect to
an $S^1$-action,
Theorem~\ref{thm;19.2.16.20}
is gene\-ra\-li\-zed for $K$-homogeneous good wild harmonic bundles
as explained in Section~\ref{section;19.2.16.30},
where $K$ is any compact Lie group.

\section[Good filtered lambda-flat bundles and wild harmonic bundles]
{Good filtered $\boldsymbol\lambda$-flat bundles and wild harmonic bundles}\label{section;19.2.18.1}

\subsection[Filtered sheaves and filtered lambda-flat sheaves]
{Filtered sheaves and filtered $\boldsymbol\lambda$-flat sheaves}

\subsubsection{Filtered sheaves}\label{subsection;19.2.18.50}

Let $X$ denote a complex manifold
with a simple normal crossing hypersurface $H$.
Let $H=\bigcup_{i\in\Lambda} H_{i}$
denote a decomposition
such that each $H_i$ is smooth.
Note that $H_i$ are not necessarily connected.
For any $P\in H$,
a holomorphic coordinate neighbourhood
$(X_P,z_1,\ldots,z_n)$ around~$P$
is called admissible
if $H_P:=H\cap X_P=\bigcup_{i=1}^{\ell(P)}\{z_i=0\}$.
For such an admissible coordinate neighbourhood,
there exists the map
$\rho_P\colon \{1,\ldots,\ell(P)\}\lrarr \Lambda$
determined by $H_{\rho_P(i)}\cap X_P=\{z_i=0\}$. We~obtain the map
$\kappa_P\colon \real^{\Lambda}\lrarr \real^{\ell(P)}$
by $\kappa_P(\veca)=\big(a_{\rho(1)},\ldots,a_{\rho(\ell(P))}\big)$.

Let $\nbigo_X(\ast H)$ denote the sheaf of
meromorphic functions which may have poles along $H$.
Let~$\nbige$ be any coherent torsion free $\nbigo_X(\ast H)$-module.
A filtered sheaf over $\nbige$
is defined to be a tuple of~coherent $\nbigo_X$-submodules
$\nbigp_{\veca}\nbige\subset\nbige$ $(\veca\in\real^{\Lambda})$
satisfying the following conditions:
\begin{itemize}\itemsep=0pt
\item
 $\nbigp_{\veca}\nbige\subset\nbigp_{\vecb}\nbige$
 if $\veca\leq\vecb$,
 i.e., $a_i\leq b_i$ for any $i\in\Lambda$.
\item
 $\nbigp_{\veca}\nbige(\ast H)=\nbige$ for any
 $\veca\in\real^{\Lambda}$.
\item
 $\nbigp_{\veca+\vecn}\nbige
 =\nbigp_{\veca}\nbige\bigl(
 \sum_{i\in\Lambda} n_iH_i
 \bigr)$
 for any $\veca\in\real^{\Lambda}$
 and $\vecn\in\seisuu^{\Lambda}$.
\item
 For any $\veca\in\real^{\Lambda}$
 there exists $\vecepsilon\in\real_{>0}^{\Lambda}$
 such that
 $\nbigp_{\veca+\vecepsilon}\nbige=\nbigp_{\veca}\nbige$.
\item
 For any $P\!\in H$,
 we take an admissible coordinate neighbourhood
 $(X_P,z_1,\ldots,z_n)$ around~$P$.
 Then, for any $\veca\in\real^{\Lambda}$,
 $\nbigp_{\veca}\nbige_{|X_P}$
 depends only on $\kappa_P(\veca)$.
\end{itemize}
For any coherent $\nbigo_X(\ast H)$-submodule $\nbige'\subset\nbige$,
we obtain a filtered sheaf
$\nbigp_{\ast}\nbige'$ over $\nbige'$
by $\nbigp_{\veca}\nbige':=\nbigp_{\veca}\nbige\cap\nbige'$.
If $\nbige'$ is saturated, i.e.,
$\nbige'':=\nbige/\nbige'$ is torsion-free,
we obtain a filtered sheaf
$\nbigp_{\ast}\nbige''$ over $\nbige''$
by $\nbigp_{\veca}\nbige'':=
 \Image\bigl(
 \nbigp_{\veca}\nbige\lrarr\nbige''
 \bigr)$.

A morphism of filtered sheaves
$f\colon \nbigp_{\ast}\nbige_1\lrarr\nbigp_{\ast}\nbige_2$
is defined to be a morphism
$f\colon \nbige_1\lrarr\nbige_2$ of
$\nbigo_X(\ast H)$-modules
such that
$f(\nbigp_{\veca}\nbige_1)\subset\nbigp_{\veca}\nbige_2$
for any $\veca\in\real^{\Lambda}$.

\begin{Remark}
The concept of filtered bundles on curves
was introduced by
Mehta and Sesha\-dri~\cite{mehta-seshadri}
and Simpson~\cite{s1, s2}.
A higher dimensional version
was first studied by Maruyama and~Yokogawa~\cite{my}
for the purpose of the construction of the moduli spaces.
\end{Remark}

\subsubsection{Restriction and gluing}

Let $U\subset X$ be any open subset. We~set $H_U=H\cap U$.
Let $H_U=\bigcup_{j\in \Lambda_U}H_{U,j}$
be the irreducible decomposition.
For any $j\in\Lambda_U$,
we have $i(j)\in \Lambda$ such that
$H_{U,j}$ is a connected component of~$H_{i(j)}\cap U$.
For any $P\in H_U$,
we set $\Lambda_U(P):=\{j\in\Lambda_U\mid P\in H_{U,j}\}$.

Let $\nbigp_{\ast}\nbige$ be a filtered sheaf over $\nbige$. We~shall define a filtered sheaf
over the $\nbigo_U(\ast H_U)$-module~$\nbige_{|U}$.
Let $\vecb\in\real^{\Lambda_U}$.
For any $P\in H_U$,
we choose $\veca(P,\vecb)\in\real^{\Lambda}$
such that
$a(P,\vecb)_{i(j)}=b_j$
for any $j\in \Lambda_U(P)$,
and we obtain the following
$\nbigo_{U,P}$-submodule of
the stalk $(\nbige_{|U})_P$:
\begin{gather*}
 \nbigp_{\vecb}(\nbige_{|U})_P:=
 \nbigp_{\veca(P,\vecb)}(\nbige)_P.
\end{gather*}
It~is independent of the choice of
$\veca(P,\vecb)$ as above.
There uniquely exists
a coherent $\nbigo_{U}$-submodule
$\nbigp_{\vecb}(\nbige_{|U})$
of $\nbige_{|U}$
such that
$(i)$ $\nbigp_{\vecb}(\nbige_{|U})(\ast H_U)=\nbige_{|U}$,
and
$(ii)$ for any $P\in H_U$,
the stalk of $\nbigp_{\vecb}(\nbige_{|U})$ at $P$ is
equal to $\nbigp_{\vecb}(\nbige_{|U})_P$.
Thus, we obtain a filtered sheaf
$\nbigp_{\ast}(\nbige_{|U})$ over $\nbige_{|U}$,
which is denoted as
$\nbigp_{\ast}\nbige_{|U}$.

Let $X=\bigcup_{k\in \Gamma}X^{(k)}$ be an open covering. We~set $H^{(k)}=H\cap X^{(k)}$.
For any filtered sheaf $\nbigp_{\ast}\nbige$ over $\nbige$,
we obtain
filtered sheaves
$\nbigp_{\ast}\nbige_{|X^{(k)}}$
over $\nbige_{|X^{(k)}}$
as the restriction.
Conversely,
let $\nbigp_{\ast}\big(\nbige_{|X^{(k)}}\big)$ $(k\in\Gamma)$
be filtered sheaves over $\nbige_{|X^{(k)}}$
such that
$\nbigp_{\ast}\big(\nbige_{|X^{(k)}}\big)_{|X^{(k)}\cap X^{(\ell)}}
=\nbigp_{\ast}\big(\nbige_{|X^{(\ell)}}\big)_{|X^{(k)}\cap X^{(\ell)}}$
for any $k,\ell\in\Gamma$.
\begin{Lemma}
 There uniquely exists
 a filtered sheaf $\nbigp_{\ast}\nbige$ over $\nbige$
 such that
 $\nbigp_{\ast}\nbige_{|X^{(k)}}
=\nbigp_{\ast}\big(\nbige_{|X^{(k)}}\big)$
for any $k\in\Gamma$.
\end{Lemma}
\begin{proof}
Let $\veca\in\Lambda$.
For any $P\in H$,
there exists $k\in\Gamma$
such that $P\in X^{(k)}$.
Let $H^{(k)}=\bigcup_{j\in \Lambda^{(k)}}H^{(k)}_j$
be the irreducible decomposition.
For any $j\in\Lambda^{(k)}$,
we have $i(k,j)\in \Lambda$ such that~$H^{(k)}_j$ is a connected component of
$H_{i(k,j)}\cap X^{(k)}$.
Thus, we obtain a map
$\Lambda^{(k)}\lrarr \Lambda$.
For any $\veca\in\real^{\Lambda}$,
let $\veca^{(k)}$ be the image of $\veca$
by the induced map
$\real^{\Lambda}\lrarr \real^{\Lambda^{(k)}}$,
and we obtain the following
$\nbigo_{X,P}$-submodule of $\nbige_P$:
\begin{gather*}
 \nbigp_{\veca}(\nbige)_{P}:=
 \nbigp_{\veca^{(k)}}\big(\nbige_{|X^{(k)}}\big)_{P}.
\end{gather*}
There uniquely exists a coherent $\nbigo_X$-submodule
$\nbigp_{\veca}\nbige$ of $\nbige$
such that
$(i)$ $\nbigp_{\veca}\nbige(\ast H)=\nbige$,
and~$(ii)$~for any $P\in H$,
the stalk of $\nbigp_{\veca}(\nbige)$ at $P$
is equal to $\nbigp_{\veca}(\nbige)_{P}$.
Thus, we obtain a filtered sheaf~$\nbigp_{\ast}\nbige$ over $\nbige$
with the desired property.
The uniqueness is also clear.
\end{proof}

\subsubsection{Reflexive filtered sheaves}

A filtered sheaf $\nbigp_{\ast}\nbige$ on $(X,H)$
is called reflexive
if each $\nbigp_{\veca}\nbige$ is a reflexive $\nbigo_X$-module.
Note that it is equivalent to
the ``reflexive and saturated'' condition
in~\cite[Definition 3.17]{mochi4}
by the following lemma.

\begin{Lemma}
Suppose that
$\nbigp_{\ast}\nbige$ is reflexive.
Let $\veca\in\real^{\Lambda}$. We~take $a_i-1<b\leq a_i$,
and let $\veca'\in\real^{\Lambda}$
be determined by
$a_j'=a_j$ $(j\neq i)$ and $a'_i=b$.
Then,
$\nbigp_{\veca}\nbige/\nbigp_{\veca'}\nbige$
is a torsion-free $\nbigo_{H_i}$-module.
\end{Lemma}

\begin{proof}
Let $s$ be a section of
$\nbigp_{\veca}\nbige/\nbigp_{\veca'}\nbige$
on an open set $U\subset D_i$.
There exists an open subset $\Utilde\subset X$
and a section $\stilde$ of $\nbigp_{\veca}\nbige$ on $\Utilde$
such that $\Utilde\cap D_i=U$
and that $\stilde$ induces $s$.
Note that there exists $Z\subset \Utilde$
of codimension $2$
such that $\stilde_{|\Utilde\setminus Z}$
is a section of
$\nbigp_{\veca'}\nbige_{|\Utilde\setminus Z}$.
Because $\nbigp_{\veca'}\nbige$ is~reflexive,
there exists a section $\stilde'$ of
$\nbigp_{\veca'}\nbige$ on $\Utilde$
such that
$\stilde'_{|\Utilde\setminus Z}=\stilde_{|\Utilde\setminus Z}$.
Hence, we obtain that~$\stilde$ is a section of
$\nbigp_{\veca'}\nbige$,
i.e., $s=0$.
\end{proof}

The following lemma is clear.
\begin{Lemma}
Let $\nbigp_{\ast}\nbige$ be a reflexive filtered sheaf on $(X,H)$.
Then a coherent $\nbigo_X(\ast H)$-submodule $\nbige'\subset\nbige$
is saturated if and only if
the induced filtered sheaf $\nbigp_{\ast}\nbige'$ is reflexive.
\end{Lemma}

\subsubsection[Filtered lambda-flat sheaves]{Filtered $\boldsymbol\lambda$-flat sheaves}

Let $\lambda$ be any complex number.
Let $\nbige$ be a coherent torsion-free $\nbigo_X(\ast H)$-module.
A $\lambda$-connection
$\DDlambda\colon \nbige\lrarr \Omega_X^1\otimes\nbige$
is a $\cnum$-linear morphism of sheaves such that
$\DDlambda(fs)=f\DDlambda(s)+\lambda {\rm d}f\otimes s$
for any local sections $f$ and $s$
of $\nbigo_X$ and $\nbige$, respectively.
Note that
an $\nbigo_X$-morphism
$\DDlambda\circ\DDlambda\colon \nbige\lrarr
 \Omega_X^2\otimes\nbige$ is induced.
If $\DDlambda\circ\DDlambda=0$,
it is called a flat $\lambda$-connection.
When $\nbige$ is equipped with a flat $\lambda$-connection,
a $\lambda$-flat subsheaf of $\nbige$
means a coherent $\nbigo_X$-submodule $\nbige'\subset\nbige$
such that
$\DDlambda(\nbige')\subset
 \Omega_X^1\otimes\nbige'$.
A pair of a filtered sheaf $\nbigp_{\ast}\nbige$ over $\nbige$
and a flat $\lambda$-connection
$\DDlambda$ of $\nbige$ is called
a filtered $\lambda$-flat connection.
It~is called reflexive
if $\nbigp_{\ast}\nbige$ is reflexive.

\subsection[mu L-stability condition for filtered lambda-flat sheaves]
{$\boldsymbol{\mu_L}$-stability condition for filtered $\boldsymbol\lambda$-flat sheaves}

Let $X$ be a connected projective manifold with a simple normal crossing
hypersurface $H=\bigcup_{i\in\Lambda}H_i$.
Let $L$ be an ample line bundle.

\subsubsection{Slope of filtered sheaves}
\label{subsection;19.2.18.100}

Let $\nbigp_{\ast}\nbige$ be a filtered sheaf on $(X,H)$.
Recall the definition of
the parabolic first Chern class~$c_1(\nbigp_{\ast}\nbige)$.
Let $\eta_i$ be the generic point of $H_i$.
Note that
$\nbigo_{X,\eta_i}$-modules
$(\nbigp_{\veca}\nbige)_{\eta_i}$
depends only on $a_i$,
which is denoted by $\nbigp_{a_i}(\nbige_{\eta_i})$. We~obtain $\nbigo_{H_i,\eta_i}$-modules
$\Gr^{\nbigp}_{a}(\nbige_{\eta_i}):=
 \nbigp_{a}(\nbige_{\eta_i})\big/
 \nbigp_{<a}(\nbige_{\eta_i})$.
Then, we set
\begin{gather}
\label{eq;20.7.8.1}
 c_1(\nbigp_{\ast}\nbige):=
 c_1(\nbigp_{\veca}\nbige)
-\sum_{i\in\Lambda}
 \sum_{a_i-1<a\leq a_i}
 a\rank \Gr^{\nbigp}_a(\nbige_{\eta_i})[H_i]
\in H^2(X,\real).
\end{gather}
Here, $[H_i]$ denote the cohomology class induced by
$H_i$.
It~is easy to see that
$c_1(\nbigp_{\ast}\nbigv)$
is independent of the choice of $\veca\in\real^{\Lambda}$. We~set
\begin{gather*}
 \mu_L(\nbigp_{\ast}\nbige):=
 \frac{1}{\rank\nbige}
 \int_Xc_1(\nbigp_{\ast}\nbige)
 \cdot c_1(L)^{n-1}.
\end{gather*}
It~is called the slope of $\nbigp_{\ast}\nbige$
with respect to $L$.
The following is proved in
\cite[Lemma 3.7]{mochi4}.
\begin{Lemma}
\label{lem;19.2.18.120}
Let $f\colon \nbigp_{\ast}\nbige^{(1)}\lrarr\nbigp_{\ast}\nbige^{(2)}$
be a morphism of filtered sheaves
which is generically an isomorphism,
i.e.,
the induced morphism
$\nbige^{(1)}_{\eta(X)}\lrarr\nbige^{(2)}_{\eta(X)}$
at the generic point of $X$ is an~iso\-mor\-phism.
Then, $\mu_L\big(\nbigp_{\ast}\nbige^{(1)}\big)
\leq
 \mu_L\big(\nbigp_{\ast}\nbige^{(2)}\big)$ holds.
If the equality holds,
$f$ is an isomorphism in codimension one,
i.e.,
there exists an algebraic subset $Z\subset X$
such that
$(i)$~the codimension of~$Z$ is larger than $2$,
$(ii)$~$f_{|X\setminus Z}\colon
 \nbigp_{\ast}\nbige^{(1)}_{|X\setminus Z}
\lrarr
 \nbigp_{\ast}\nbige^{(2)}_{|X\setminus Z}$
is an isomorphism.
\end{Lemma}

\subsubsection[mu L-stability condition]{$\boldsymbol{\mu_L}$-stability condition}

A filtered $\lambda$-flat sheaf
$(\nbigp_{\ast}\nbige,\DDlambda)$ on $(X,H)$
is called $\mu_L$-stable
(resp.~$\mu_L$-semistable)
if the following holds:
\begin{itemize}\itemsep=0pt
\item
 Let $\nbige'\subset\nbige$
 be any $\lambda$-flat $\nbigo_X(\ast H)$-submodule
 such that
$0<\rank(\nbige')<\rank(\nbige)$.
 Then,
 $\mu_L(\nbigp_{\ast}\nbige')<\mu_L(\nbigp_{\ast}\nbige)$
(resp.~$\mu_L(\nbigp_{\ast}\nbige')\leq\mu_L(\nbigp_{\ast}\nbige)$)
 holds.
\end{itemize}
A filtered $\lambda$-flat sheaf
$\big(\nbigp_{\ast}\nbige,\DDlambda\big)$
is called $\mu_L$-polystable
if the following holds:
\begin{itemize}\itemsep=0pt
\item $\big(\nbigp_{\ast}\nbige,\DDlambda\big)$ is $\mu_L$-semistable.
\item
 $\big(\nbigp_{\ast}\nbige,\DDlambda\big)
 =\bigoplus \big(\nbigp_{\ast}\nbige_i,\DDlambda_i\big)$,
 where each $\big(\nbigp_{\ast}\nbige_i,\DDlambda_i\big)$
 is $\mu_L$-stable.
\end{itemize}
The following is standard.
(See~\cite[Section~3.1.3]{mochi4} and~\cite[Section~2.1.4]{mochi5}.)
\begin{Lemma}
Suppose that
$(\nbigp_{\ast}\nbige,\DDlambda)$ is
a $\mu_L$-polystable reflexive filtered $\lambda$-flat sheaf.
Then, there exists a unique decomposition
$\big(\nbigp_{\ast}\nbige,\DDlambda\big)
=\bigoplus_{i=1}^N
 \big(\nbigp_{\ast}\nbige_i,\DDlambda_i\big)\otimes\cnum^{m(i)}$
such that
$(i)$ $\big(\nbigp_{\ast}\nbige_i,\DDlambda_i\big)$
are $\mu_L$-stable,
$(ii)$ $\mu_L(\nbigp_{\ast}\nbige_i)=\mu_L(\nbigp_{\ast}\nbige)$,
$(iii)$ $\big(\nbigp_{\ast}\nbige_i,\DDlambda_i\big)\not\simeq
 \big(\nbigp_{\ast}\nbige_j,\DDlambda_j\big)$ $(i\neq j)$.
\end{Lemma}

\begin{Remark}
In {\rm~\cite[Section~3.1.3]{mochi4}},
``the inequality
$\pardeg_L(\nbige'_{\ast})
<\pardeg_L(\nbige_{\ast})$''
should be corrected to
``the inequality
$\mu_L(\nbige'_{\ast})
<\mu_L(\nbige_{\ast})$''.
\end{Remark}

\subsection{Filtered bundles}

\subsubsection{Filtered bundles in the local case}\label{subsection;19.1.19.30}

We recall the notion of filtered bundle
in the local case. We~shall explain it in the global case
in Section~\ref{subsection;19.1.19.40}.
Let $U$ be a neighbourhood of
$(0,\ldots,0)$ in $\cnum^n$. We~set $H_{U,i}:=U\cap\{z_i=0\}$,
and
$H_U:=\bigcup_{i=1}^{\ell}H_{U,i}$
for some $0\leq \ell\leq n$.
Let $\nbigv$ be a locally free
$\nbigo_U(\ast H_U)$-module.
A filtered bundle
$\nbigp_{\ast}\nbigv$ over $\nbigv$
is a tuple of locally free $\nbigo_U$-submodules
$\nbigp_{\veca}\nbigv$ $(\veca\in\real^{\ell})$
such that the following holds:
\begin{itemize}\itemsep=0pt
\item
 $\nbigp_{\veca}\nbigv\subset\nbigp_{\vecb}\nbigv$
 if $\veca\leq\vecb$,
 i.e., $a_i\leq b_i$ for any $i=1,\ldots,\ell$.
\item
 There exists a frame
 $\vecv=(v_1,\ldots,v_r)$
 of $\nbigv$
 and tuples $\veca(v_j)\in\real^{\ell}$ $(j=1,\ldots,r)$
 such that
\begin{gather}
\label{eq;21.5.4.1}
 \nbigp_{\vecb}\nbigv
=\bigoplus_{j=1}^r
 \nbigo_{U}\bigg(
 \sum_i \bigl[b_i-a_i(v_j)\bigr]H_{U,i}
 \bigg)\cdot v_j,
\end{gather}
where we set $[c]:=\max\{p\in\seisuu\mid p\leq c\}$
for any $c\in\real$.
\end{itemize}
Clearly,
a filtered bundle over $\nbigv$
is a filtered sheaf over $\nbigv$.

\begin{Remark}
\label{rem;20.2.14.2}
We set $\nbigr:=\cnum[\![z_1,\ldots,z_n]\!]$
and $\nbigrtilde:=\nbigr\big[z_1^{-1},\ldots,z_{\ell}^{-1}\big]$.
For a free $\nbigrtilde$-module $\nbigvhat$,
a~filtered bundle over $\nbigvhat$
is defined to be a tuple $\nbigp_{\ast}\nbigvhat:=
 \bigl( \nbigp_{\veca}\nbigvhat\mid \veca\in\real^{\ell(P)}\bigr)$
of free $\nbigr$-submodules
satisfying similar conditions as above.
\end{Remark}

\subsubsection[Pull back, push-forward and descent with respect to ramified coverings in the local case]
{Pull back, push-forward and descent with respect to ramified coverings\\ in the local case}
\label{subsection;20.2.14.1}

Let $\varphi\colon \cnum^n\lrarr\cnum^n$
be given by
$\varphi(\zeta_1,\ldots,\zeta_n)=
 \big(\zeta_1^{m_1},\ldots,\zeta_{\ell}^{m_{\ell}},
 \zeta_{\ell+1},\ldots,\zeta_n\big)$. We~set $U':=\varphi^{-1}(U)$,
$H_{U',i}:=\varphi^{-1}(H_{U,i})$
and $H_{U'}:=\varphi^{-1}(H_U)$.
The induced ramified covering
$U'\lrarr U$ is also denoted by $\varphi$.

For any $\vecb\in\real^{\ell}$,
we set
$\varphi^{\ast}(\vecb)=(m_ib_i)\in\real^{\ell}$.
For any filtered bundle $\nbigp_{\ast}\nbigv_1$ on
$(U,H_U)$,
we~define a filtered bundle $\nbigp_{\ast}\nbigv'_1$
on $(U',H_{U'})$ as follows:
\begin{gather*}
 \nbigp_{\veca}\nbigv'_1
=\sum_{\varphi^{\ast}(\vecb)+\vecn\leq\veca}
 \varphi^{\ast}\bigl(\nbigp_{\vecb}\nbigv_1\bigr)
 \Bigl(\sum n_iH_{U',i}\Bigr).
\end{gather*}
We set $\varphi^{\ast}(\nbigp_{\ast}\nbigv_1):=\nbigp_{\ast}\nbigv'_1$.
Thus,
we obtain the pull back functor $\varphi^{\ast}$
from the category of filtered bundles on $(U,H_U)$
to the category of filtered bundles on $(U',H_{U'})$.

For any $\vecb\in\real^{\ell}$,
we set $\varphi_{\ast}(\vecb)=(m_i^{-1}b_i)$.
For any filtered bundle $\nbigp_{\ast}\nbigv_2$ on $(U',H_{U'})$,
we obtain the following filtered bundle
\begin{gather*}
 \nbigp_{\vecb}\varphi_{\ast}(\nbigv_2):=
 \varphi_{\ast}\nbigp_{\varphi_{\ast}(\vecb)}\nbigv_2.
\end{gather*}
In this way, we obtain a functor $\varphi_{\ast}$
from the category of filtered bundles on $(U',H_{U'})$
to the category of filtered bundles on $(U,H_U)$.

We set $G:=\prod_{i=1}^{\ell}\big\{\mu_i\in\cnum^{\ast}\mid \mu_i^{m_i}=1\big\}$. We~define the action of $G$ on $U'$
by
\begin{gather*}
 (\mu_1,\ldots,\mu_{\ell})(\zeta_1,\ldots,\zeta_n)
=(\mu_1\zeta_1,\ldots,\mu_{\ell}\zeta_{\ell},\zeta_{\ell+1},\ldots,\zeta_n).
\end{gather*}
We identify $G$ as the Galois group of the ramified covering
$U'\lrarr U$.
Let $\nbigp_{\ast}\nbigv_3$ be a $G$-equivariant
filtered bundles on $(U',H_{U'})$.
Then,
$\nbigp_{\ast}\varphi_{\ast}\nbigv_3$
is equipped with an induced $G$-action. We~obtain a filtered bundle
$(\nbigp_{\ast}\varphi_{\ast}\nbigv_3)^{G}$
on $(U,H_U)$
as the $G$-invariant part of $\nbigp_{\ast}\varphi_{\ast}\nbigv_3$,
which is~called the descent of $\nbigp_{\ast}\nbigv_3$
with respect to the $G$-action. In~this way,
we obtain a functor
from the category of $G$-equivariant
filtered bundles on $(U',H_{U'})$
to the category of filtered bundles on $(U,H_U)$.

For a filtered bundle $\nbigp_{\ast}\nbigv_1$ on $(U,H_U)$,
the pull back
$\varphi^{\ast}(\nbigp_{\ast}\nbigv_1)$ is
a $G$-equivariant filtered bundle on $(U',H_{U'})$,
and its descent is naturally isomorphic to
$\nbigp_{\ast}\nbigv_1$.

\subsubsection{Filtered bundles in the global case}
\label{subsection;19.1.19.40}

We use the notation in Section~\ref{subsection;19.2.18.50}.
Let $\nbigv$ be a locally free $\nbigo_X(\ast H)$-module.
A filtered bundle
$\nbigp_{\ast}\nbigv=\bigl(
\nbigp_{\veca}\nbigv\mid
 \veca\in\real^{\Lambda}
 \bigr)$ over $\nbigv$
be a sequence of locally free
$\nbigo_X$-submodules
$\nbigp_{\veca}\nbigv$ of $\nbigv$
such that the following holds:
\begin{itemize}\itemsep=0pt
\item
For any $P\!\in H$,
we take an admissible coordinate neighbourhood $(X_P,z_1,\ldots,z_n)$
around~$P$.
Then, for any $\veca\in\real^{\Lambda}$,
$\nbigp_{\veca}\nbigv_{|X_P}$
depends only on $\kappa_P(\veca)$,
denoted as
$\nbigp^{(P)}_{\kappa_P(\veca)}(\nbigv_{|X_P})$.
\item
The sequence
$\big(\nbigp^{(P)}_{\vecb}(\nbigv_{|X_P})\mid \vecb\in\real^{\ell(P)}\big)$
is a filtered bundle over $\nbigv_{|X_P}$
in the sense of~Sec\-tion~\ref{subsection;19.1.19.30}.
\end{itemize}
In other words,
a filtered bundle is a filtered sheaf
(see Section~\ref{subsection;19.2.18.50})
satisfying the condition in~Section~\ref{subsection;19.1.19.30}
locally around any point of $H$.

\begin{Remark}
The higher dimensional version of filtered bundles
was introduced in~\cite{mochi2-I,mochi2-II}
with a different formulation.
See also~\cite{borne1, borne2}. In~this paper, we essentially follow Iyer and Simpson~\cite{i-s}.
\end{Remark}

\subsubsection{The induced bundles and filtrations}

For any $I\subset\Lambda$,
let $\vecdelta_I\in\real^{\Lambda}$
be the element
whose $j$-th component is $0$ $(j\in\Lambda\setminus I)$
or $1$ $(j\in I)$. We~also set $H_I:=\bigcap_{i\in I}H_i$
and
$\del H_I:=H_I\cap\bigl(
 \bigcup_{j\in \Lambda\setminus I}H_j\bigr)$.

Let $\nbigp_{\ast}\nbigv$ be a filtered bundle on $(X,H)$.
Take $i\in \Lambda$.
Let $\veca\in\real^{\Lambda}$.
For any $a_i-1\leq b\leq a_i$,
we set $\veca(b,i):=\veca+(b-a_i)\vecdelta_i$. We~set
\begin{gather*}
 \lefttop{i}F_b\bigl(
 \nbigp_{\veca}(\nbigv)_{|H_i}
 \bigr)
:=\nbigp_{\veca(b,i)}\nbigv/\nbigp_{\veca(a_i-1,i)}\nbigv.
\end{gather*}
It~is naturally regarded as
a locally free $\nbigo_{H_i}$-module.
Moreover, it is a subbundle of
$\nbigp_{\veca}(\nbigv)_{|H_i}$. In~this way,
we obtain a filtration
$\lefttop{i}F$ of
$\nbigp_{\veca}(\nbigv)_{|H_i}$
indexed by
$\openclosed{a_i-1}{a_i}$. We~shall also denote it~as just $F$
if there is no risk of confusion.

We obtain the induced filtrations
$\lefttop{i}F$
of $\nbigp_{\veca}\nbigv_{|H_I}$
if $i\in I$.
Let $\veca_I\in\real^I$ denote the image of $\veca$
by the projection $\real^{\Lambda}\lrarr\real^I$.
Set $\openclosed{\veca_I-\vecdelta_I}{\veca_I}:=
 \prod_{i\in I}\openclosed{a_i-1}{a_i}$.
For any $\vecb\in\openclosed{\veca_I-\vecdelta_I}{\veca_I}$,
we~set
\begin{gather*}
 \lefttop{I}F_{\vecb}\big(\nbigp_{\veca}\nbigv_{|H_I}\big):=
 \bigcap_{i\in I}
 \lefttop{i}F_{b_i}\big(\nbigp_{\veca}\nbigv_{|H_I}\big).
\end{gather*}
By the condition of filtered bundles,
the following compatibility condition holds.
\begin{itemize}\itemsep=0pt
\item
 Let $P$ be any point of $H_I$.
 There exist a neighbourhood $X_P$ of $P$ in $X$
 and a non-canonical decomposition
\begin{gather*}
 \nbigp_{\veca}\nbigv_{|X_P\cap H_I}
=\bigoplus_{\vecb\in \openclosed{\veca_I-\vecdelta_I}{\veca_I}}
 \nbigg_{P,\vecb}
\end{gather*}
such that the following holds
for any $\vecc\in \openclosed{\veca_I-\vecdelta_I}{\veca_I}$:
\begin{gather}
\label{eq;21.5.4.2}
 \lefttop{I}F_{\vecc}(\nbigp_{\veca}\nbigv_{|H_I\cap X_P})
=\bigoplus_{\vecb\leq\vecc}
 \nbigg_{P,\vecb}.
\end{gather}
 Indeed, there exists a frame $\vecv=(v_1,\ldots,v_r)$
 of $\nbigp_{\veca}\nbigv$ around $P$
 with tuples $\veca(v_i)\in\real^{\ell(P)}$ of real numbers
 satisfying (\ref{eq;21.5.4.1}),
 where $\vecb$ is replaced with $\veca$.
 There exists the bijection
 $\kappa\colon I\simeq \{1,\ldots,\ell(P)\}$
 determined by $H_i\cap X_P=\{z_{\kappa(i)}=0\}$,
 by which we identify $I$ with $\{1,\ldots,\ell(P)\}$.
 Let $\nbigg_{P,\vecb}$
 be the subbundle of
 $\nbigp_{\veca}\nbigv_{|X_P\cap H_I}$
 generated by $v_{j|X_P\cap H_I}$
 satisfying
 $\veca(v_j)=\vecb$.
 Then, we obtain the decomposition (\ref{eq;21.5.4.2}).
\end{itemize}
For any $\vecc\in\openclosed{\veca_I-\vecdelta_I}{\veca_I}$,
we obtain the following locally free $\nbigo_{H_I}$-modules:
\begin{gather*}
 \lefttop{I}\Gr^{F}_{\vecc}(\nbigp_{\veca}\nbigv)
:=\frac{\lefttop{I}F_{\vecc}\big(\nbigp_{\veca}\nbigv_{|H_I}\big)}{
 \sum_{\vecb\lneq\vecc}\lefttop{I}F_{\vecb}\big(\nbigp_{\veca}\nbigv_{|H_I}\big)}.
\end{gather*}
Here, $\vecb=(b_i)\lneq\vecc=(c_i)$ means
that $b_i\leq c_i$ for any $i$
and that $\vecb\neq\vecc$.
Clearly,
if $\vecc\in\openclosed{\veca'_I-\vecdelta_I}{\veca'_I}$
and $\veca'_{\Lambda\setminus I}=\veca_{\Lambda\setminus I}$,
we obtain
$\lefttop{I}\Gr^F_{\vecc}(\nbigp_{\veca}\nbigv)
=\lefttop{I}\Gr^F_{\vecc}(\nbigp_{\veca'}\nbigv)$.

\subsubsection{The induced filtered bundles}

For $\vecc\in\real^I$,
we choose
$\veca\in\real^{\Lambda}$
such that $\vecc\in\openclosed{\veca_I-\vecdelta_I}{\veca_I}$,
and we obtain the following
$\nbigo_{H_I}(\ast \del H_I)$-module:
\begin{gather*}
 \lefttop{I}\Gr^{F}_{\vecc}(\nbigv):=
 \lefttop{I}\Gr^{F}_{\vecc}(\nbigp_{\veca}\nbigv)(\ast \del H_I).
\end{gather*}
It~is independent of the choice of $\veca$ as above. We~obtain the irreducible decomposition
$\del H_I=\bigcup_{i\in\Lambda(I)}H_{I,i}$.
For any $i\in\Lambda(I)$,
there exists $j(i)\in \Lambda\setminus I$
such that
$H_{I,i}$ is a connected component of~$H_I\cap H_{j(i)}$.
Let $\vecd\in \real^{\Lambda(I)}$.
For $P\in \del H_I$,
there exists $\veca(\vecc,\vecd,P)\in\real^{\Lambda}$
such that
$(i)$~$a(\vecc,\vecd,P)_j=c_j$ $(j\in I)$,
$(ii)$~$a(\vecc,\vecd,P)_{j(i)}=d_{i}$ $(i\in \Lambda(I),\,P\in H_{I,i})$. We~obtain an $\nbigo_{X,P}$-submodule
\begin{gather*}
 \nbigp_{\vecd}\bigl( \lefttop{I}\Gr^{F}_{\vecc}(\nbigv) \bigr)_P:=
 \lefttop{I}\Gr^F_{\vecc}\bigl(\nbigp_{\veca(\vecc,\vecd,P)}\nbigv\bigl)_P
 \subset \lefttop{I}\Gr^{F}_{\vecc}(\nbigv)_P.
\end{gather*}
Note that
$\nbigp_{\vecd}\bigl(
 \lefttop{I}\Gr^{F}_{\vecc}(\nbigv)
 \bigr)_P$
is independent of the choice of $\veca(P)$.
There uniquely exists
an $\nbigo_{H_I}$-submodule
$\nbigp_{\vecd}\bigl(
 \lefttop{I}\Gr^{F}_{\vecc}(\nbigv)
 \bigr)
 \subset
\lefttop{I}\Gr^{F}_{\vecc}(\nbigv)$
whose stalk at $P$ $(P\in \del H_I)$
 are equal to
 $\nbigp_{\vecd}\bigl(
 \lefttop{I}\Gr^{F}_{\vecc}(\nbigv)
 \bigr)_P$.
 Thus, we obtain the following filtered bundle
 over $\lefttop{I}\Gr^F_{\vecc}(\nbigv)$
 on $(H_I,\del H_I)$:
\begin{gather*}
\lefttop{I}\Gr^{F}_{\vecc}(\nbigp_{\ast}\nbigv):=
\bigl( \nbigp_{\vecd}\bigl( \lefttop{I}\Gr^{F}_{\vecc}(\nbigv)
 \bigr)\mid\vecd\in\real^{\Lambda(I)}\bigr).
\end{gather*}

\subsubsection{First and second Chern characters for filtered bundles}

Let $\nbigp_{\ast}\nbigv$
be a filtered bundle over $(X,H)$.
Take any $\veca\in\real^{\Lambda}$.
As recalled in Section~\ref{subsection;19.2.18.100},
we~obtain the parabolic first Chern class:
\begin{gather*}
 c_1(\nbigp_{\ast}\nbigv)=
 c_1(\nbigp_{\veca}\nbigv)
-\sum_{i\in\Lambda}
 \sum_{a_i-1<b\leq a_i}
 a_i\rank\lefttop{i}\Gr^{F}_{b}(\nbigp_{\veca}\nbige_{|H_i})
 \cdot [H_i]
\in H^2(X,\real).
\end{gather*}

To explain the second parabolic Chern character in $H^4(X,\real)$,
let us introduce some notation.
Let $\Irr(H_i\cap H_j)$ be the set of the irreducible components
of $H_i\cap H_j$.
For $C\in \Irr(H_I)$,
let
$[C]\in H^4(X,\real)$ denotes the induced cohomology class,
and let $\lefttop{C}\Gr^F_{\vecc}(\nbigp_{\veca}\nbigv)$
denote the restriction of
$\lefttop{I}\Gr^F_{\vecc}(\nbigp_{\veca}\nbigv)$
to $C$.
Moreover,
$\iota_{i\ast}\colon H^2(H_i,\real)\lrarr H^4(X,\real)$
denotes the Gysin map
induced by~$\iota_i\colon H_i\lrarr X$.
Then, the second parabolic Chern character is given as follows.
\begin{align*}
 \ch_2(\nbigp_{\ast}\nbigv):={}& \ch_2(\nbigp_{\veca}\nbigv)
-\sum_{i\in \Lambda} \sum_{a_i-1<b\leq a_i}
 b\cdot \iota_{i\ast}\bigl( c_1\bigl(
 \lefttop{i}\Gr^F_b\big(\nbigp_{\veca}\nbigv_{|H_i}\big) \bigr) \bigr)
 \\
&+\frac{1}{2} \sum_{i\in\Lambda} \sum_{a_i-1<b\leq a_i}
 b^2\rank\big( \lefttop{i}\Gr^F_{b}(\nbigp_{\veca}\nbigv) \big)
\cdot [H_i]^2
\\
&+\frac{1}{2}
 \sum_{\substack{(i,j)\in \Lambda^2\\ i\neq j}}
 \sum_{C\in\Irr(H_i\cap H_j)}
 \sum_{\substack{a_i-1<c_i\leq a_i\\ a_j-1<c_j\leq a_j}}
 c_i\cdot c_j \rank
 \lefttop{C}\Gr^F_{(c_i,c_j)}(\nbigp_{\veca}\nbigv) \cdot [C].
\end{align*}

\begin{Remark}
The higher Chern character
for filtered sheaves was defined by
Iyer and Simpson~\cite{i-s}
in a systematic way. In~this paper, we adopt the definition
of $\ch_2(\nbigp_{\ast}\nbigv)$
in~\cite{mochi4}.
\end{Remark}

\subsection[Good filtered lambda-flat bundles]{Good filtered $\boldsymbol\lambda$-flat bundles}

Let $X$ be a complex manifold with a simple normal crossing
hypersurface $H=\bigcup_{i\in\Lambda}H_i$.

\subsubsection[Good set of irregular values at P]{Good set of irregular values at $\boldsymbol P$}

Let $P$ be any point of $H$. We~take an admissible holomorphic coordinate neighbourhood
$(X_P,z_1,\ldots,z_n)$ around $P$.
Let $f\in \nbigo_X(\ast H)_P$. If $f\in\nbigo_{X,P}$,
we set $\ord(f):=(0,\ldots,0)\in\real^{\ell(P)}$. If there exists
$\vecn\in\seisuu_{\leq 0}^{\ell(P)}\setminus\{(0,\ldots,0)\}$
such that $(i)$~$g:=f\prod z_i^{-n_i}\in \nbigo_{X,P}$,
$(ii)$~$g(P)\neq 0$, then we set $\ord(f):=\vecn$. Otherwise,
$\ord(f)$ is not defined.

For any $\gminia\in\nbigo_X(\ast H)_P/\nbigo_{X,P}$,
we take a lift $\gminiatilde\in\nbigo_X(\ast H)_P$.
If $\ord(\gminiatilde)$ is defined,
we set
$\ord(\gminia):=\ord(\gminiatilde)$.
Otherwise, $\ord(\gminia)$ is not defined.
Note that it is independent of the choice of
a lift $\gminiatilde$.

Let $\nbigi_P\subset \nbigo_{X}(\ast H)_P/\nbigo_{X,P}$
be a finite subset. We~say that $\nbigi_P$ is a good set of irregular values
if the following conditions are satisfied:
\begin{itemize}\itemsep=0pt
\item
 $\ord(\gminia)$ is defined for any $\gminia\in \nbigi_P$.
\item
 $\ord(\gminia-\gminib)$ is defined for any
 $\gminia,\gminib\in\nbigi_P$.
\item
$\{\ord(\gminia-\gminib)\mid \gminia,\gminib\in\nbigi_P\}$
is totally ordered
with respect to
the order $\leq_{\seisuu^{\ell(P)}}$.
Here, we define
$\vecn\leq_{\seisuu^{\ell(P)}}\vecn'$
if $n_i\leq n_i'$ for any $i$.
\end{itemize}

\subsubsection[Good filtered lambda-flat bundles]{Good filtered $\boldsymbol\lambda$-flat bundles}

Let $\nbigv$ be a locally free $\nbigo_X(\ast H)$-module
with a flat $\lambda$-connection.
Let $\nbigp_{\ast}\nbigv$ be a filtered bundle over $\nbigv$.
For any $P\in X$,
let $\nbigo_{X,\Phat}$ denote the completion of
the local ring $\nbigo_{X,P}$
with respect to the maximal ideal.
Note that Remark~\ref{rem;20.2.14.2}
has a natural generalization
to filtered $\lambda$-flat bundles. We~say that $\big(\nbigp_{\ast}\nbigv,\DDlambda\big)$
is unramifiedly good at $P$ if the following holds:
\begin{itemize}\itemsep=0pt
\item
 There exist a good set of irregular values
 $\nbigi_P\subset \nbigo_{X}(\ast H)_P/\nbigo_{X,P}$
 and a decomposition of~filtered $\lambda$-flat bundles
\begin{gather}
\label{eq;20.2.1.1}
 \big(\nbigp_{\ast}\nbigv,\DDlambda\big)\otimes\nbigo_{X,\Phat}
=\bigoplus_{\gminia\in\nbigi_P}
 \big(\nbigp_{\ast}\nbigv_{\gminia},\DDlambda_{\gminia}\big)
\end{gather}
such that
$\DDlambda_{\gminia}-{\rm d}\gminiatilde\id_{\nbigv_{\gminia}}$
are logarithmic with respect to the lattices
$\nbigp_{\veca}\nbigv_{\gminia}$
for any $\veca\in\real^{\ell(P)}$
and $\gminia\in\nbigi_P$,
i.e.,
\begin{gather}
\label{eq;20.2.8.1}
 \big(\DDlambda_{\gminia}-{\rm d}\gminiatilde\id_{\nbigv_{\gminia}}\big)
 \nbigp_{\veca}\nbigv_{\gminia}
\subset
 \nbigp_{\veca}\nbigv_{\gminia}
 \otimes\Omega^1_X(\log H).
\end{gather}
Here,
$\gminiatilde$ denote lifts of $\gminia$
to $\nbigo_X(\ast H)_P$.
\end{itemize}
We say that $\big(\nbigp_{\ast}\nbigv,\DDlambda\big)$
is good at $P$ if the following holds:
\begin{itemize}\itemsep=0pt
\item
There exist a neighbourhood $X_P$ of $P$ in $X$
and a covering map $\varphi_P\colon X_P'\lrarr X_P$
rami\-fied over $H_P=H\cap X_P$
such that
$\varphi_P^{\ast}\big(\nbigp_{\ast}\nbigv,\DDlambda\big)$
is unramifiedly good at $\varphi_P^{-1}(P)$.
(See Section~\ref{subsection;20.2.14.1}
for the pull back of filtered bundles.)
\end{itemize}
We say that
$\big(\nbigp_{\ast}\nbigv,\DDlambda\big)$
is good (resp.~unramifiedly good)
if it is good (resp.~unramifiedly good)
at any point of $H$.

\subsection{Prolongation of holomorphic vector bundles
with a Hermitian metric}

Let $X$ be any complex manifold
with a simple normal crossing hypersurface
$H=\bigcup_{i\in\Lambda} H_i$.
Let~$\big(E,\delbar_E\big)$ be a holomorphic vector bundle
on $X\setminus H$
with a Hermitian metric $h$.
Let us recall
the construction of
$\nbigo_X(\ast H)$-module $\nbigp^hE$
and $\nbigo_X$-modules $\nbigp^h_{\veca}E$ $(\veca\in\real^{\Lambda})$.

Let $\veca\in\real^{\Lambda}$.
For any open subset $\nbigu\subset X$,
let $\nbigp_{\veca}^hE(\nbigu)$
be the space of holomorphic sections~$s$ of
$E_{|\nbigu\setminus H}$
satisfying the following condition:
\begin{itemize}\itemsep=0pt
\item
For any point $P$ of $\nbigu\cap H$,
let $(X_P,z_1,\ldots,z_n)$
be an admissible holomorphic coordinate neighbourhood
around $P$
such that $X_P$ is relatively compact in $\nbigu$.
Set $\vecc=\kappa_P(\veca)$.
(See Section~\ref{subsection;19.2.18.50}.)
Then,
\begin{gather*}
 \bigl|s\bigr|_{h}=
 O\Biggl(\prod_{i=1}^{\ell(P)}|z_i|^{-c_i-\epsilon}\Biggr)
\end{gather*}
holds on $X_P\setminus H$ for any $\epsilon>0$.
\end{itemize}
We obtain an $\nbigo_X$-module
$\nbigp^h_{\veca}E$. We~set
$\nbigp^h E:=\bigcup_{\veca\in\real^{\Lambda}}\nbigp^h_{\veca}E$
which is an $\nbigo_X(\ast H)$-module.
Note that in general,
$\nbigp^h_{\veca}E$
are not necessarily coherent $\nbigo_X$-modules.

\begin{Definition}
Let $\nbigp_{\ast}\nbigv$ be a filtered bundle over $(X,H)$.
Let $\big(E,\delbar_E\big)$ be the holomorphic vector bundle
obtained as the restriction of $\nbigv$ to $X\setminus H$.
A Hermitian metric $h$ is called adapted to
$\nbigp_{\ast}\nbigv$
if $\nbigp^h_{\ast}E=\nbigp_{\ast}\nbigv$
in $\iota_{\ast}(E)=\iota_{\ast}\big(\nbigv_{|X\setminus H}\big)$,
where $\iota\colon X\setminus H\lrarr X$ denotes
the inclusion.
\end{Definition}

\subsubsection{A sufficient condition}

We mention a useful sufficient condition
for $\nbigp^h_{\ast}E$ to be a filtered bundle,
although we do not use it in this paper.
Let $g_{X\setminus H}$ be a K\"ahler metric
satisfying the following condition
\cite{cg}:
\begin{itemize}\itemsep=0pt
\item
For any $P\in H$,
take an admissible holomorphic coordinate neighbourhood
$(X_P,z_1,\ldots,z_n)$ around $P$
such that
$X_P$ is isomorphic to
$\prod_{i=1}^{n}\{|z_i|<1\}$
by the coordinate system.
Set $X_P':=\prod_{i=1}^n\{|z_i|<1/2\}$.
Then, $g_{|X_P'\setminus H}$ is mutually bounded
with the restriction of the Poincar\'e metric
\begin{gather*}
 \sum_{i=1}^{\ell(P)}\frac{{\rm d}z_i\,{\rm d}\zbar_i}{|z_i|^2(\log|z_i|^2)^2}
+\sum_{i=\ell(P)+1}^n{\rm d}z_i\,{\rm d}\zbar_i.
\end{gather*}
\end{itemize}
A Hermitian metric $h$ of $\big(E,\delbar_E\big)$ is called acceptable
if the curvature of the Chern connection
is bounded with respect to $h$ and $g_{X\setminus H}$.
The following theorem is
proved in
\cite[Theorem~21.3.1]{Mochizuki-wild}.
\begin{Theorem}
If $h$ is acceptable,
then $\nbigp^h_{\ast}E$ is a filtered bundle,
and $\nbigp^hE$ is a locally free
$\nbigo_X(\ast H)$-module.
\end{Theorem}

\subsection{Harmonic bundles}

\subsubsection[Pluri-harmonic metrics for lambda-flat bundles]
{Pluri-harmonic metrics for $\boldsymbol\lambda$-flat bundles}

Let $Y$ be any complex manifold.
Let $E$ be a $C^{\infty}$-vector bundle on $Y$.
Let $A^{p,q}(E)$ denote the space of
$C^{\infty}$-sections of
$\Omega^{p,q}\otimes E$. We~set $A^{\ell}(E):=\bigoplus_{p+q=\ell}A^{p,q}(E)$. In~this context,
a~$\lambda$-connection of $E$
is a differential operator
$\DDlambda\colon A^0(E)\lrarr A^1(E)$
such that
$\DDlambda(fs)=
 f\DDlambda(s)
+\big(\lambda \del_Y+\delbar_Y\big)f\otimes s$
for any
$f\in C^{\infty}(Y)$ and $s\in A^0(E)$. We~obtain a section
$\DDlambda\circ\DDlambda \in A^2(\End(E))$.
A $\lambda$-connection is called flat
if $\DDlambda\circ\DDlambda=0$.

Let $\big(E,\DDlambda\big)$ be a $\lambda$-flat bundle on $Y$. We~decompose $\DDlambda={\rm d}_E''+{\rm d}_E'$ into the $(0,1)$-part
and the $(1,0)$-part.
Then, $\big(E,{\rm d}_E''\big)$ is a holomorphic vector bundle.
Let $h$ be a Hermitian metric of $E$.
From $h$ and ${\rm d}''_E$,
we obtain the differential operator $\delta'_{E,h}$
such that
${\rm d}''_E+\delta'_{E,h}$ is a Chern connection.
From ${\rm d}'_E$ and $h$,
we obtain the $(0,1)$-operator $\delta''_{E,h}$
determined by
$\lambda\del h(u,v)=h({\rm d}'_Eu,v)+h(u,\delta''_{E,h}v)$.
As in~\cite[Section~2.2.1]{mochi5},
we obtain the operators
\begin{gather*}
 \delbar_{E,h}:= \frac{1}{1+|\lambda|^2}
 \big({\rm d}''_E+\lambda\delta'_{E,h}\big),\qquad
 \del_{E,h}:= \frac{1}{1+|\lambda|^2}
 \big(\lambdabar {\rm d}'_E+\delta'_{E,h}\big),
\\
 \theta^{\dagger}_{E,h}:= \frac{1}{1+|\lambda|^2}
 \big(\lambdabar {\rm d}''_E-\delta''_{E,h}\big),\qquad
 \theta_{E,h}:= \frac{1}{1+|\lambda|^2}
 \big({\rm d}'_E-\lambda\del'_{E,h}\big).
\end{gather*}
Note that
$\DDlambda=\delbar_{E,h}+\theta_{E,h}+\lambda\big(\del_{E,h}+\theta^{\dagger}_{E,h}\big)$. We~set
$\DD^{\lambda\star}_{E,h}:=\delta'_{E,h}-\delta''_{E,h}
=\del_{E,h}+\theta_{E,h}^{\dagger}\allowbreak-\lambdabar\big(\delbar_{E,h}+\theta_{E,h}\big)$,
and~$G(h):=\bigl[\DDlambda,\DD^{\lambda\star}_{E,h}\bigr]$.

\begin{Definition}
$h$ is called a pluri-harmonic metric of $\big(E,\DDlambda\big)$
if $G(h)=0$. Such a tuple $\big(E,\DDlambda,h\big)$ is called a harmonic bundle.
\end{Definition}

If $\lambda\neq 0$, because
$\big(1+|\lambda|^2\big)\big(\delbar_{E,h}+\theta_{E,h}\big)
=\DDlambda-\lambda\DD_{E,h}^{\lambda\star}$,
and $\big(\DDlambda\big)^2=\big(\DD_{E,h}^{\lambda\star}\big)^2=0$,
we obtain
\begin{gather*}
G(h)=-\frac{(1+|\lambda|^2)^2}{\lambda}\big(\delbar_{E,h}+\theta_{E,h}\big)^2
=-\frac{(1+|\lambda|^2)^2}{\lambda}
 \bigl( \delbar_{E,h}^2+\delbar_{E,h}\theta_{E,h} +\theta_{E,h}^2 \bigr).
\end{gather*}
Hence,
$G(h)=0$ implies that
$\big(E,\delbar_{E,h},\theta_{E,h}\big)$ is a Higgs bundle.
The metric $h$ is a pluri-harmonic metric for
$\big(E,\delbar_{E,h},\theta_{E,h}\big)$.
Conversely,
if $h$ is a pluri-harmonic metric for
a Higgs bundle $\big(E,\delbar_E,\theta\big)$,
we obtain
the Chern connection
$\delbar_E+\del_{E,h}$
associated with $\delbar_E$ and $h$,
and the adjoint $\theta^{\dagger}_h$
of $\theta$ with respect to $h$. We~obtain a flat $\lambda$-connection
$\DDlambda_h=
\delbar_E+\lambda\theta_h
+\lambda\del_{E,h}+\theta$.
The metric $h$ is a pluri-harmonic metric
for $(E,\DDlambda_h)$.

\begin{Remark}
If $\lambda=0$,
a flat $0$-connection is equivalent to
a Higgs bundle
$\big(E,\delbar_E,\theta\big)$
by the relation $\DD^1=\delbar_E+\theta$. In~this case,
we obtain
$G(h)=\bigl[ \delbar_E+\theta,\del_{E,h}+\theta_h^{\dagger} \bigr]$,
and hence
$2G(h)$ is equal to the curvature $F(h)$ of
the connection
$\DD^1_h=\delbar_E+\theta_h^{\dagger}+\del_{E,h}+\theta$.
\end{Remark}

\subsubsection[The case lambda neq 0]{The case $\boldsymbol{\lambda\neq 0}$}

Let $G(h)=G(h)^{2,0}+G(h)^{1,1}+G(h)^{0,2}$
denote the decomposition into $(p,q)$-parts.
If $G(h)=0$,
we clearly obtain $G(h)^{1,1}=0$.
If $\lambda\neq 0$,
we obtain the converse.

\begin{Proposition}
 \label{prop;20.1.29.10}
Suppose $\lambda\neq 0$. If $G(h)^{1,1}=0$,
we obtain $G(h)=0$, i.e.,
$h$ is a pluri-harmonic metric of $\big(E,\DDlambda\big)$.
\end{Proposition}
\begin{proof}
As in~\cite[Lemma 2.28]{mochi5},
the following holds:
\begin{gather}
\label{eq;20.1.31.1}
 \lambdabar^{-1}\del_{E,h}^2+\lambda^{-1}\theta_{E,h}^2=0,\qquad
 \lambda^{-1}\delbar_{E,h}^2+\lambdabar^{-1}\big(\theta_{E,h}^{\dagger}\big)^2=0.
\end{gather}
It~is easy to check that
$\delbar_{E,h}^2=-(\del_{E,h}^2)^{\dagger}$,
$\big(\theta_{E,h}^{\dagger}\big)^2=-(\theta_{E,h}^2)^{\dagger}$
and
$\big(\delbar_{E,h}\theta_{E,h}\big)^{\dagger}=\del_{E,h}\theta^{\dagger}_{E,h}$.

From the flatness $\DDlambda\circ\DDlambda=0$, we obtain
\begin{gather}
\label{eq;20.1.31.2}
 (\lambda\del_{E,h}+\theta_{E,h})^2=
 \lambda^2\del_{E,h}^2+\lambda\del_{E,h}\theta_{E,h}
 +\theta_{E,h}^2=0,
\\
 \bigl[\delbar_{E,h}+\lambda\theta^{\dagger}_{E,h},
 \lambda\del_{E,h}+\theta_{E,h} \bigr]\nonumber
 \\ \qquad
{}=\lambda\bigl(\bigl[\delbar_{E,h},\del_{E,h}\bigl]
+\bigl[\theta_{E,h}\theta^{\dagger}_{E,h}\bigl]\bigr)
+\delbar_{E,h}\theta_{E,h}+\lambda^2\del_{E,h}\theta^{\dagger}_{E,h}=0.
\label{eq;20.1.31.3}
\end{gather}
From (\ref{eq;20.1.31.1}) and (\ref{eq;20.1.31.2}),
we obtain
\begin{gather*}
 \del_{E,h}\theta_{E,h}=-\lambda^{-1}\big(1-|\lambda|^2\big) \theta_{E,h}^2.
\end{gather*}
Because
$(\delbar_{E,h}\theta^{\dagger}_{E,h})^{\dagger}
=\del_{E,h}\theta_{E,h}$,
we obtain
\begin{gather}
\label{eq;20.1.29.2}
 \delbar_{E,h}\theta_{E,h}^{\dagger}
=\lambdabar^{-1}\big(1-|\lambda|^2\big)\big(\theta^{\dagger}_{E,h}\big)^2.
\end{gather}

Note that $G(h)^{1,1}=0$ is equivalent to
$\delbar_{E,h}\theta_{E,h}=0$ and
$\del_{E,h}\theta_{E,h}^{\dagger}=0$.
To obtain $G(h)=0$, it is enough to prove
$\Tr\big(\theta_{E,h}^2\big(\theta_{E,h}^{\dagger}\big)^2\big)=0$.
Indeed, there exists $C\neq 0$ depending on $\dim Y$
such that for any K\"ahler form $\omega$ of $Y$
we obtain
$\Tr\big(\theta_{E,h}^2\big(\theta_{E,h}^{\dagger}\big)^2\big)\omega^{\dim Y-2}
=C\bigl|\theta_{E,h}^2\bigr|^2_{h,\omega}\omega^{\dim Y}$.
Hence, the vanishing
$\Tr\big(\theta_{E,h}^2\big(\theta_{E,h}^{\dagger}\big)^2\big)=0$
implies $\theta_{E,h}^2=\big(\theta_{E,h}^{\dagger}\big)^2=0$
and $\del_{E,h}^2=\delbar_{E,h}^2=0$.

From (\ref{eq;20.1.29.2}) and $\del_{E,h}\theta^{\dagger}_{E,h}=0$,
we obtain $\del_{E,h}\delbar_{E,h}\theta_{E,h}^{\dagger}
=\lambdabar^{-1}\big(1-|\lambda|^2\big) \del_{E,h}\bigl((\theta^{\dagger}_{E,h})^2\bigr) =0$.
 We~also have
 $\delbar_{E,h}\del_{E,h}\theta_{E,h}^{\dagger}=0$.
Hence, we obtain the following equality:
\begin{gather}
\label{eq;20.1.31.4}
 0=\Tr\big(\theta_{E,h}\del_{E,h}\delbar_{E,h}\theta^{\dagger}_{E,h}\big)
=\Tr\bigl( \theta_{E,h}\cdot\big[\delbar_{E,h}\del_{E,h} +\del_{E,h}\delbar_{E,h},\theta_{E,h}^{\dagger}\big] \bigr).
\end{gather}
From (\ref{eq;20.1.31.3}),
$\delbar_{E,h}\theta_{E,h}=0$
and
$\del_{E,h}\theta_{E,h}^{\dagger}=0$,
we obtain
\begin{gather*}
 \big[\delbar_{E,h},\del_{E,h}\big]+\big[\theta_{E,h},\theta_{E,h}^{\dagger}\big]=0.
\end{gather*}
Hence, we obtain the following:
\begin{align}
 \Tr\bigl( \theta_{E,h}\cdot
 \big[\delbar_{E,h}\del_{E,h}
 +\del_{E,h}\delbar_{E,h},\theta_{E,h}^{\dagger}\big] \bigr)
&=\Tr\bigl( \theta_{E,h}\cdot\bigl[ {-}\big[\theta_{E,h},\theta_{E,h}^{\dagger}\big],\theta_{E,h}^{\dagger} \bigr] \bigr)\nonumber
 \\
&=-2\Tr\bigl(\theta_{E,h}^2\big(\theta^{\dagger}_{E,h}\big)^2\bigr).
\label{eq;20.1.31.5}
\end{align}
We obtain the claim of the proposition
from (\ref{eq;20.1.31.4}) and (\ref{eq;20.1.31.5}).
\end{proof}

By using Proposition~\ref{prop;20.1.29.10}
we can improve~\cite[Corollary 2.30]{mochi5}
as follows.
\begin{Corollary}
\label{cor;20.2.3.1}
If $\lambda\neq 0$,
the pluri-harmonicity of the metric $h$ is
equivalent to
the vanishing
$G(h)^{1,1}=0$,
i.e., $\delbar_{E,h}\theta_{E,h}=0$.
\end{Corollary}

\begin{Remark}
In~\cite[Lemma 2.29]{mochi5}, the claim
$\big[\,\delbar_{V,h},\del_{V,h}\big]+\big[\theta_{V,h},\theta_{V,h}^{\dagger}\big]=0$
is incorrect, in~gene\-ral. The author thanks Pengfei Huang for pointing out it.
\end{Remark}

\begin{Remark}
If $\lambda=1$,
the claim of Proposition~\ref{prop;20.1.29.10}
also follows from
a Bochner type formula~\cite[Proposition~21.39]{mochi2-II},
which originally goes back to
the study of Simpson~\cite{s5}
in the context of harmonic bundles,
the study of Corlette~\cite{corlette}
in the context of harmonic metrics
for flat bundles on Riemannian manifolds,
and the study of Siu~\cite{siu-Bochner}
in the context of harmonic maps.
\end{Remark}

\subsection{Wild harmonic bundles}

\subsubsection{Higgs case}\label{subsection;19.2.21.1}

Let $X$ be a complex manifold
with a simple normal crossing hypersurface
$H=\bigcup_{i\in\Lambda}H_i$.
Let $\big(E,\delbar_E,\theta,h\big)$ be a harmonic bundle
on $X\setminus H$.
It~is called wild on $(X,H)$
if the following holds:
\begin{itemize}\itemsep=0pt
\item
 Let $\Sigma_{\theta}\subset T^{\ast}(X\setminus H)$
 denote the spectral cover of~$\theta$,
 i.e., $\Sigma_{\theta}$ denotes the support
 of the coherent $\nbigo_{T^{\ast}(X\setminus H)}$-module
 induced by $\big(E,\delbar_E,\theta\big)$.
 Then, the closure of~$\Sigma_{\theta}$
 in the relatively projective completion of
 $T^{\ast}X$ with respect to $X$ is complex analytic.
\end{itemize}
A wild harmonic bundle $\big(E,\delbar_E,\theta,h\big)$
is called unramifiedly good at $P\in H$ if the following holds:
\begin{itemize}\itemsep=0pt
\item
 There exists a good set of irregular values
 $\nbigi_P\subset\nbigo_{X}(\ast H)_P/\nbigo_{X,P}$,
 a neighbourhood $X_P$,
 and a decomposition
\begin{gather*}
 \big(E,\delbar_E,\theta\big)_{|X_P\setminus H}
=\bigoplus_{\gminia\in\nbigi_P}
 \big(E_{\gminia},\delbar_{E_{\gminia}},\theta_{\gminia}\big)
\end{gather*}
such that the closure of the spectral cover $\Sigma_{\gminia}$
of $\theta_{\gminia}-{\rm d}\gminiatilde\id_{E_{\gminia}}$
in $T^{\ast}X_P(\log (X_P\cap H))$
is proper over $X_P$,
where $\gminiatilde$ denote lifts of $\gminia$
to $\nbigo_X(\ast H)_P$.
\end{itemize}
A wild harmonic bundle
$\big(E,\delbar_E,\theta,h\big)$
is called good at $P\in H$
if the following holds:
\begin{itemize}\itemsep=0pt
\item
 There exist a neighbourhood $X_P$
 and a covering $\varphi_P\colon X_P'\lrarr X_P$
 ramified along $H_P'$
 such that the pull back
 $\varphi_P^{-1}\big(E,\delbar_E,\theta,h\big)_{|X_P}$
 is unramifiedly good wild at any point of
 $\varphi_P^{-1}(H)$.
\end{itemize}
We say that
$\big(E,\delbar_E,\theta,h\big)$ is good wild
(resp.~unramifiedly good wild) on $(X,H)$
if it is good wild (resp.~unramifiedly good wild)
at any point of $H$.

Note that not every wild harmonic bundle on $(X,H)$
is necessarily good on $(X,H)$.
But, the following is known
\cite[Corollary 15.2.8]{Mochizuki-MTM}.
\begin{Theorem}
Let $\big(E,\delbar_E,\theta,h\big)$ be a wild harmonic bundle
on $(X,H)$.
Then, there exists a~proper birational morphism
$\varphi\colon X'\lrarr X$ of complex manifolds such that
$(i)$~$H':=\varphi^{-1}(H)$ is simple normal crossing,
$(ii)$~$X'\setminus H'\simeq X\setminus H$,
$(iii)$~$\varphi^{-1}\big(E,\delbar_E,\theta,h\big)$
 is good wild on $(X',H')$.
\end{Theorem}

\subsubsection[The case of lambda-flat bundles]{The case of $\boldsymbol\lambda$-flat bundles}
A $\lambda$-flat bundle $\big(E,\DDlambda\big)$
with a pluri-harmonic metric $h$ on $X\setminus H$
is called (good, unramifiedly good) wild
if the associated Higgs bundle with a pluri-harmonic metric
 $\big(E,\delbar_{E,h},\theta_{E,h},h\big)$ is a~(good, unramifiedly good) wild harmonic bundle.

\subsubsection[Prolongation of good wild harmonic bundles to good filtered lambda-flat bundles]
{Prolongation of good wild harmonic bundles to good filtered $\boldsymbol\lambda$-flat bundles}

The following is one of the fundamental theorems
in the study of wild harmonic bundles
\cite[Theorem~7.4.3]{Mochizuki-wild}.
\begin{Theorem}
If $\big(E,\DDlambda,h\big)$ is a good wild harmonic bundle
on $(X,H)$,
then $\big(\nbigp^h_{\ast}E,\DDlambda\big)$
is a~good filtered $\lambda$-flat bundle on $(X,H)$.
\end{Theorem}

The following is a consequence of
the norm estimate for good wild harmonic bundles
\cite[Theorem~11.7.2]{Mochizuki-wild}.
\begin{Theorem}
Let $\big(E,\DDlambda,h_i\big)$ $(i=1,2)$ be
good wild harmonic bundles on $X\setminus H$
such that
$\nbigp^{h_1}_{\ast}E=\nbigp^{h_2}_{\ast}E$.
 Then, $h_i$ are mutually bounded
 around any point of $H$.
\end{Theorem}

\subsubsection{Prolongation of good wild harmonic bundles in the projective case}

Suppose that $X$ is projective and connected.
Let $L$ be any ample line bundle on $X$.
The~fol\-lo\-wing is proved in~\cite[Propositions 13.6.1 and~13.6.4]{Mochizuki-wild}.
\begin{Proposition}\label{prop;19.2.12.200}
Let $\big(E,\DDlambda,h\big)$ be a good wild harmonic bundle
on $(X,H)$.
\begin{itemize}\itemsep=0pt
\item
 $\big(\nbigp^h_{\ast}E,\DDlambda\big)$ is
 $\mu_L$-polystable
 with $\mu_L\big(\nbigp^h_{\ast}E\big)=0$.
\item
We obtain
$c_1(\nbigp_{\ast}E)=0$
and
$\int_{X}\ch_2(\nbigp_{\ast}E)c_1(L)^{\dim X-2}=0$.
\item
 Let $h'$ be another pluri-harmonic metric of
 $\big(E,\DDlambda,h\big)$
 such that $\nbigp^{h'}_{\ast}E=\nbigp^h_{\ast}E$.
 Then, there exists a decomposition of the $\lambda$-flat bundle
 $\big(E,\DDlambda\big)=
 \bigoplus
 \big(E_j,\DDlambda_j\big)$
 such that
 $(i)$ the decomposition is orthogonal
 with respect to both $h$ and $h'$,
 $(ii)$ $h_{|E_i}=a_i\cdot h'_{|E_i}$ for some $a_i>0$.
\item
Let $\big(\nbigp_{\ast}\nbigv_1,\DDlambda_1\big)$
be any direct summand
of $\big(\nbigp^h_{\ast}E,\DDlambda\big)$.
Let $\big(E_1,\DDlambda_1\big)$ be
the $\lambda$-flat bundle on $X\setminus H$
obtained as the restriction of $\big(\nbigv_1,\DDlambda_1\big)$,
and let $h_1$ be the metric of $E_1$ induced by~$h$.
Then, $\big(E_1,\DDlambda_1,h_1\big)$ is a harmonic bundle. In~particular, we obtain
$c_1(\nbigp_{\ast}\nbigv_1)=0$
and
$\int_{X}\ch_2(\nbigp_{\ast}\nbigv_1)c_1(L)^{\dim X-2}=0$.
\end{itemize}
\end{Proposition}

\subsection{Main existence theorem in this paper}

Let $X$ be a smooth connected projective complex manifold
with a simple normal crossing hypersurface $H$.
Let $L$ be any ample line bundle on $X$.
Let $\big(\nbigp_{\ast}\nbigv,\DDlambda\big)$ be a good filtered
$\lambda$-flat bundle on $(X,H)$.
Let $\big(E,\delbar_E,\DDlambda\big)$ be the $\lambda$-flat bundle
obtained as the restriction of $\big(\nbigp_{\ast}\nbigv,\DDlambda\big)$
to $X\setminus H$.

\begin{Theorem}
\label{thm;19.2.11.20}
Suppose that
$\big(\nbigp_{\ast}\nbigv,\DDlambda\big)$ is $\mu_L$-polystable,
and that the following vanishing holds:
\begin{gather}
\label{eq;19.2.18.30}
 \mu_L(\nbigp_{\ast}\nbigv)=0,\qquad
 \int_X\ch_{2}(\nbigp_{\ast}\nbigv)c_1(L)^{\dim X-2}=0.
\end{gather}
Then, there exists a pluri-harmonic metric $h$
of $\big(E,\delbar_E,\DDlambda\big)$ such that
$\big(\nbigv,\DDlambda\big)_{|X\setminus H}\simeq \big(E,\DDlambda\big)$
extends to
$\big(\nbigp_{\ast}\nbigv,\DDlambda\big)\simeq \big(\nbigp^h_{\ast}E,\DDlambda\big)$.
\end{Theorem}

We proved the claim of the theorem in the case $\lambda=1$
in~\cite[Theorem~16.1.1]{Mochizuki-wild}. We~shall explain the proof in
Sections~\ref{section;19.2.11.10}--\ref{section;19.2.11.11}.
Note that the one dimensional case is due to
Biquard--Boalch~\cite{biquard-boalch}.

\begin{Corollary}\label{cor;19.2.18.3}
There exists the equivalence of the following objects
for each $\lambda$:
\begin{itemize}\itemsep=0pt
\item
 Good wild harmonic bundles on $(X,H)$.
\item
 $\mu_L$-polystable good filtered $\lambda$-flat bundles
 $\big(\nbigp_{\ast}\nbigv,\DDlambda\big)$
 satisfying the condition~$(\ref{eq;19.2.18.30})$.
\end{itemize}
\end{Corollary}

\begin{Remark}
One of the referees raised the following interesting question.
Let $L$ be a big and nef line bundle on $X$
such that there exists a positive current
$T$ representing $c_1(L)$
whose restriction to $X\setminus H$
is a smooth K\"ahler form
with at most Poincar\'{e} growth near $H$. We~do not assume that $L$ is ample. We~can define the slope $\mu_L(\nbigp_{\ast}\nbigv)$
for a filtered sheaves,
by using which we can introduce a stability condition
for good filtered $\lambda$-flat bundles.
Then, we may ask whether a statement
similar to Theorem~\ref{thm;19.2.11.20} holds.
This question might also be related with a~generalization
of Kobayashi--Hitchin correspondence
to the context of $\nbigd$-modules.
\end{Remark}

\subsubsection{Outline of the proof
Theorem~\ref{thm;19.2.11.20}}

Let us explain a rough outline of the proof. We~shall omit some technical details.
Let $\big(\nbigp_{\ast}\nbigv,\DDlambda\big)$
be a $\mu_L$-stable good filtered $\lambda$-flat bundle on $(X,H)$
such that
$\ch_i(\nbigp_{\ast}\nbigv)=0$ $(i=1,2)$.
Let~$\big(V,\DDlambda\big)$ be the $\lambda$-flat bundle
on $X\setminus H$
obtained as the restriction of $\big(\nbigp_{\ast}\nbigv,\DDlambda\big)$
to $X\setminus H$.

In the case $\dim X=1$,
we shall apply the argument in~\cite{s2}
as follows.
For each $P\in H$,
we take a holomorphic coordinate neighbourhood
$\big(X_P,z_P\big)$ around $P$. We~take a K\"ahler metric $g_{X\setminus H}$
such that $g_{X\setminus H|X_P\setminus\{P\}}$ $(P\in H)$
are mutually bounded with
$|z_P|^{2\eta-2}{\rm d}z_P\,{\rm d}\zbar_P$
for some $\eta>0$.
If~$\eta$ is sufficiently small,
there exists a Hermitian metric $h_{\inn}$ of $V$
such that
$(i)$~$\nbigp^{h_{\inn}}_{\ast}V=\nbigp_{\ast}\nbigv$,
and
$(ii)$~$G(h_{\inn})$ is bounded with respect to
$h_{\inn}$ and $g_{X\setminus H}$,
$(iii)$~$\det(h_{\inn})$ is flat.
(See Corollary~\ref{cor;19.1.20.2}.
Though we state it as a corollary of Proposition~\ref{prop;20.2.6.10},
which also deals with a perturbation,
it is easy to deduce it directly from
the estimate in the tame case~\cite{s2}.)
Moreover for any filtered $\lambda$-flat subsheaf
$\nbigp_{\ast}\nbigv'\subset\nbigp_{\ast}\nbigv$,
$\deg(\nbigp_{\ast}\nbigv')$ is equal
to the analytic degree of $\big(\nbigv',\DDlambda\big)_{|X\setminus H}$
with respect to $h_{\inn}$ and $g_{X\setminus H}$.
Then, by~\cite[Theorem~1]{s1},
if $\big(\nbigp_{\ast}\nbigv,\DDlambda\big)$ is stable
of degree $0$,
there exists a harmonic metric $h$ of
$\big(\nbigp_{\ast}\nbigv,\DDlambda\big)$
such that
$h$ and $h_{\inn}$ are mutually bounded
(Theorem~\ref{thm;17.11.1.31}).
Let us note that the proof allows us to
obtain the inequality for the Donaldson functional
$M(h_{\inn},h)\leq 0$ (Proposition~\ref{prop;19.1.28.10}).
This inequality is useful for the study of
the continuity of the family of harmonic metrics
of some family of good filtered $\lambda$-flat bundles
(Proposition~\ref{prop;19.1.29.10}).

For the higher dimensional case,
we use the same strategy in~\cite{mochi4,mochi5} and~\cite{Mochizuki-wild}.
It~is a key step to study the case $\dim X=2$.
There are two naive ideas which are not available as they are.

One is to apply~\cite[Theorem~1]{s1}
by constructing a Hermitian metric $h_{\inn}$ of $V$
such that
$(i)$~$\nbigp_{\ast}^{h_{\inn}}V=\nbigp_{\ast}\nbigv$,
$(ii)$~$G(h_{\inn})$ is dominated in an appropriate way,
$(iii)$~$\det(h_{\inn})$ is flat.
For the construction of such a Hermitian metric $h_{\inn}$,
a compatibility condition seems necessary
between the nilpotent parts of the induced endomorphisms
$\Res_i(\DDlambda)$ and $\Res_j(\DDlambda)$
on $\lefttop{i,j}\Gr^F(\nbigp_{\veca}\nbigv)$.
(See Section~\ref{subsection;21.5.5.10} for
the endomorphisms $\Res_i\big(\DDlambda\big)$.)
Once we prove the existence of a pluri-harmonic,
it turns out that such a compatibility condition is satisfied.
However, before proving the existence,
it is not clear whether such a compatibility condition
is satisfied. As a result, it~is~difficult to construct
a Hermitian metric $h_{\inn}$ with the desired property,
in general.

The other is to use Mehta--Ramanathan type theorem
(Proposition~\ref{prop;19.1.30.100}),
according to which
there exists $m>0$ such that
for the $0$-set $Y\subset X$ of a generic section of
$H^0\big(X,L^{\otimes m}\big)$,
the restriction
$\big(\nbigp_{\ast}\nbigv,\DDlambda\big)_{|Y}$ is also stable.
Hence,
if we fix a flat metric $h_{\det(\nbigv)}$
of $\det(\nbigv)_{|X\setminus H}$
adapted to $\det(\nbigp_{\ast}\nbigv)$,
there exists a harmonic metric
$h_Y$ of
$\big(V,\DDlambda\big)_{|Y\setminus H}$
adapted to $\nbigp_{\ast}\nbigv_{|Y}$
such that
$\det(h_Y)=h_{\det(\nbigv)|Y\setminus H}$.
If we can prove that there exists a Hermitian metric $h$ of $V$
such that~$h_{|Y}=h_Y$ for such generic hypersurfaces $Y$,
then $h$ should be the desired pluri-harmonic metric
for $\big(\nbigp_{\ast}\nbigv,\DDlambda\big)$.
But, the existence of such $h$ is not clear.

Roughly speaking, we combine these two ideas
as follows.
For any small $\epsilon>0$,
there exists a~filtered bundle
$\nbigp^{(\epsilon)}_{\ast}\nbigv$ over $\nbigv$
such that
$(i)$~$\big(\nbigp^{(\epsilon)}_{\ast}\nbigv,\DDlambda\big)$
is a $\mu_L$-stable good filtered $\lambda$-flat bundle,
$(ii)$~$\Res_i\big(\DDlambda\big)$ are semisimple
for $\nbigp^{(\epsilon)}_{\ast}\nbigv$,
$(iii)$~$\det\big(\nbigp^{(\epsilon)}_{\ast}\nbigv\big)=\det(\nbigp_{\ast}\nbigv)$,
$(iv)$~the difference of
$\nbigp^{(\epsilon)}_{\ast}\nbigv$
and $\nbigp_{\ast}\nbigv$ are dominated by $\epsilon$.
(See Section~\ref{subsection;19.1.30.101}
for more precise conditions.)
The~last condition implies that
$\lim\limits_{\epsilon\to 0}\int\ch_2\big(\nbigp^{(\epsilon)}_{\ast}\big)=0$.
For $\big(\nbigp_{\ast}^{(\epsilon)}\nbigv,\DDlambda\big)$,
we can construct $h^{(\epsilon)}_{\inn}$
such that
$(i)$~$\nbigp_{\ast}^{h^{(\epsilon)}_{\inn}}V=\nbigp_{\ast}\nbigv$,
$(ii)$~$G\big(h^{(\epsilon)}_{\inn}\big)$ is dominated in an appropriate way,
$(iii)$~$\det\big(h^{(\epsilon)}_{\inn}\big)=h_{\det(\nbigv)}$. By~\cite[Theorem~1]{s1},
there exists a Hermitian--Einstein metric $h^{(\epsilon)}_{\HE}$
of $\big(V,\DDlambda\big)$ such that
$h^{(\epsilon)}_{\HE}$ and $h^{(\epsilon)}_{\inn}$
are mutually bounded,
and that $\det\big(h^{(\epsilon)}_{\HE}\big)=h_{\det(\nbigv)}$.
(See Section~\ref{subsection;21.5.5.11} for Hermitian--Einstein metrics
of Higgs bundles.)
Moreover, $G\big(h^{(\epsilon)}_{\HE}\big)\to 0$ in $L^2$ as $\epsilon\to 0$.
Hence, we would like to construct
the desired pluri-harmonic metric
as $\lim\limits_{\epsilon \to 0}h^{(\epsilon)}_{\HE}$.
If $\epsilon$ is sufficiently small,
$\big(\nbigp^{(\epsilon)}_{\ast}\nbigv,\DDlambda\big)_{|Y}$
is also stable for the $0$-set $Y$
of a generic section of $H^0\big(X,L^{\otimes m}\big)$,
and hence
$\big(\nbigp^{(\epsilon)}_{\ast}\nbigv,\DDlambda\big)_{|Y}$
has a harmonic metric $h^{(\epsilon)}_Y$
such that $\det\big(h^{(\epsilon)}_Y\big)=h_{\det(\nbigv)|Y\setminus H}$. By~the continuity of a family of harmonic metrics mentioned above,
the sequence $h^{(\epsilon)}_Y$ is convergent to $h_Y$
as $\epsilon\to 0$
(Proposition~\ref{prop;19.1.29.10}).
Because $h^{(\epsilon)}_{\HE|Y}$ is not necessarily
a harmonic metric of
$\big(\nbigp^{(\epsilon)}_{\ast}\nbigv,\DDlambda\big)_{|Y}$,
it is not necessarily equal to $h^{(\epsilon)}_{Y}$.
But, because the $L^2$-norm of $G\big(h^{(\epsilon)}_{\HE|Y}\big)$
is dominated by $\epsilon$,
we~can deduce the convergence of
the sequence $h^{(\epsilon)}_{\HE|Y}$ to $h_Y$ as $\epsilon\to 0$
(Proposition~\ref{prop;13.1.7.20}).
Hence, we~obtain the convergence of $h^{(\epsilon)}_{\HE}$
almost everywhere, and the limit satisfies
$h_{|Y}=h_Y$ for the $0$-set of generic section $s$
of $H^0\big(X,L^{\otimes m}\big)$.
Thus, we can prove the theorem in the case $\dim X=2$.
(See Section~\ref{subsection;19.2.1.20} for a~more precise argument.)

In the case $\dim X\geq 3$,
we use an induction on $\dim X$. By~the Mehta--Ramanathan type theorem,
there exists $m>0$ such that
for the $0$-sets $Y_i$ $(i=1,2)$ of generic sections of
$H^0\big(X,L^{\otimes m}\big)$,
$\big(\nbigp_{\ast}\nbigv,\DDlambda\big)_{|Y_i}$
and $\big(\nbigp_{\ast}\nbigv,\DDlambda\big)_{|Y_1\cap Y_2}$
are $\mu_L$-stable. By~fixing a flat metric
$h_{\det(\nbigv)}$ for $\det\big(\nbigp_{\ast}\nbigv,\DDlambda\big)$,
there exist pluri-harmonic metric $h_{Y_i}$ of
$\big(\nbigp_{\ast}\nbigv,\DDlambda\big)_{|Y_i}$
such that $\det(h_{Y_i})=h_{\det(\nbigv)|Y_i\setminus H}$.
Because $\big(\nbigp_{\ast}\nbigv,\DDlambda\big)_{|Y_1\cap Y_2}$
is also $\mu_L$-stable,
we obtain that
$h_{Y_1|(Y_1\cap Y_2)\setminus H}
=h_{Y_2|(Y_1\cap Y_2)\setminus H}$.
Hence,
there exists a Hermitian metric $h$ of
$V_{|X\setminus (H\cup W)}$ for a finite subset $W$,
such that
$h_{|Y\setminus H}=h_{Y}$
for the $0$-set $Y$ of a generic section of
$H^0\big(X,L^{\otimes m}\big)$.
It~is easy to see that $h$ is the desired pluri-harmonic metric.
(See Section~\ref{subsection;19.2.1.21} for a more precise argument.)

\section{Preliminaries}
\label{section;19.2.11.10}

\subsection[Hermitian--Einstein metrics of lambda-flat bundles]
{Hermitian--Einstein metrics of $\boldsymbol\lambda$-flat bundles}\label{subsection;21.5.5.11}

Let $Y$ be a K\"ahler manifold with a K\"ahler form $\omega$.
Let $\big(E,\DDlambda\big)$ be a $\lambda$-flat bundle on $Y$
with a Hermitian metric.
Recall that $h$ is called a Hermitian--Einstein metric
of the $\lambda$-flat bundle
if~$\Lambda_{\omega} G(h)^{\bot}=0$,
where $G(h)^{\bot}$ denote the trace-free part of
$G(h)$,
and $\Lambda_{\omega}$ denote the adjoint
of the multiplication by $\omega$ (see~\cite[Section~3.2]{koba}).
The following is a generalization of
Kobayashi--L\"ubke inequality
to the context of $\lambda$-flat bundles
due to Simpson~\cite[Proposition~3.4]{s1}.
\begin{Proposition}[Simpson]
\label{prop;19.1.30.110}
If $h$ is a Hermitian--Einstein metric,
there exists $C>0$ depending only on $n=\dim Y$
such that the following holds:
\begin{gather*}
 \Tr\bigl(\big(G(h)^{\bot}\big)^2 \bigr)\omega^{n-2}
=C\bigl|G(h)^{\bot}\bigr|_{h,\omega}^2\omega^{n}.
\end{gather*}
As a result,
if $\Tr\bigl(\big(G(h)^{\bot}\big)^2\bigr)\omega^{n-2}=0$,
then we obtain $G(h)^{\bot}=0$.
\end{Proposition}

\subsection{Rank one case}\label{subsection;20.2.10.1}

Let $X$ be an $n$ dimensional smooth connected
projective variety
with a simple normal crossing hypersurface $H$.
Let $\omega$ be a K\"ahler form.
Let $H=\bigcup_{i\in\Lambda} H_i$ be
the irreducible decomposition.
Let $g_i$ be a $C^{\infty}$-Hermitian metric
of the line bundle $\nbigo(H_i)$.
Let $\sigma_i$ denote the section of $\nbigo_X(H_i)$
induced by
the inclusion $\nbigo_X\lrarr\nbigo_X(H_i)$.

Let $\big(\nbigp_{\ast}\nbigv,\DDlambda\big)$ be a good filtered
$\lambda$-flat bundle on $(X,H)$ of rank one.
For each $i\in\Lambda$, there uniquely exists
$a_i\in\openclosed{{-}1}{0}$ such that
$\lefttop{i}\Gr^F_{a_i}(\nbigp_{\ast}\nbigv)\neq 0$.
Let $A$ be the constant determined by
\begin{gather*}
 A\int_X\omega^n=2\pi n\big(1+|\lambda|^2\big)\int_Xc_1(\nbigp_{\ast}\nbigv)
 \omega^{n-1}.
\end{gather*}
The following proposition is standard.

\begin{Proposition}\label{prop;19.1.30.31}
There exists a Hermitian metric $h$
of the line bundle $E:=\nbigv_{|X\setminus D}$
such that
$(i)$~$\sqrt{-1}\Lambda_{\omega} G(h)=A$,
$(ii)$~$h\prod_{i\in\Lambda}|\sigma_i|_{g_i}^{2a_i}$
is a Hermitian metric of $\nbigp_{\veca}(\nbigv)$
of $C^{\infty}$-class.
Such a~metric is unique up to
the multiplication by a positive constant.
Moreover,
if $c_1(\nbigp_{\ast}\nbigv)=0$,
then $R(h)=0$ holds,
and hence
$h$ is a pluri-harmonic metric
of $\big(E,\DDlambda\big)$.
\end{Proposition}
\begin{proof}
Note that $G(h)=\big(1+|\lambda|^2\big)R(h)$ holds in the rank one case.
(See~\cite[Lemma~2.31]{mochi5}.)
Let $h_0'$ be a $C^{\infty}$-metric of
$\nbigp_{\veca}E$. We~obtain the metric
$h_0:=h_0'\cdot \prod_{i\in\Lambda}|\sigma_i|^{-2a_i}_{g_{i}}$
of $E$ on $X\setminus H$.
It~is well known that
$\frac{\sqrt{-1}}{2\pi}R(h_0)$ naturally extends to
a closed $(1,1)$-form on $X$ of $C^{\infty}$-class
which represents $c_1(\nbigp_{\ast}E)$. By~the condition of $A$,
we obtain
$\int_X
 \bigl(
 \sqrt{-1}\Lambda_{\omega} R(h_0)-
 \big(1+|\lambda|^2\big)^{-1}A
 \bigr)\omega^{n}=0$.
Note that
$\sqrt{-1}\Lambda_{\omega} R(h_0{\rm e}^{\varphi})
=\sqrt{-1}\Lambda_{\omega} R(h_0)
+\sqrt{-1}\Lambda_{\omega}\delbar\del\varphi$.
Hence,
there exists
an $\real$-valued $C^{\infty}$-function $\varphi_0$
such that
$\sqrt{-1}\Lambda_{\omega} R(h_0{\rm e}^{\varphi_0})-\big(1+|\lambda|^2\big)^{-1}A=0$.
The metric $h=h_0{\rm e}^{\varphi_0}$ has the desired property.
The uniqueness is clear.

Suppose that $c_1(\nbigp_{\ast}E)=0$. In~the rank one case,
a Hermitian metric of $E$ is
a pluri-harmonic metric of $\big(E,\DDlambda\big)$,
if and only if $R(h)=0$.
Because the cohomology class of $R(h_0)$ is $0$,
there exists
an $\real$-valued $C^{\infty}$-function $\varphi_0$
such that
$R(h_0{\rm e}^{\varphi_0})=0$
by the standard $\del\delbar$-lemma. By~the uniqueness, we obtain the second claim of the lemma.
\end{proof}

For the metric $h$ in Proposition~\ref{prop;19.1.30.31},
$\frac{\sqrt{-1}}{2\pi}R(h)$ induces
a closed $(1,1)$-form on $X$ of $C^{\infty}$-class
which represents $c_1(\nbigp_{\ast}E)$.

\subsection[beta-subobject and socle for reflexive filtered lambda-flat sheaves]
{$\boldsymbol\beta$-subobject and socle for reflexive filtered $\boldsymbol\lambda$-flat sheaves}

Let $X$ and $H$ be as in Section~\ref{subsection;20.2.10.1}.
Let $L$ be an ample line bundle $L$ on $X$.
For any coherent $\nbigo_X$-module $\nbigm$,
we set $\deg_L(\nbigm):=\int_Xc_1(\nbigm)c_1(L)^{\dim X-1}$.

\subsubsection[beta-subobjects]{$\boldsymbol\beta$-subobjects}

Let $\big(\nbigp_{\ast}\nbigv,\DDlambda\big)$
be a reflexive filtered $\lambda$-flat sheaf on $(X,H)$.
For any $A\in\real$,
let $\nbigs(\nbigp_{\veczero}\nbigv,A)$
denote the family of saturated coherent subsheaves $\nbigf$
of $\nbigp_{\veczero}\nbigv$
such that $\deg_L(\nbigf)\geq -A$
and that
$\nbigf(\ast H)$ is a $\lambda$-flat subsheaf of $\nbigv$.
Any $\nbigf\in\nbigs(\nbigp_{\veczero}\nbigv,A)$
induces a reflexive filtered sheaf
$\nbigp_{\ast}(\nbigf(\ast H))$
by $\nbigp_{\vecc}(\nbigf(\ast H)):=
 \nbigp_{\vecc}\nbigv
 \cap
 \nbigf(\ast H)$
for any $\vecc\in\real^{\Lambda}$. We~set $f_A(\nbigf):=\mu_L(\nbigp_{\ast}(\nbigf(\ast H)))$.
Thus, we obtain a function
$f_A$ on $\nbigs(\nbigp_{\veczero}\nbigv,A)$.
\begin{Lemma}
The image
$f_A\bigl(
 \nbigs(\nbigp_{\veczero}\nbigv,A)
 \bigr)$
is a finite subset of $\real$. In~particular,
$f_A$ has the maximum.
\end{Lemma}
\begin{proof}
According to~\cite[Lemma 2.5]{grothendieck},
the family $\nbigs(\nbigp_{\veczero}\nbigv,A)$
is bounded.
Hence, by using the flattening stratifications~\cite[Section~8]{Mumford-curves-on-surface},
it is easy to see that
there exists a finite decomposition
$\nbigs(\nbigp_{\veczero}\nbigv,A)
=\coprod_{i=1}^N \nbigs_i(\nbigp_{\veczero}\nbigv,A)$
such that
$f_A$ is constant on each $\nbigs_i(\nbigp_{\veczero}\nbigv,A)$.
\end{proof}

It~is standard that any reflexive filtered $\lambda$-flat sheaf
has a $\beta$-subobject,
i.e., the following holds.
\begin{Proposition}
\label{prop;19.2.18.110}
For any reflexive filtered $\lambda$-flat sheaf
$\big(\nbigp_{\ast}\nbigv,\DDlambda\big)$,
there uniquely exists a non-zero $\lambda$-flat
subsheaf $\nbigv_0\subset\nbigv$
such that the following holds
for any non-zero reflexive $\lambda$-flat sub\-sheaf~$\nbigv'\subset\nbigv$:
\begin{itemize}\itemsep=0pt
\item
$\mu_L(\nbigp_{\ast}\nbigv')\leq\mu_L(\nbigp_{\ast}\nbigv_0)$
holds.
\item
If $\mu_L(\nbigp_{\ast}\nbigv')=\mu_L(\nbigp_{\ast}\nbigv_0)$
holds,
then we obtain $\nbigv'\subset\nbigv_0$.
\end{itemize}
\end{Proposition}
\begin{proof}
By the formula (\ref{eq;20.7.8.1}),
there exists $N>0$ such that
the following holds for any
saturated subsheaf $\nbigf\subset\nbigp_{\veczero}\nbigv$:
\begin{gather*}
\left| \deg_L(\nbigf)-\rank(\nbigf)\mu_L(\nbigp_{\ast}\nbigf(\ast H))\right|<N.
\end{gather*}
We set $A_0:=|\deg_L(\nbigp_{\veczero}\nbigv)|+10N$.
Let $B_0$ denote the maximum of $f_{A_0}$.
Then, it is easy to see that
$\mu_L(\nbigp_{\ast}\nbigv')\leq B_0$
for any saturated $\lambda$-flat subsheaf $\nbigv'\subset\nbigv$.
Moreover,
if $\mu_L(\nbigp_{\ast}\nbigv')=B_0$,
then $\big(\nbigp_{\ast}\nbigv',\DDlambda_{\nbigv'}\big)$
is $\mu_L$-semistable,
where $\DDlambda_{\nbigv'}$ denote
the flat $\lambda$-connection
induced by $\DDlambda$.

Suppose that the $\lambda$-flat subsheaves
$\nbigv_i\subset\nbigv$ $(i=1,2)$ satisfy
$\mu_L(\nbigp_{\ast}\nbigv_i)=B_0$. We~obtain the subsheaf
$\nbigv_1+\nbigv_2\subset\nbigv$.
Because $\nbigv_1+\nbigv_2$ is a quotient of
$\nbigv_1\oplus\nbigv_2$,
we obtain a filtered sheaf $\nbigp_{\ast}(\nbigv_1+\nbigv_2)$
over $\nbigv_1+\nbigv_2$.
induced by $\nbigp_{\ast}\nbigv_1\oplus\nbigp_{\ast}\nbigv_2$.
Then, by the $\mu_L$-semistability of
$\big(\nbigp_{\ast}\nbigv_i,\DDlambda_i\big)$,
we obtain that
$B_0=\mu_L(\nbigp_{\ast}\nbigv_1\oplus
 \nbigp_{\ast}\nbigv_2)
\leq
 \mu_L(\nbigp_{\ast}(\nbigv_1+\nbigv_2))$.
Let $\nbigv_3$ denote the saturated subsheaf of $\nbigv$
generated by $\nbigv_1+\nbigv_2$. We~obtain a filtered sheaf $\nbigp_{\ast}\nbigv_3$
by $\nbigp_{\veca}\nbigv_3=\nbigp_{\veca}(\nbigv)\cap\nbigv_3$.
Because the natural morphism
$\nbigp_{\ast}(\nbigv_1+\nbigv_2)
\lrarr
 \nbigp_{\ast}\nbigv_3$
is generically an isomorphism,
we obtain
$\mu_L(\nbigp_{\ast}(\nbigv_1+\nbigv_2))
\leq
 \mu_L(\nbigp_{\ast}\nbigv_3)\leq B_0$
by Lemma~\ref{lem;19.2.18.120}.
Hence, we obtain $\mu_L(\nbigp_{\ast}\nbigv_3)=B_0$.
Then, the claim of the lemma is clear.
\end{proof}

\subsubsection{Socle}

Let $\big(\nbigp_{\ast}\nbigv,\DDlambda\big)$
be a $\mu_L$-semistable reflexive filtered $\lambda$-flat sheaf
on $(X,H)$.
Let $\nbigt$
denote the fa\-mily of saturated $\lambda$-flat
subsheaves $\nbigv'\subset\nbigv$
such that
the induced filtered $\lambda$-flat sheaf
$\big(\nbigp_{\ast}\nbigv',\DDlambda_{\nbigv'}\big)$ is $\mu_L$-stable
with $\mu_L(\nbigp_{\ast}\nbigv')=\mu_L(\nbigp_{\ast}\nbigv)$.
Let $\nbigv_1$ be the saturated $\nbigo_X(\ast H)$-submodule of
$\nbigv$
generated by $\sum_{\nbigv'\in\nbigt}\nbigv'$.
It~is a $\lambda$-flat subsheaf of $\nbigv$.

\begin{Proposition}
\label{prop;19.2.18.140}
$\big(\nbigp_{\ast}\nbigv_1,\DDlambda_{\nbigv_1}\big)$
is equal to the direct sum
$\bigoplus_{k=1}^{\ell}\big(\nbigp_{\ast}\nbigv^{(k)},\DDlambda_{\nbigv^{(k)}}\big)$
of a tuple of $\mu_L$-stable filtered $\lambda$-flat subsheaves
of $\big(\nbigp_{\ast}\nbigv,\DDlambda\big)$. In~particular,
$\big(\nbigp_{\ast}\nbigv_1,\DDlambda_1\big)$
is $\mu_L$-polystable.
The filtered $\lambda$-flat subsheaf
$\big(\nbigp_{\ast}\nbigv_1,\DDlambda_{\nbigv_1}\big)$
is called the socle of $\big(\nbigp_{\ast}\nbigv,\DDlambda\big)$.
\end{Proposition}
\begin{proof}
Let $\nbigv^{(i)}$ $(i=1,2)$ be saturated $\lambda$-flat
subsheaves of $\nbigv$ such that
$(i)$~$\mu_L\big(\nbigp_{\ast}\nbigv^{(i)}\big)=\mu_L(\nbigp_{\ast}\nbigv)$,
$(ii)$~$\big(\nbigp_{\ast}\nbigv^{(1)},\DDlambda_{\nbigv^{(1)}}\big)$ is $\mu_L$-semistable,
$(iii)$~$\big(\nbigp_{\ast}\nbigv^{(2)},\DDlambda_{\nbigv^{(2)}}\big)$ is $\mu_L$-stable.

\begin{Lemma}
Either $\nbigv^{(2)}\subset\nbigv^{(1)}$
or $\nbigv^{(1)}\cap\nbigv^{(2)}=0$ holds.
\end{Lemma}
\begin{proof}
Let us consider the morphism
$\iota_1-\iota_2\colon
 \nbigv^{(1)}\oplus\nbigv^{(2)}\lrarr \nbigv$,
where $\iota_i\colon \nbigv^{(i)}\lrarr\nbigv$
denote the inclusions.
Let $\nbigk$ denote the kernel. We~obtain a filtered sheaf $\nbigp_{\ast}\nbigk$
over $\nbigk$
by $\nbigp_{\veca}\nbigk:=
 \nbigk\cap\nbigp_{\veca}(\nbigv_1\oplus\nbigv_2)$.
The projection
$\nbigv^{(1)}\oplus\nbigv^{(2)}\lrarr\nbigv^{(2)}$
induces
$\nbigk\simeq\nbigv^{(1)}\cap\nbigv^{(2)}=:\nbigi$.
It~induces a morphism of filtered $\lambda$-flat sheaves
$g\colon \big(\nbigp_{\ast}\nbigk,\DDlambda_{\nbigk}\big)
\lrarr\big(\nbigp_{\ast}\nbigv^{(2)},\DDlambda_{\nbigv^{(2)}}\big)$. We~set
$\mu_0:=\mu_L(\nbigp_{\ast}\nbigv)$.
Because
$\bigoplus_{i=1,2}\big(\nbigp_{\ast}\nbigv^{(i)},\DDlambda_{\nbigv^{(i)}}\big)$
and $\big(\nbigp_{\ast}\nbigv,\DDlambda\big)$ are $\mu_L$-semistable
with the same slope $\mu_0$,
we~obtain that
$\big(\nbigp_{\ast}\nbigk,\DDlambda_{\nbigk}\big)$
is also $\mu_L$-semistable
with $\mu_L(\nbigp_{\ast}\nbigk)=\mu_0$.

Suppose that $\nbigk\neq 0$,
i.e., $\nbigi\neq 0$.
Because $\nbigi$ is a subsheaf of $\nbigv^{(2)}$,
we also obtain a filtered sheaf $\nbigp_{\ast}\nbigi$
induced by~$\nbigp_{\ast}\nbigv^{(2)}$.
Because $\nbigi\simeq\nbigk$,
we obtain a filtered sheaf $\nbigp'_{\ast}\nbigi$ over $\nbigi$
induced by~$\nbigp_{\ast}\nbigk$.
Then, we obtain
\begin{gather*}
\mu_0=\mu_L(\nbigp_{\ast}\nbigk)
=\mu_L(\nbigp'_{\ast}\nbigi)
\leq \mu_L(\nbigp_{\ast}\nbigi)
\leq \mu_L\big(\nbigp_{\ast}\nbigv^{(2)}\big)=\mu_0.
\end{gather*}
Because $\big(\nbigp_{\ast}\nbigv^{(2)},\DDlambda_{\nbigv^{(2)}}\big)$
is $\mu_L$-stable
and because $\nbigi\neq 0$,
we obtain that
$\rank(\nbigi)=\rank\nbigv^{(2)}$,
i.e.,~$\nbigi$ and~$\nbigv^{(2)}$
are generically isomorphic.
Because $\mu_L(\nbigp_{\ast}\nbigi)=\mu_L\big(\nbigp_{\ast}\nbigv^{(2)}\big)$,
Lemma~\ref{lem;19.2.18.120} implies that
$\nbigp_{\ast}\nbigi\lrarr\nbigp_{\ast}\nbigv^{(2)}$
is an isomorphism in codimension $1$.
Hence, there exists a closed algebraic subset $Z\subset X$
such that
$(i)$ the codimension of $Z$ is larger than $2$,
$(ii)$ $\nbigv^{(2)}_{|X\setminus Z}\subset \nbigv^{(1)}_{|X\setminus Z}$.
Because~$\nbigv^{(1)}$ is reflexive
we obtain that $\nbigv^{(2)}\subset\nbigv^{(1)}$.
\end{proof}

Let us study the case where $\nbigv^{(1)}\cap\nbigv^{(2)}=0$.
Let $\nbigv^{(3)}$ denote
the saturated $\lambda$-flat subsheaf of~$\nbigv$ generated by
$\nbigv^{(1)}+\nbigv^{(2)}$.
Let $\nbigp_{\ast}\nbigv^{(3)}$
denote the filtered sheaf over $\nbigv^{(3)}$
induced by $\nbigp_{\ast}\nbigv$.

\begin{Lemma}
\label{lem;19.2.18.130}
$\big(\nbigp_{\ast}\nbigv^{(3)},\DDlambda_{\nbigv^{(3)}}\big)$
is $\mu_L$-semistable,
and the induced morphism
$g\colon \nbigp_{\ast}\nbigv^{(1)}\oplus
 \nbigp_{\ast}\nbigv^{(2)}\lrarr\nbigp_{\ast}\nbigv^{(3)}$
is an isomorphism in codimension one.
\end{Lemma}

\begin{proof}
We obtain
$\mu_0=\mu_L\big(\nbigp_{\ast}\big(\nbigv^{(1)}\oplus\nbigv^{(2)}\big)\big)
\leq \mu_L\big(\nbigp_{\ast}\nbigv^{(3)}\big)
\leq \mu_L(\nbigp_{\ast}\nbigv)=\mu_0$.
Hence, we obtain that
$\mu_L\big(\nbigp_{\ast}\nbigv^{(3)}\big)=\mu_0$ and that
$\big(\nbigp_{\ast}\nbigv^{(3)},\DDlambda_{\nbigv^{(3)}}\big)$
is $\mu_L$-semistable.
Because
$g\colon \nbigp_{\ast}\nbigv^{(1)}\oplus
 \nbigp_{\ast}\nbigv^{(2)}\lrarr\nbigp_{\ast}\nbigv^{(3)}$
is generically an isomorphism,
and because they have the same slope,
$g$ is an isomorphism in codimension one
by Lemma~\ref{lem;19.2.18.120}.
\end{proof}

By Lemma~\ref{lem;19.2.18.130},
it is easy to observe that
there exists a finite sequence
of reflexive $\lambda$-flat subsheaves $\nbigv'_j$ $(j=1,\ldots,m)$
such that
$(i)$ the induced filtered $\lambda$-flat sheaves
 $\big(\nbigp_{\ast}\nbigv'_j,\DDlambda_{\nbigv_j'}\big)$ are~$\mu_L$-stable,
$(ii)$ the image of the induced morphism
$g\colon \nbigvtilde:=\bigoplus\nbigv'_j\lrarr\nbigv_1$
is generically an isomorphism.
Because
$\mu_0=\mu_L\bigl(\nbigp_{\ast}\nbigvtilde\bigr)
\leq\mu_L(\nbigp_{\ast}\nbigv_1)
\leq\mu_L(\nbigp_{\ast}\nbigv)=\mu_0$,
we obtain that
$\mu_L\bigl(\nbigp_{\ast}\nbigvtilde\bigr)
=\mu_L(\nbigp_{\ast}\nbigv_1)
=\mu_L(\nbigp_{\ast}\nbigv)$.
Hence,
$g$ is an isomorphism in codimension one
by Lemma~\ref{lem;19.2.18.120}.
Because both $\nbigp_{\ast}\nbigvtilde$
and $\nbigp_{\ast}\nbigv_1$
are reflexive,
we obtain that
$\nbigp_{\ast}\nbigvtilde\simeq\nbigp_{\ast}\nbigv_1$.
Thus, we obtain Proposition~\ref{prop;19.2.18.140}.
\end{proof}

\subsection{Mehta--Ramanathan type theorems}
\label{subsection;20.7.7.21}

Let $X$ be a smooth connected $n$-dimensional
projective variety
with a simple normal crossing hypersurface $H$.
Let $H=\bigcup_{i\in\Lambda}H_i$ be the irreducible decomposition.
Let $L$ be an ample line bundle on $X$.

\subsubsection{Restriction to general curves}

Let $\big(\nbigp_{\ast}\nbigv,\DDlambda\big)$ be
a reflexive filtered $\lambda$-flat sheaf
on $(X,H)$.
There exists a Zariski closed subset
$W\subset X$ with $\dim W<\dim H$
such that
$(i)$~the singular locus of $H$ is contained in $W$,
$(ii)$~$\nbigp_{\ast}\nbigv_{|X\setminus W}$
is a filtered bundle on $(X\setminus W,H\setminus W)$.

Let $Y$ be a smooth curve in $X$ such that
$(i)$~$Y\cap W=\varnothing$,
$(ii)$~$Y$ intersects with the smooth part of $H$ transversally.
Set $H_Y:=H\cap Y$. We~obtain a locally free $\nbigo_Y(\ast H_Y)$-module
$\nbigv_{|Y}$.
It~is equipped with the induced flat $\lambda$-connection
$\DDlambda_{|Y}$.
Let $\vecb\in\real^{H_Y}$.
For any $P\in H_Y$,
there exists $i\in\Lambda$
such that $P\in H_i$. We~choose $\veca(P,\vecb)\in\real^{\Lambda}$
such that $a(P,\vecb)_i=b(P)$,
and we obtain
an $\nbigo_{Y,P}$-submodule
$\nbigp_{\vecb}\big(\nbigv_{|Y}\big)_P:=
\nbigp_{\veca(P,\vecb)}(\nbigv)_P$
of $\big(\nbigv_{|Y}\big)_P$,
which is independent of the choice of $\veca(P,\vecb)$
as above.
There exists a locally free $\nbigo_Y$-module
$\nbigp_{\vecb}\big(\nbigv_{|Y}\big)
\subset\nbigv_{|Y}$
whose stalk at $P$ is $\nbigp_{\vecb}\big(\nbigv_{|Y}\big)_P$.
Thus, we obtain a filtered $\lambda$-flat bundle
$\big(\nbigp_{\ast}\big(\nbigv_{|Y}\big),\DDlambda_{|Y}\big)$
which is denoted by~$\big(\nbigp_{\ast}\nbigv,\DDlambda\big)_{|Y}$.

\subsubsection{The stability condition}

\begin{Proposition}
\label{prop;19.1.30.100}
 A reflexive filtered $\lambda$-flat sheaf
 $\big(\nbigp_{\ast}\nbigv,\DDlambda\big)$
 on $(X,H)$
is $\mu_L$-stable (resp.~$\mu_L$-semistable)
if and only if the following holds:
\begin{itemize}\itemsep=0pt
\item
 For any $m_1>0$,
 there exists $m>m_1$
 such that
 $\big(\nbigp_{\ast}\nbigv,\DDlambda\big)_{|Y}$
 is $\mu_L$-stable $($resp.~$\mu_L$-semistable$)$,
 where $Y$ denotes a generic $1$-dimensional
 complete intersection of
 hypersurfaces of $L^{\otimes m}$.
\end{itemize}
\end{Proposition}
\begin{proof}
The case $\lambda=1$ is already studied in~\cite[Section~13.2]{Mochizuki-wild}.
The case $\lambda\neq 0$ is reduced to the case $\lambda=1$.
As for the case $\lambda=0$,
we can prove the claim of the proposition
by the argument in~\cite[Section~3.4]{mochi4},
which closely follows
the arguments of Mehta--Ramanathan~\cite{mehta-ramanathan1, mehta-ramanathan2}
and Simpson~\cite{s5}. We~use $\nbigw=\Omega^1(ND)$
for a large $N$ instead of $\Omega^1(\log D)$
in~\cite[Section~3.4]{mochi4}.
(See~also~\cite[Section~13.2]{Mochizuki-wild}.)
\end{proof}

\subsubsection{Restrictions of morphisms and the polystability condition}

Let us give a complement on
the restriction of
morphisms of reflexive filtered $\lambda$-flat sheaves
to generic complete intersection curves,
which is a variant of~\cite[Lemma 3.9]{s5}.
 Let $\big(\nbigp_{\ast}\nbigv_i,\DDlambda\big)$ $(i=1,2)$
 be reflexive filtered $\lambda$-flat sheaves
 on $(X,H)$.
 Let $\Hom\bigl(\big(\nbigp_{\ast}\nbigv_1,\DDlambda_1\big),
 \big(\nbigp_{\ast}\nbigv_2,\DDlambda_2\big)\bigr)$
 denote the vector space of
 morphisms of filtered $\lambda$-flat sheaves
 $\big(\nbigp_{\ast}\nbigv_1,\DDlambda_1\big)
 \lrarr
 \big(\nbigp_{\ast}\nbigv_2,\DDlambda_2\big)$. We~shall prove
a refined claim (Proposition~\ref{prop;20.7.7.5})
of the following proposition
in Sections~\ref{subsection;20.7.7.2}--\ref{subsection;20.7.7.4}.

\begin{Proposition}\label{prop;20.7.7.1}
 There exists a positive integer $m_0>0$
 such that the restriction
\begin{gather*}
 \Hom\bigl(\big(\nbigp_{\ast}\nbigv_1,\DDlambda_1\big),
 \bigl(\nbigp_{\ast}\nbigv_2,\DDlambda_2\bigl)\bigr)
 \lrarr
 \Hom\bigl(\big(\nbigp_{\ast}\nbigv_1,\DDlambda_1\big)_{|Y},
 \bigl(\nbigp_{\ast}\nbigv_2,\DDlambda_2\bigl)_{|Y}\bigr)
\end{gather*}
is an isomorphism
for a generic $1$-dimensional complete intersection $Y$
of hypersurfaces of $L^{\otimes m}$ $(m\geq m_0)$.
 \end{Proposition}

Before going to the proof of Proposition~\ref{prop;20.7.7.1},
we state a variant of Proposition~\ref{prop;19.1.30.100}
on the $\mu_L$-polystability condition.
\begin{Corollary}
\label{cor;20.7.7.20}
 A reflexive filtered $\lambda$-flat sheaf
 $\big(\nbigp_{\ast}\nbigv,\DDlambda\big)$
 on $(X,H)$
is $\mu_L$-polystable
if and only if the following holds:
\begin{itemize}\itemsep=0pt
\item
 For any $m_1>0$,
 there exists $m>m_1$
 such that
 $\big(\nbigp_{\ast}\nbigv,\DDlambda\big)_{|Y}$
 is $\mu_L$-polystable,
 where $Y$ denotes the $1$-dimensional
 complete intersection of
 generic hypersurfaces of $L^{\otimes m}$.
\end{itemize}
\end{Corollary}
\begin{proof}
If $\big(\nbigp_{\ast}\nbigv,\DDlambda\big)$ is
$\mu_L$-polystable,
we obtain a decomposition
$\big(\nbigp_{\ast}\nbigv,\DDlambda\big)
=\bigoplus \big(\nbigp_{\ast}\nbigv_i,\DDlambda_i\big)$
into $\mu_L$-stable filtered $\lambda$-flat sheaves.
Applying Proposition~\ref{prop;19.1.30.100}
to each stable component,
we obtain the ``only if'' claim.

Let $m_1$ be an integer larger than $m_0$
in Proposition~\ref{prop;20.7.7.1}
for
$\Hom\bigl(
\big(\nbigp_{\ast}\nbigv,\DDlambda\big),
\big(\nbigp_{\ast}\nbigv,\DDlambda\big)
\bigr)$.
Suppose that there exists $m>m_1$
such that
 $\big(\nbigp_{\ast}\nbigv,\DDlambda\big)_{|Y}$
 is $\mu_L$-polystable
 for a generic $1$-dimensional
 complete intersection $Y$ of
 hypersurfaces of $L^{\otimes m}$.
 We obtain the decomposition
 \begin{gather}
\label{eq;20.7.7.11}
 \big(\nbigp_{\ast}\nbigv,\DDlambda\big)_{|Y}
 =\bigoplus_{i=1}^{\ell}
 \big(\nbigp_{\ast}\nbigv_{Y,i},\DDlambda_{Y,i}\big)
 \end{gather}
 into stable filtered Higgs bundles.
 Let $\pi_{Y,i}$ denote the endomorphisms of
 $\big(\nbigp_{\ast}\nbigv,\DDlambda\big)_{|Y}$
obtai\-ned by composing the projection
$\big(\nbigp_{\ast}\nbigv,\DDlambda\big)_{|Y}
\lrarr
\big(\nbigp_{\ast}\nbigv_{Y,i},\DDlambda_{Y,i}\big)$
with respect to
the decomposition~(\ref{eq;20.7.7.11}),
with the inclusion
$\big(\nbigp_{\ast}\nbigv_{Y,i},\DDlambda_{Y,i}\big)
 \lrarr
 \big(\nbigp_{\ast}\nbigv,\DDlambda\big)_{|Y}$.
 Note that they satisfy
 $\pi_{Y,i}\circ\pi_{Y,i}=\pi_{Y,i}$,
 $\pi_{Y,i}\circ\pi_{Y,j}=0$ $(i\neq j)$
 and
 $\sum \pi_{Y,i}=\id$. By~Proposition~\ref{prop;20.7.7.1},
there uniquely exist the endomorphisms
$\pi_i$ of $\big(\nbigp_{\ast}\nbigv,\DDlambda\big)$
such that $\pi_{i|Y}=\pi_{Y,i}$. By~Proposition~\ref{prop;20.7.7.1} again,
they satisfy
 $\pi_{i}\circ\pi_{i}=\pi_{i}$,
 $\pi_{i}\circ\pi_{j}=0$ $(i\neq j)$
 and
 $\sum \pi_{i}=\id$.
Let $\nbigv_i\subset\nbigv$ denote the image of $\pi_i$. We~define
$\nbigp_{\veca}\nbigv_i
=\nbigv_i\cap\nbigp_{\veca}\nbigv$
for any $\veca\in\real^{\Lambda}$.
Because $\pi_i$ are compatible with $\DDlambda$
and the filtration $\nbigp_{\ast}\nbigv$,
we obtain the decomposition
$\big(\nbigp_{\ast}\nbigv,\DDlambda\big)
 =\bigoplus \big(\nbigp_{\ast}\nbigv_i,\DDlambda_i\big)$.
 By the construction,
 $\big(\nbigp_{\ast}\nbigv_i,\DDlambda_i\big) _{|Y}
 =\big(\nbigp_{\ast}\nbigv_{Y,i},\DDlambda_{Y,i}\big)$
 are stable.
 Hence,
 $\big(\nbigp_{\ast}\nbigv_i,\DDlambda_i\big)$
 are
 $\mu_L$-stable with
 $\mu_L\big(\nbigp_{\ast}\nbigv_i,\DDlambda_i\big)
 =\mu_L\big(\nbigp_{\ast}\nbigv_j,\DDlambda_j\big)$
 $(i\neq j)$, i.e.,
 $\big(\nbigp_{\ast}\nbigv,\DDlambda\big)$
 is $\mu_L$-polystable.
\end{proof}

\subsubsection{General Enriques--Severi lemma due to Mehta--Ramanathan}
\label{subsection;20.7.7.2}

To prove Proposition~\ref{prop;20.7.7.1},
we recall the general Enriques--Severi lemma
in~\cite{mehta-ramanathan1}.
Recall $n=\dim X$.
For a positive integer $m$,
let $S_m$ denote the projective space of
lines in $H^0\big(X,L^{\otimes m}\big)$.
For~sequ\-ences $\vecm=(m_1,\ldots,m_t)\in\seisuu_{>0}^t$
with $t\leq n-1$,
we set $S_{\vecm}:=\prod S_{m_i}$.
There exists the~cor\-res\-pon\-dence variety
$Z_{\vecm}\subset X\times S_{\vecm}$,
i.e.,
$Z_{\vecm}= \bigl\{(x,s_1,\ldots,s_t)\in X\times S_{\vecm}\mid
 s_i(x)=0,\,1\leq i\leq t\bigr\}$.
For any $s\in S_{\vecm}$, we set
$X_s:=Z_{\vecm}\times_{S_{\vecm}}\times\{s\}$.

Let $F$ be a coherent reflexive $\nbigo_X$-module on $X$.
For any $\vecm=(m_1,\ldots,m_t)\in\seisuu_{>0}^t$
with $t\leq n-1$,
and for any $s\in S_{\vecm}$,
we set $F_s:=F\otimes_{\nbigo_X}\nbigo_{X_s}$.
For any integer $m$,
let $F_s(m)=F_s\otimes_{\nbigo_X}L^{\otimes m}$.

According to~\cite[Proposition~1.5]{mehta-ramanathan1},
there exists a non-empty Zariski open subset
$\Stilde_{\vecm}\subset S_{\vecm}$
such that the following holds.
\begin{itemize}\itemsep=0pt
 \item $\Stilde_{\vecm}\times_{S_{\vecm}}Z_{\vecm}
 \lrarr \Stilde_{\vecm}$ is smooth.
 \item For any $s\in S_{\vecm}$,
 $F_s$ is a reflexive $\nbigo_{X_s}$-module.
\end{itemize}

In the proof of~\cite[Proposition~3.2]{mehta-ramanathan1},
the following proposition is proved.
 \begin{Proposition}
\label{prop;20.5.19.1}
 Let $t\leq n-2$.
 There exists a positive integer $m_0$
 depending only on $F$ such that the following holds:
 \begin{itemize}\itemsep=0pt
 \item For any $\vecm=(m_1,\ldots,m_t)\in\seisuu_{>0}^{t}$
	 with $m_i\geq m_0$,
	 there exists a non-empty Zariski open subset
	 $U\subset \Stilde_{\vecm}$
	 such that $H^1(X_s,F_s(-\ell))=0$
	 for any $s\in U$ and any $\ell\geq m_0$.
 \end{itemize}
 \end{Proposition}

\begin{Corollary}
\label{cor;20.5.19.2}
 Let $t\leq n-1$.
 There exists a positive integer $m_0$
 depending only on $F$ such that the following holds:
 \begin{itemize}\itemsep=0pt
 \item For any $\vecm=(m_1,\ldots,m_t)\in\seisuu_{>0}^{t}$
 with $m_i\geq m_0$,
	 there exists a non-empty Zariski open subset
	 $U\subset \Stilde_{\vecm}$ such that
	 $H^0(X_s,F_s(-\ell))=0$
	 for any $s\in U$
	 and any $\ell\geq m_0$.
 \end{itemize}
\end{Corollary}

\begin{proof}
Let $m_0$ be a positive integer
as in Proposition~\ref{prop;20.5.19.1}. We~also assume $H^0(X,F(-\ell))=0$
for any $\ell\geq m_0$. We~use an induction on $t$.
For $\vecm=(m_1,\ldots,m_t)$ with $m_i\geq m_0$,
we set
$\vecm'=(m_1,\ldots,m_{t-1})$. By~the assumption of the induction
and Proposition~\ref{prop;20.5.19.1},
there exists a non-empty Zariski open subset
$U_1'\subset \Stilde_{\vecm'}$
such that
$H^i(X_{s'},F_{s'}(-\ell))=0$
$(i=0,1)$
for any~$s'\in U_1'$
and any $\ell\geq m_0$.

For any $s\in S_{\vecm}$,
let $s'$ denote the image $s$ in $S_{\vecm'}$
by the projection $S_{\vecm}\lrarr S_{\vecm'}$.
There exists the exact sequence
\begin{gather}
\label{eq;20.5.19.30}
0\lrarr \nbigo_{X_{s'}}(-m_t)
\lrarr \nbigo_{X_{s'}}\lrarr \nbigo_{X_s}\lrarr 0.
\end{gather}
By~\cite[Proposition~1.5]{mehta-ramanathan1},
there exists a Zariski open subset
$U_1\subset \Stilde_{\vecm}$
such that if $s\in U_1$
then we obtain the following exact sequence
from (\ref{eq;20.5.19.30})
by taking the tensor product with $F$:
\begin{gather*}
 0\lrarr F_{s'}(-m_t)
\lrarr F_{s'}\lrarr F_s\lrarr 0.
\end{gather*}
We shrink $U_1$ so that
$U_1'$ contains the image of $U_1$
by the projection $S_{\vecm}\lrarr S_{\vecm'}$.
Let $\ell\geq m_0$.
For any $s\in U_1$,
we obtain
$H^0(X_{s'},F_{s'}(-\ell))=0$
and
$H^1(X_{s'},F_{s'}(-\ell-m_t))=0$
because $s'\in U_1'$.
Hence, we obtain
$H^0(X_s,F_s(-\ell))=0$.
\end{proof}

\begin{Corollary}
 \label{cor;20.5.19.3}
Let $t\leq n-1$.
 There exists $m_0$ depending only on $F$
 such that the following holds:
 \begin{itemize}\itemsep=0pt
 \item For any $\vecm=(m_1,\ldots,m_t)\in\seisuu_{>0}^{t}$
	with $m_i\geq m_0$,
	there exists a non-empty Zariski open subset
	$U\subset \Stilde_{\vecm}$
	such that
	the natural morphism
	$H^0(X,F)\lrarr H^0(X_{s},F_s)$
	is an isomorphism
	for any $s\in U$.
 \end{itemize}
\end{Corollary}

\begin{proof}
It~is enough to apply the argument
in the first paragraph of the proof of~\cite[Propo\-si\-tion~3.2]{mehta-ramanathan1}
with Proposition~\ref{prop;20.5.19.1}
and Corollary~\ref{cor;20.5.19.2}.
(This is essentially pointed out
in~\cite[Lem\-ma~3.9]{s5}.)
\end{proof}

Let $T^{\ast}_{X_s}X$ denote
the conormal bundle of $X_s$ in $X$
for any $s\in \Stilde_{\vecm}$.

 \begin{Corollary}
 \label{cor;20.5.19.4}
Let $t\leq n-1$.
 There exists $m_0$ depending only on $F$
 such that the following holds:
 \begin{itemize}\itemsep=0pt
 \item For any $\vecm=(m_1,\ldots,m_t)\in\seisuu_{>0}^{t}$
	 with $m_i\geq m_0$,
	 there exists a non-empty Zariski open subset
	 $U\subset \Stilde_{\vecm}$ such that
	 $H^0\big(X_s,T^{\ast}_{X_s}X\otimes F_s\big)=0$
 	 for any $s\in U$.
 \end{itemize}
 \end{Corollary}

\begin{proof}
Because
 $T^{\ast}_{X_s}X
 \simeq
 \bigoplus_{j=1}^t\nbigo_{X_s}(-m_j)$,
 the claim follows from Corollary~\ref{cor;20.5.19.2}.
 \end{proof}

\subsubsection[Flat sections of reflexive filtered lambda-flat sheaves]
{Flat sections of reflexive filtered $\boldsymbol\lambda$-flat sheaves}

Let $H$ be a simple normal crossing hypersurface of $X$
with the irreducible decomposition
$H=\bigcup_{i\in\Lambda}H_i$.
Let $\big(\nbigp_{\ast}\nbigv,\DDlambda\big)$
be a reflexive filtered $\lambda$-flat sheaf
on $(X,H)$.
Note that there exists a~Zariski closed subset
$W\subset X$ with $\dim W<\dim H$
such that
$(i)$~$W$ contains the singular locus of $H$,
$(ii)$~$\nbigp_{\ast}\nbigv_{|X\setminus W}$
is a filtered bundle on $(X\setminus W,H\setminus W)$.
For any $\vecm\in\seisuu^{n-1}$
and for any $s\in S_{\vecm}$,
we set $H_s:=X_s\times_XH$.

According to~\cite[Proposition~1.5]{mehta-ramanathan1},
there exists a non-empty Zariski open subset
$S^{\circ}_{\vecm}\subset S_{\vecm}$
such that the following holds:
\begin{itemize}\itemsep=0pt
 \item $Z_{\vecm}\times_{S_{\vecm}}S^{\circ}_{\vecm}
 \lrarr S^{\circ}_{\vecm}$
 is smooth.
 \item For any $s\in S^{\circ}_{\vecm}$,
 $X_s\cap W=\varnothing$ holds,
 and $X_s$ intersects with $H$
 in $H\setminus W$ transversally.
 Moreover,
 $\nbigp_{\veca}\nbigv_s:=
 \nbigp_{\veca}\nbigv_{|X_s}$
 $(\veca\in\real^{\Lambda})$
 are locally free $\nbigo_{X_s}$-modules.
 \end{itemize}

There exists a non-negative integer $N$ such that
$\DDlambda$ induces a morphism of sheaves
$\DDlambda$: $\nbigp_{\veca}\nbigv\lrarr
\nbigp_{\veca}\nbigv
\otimes\Omega^1_X(\log H)\otimes\nbigo_X(NH)$.
For $j=0,1,\ldots,n$, we set
\begin{gather*}
 \nbigc^j_N(\nbigp_{\veca}\nbigv)
 =\nbigp_{\veca}\nbigv
 \otimes\Omega^j_X(\log H)\otimes\nbigo_X(jNH).
\end{gather*}
The flat $\lambda$-connection $\DDlambda$ induces
$\DDlambda\colon \nbigc^j_{N}(\nbigp_{\veca}\nbigv)
\lrarr
\nbigc^{j+1}_{N}(\nbigp_{\veca}\nbigv)$
such that
$\DDlambda\wedge\DDlambda=0$.
Thus,
we obtain a complex of sheaves
$\nbigc^{\bullet}_{N}\big(\nbigp_{\veca}\nbigv,\DDlambda\big)$ on $X$.
 Clearly, the following holds:
\begin{gather*}
 \hyperh^0\bigl(X, \nbigc^{\bullet}_{N}\big(\nbigp_{\veca}\nbigv,\DDlambda\big) \bigr)
 = \Ker\Bigl(
 H^0(X,\nbigp_{\veca}\nbigv)
 \stackrel{\DDlambda}{\lrarr}
 H^0\bigl(X,
 \nbigp_{\veca}\nbigv\otimes
 \Omega_X^1(\log H)\otimes\nbigo_X(NH)\bigr)
 \Bigr).
\end{gather*}
For any $s\in S^{\circ}_{\vecm}$,
we obtain the filtered $\lambda$-flat bundle
$\big(\nbigp_{\ast}\nbigv_s,\DDlambda_s\big):=
\big(\nbigp_{\ast}\nbigv,\DDlambda\big)_{|X_s}$.
Let $\veca(s)$ denote the image of
$\real^{\Lambda}\lrarr\real^{H_s}$
induced by the natural map $H_s\lrarr\Lambda$.
Let $\iota_s\colon X_s\lrarr X$ denote the inclusion. We~obtain the natural morphism of complexes of sheaves
$\nbigc_N^{\bullet}\big(\nbigp_{\veca}\nbigv,\DDlambda\big)
\lrarr\iota_{s\ast}
 \nbigc_N^{\bullet}\big(\nbigp_{\veca(s)}\nbigv_s,\DDlambda_s\big)$,
which induces
\begin{gather}
\label{eq;20.7.8.10}
	 \hyperh^0\big(X,\nbigc^{\bullet}_{N}
	 \big(\nbigp_{\veca}\nbigv,\DDlambda\big)\big)
	 \lrarr	 \hyperh^0\big(X_s,\nbigc^{\bullet}_{N}
	 \big(\nbigp_{\veca(s)}\nbigv_s,\DDlambda_s\big)\big).
\end{gather}

 The following proposition is essentially~\cite[Lemma 3.9]{s5}.
 \begin{Proposition}
\label{prop;20.5.19.10}
 There exists a positive integer $m_0>0$
 such that the following claim holds
 for any $\vecm=(m_1,\ldots,m_{n-1})$ with $m_i\geq m_0$
 and a non-empty Zariski open subset $U\subset S^{\circ}_{\vecm}$.
 \begin{itemize}\itemsep=0pt
 \item For any $s\in U$,
	 the natural morphism {\rm(\ref{eq;20.7.8.10})}
 is an isomorphism.
 \end{itemize}
 \end{Proposition}

 \begin{proof}
 According to Corollary~\ref{cor;20.5.19.3},
 if $m_0$ is sufficiently large,
 there exists a non-empty Zariski open subset
 $U_1\subset S^{\circ}_{\vecm}$
 such that
 the following natural morphisms are isomorphisms
 for any~$s\in U$:
 \begin{gather*}
 H^0(X,\nbigp_{\veca}\nbigv) \lrarr
 H^0(X_s,\nbigp_{\veca(s)}\nbigv_s),
\\
 H^0\bigl(X, \nbigp_{\veca}\nbigv\otimes_{\nbigo_X}\Omega^1_X(\log H)
 \otimes\nbigo_X(NH) \bigr)
 \\ \qquad
 {}\lrarr
 H^0\bigl(X_s,\nbigp_{\veca(s)}\nbigv_s\otimes_{\nbigo_{X_s}}
 \bigl(\Omega^1_{X}(\log H) \otimes\nbigo_{X_s}(NH_s)\bigr) \bigr).
 \end{gather*}
 There exists the following exact sequence:
 \begin{align*}
 0&\lrarr
 T_{X_s}^{\ast}X\otimes\nbigp_{\veca(s)}\nbigv_s\otimes\nbigo_{X_s}(NH_s)
 \lrarr
 \bigl(\Omega^1_X(\log H)\otimes_{\nbigo_X}\nbigo_{X_s}(NH_s)\bigr)
 \otimes\nbigp_{\veca(s)}\nbigv_s
 \\
&\lrarr \Omega^1_{X_s}(\log H_s)\otimes\nbigp_{\veca(s)}\nbigv_s
 \otimes\nbigo_{X_s}(NH_s)
 \lrarr 0.
 \end{align*}
 According to Corollary~\ref{cor;20.5.19.4},
 if $m_0$ is sufficiently large,
 there exists a non-empty Zariski open subset
 $U_2\subset U_1$ such that the following holds
 for any $s\in U_2$:
 \begin{gather*}
 H^0\bigl(X_s,
 T_{X_s}^{\ast}X\otimes\nbigp_{\veca(s)}\nbigv_s
 \otimes\nbigo_{X_s}(NH_s) \bigr)=0.
 \end{gather*}
 Hence, the natural morphism
\begin{gather*}
 H^0\bigl( X_s,(\Omega^1_X(\log H)\otimes_{\nbigo_X}\nbigo_{X_s}(NH_s))
 \otimes\nbigp_{\veca(s)}\nbigv_s \bigr)
 \\ \qquad
\lrarr H^0\bigl(X_s,
 \Omega^1_{X_s}(\log H_s)\otimes\nbigp_{\veca(s)}\nbigv_s
 \otimes\nbigo_{X_s}(NH_s) \bigr)
\end{gather*}
 is injective for any $s\in U_2$. We~obtain the injectivity of
 the following natural morphism
 for any~$s\in U_2$:
\begin{gather*}
 H^0\bigl(X, \Omega^1_X(\log H)\otimes\nbigp_{\veca}\nbigv\otimes\nbigo_{X}(NH) \bigr)
\lrarr H^0\bigl(X_s, \Omega^1_{X_s}(\log H_s)\otimes\nbigp_{\veca(s)}\nbigv_s
 \otimes\nbigo_{X_s}(NH_s) \bigr).
\end{gather*}
 Then, we obtain the claim of the proposition.
 \end{proof}

\subsubsection[Morphisms of reflexive filtered lambda-flat sheaves]
{Morphisms of reflexive filtered $\boldsymbol\lambda$-flat sheaves}
\label{subsection;20.7.7.4}

Let $\nbigp_{\ast}\nbigv_i$ $(i=1,2)$
 be reflexive filtered sheaves
 with meromorphic flat $\lambda$-connection $\DDlambda_i$
 on $(X,H)$.
 Let $\Hom\bigl(\big(\nbigp_{\ast}\nbigv_1,\DDlambda_1\big),
 \big(\nbigp_{\ast}\nbigv_2,\DDlambda_2\big)\bigr)$
 denote the vector space of
 morphisms of filtered $\lambda$-flat sheaves
 $\big(\nbigp_{\ast}\nbigv_1,\DDlambda_1\big)
 \lrarr
 \big(\nbigp_{\ast}\nbigv_2,\DDlambda_2\big)$.
 \begin{Proposition}
\label{prop;20.7.7.5}
 There exists a positive integer $m_0>0$
 such that the following claim holds
 for any $\vecm=(m_1,\ldots,m_{n-1})$ with $m_i\geq m_0$
 and for a non-empty Zariski open subset $U\subset S_{\vecm}$.
 \begin{itemize}\itemsep=0pt
 \item For any $s\in U$,
	let $\big(\nbigp_{\ast}\nbigv_{i,s},\DDlambda_{i,s}\big)$
	 denote the induced filtered $\lambda$-flat bundles
	 on $(X_s,H_s)$.
	 Then, the natural morphism
	 \begin{gather*}
	 \Hom\bigl(\big(\nbigp_{\ast}\nbigv_1,\DDlambda_1\big),
	 \big(\nbigp_{\ast}\nbigv_2,\DDlambda_2\big) \bigr)
	 \lrarr 	 \Hom\bigl(\big(\nbigp_{\ast}\nbigv_{1,s},\DDlambda_{1,s}\big),
	 \big(\nbigp_{\ast}\nbigv_{2,s},\DDlambda_{2,s}\big) \bigr)
	 \end{gather*}
 is an isomorphism.
 \end{itemize}
 \end{Proposition}
\begin{proof}
For any $\veca\in\real^{\Lambda}$,
let $\nbigp_{\veca}\nhom(\nbigv_1,\nbigv_2)$
denote the subsheaf of
 the $\nbigo_X(\ast H)$-module
 $\nhom_{\nbigo_X(\ast H)}(\nbigv_1,\nbigv_2)$
 determined as follows
 for any open subset $U\subset X$:
\begin{gather*}
 H^0\bigl( U,\nbigp_{\veca}(\nhom(\nbigv_1,\nbigv_2)) \bigr)
\! =
 \bigl\{ f\!\in\! H^0(U,\nhom(\nbigv_1,\nbigv_2))\mid
 f\big(\nbigp_{\vecb}\nbigv_{1|U}\big)
\! \subset\nbigp_{\veca+\vecb}\big(\nbigv_{2|U}\big)
 \,\forall\vecb\in\real^{\Lambda}
 \bigr\}.
\end{gather*}
It~is easy to see that
$\nbigp_{\veca}\nhom(\nbigv_1,\nbigv_2)$ are reflexive
$\nbigo_X$-modules.
Thus, we obtain a reflexive filtered sheaf
$\nbigp_{\ast}\nhom(\nbigv_1,\nbigv_2)$
with the induced flat $\lambda$-connection $\DDlambdatilde$. We~can easily observe that
\begin{gather*}
 \Hom\bigl( \big(\nbigp_{\ast}\nbigv_1,\DDlambda_1\big),
 \big(\nbigp_{\ast}\nbigv_2,\DDlambda_2\big) \bigr)
 =\hyperh^0\bigl( X,\nbigc^{\bullet}_N\bigl(
 \nbigp_{\veczero}\nhom(\nbigv_1,\nbigv_2),\DDlambdatilde \bigr) \bigr)
\end{gather*}
for any large $N$,
where $\veczero=(0,\ldots,0)\in\real^{\Lambda}$.
Then, the claim follows from Proposition~\ref{prop;20.5.19.10}.
\end{proof}

\subsection[Good filtered lambda-flat bundles and ramified coverings]
{Good filtered $\boldsymbol\lambda$-flat bundles and ramified coverings}

\subsubsection{Pull back}\label{subsection;20.2.14.11}

Let $X$ be any complex manifold with a simple normal crossing
hypersurface $H$.
Let $\big(\nbigp_{\ast}\nbigv,\DDlambda\big)$ be a good
filtered $\lambda$-flat bundle. We~set $e:=\rank(\nbigv)!$.

For any point $P\in X$,
let $(X_P,z_1,\ldots,z_n)$ denote
an admissible holomorphic coordinate neighbourhood around $P$. We~set $H_P:=X_P\cap H$
and $H_{P,i}:=H_P\cap\{z_i=0\}$. We~set $\big(\nbigp_{\ast}\nbigv_P,\DDlambda_P\big):=
 \big(\nbigp_{\ast}\nbigv,\DDlambda\big)_{|X_P}$.

By using the coordinate system,
we may regard $X_P$
as an open subset of $\cnum^n$.
Let $\varphi_P$: $\cnum^n\lrarr \cnum^n$
be given by
$\varphi_P(\zeta_1,\ldots,\zeta_n)
=\big(\zeta_1^{e},\ldots,\zeta_{\ell(P)}^{e},\zeta_{\ell(P)+1},
\ldots,\zeta_{n}\big)$. We~set
$\Xtilde_P:=
 \varphi_P^{-1}(X_P)$,
$\Htilde_P:=\varphi_{P}^{-1}(H_P)$
and
$\Htilde_{P,i}:=\varphi_P^{-1}(H_{P,i})$. We~set
$G_P:=\prod_{i=1}^{\ell(P)}\{\alpha_i\in\cnum\mid \alpha_i^{e}=1\}$.
It~is identified with the Galois group of the ramified covering
$\varphi_P$
by the action as in Section~\ref{subsection;20.2.14.1}.

We obtain the $G_P$-equivariant
good filtered $\lambda$-flat bundle
$\big(\nbigp_{\ast}\nbigvtilde_P,\DDtilde^{\lambda}_P\big)
:=\varphi_P^{\ast}\big(\nbigp_{\ast}\nbigv_P,\DDlambda_P\big)$
on~$\big(\Xtilde_P,\Htilde_P\big)$.

\begin{Lemma}
$\big(\nbigp_{\ast}\nbigvtilde_P,\DDtilde^{\lambda}_P\big)$
is unramifiedly good.
\end{Lemma}
\begin{proof}
See~\cite[Lemma 2.2.7]{Mochizuki-wild}.
\end{proof}

\subsubsection{The associated graded bundles}

We obtain
the $G_P$-equivariant filtered bundles
$\lefttop{1}\Gr^F_c\big(\nbigp_{\ast}\nbigvtilde_P\big)$ $(c\in\real)$
on $\big(\Htilde_{P,1},\del \Htilde_{P,1}\big)$.
There exists the $G_P$-equivariant decomposition
\begin{gather*}
 \lefttop{1}\Gr^F_c\big(\nbigp_{\ast}\nbigvtilde_P\big)
=\bigoplus_{\ell=0}^{e-1}
 \GG_{\ell}
 \lefttop{1}\Gr^F_c\big(\nbigp_{\ast}\nbigvtilde_P\big),
\end{gather*}
where $(\alpha,1,\ldots,1)$
acts on $\GG_{\ell}\lefttop{1}\Gr^F_c\big(\nbigp_{\ast}\nbigvtilde_P\big)$
as the multiplication by $\alpha^{\ell}$.

\begin{Lemma}
\label{lem;20.2.14.10}
The pull back naturally induces the isomorphism
$ (\varphi_{|\Htilde_{P,1}})^{\ast}
 \bigl(
 \lefttop{1}\Gr^F_{c/e}(\nbigp_{\ast}\nbigv_P)
 \bigr)
\simeq
\GG_0\lefttop{1}\Gr^F_c\big(\nbigp_{\ast}\nbigvtilde_P\big)$.
As a result,
$\lefttop{1}\Gr^F_{c/e}(\nbigp_{\ast}\nbigv_P)$
is the descent of
$\GG_0\lefttop{1}\Gr^F_c\big(\nbigp_{\ast}\nbigvtilde_P\big)$.
More generally,
the pull back and the multiplication by $\zeta_1^{\ell}$
induces an isomorphism
$ (\varphi_{|\Htilde_{P,1}})^{\ast}
 \bigl(
 \lefttop{1}\Gr^F_{(c+\ell)/e}(\nbigp_{\ast}\nbigv_P)
 \bigr)
\simeq
\GG_{\ell}\lefttop{1}\Gr^F_c\big(\nbigp_{\ast}\nbigvtilde_P\big)$.
\end{Lemma}
Clearly, there exist a similar decomposition
 $\lefttop{i}\Gr^F_c\big(\nbigp_{\ast}\nbigvtilde_P\big)
=\bigoplus_{\ell=0}^{e-1}
 \GG_{\ell}
 \lefttop{i}\Gr^F_c\big(\nbigp_{\ast}\nbigvtilde_P\big)$
and isomorphisms for any $i=1,\ldots,\ell(P)$.

\subsubsection{Residues}
\label{subsection;21.5.5.10}

Let us recall that we obtain the endomorphisms
$\Res_j\big(\DDlambda\big)$ $(j\in\Lambda)$
on $\lefttop{j}\Gr_c^F(\nbigp_{\ast}\nbigv)$
by using Lemma~\ref{lem;20.2.14.10}.
(See~\cite[Section~2.5.2]{Mochizuki-wild}
for more detailed explanations.)

Let $P$ be any point of $H$.
First, let us construct the residues
$\Res_{1}\big(\DDlambda_P\big)$
on
$\lefttop{1}\Gr^F(\nbigp_{\ast}\nbigv_P)$.
At~any $Q\in \Htilde_{P,1}$,
we obtain the formal decomposition
$\big(\nbigp_{\veca}\nbigvtilde_P,\DDlambda_P\big)
 \otimes\nbigo_{\Xtilde_P,\Qhat}
=\bigoplus
 \big(\nbigp_{\veca}\nbigvtilde_{\gminia},\DDtilde^{\lambda}_{\gminia}\big)$
as in~(\ref{eq;20.2.1.1}).
For $a_1-1<c\leq a_1$,
we obtain the endomorphisms
$\Res_1\big(\DDtilde_P^{\lambda}\big)_Q$
of~$\lefttop{1}\Gr^F_{c}\big(\nbigp_{\veca}\nbigvtilde\big)_{|Q}$
as the residue of
$\bigoplus \bigl( \DDtilde^{\lambda}_{\gminia}
-{\rm d}\gminiatilde\id_{\nbigvtilde_{\gminia}} \bigr)$ at $Q$.
According to
\cite[Lemma 2.5.2]{Mochizuki-wild},
by varying $Q\in \Htilde_{P,1}$,
we obtain the endomorphism
$\Res_1\big(\DDtilde_P^{\lambda}\big)$
of the filtered bundle
$\lefttop{1}\Gr^F_{c}\big(\nbigp_{\ast}\nbigvtilde_P\big)$.
It~is $G_P$-equivariant.
Hence, we obtain
$\Res_1\big(\DDlambda_P\big)$
on $\lefttop{1}\Gr^F_{c}(\nbigp_{\ast}\nbigv_P)$
as the descent of
$\frac{1}{e}\Res_1(\DDtilde^{\lambda}_P)$
on
$\GG_0\lefttop{1}\Gr^F_{ec}\big(\nbigp_{\ast}\nbigvtilde_P\big)$.
The factor $\frac{1}{e}$
comes from the relation
${e}\,{\rm d}\zeta_1/\zeta_1={\rm d}z_1/z_1$.
Similarly, we obtain
$\Res_i\big(\DDlambda_P\big)$
$\lefttop{i}\Gr^F_{c}(\nbigp_{\ast}\nbigv_P)$
for $i=1,\ldots,\ell(P)$.

It~is easy to see that
there exists a globally defined endomorphism
$\Res_j\big(\DDlambda\big)$
on $\lefttop{j}\Gr^F_c(\nbigp_{\ast}\nbigv)$
which is equal to the endomorphisms
constructed locally around $P\in H_j$ as above.

\subsubsection{Parabolic weights}
\label{subsection;20.7.9.1}

We introduce some notation. We~set
$\Par(\nbigp_{\ast}\nbigv,j):=\bigl\{ b\in\real\mid
 \lefttop{j}\Gr^F_{b}(\nbigp_{\ast}\nbigv)\neq 0 \bigr\}$
for $j\in\Lambda$.
Because we shall often use the pull back by a ramified covering
as in Section~\ref{subsection;20.2.14.11},
for a fixed $e$,
it is convenient to consider
\begin{gather*}
 \Partilde(\nbigp_{\ast}\nbigv,j):=
 \bigl\{ c+m/e\mid c\in\Par(\nbigp_{\ast}\nbigv,j),\, m\in\seisuu \bigr\}.
\end{gather*}
Note that
$\Par\big(\nbigp_{\ast}\nbigvtilde_P,j\big)
=\bigl\{ eb\mid b\in\Partilde(\nbigp_{\ast}\nbigv_P,j) \bigr\}$
for $P\in H$.
(See Lemma~\ref{lem;20.2.14.10}.)
We define
\begin{gather*}
 \gaptilde(\nbigp_{\ast}\nbigv,j):= \min \bigl\{|b_1-b_2|\mid
 b_1,b_2\in\Partilde(\nbigp_{\ast}\nbigv,j),\,b_1\neq b_2 \bigr\}.
\end{gather*}
If $\Lambda$ is finite, we also set
$\gaptilde(\nbigp_{\ast}\nbigv):=
\min\limits_{j\in\Lambda}\gaptilde(\nbigp_{\ast}\nbigv,j)$.

For each $\veca\in\real^{\Lambda}$, we set
\begin{gather*}
\Par(\nbigp_{\ast}\nbigv,\veca,i):=
\Par(\nbigp_{\ast}\nbigv,i)\cap\openclosed{a_i-1}{a_i},\qquad
\Partilde(\nbigp_{\ast}\nbigv,\veca,i):=
\Partilde(\nbigp_{\ast}\nbigv,i)\cap\openclosed{a_i-1}{a_i}.
\end{gather*}

We remark the following obvious lemma.
\begin{Lemma}
\label{lem;20.2.12.20}
For each $j$,
there exists $a_j\in\real$
such that
$|a_j-b|>(4e\rank\nbigv)^{-1}$
for any $b\in\Partilde(\nbigp_{\ast}\nbigv,j)$.
\end{Lemma}

\subsection[Approximation by model filteredlambda-flat bundles]
{Approximation by model filtered $\boldsymbol\lambda$-flat bundles}
\label{subsection;20.2.6.20}

\subsubsection[Model filtered lambda-flat bundles]
{Model filtered $\boldsymbol\lambda$-flat bundles}

\label{subsection;20.2.9.1}

Let $Z$ be a complex manifold.
Let $Y$ be a neighbourhood of
$\{0\}\times Z$ in $\cnum\times Z$. We~set $H_Y:=\{0\}\times Z$.
Let $e$ be a positive integer.
Let $\zeta$ be the standard complex coordinate of $\cnum$.
Consider
$\varphi\colon \cnum\times Z\lrarr \cnum\times Z$
induced by $\zeta\longmapsto \zeta^e$. We~set
$\Ytilde:=\varphi^{-1}(Y)$
and $\Htilde:=\varphi^{-1}(H)$.
The induced morphism $\Ytilde\lrarr Y$
is also denoted by $\varphi$.
Let $G$ denote the group of the $e$-th roots of $1$,
which is naturally identified with the Galois group of
the ramified covering $\varphi$.

Let $\nbigi$ be a finite subset
of $H^0\big(\Ytilde,\nbigo_{\Ytilde}\big({\ast} \Htilde\big)\big)$
which is preserved by the $G$-action.
Let $S_1$ and~$S_2$ be finite subsets of
$\openclosed{{-}1}{0}$
and $\cnum$, respectively.
Let $V_{\gminia,a,\alpha}$
$((\gminia,a,\alpha)\in \nbigi\times S_1\times S_2)$
be finite dimensional $\cnum$-vector spaces
equipped with a nilpotent endomorphism
$f_{\gminia,a,\alpha}$.
Note that
$V_{\gminia,a,\alpha}$ may be $0$. We~suppose that
$\bigoplus_{\gminia,a,\alpha}V_{\gminia,a,\alpha}$
is a $G$-representation
such that
$(i)$ it is $G$-equivariant as a vector bundle
over $\nbigi$,
$(ii)$ $\bigoplus f_{\gminia,a,\alpha}$ commutes with
the $G$-action.

We set
$\nbigvtilde_{\gminia,a,\alpha}:=
 \nbigo_{\Ytilde}\big({\ast} \Htilde\big)
 \otimes V_{\gminia,a,\alpha}$. We~define the filtered bundle
$\nbigp_{\ast}\nbigvtilde_{\gminia,a,\alpha}$
over $\nbigvtilde_{\gminia,a,\alpha}$
by setting
\begin{gather*}
 \nbigp_{b}\nbigvtilde_{\gminia,a,\alpha}:=
 \nbigo_{\Ytilde}\big([b-a]\Htilde\big)\otimes V_{\gminia,a,\alpha}
\end{gather*}
for any $b\in\real$,
where $[b-a]:=\max\{n\in\seisuu\mid n\leq b-a\}$. We~define
the flat $\lambda$-connection $\DDtilde^{\lambda}_{\gminia,a,\alpha}$
on $\nbigvtilde_{\gminia,a,\alpha}$
by setting
\begin{gather*}
 \DDtilde^{\lambda}_{\gminia,a,\alpha}(v)
={\rm d}\gminia\cdot v+(\alpha v+f_{\gminia,a,\alpha}(v))\,{\rm d}\zeta/\zeta
\end{gather*}
for any $v\in V_{\gminia,a,\alpha}$,
which we regard as a section of $\nbigvtilde_{\gminia,a,\alpha}$
in a natural way.
Thus, we obtain a $G$-equivariant
filtered $\lambda$-flat bundle
$\bigoplus_{\gminia,a,\alpha}
 \big(\nbigp_{\ast}\nbigvtilde_{\gminia,a,\alpha},
 \DDtilde^{\lambda}_{\gminia,a,\alpha}\big)$,
called a model filtered $\lambda$-flat bundle.
If $\nbigi$ induces a good set of irregular values
in $\nbigo_{\Ytilde}\big({\ast} \Htilde\big)_Q/\nbigo_{\Ytilde,Q}$
at each $Q\in \Htilde$,
then $\bigoplus_{\gminia,a,\alpha}
 \big(\nbigp_{\ast}\nbigvtilde_{\gminia,a,\alpha},
 \DDtilde^{\lambda}_{\gminia,a,\alpha}\big)$
is an unramifiedly good filtered $\lambda$-flat bundle.
It~induces a filtered $\lambda$-flat bundle
on $(Y,H)$ as the descent,
which is also called a model filtered $\lambda$-flat bundle.

\subsubsection[Approximation of good filtered lambda-flat bundles]
{Approximation of good filtered $\boldsymbol\lambda$-flat bundles}

Let $\big(\nbigp_{\ast}\nbigv,\DDlambda\big)$
be any good filtered $\lambda$-flat bundle
on $(Y,H)$.
Assume the following condition:
\begin{Condition}
\label{condition;20.2.8.2}
For each $a\in\Par(\nbigp_{\ast}\nbigv)$,
the conjugacy classes of
 $\Res\big(\DDlambda\big)$
 on
 $\Gr^F_a(\nbigp_{\ast}\nbigv)_{|P}$
are independent of $P\in H$.
Note that this condition is trivially satisfied
if $\lambda\neq 0$.
\end{Condition}

We set $e:=\rank(\nbigv)!$.
Let $\varphi$ and $\big(\Ytilde,\Htilde\big)$
be as in Section~\ref{subsection;20.2.9.1}. We~set
$\big(\nbigp_{\ast}\nbigvtilde,\DDlambdatilde\big):=
\varphi^{\ast}\big(\nbigp_{\ast}\nbigv,\DDlambda\big)$.
For each $Q\in \Htilde$,
there exists a decomposition
\begin{gather*}
\big (\nbigp_{\ast}\nbigvtilde,\DDtilde^{\lambda}\big) \otimes \nbigo_{\Ytilde,\Qhat}
=\bigoplus_{\gminia\in\nbigi(Q)}\big(\nbigp_{\ast}\nbigvtilde_{\gminia},
 \DDtilde^{\lambda}_{\gminia}\big)
\end{gather*}
as in (\ref{eq;20.2.8.1}). We~obtain the vector spaces
$\Gr^F_a\big(\nbigp_0\nbigvtilde_{\gminia}\big)_{|Q}$
$(-1<a\leq 0)$
equipped with the endomorphisms
$\Res\big(\DDtilde^{\lambda}\big)$.
Condition~\ref{condition;20.2.8.2}
is equivalent to the following.
\begin{itemize}\itemsep=0pt
\item
 The conjugacy classes of
 $\Res\big(\DDtilde^{\lambda}\big)$
 on
 $\Gr^F_a\big(\nbigp_0\nbigvtilde_{\gminia}\big)_{|Q}$
 are independent of $Q\in \Htilde$
 for any $-1<a\leq 0$.
\end{itemize}
In particular,
the condition implies that there exists a decomposition
\begin{gather*}
 \Gr^F_a\big(\nbigp_0\nbigvtilde_{\gminia}\big)
=\bigoplus_{\alpha\in\cnum}
 \EE_{\alpha}\Gr^F_a\big(\nbigp_0\nbigvtilde_{\gminia}\big)
\end{gather*}
on $\Htilde$,
where $\Res\big(\DDtilde^{\lambda}\big)-\alpha\id$ are nilpotent
on $\EE_{\alpha}\Gr^F_a\big(\nbigp_0\nbigvtilde_{\gminia}\big)$.

Fix $P\in H$.
Let $\Ptilde\in \Htilde$
be determined by $\varphi\big(\Ptilde\big)=P$. We~set
$V_{\gminia,a,\alpha}:=
\EE_{\alpha}\Gr^F_a\big(\nbigp_0\nbigvtilde_{\gminia}\big)_{|\Ptilde}$.
Let $f_{\gminia,a,\alpha}$
be the nilpotent part of
$\Res\big(\DDtilde^{\lambda}\big)$
on
$V_{\gminia,a,\alpha}$.
For a neighbourhood
$Y_P$ of $P$ in $Y$,
we set
$H_P:=Y_P\cap H$,
$\Ytilde_P:=\varphi^{-1}(Y_P)$
and
$\Htilde_P:=\varphi^{-1}(H_P)$. We~may assume that
any $\gminia\in\nbigi\big(\Ptilde\big)$
has a lift $\gminiatilde$
in $H^0\big(\Ytilde_P,\nbigo_{\Ytilde_P}\big({\ast} \Htilde_P\big)\big)$.
From the set
$\{(\gminia,a,\alpha)\}\subset
\nbigi\big(\Ptilde\big)\times\openclosed{{-}1}{0} \times\cnum$
and
the tuples $(V_{\gminia,a,\alpha},f_{\gminia,a,\alpha})$,
we obtain a model filtered $\lambda$-flat bundle
\begin{gather*}
 \bigl( \nbigp_{\ast}\nbigvtilde_0,\DDtilde_0^{\lambda} \bigr)
:=\bigoplus_{\gminia,a,\alpha}
 \bigl(\nbigp_{\ast}\nbigvtilde_{\gminia,a,\alpha}, \DDtilde^{\lambda}_{\gminia,a,\alpha}\bigl)
\end{gather*}
on $\big(\Ytilde_P,\Htilde_P\big)$.
It~is unramifiedly good and naturally $G$-equivariant.
As the descent,
we obtain a~good filtered $\lambda$-flat bundle
$\big(\nbigp_{\ast}\nbigv_0,\DD_0^{\lambda}\big)$
on $(Y_P,H_P)$.

\begin{Lemma}[{assume Condition~\ref{condition;20.2.8.2}}]
\label{lem;20.2.14.20}
For any positive integer $m$,
there exist a neighbourhood~$Y_{P}$
of $P$ in $Y$
and an isomorphism of filtered bundles
$\Phi_m\colon \nbigp_{\ast}\nbigv_0\simeq\nbigp_{\ast}\nbigv_{|Y_P}$
such that the following holds:
\begin{itemize}\itemsep=0pt
\item
We set $\Phitilde_m:=\varphi^{\ast}(\Phi_m)$
and
$A:=\big(\Phitilde_m\big)^{\ast}\big(\DDtilde^{\lambda}\big)-\DDtilde_0^{\lambda}$
on $\Ytilde_P$.
Let
$A=\sum A_{(\gminib,b,\beta),(\gminia,a,\alpha)}$
be the decomposition such that
\begin{gather*}
 A_{(\gminib,b,\beta),(\gminia,a,\alpha)}\in
 \Hom\bigl(\nbigvtilde_{\gminia,a,\alpha},
 \nbigvtilde_{\gminib,b,\beta}\bigr)
 \otimes\Omega^1_{\Ytilde_P}.
\end{gather*}
Then,
we obtain the following for any $c\in\real$:
\begin{gather*}
 A_{(\gminib,b,\beta),(\gminia,a,\alpha)}\cdot
 \nbigp_{c}\nbigvtilde_{\gminia,a,\alpha}
 \subset \nbigp_{c-10m}\nbigvtilde_{\gminib,b,\beta}
 \otimes\Omega^1_{\Ytilde_P}\qquad
 (\gminia\neq\gminib),
 \\
 A_{(\gminia,b,\beta),(\gminia,a,\alpha)}\cdot
 \nbigp_{c}\nbigvtilde_{\gminia,a,\alpha}
 \subset \nbigp_{c}\nbigvtilde_{\gminia,b,\beta}
 \otimes\Omega^1_{\Ytilde_P}\big(\log \Htilde_P\big),\qquad
 (a,\alpha)\neq(b,\beta),
 \\
 \Res A_{(\gminia,a,\alpha),(\gminia,a,\alpha)}
 \bigl( \nbigp_{c}\nbigvtilde_{\gminia,a,\alpha} \bigr)
 \subset \nbigp_{<c}\nbigvtilde_{\gminia,a,\alpha}.
\end{gather*}
Here,
$\nbigp_{<c}\nbigvtilde_{\gminib,b,\beta}
=\bigcup_{d<c}\nbigp_d\nbigvtilde_{\gminib,b,\beta}$.
\end{itemize}
\end{Lemma}

\begin{proof}
By~\cite[Proposition~2.4.4]{Mochizuki-wild},
for any large inter $N$,
there exists a $G$-equivariant decomposition of
filtered bundles
\begin{gather*}
 \nbigp_{\ast}\nbigvtilde_{|\Ytilde_P}
 =\bigoplus_{\gminia\in\nbigi(\Ptilde)}
 \nbigp_{\ast}\nbigvtilde^{(N)}_{\gminia}
\end{gather*}
on $\big(\Ytilde_P,\Htilde_P\big)$
such that the following holds.
\begin{itemize}\itemsep=0pt
 \item Let $\pi^{(N)}_{\gminia}$ denote the projection of
 $\nbigp_{\ast}\nbigvtilde_{|\Ytilde_P}$
 onto $\nbigp_{\ast}\nbigvtilde^{(N)}_{\gminia}$,
 and let $\iota^{(N)}_{\gminia}$ denote the inclusion of
 $\nbigp_{\ast}\nbigvtilde^{(N)}_{\gminia}$
 into $\nbigp_{\ast}\nbigvtilde_{|\Ytilde_P}$.
 Then, for any $\gminia\neq\gminib$ and
 for any $c\in\real$, we obtain
\begin{gather*}
 \pi^{(N)}_{\gminib}\circ\DDlambdatilde\circ\iota^{(N)}_{\gminia}
 \bigl(
 \nbigp_{c}\nbigvtilde^{(N)}_{\gminia}\bigr)
 \subset
 \nbigp_{c-2N}\nbigvtilde^{(N)}_{\gminib}
 \otimes\Omega^1_{\Ytilde_P}.
\end{gather*}
\item For $\gminia$, we set
 $\DDtilde^{\lambda(N)}_{\gminia}:=
 \pi^{(N)}_{\gminia}\circ\DDlambdatilde\circ\iota^{(N)}_{\gminia}
 -{\rm d}\gminia\id_{\nbigvtilde^{(N)}_{\gminia}}$.
 Then, for any $c\in\real$, we obtain
\begin{gather*}
 \DDtilde^{\lambda\,(N)}
 \bigl(
 \nbigp_{c}\nbigvtilde^{(N)}_{\gminia}\bigr)
 \subset
 \nbigp_{c}\nbigvtilde^{(N)}_{\gminia}
	 \otimes\Omega^1_{\Ytilde_P}\big(\log\Htilde_P\big),
\\
 \DDtilde^{\lambda\,(N)}\circ\DDtilde^{\lambda\,(N)}
 \bigl(\nbigp_{c}\nbigvtilde^{(N)}_{\gminia}\bigr)
 \subset
 \nbigp_{c-N}\nbigvtilde^{(N)}_{\gminia}
 \otimes\Omega^2_{\Ytilde_P}\big(\log\Htilde_P\big).
\end{gather*}
\end{itemize}

Then, the claim of the lemma is clear.
\end{proof}

\subsection[Perturbation of good filtered lambda-flat bundles]
{Perturbation of good filtered $\boldsymbol\lambda$-flat bundles}

\subsubsection{Curve case}
\label{subsection;19.1.28.1}

Let $C$ be a Riemann surface
with a finite subset $D\subset C$.
Let $\big(\nbigp_{\ast}\nbigv,\DDlambda\big)$
be a good filtered $\lambda$-flat bundle on $(C,D)$. We~set $e:=\rank(\nbigv)!$. We~choose
$\eta>0$ such that
$10e\eta<\gaptilde(\nbigp_{\ast}\nbigv)$.

For any $0<10\rank(\nbigv)\epsilon\leq \eta$,
and for any $P\in D$,
let $\psi_{\epsilon,P}$ be a map
$\Partilde(\nbigp_{\ast}\nbigv,P)\lrarr \real$
such that
$(i)$ $|\psi_{\epsilon,P}(b)-b|<2\epsilon$,
$(ii)$ if $e(b_1-b_2)\in\seisuu$
then $\psi_{\epsilon,P}(b_1)-b_1=\psi_{\epsilon,P}(b_2)-b_2$. We~define
$\varphi_{\epsilon,P}\colon
 \seisuu\times\Partilde(\nbigp_{\ast}\nbigv,P)
\lrarr \real$
by
\begin{gather*}
 \varphi_{\epsilon,P}(k,b):=
 \psi_{\epsilon,P}(b)+\epsilon k.
\end{gather*}

We take $\veca\in\real^{D}$
as in Lemma~\ref{lem;20.2.12.20}.
For each $P\in D$ and
$b\in\Par(\nbigp_{\ast}\nbigv,\veca,P)$,
we obtain the endomorphism
$\Res_P\big(\DDlambda\big)$
of $\Gr^F_b\big(\nbigp_{\veca}\nbigv_{|P}\big)$.
Let $W_{\bullet}
 \Gr^F_b\big(\nbigp_{\veca}\nbigv_{|P}\big)$
denote the weight filtration
associated with
the nilpotent part of $\Res_P(\DDlambda)$.
For any $(k,b)\in\seisuu\times\Par(\nbigp_{\ast}\nbigv,\veca,P)$,
we obtain the subspace
$W_k\bigl(
 F_b\big(\nbigp_{\veca}\nbigv_{|P}\big)
 \bigr)$
as the pull back of
$W_k\Gr^F_b\big(\nbigp_{\veca}\nbigv_{|P}\big)$
by the projection
$F_b\big(\nbigp_{\veca}\nbigv_{|P}\big)\lrarr
 \Gr^F_b\big(\nbigp_{\veca}\nbigv_{|P}\big)$. We~define the filtration $F^{(\epsilon)}$
on $\nbigp_{\veca}\nbigv_{|P}$
indexed by $\openclosed{a(P)-1}{a(P)}$
as follows:
\begin{gather*}
 F^{(\epsilon)}_c\bigl(
 \nbigp_{\veca}\nbigv_{|P}
 \bigr):=
 \sum_{\substack{
 (k,b)\in \seisuu\times\Par(\nbigp_{\ast}\nbigv,\veca,P)
\\
 \varphi_{\epsilon,P}(k,b)\leq c
}}
 W_kF_b\big(\nbigp_{\veca}\nbigv_{|P}\big).
\end{gather*}
We obtain the corresponding filtered bundle
$\nbigp^{(\epsilon)}_{\ast}\nbigv$.
Note the following lemma.
\begin{Lemma}
\label{lem;20.2.12.30}
$\big(\nbigp^{(\epsilon)}_{\ast}\nbigv,\DDlambda\big)$ is
a good filtered $\lambda$-flat bundle.
\end{Lemma}

\begin{proof}
It~is enough to consider the case
where $C$ is a neighbourhood of $D=\{P\}=\{0\}$ in~$\cnum$. We~obtain
$\varphi\colon \Ctilde\lrarr C$,
$\Dtilde$,
$G$,
and $\big(\nbigp_{\ast}\nbigvtilde,\DDtilde^{\lambda}\big)$
as in Section~\ref{subsection;20.2.14.11}. We~set $\atilde:=ea(P)$
and $\Ptilde:=\varphi^{-1}(P)$.
For $(k,b)\in\seisuu\times\Par(\nbigp_{\ast}\nbigvtilde,\atilde)$,
we construct
$W_kF_b\big(\nbigp_{\atilde}\nbigvtilde_{|\Ptilde}\big)$
as above.

For $b\in\Par\big(\nbigp_{\ast}\nbigvtilde,\atilde\big)$,
note that
$b/e\in\Partilde(\nbigp_{\ast}\nbigv)$,
and we set
\begin{gather*}
 \varphitilde_{\epsilon,P}(k,b):=
e\psi_{\epsilon,P}(k,b/e)+ke\epsilon.
\end{gather*}
We set
\begin{gather*}
 F^{(\epsilon)}_d\big(\nbigp_{\atilde}\nbigvtilde_{|\Ptilde}\big):=
 \sum_{\varphitilde_{\epsilon,P}(k,b)\leq d}
 W_kF_b\big(\nbigp_{\atilde}\nbigvtilde_{|\Ptilde}\big).
\end{gather*}
We obtain the corresponding
$G$-equivariant filtered bundle
$\nbigp^{(\epsilon)}_{\ast}\nbigvtilde$. We~can easily observe that
$\varphi^{\ast}\nbigp_{\ast}^{(\epsilon)}\nbigv
\simeq
 \nbigp_{\ast}^{(\epsilon)}\nbigvtilde$
by using
$10^2e\rank(\nbigv)\epsilon<\gaptilde(\nbigp_{\ast}\nbigv)$.

There exists the decomposition
\begin{gather*}
\big (\nbigp_{\ast}\nbigvtilde,\DDtilde^{\lambda}\big) \otimes\cnum[\![\zeta]\!]
=\bigoplus_{\gminia\in\nbigi}
 \bigl( \nbigp_{\ast}\nbigvtilde_{\gminia}, \DDtilde^{\lambda}_{\gminia} \bigr),
\end{gather*}
as in (\ref{eq;20.2.1.1}). We~apply the same procedure to
each
$\nbigp_{\ast}\nbigvtilde_{\gminia}$
by using $\varphitilde_{\epsilon,P}$,
and we obtain filtered bundles
$\nbigp^{(\epsilon)}_{\ast}\nbigvtilde_{\gminia}$
for which $\DDtilde^{\lambda}_{\gminia}-{\rm d}\gminia\id$
are logarithmic.
Because
$\nbigp^{(\epsilon)}_{\ast}\nbigvtilde\otimes\cnum[\![\zeta]\!]
=\bigoplus
 \nbigp^{(\epsilon)}_{\ast}\nbigvtilde_{\gminia}$,
we obtain that
$\big(\nbigp^{(\epsilon)}_{\ast}\nbigvtilde,\DDtilde^{\lambda}\big)$
is unramifiedly good.
Hence, we obtain that
$\big(\nbigp^{(\epsilon)}_{\ast}\nbigv,\DD^{\lambda}\big)$
is good.
\end{proof}

Suppose that $C$ is compact and connected. We~clearly obtain
$\lim\limits_{\epsilon\to 0}
 c_1\big(\nbigp^{(\epsilon)}_{\ast}\nbigv\big)=
 c_1(\nbigp_{\ast}\nbigv)$.
The following is also standard.
\begin{Lemma}
\label{lem;19.1.28.30}
Suppose that $\big(\nbigp_{\ast}\nbigv,\DDlambda\big)$ is stable.
Then, if $\epsilon>0$ is sufficiently small,
$\big(\nbigp^{(\epsilon)}_{\ast}\nbigv,\DDlambda\big)$
is also stable.
\end{Lemma}

\begin{proof}
For any positive integer $s\leq \rank(\nbigv)$,
and for any $P\in D$,
let $\nbigt(s,P)$ be the set of~real numbers
expressed as
\begin{gather*}
 \sum_{b\in \Par(\nbigp_{\ast}\nbigv,\veca,P)}
 b\,s_b,
\end{gather*}
where $s_b$ are non-negative integers such that
$\sum s_b=s$.
Let $\nbigs(s)$ denote the set of
real numbers expressed as
$\frac{1}{s}\bigl( m+\sum_{P\in D} c_P\bigr)$,
where $m\in\seisuu$ and $c_P\in\nbigt(s,P)$.
Then,
$\bigcup_{0<s\leq\rank\nbigv}\nbigs(s)$
is discrete in $\real$.
Hence, there exists $\delta>0$
such that
if $t\in\bigcup_{0<s\leq\rank\nbigv}\nbigs(s)$
satisfies $t<\mu(\nbigp_{\ast}\nbigv)$,
we obtain
$t<\mu(\nbigp_{\ast}\nbigv)-\delta$.
Then, the claim of the lemma is clear.
\end{proof}

\subsubsection{Surface case}
\label{subsection;19.1.30.101}

Let $X$ be a complex projective surface
with a simple normal crossing hypersurface
$H=\bigcup_{i\in\Lambda}H_i$.
Let $\big(\nbigp_{\ast}\nbigv,\DDlambda\big)$ be
a good filtered $\lambda$-flat bundle
on $(X,H)$. We~shall explain a similar perturbation of
good filtered $\lambda$-flat bundles.
Set $e:=\rank(\nbigv)!$. We~choose $\eta>0$ such that
$0<10e\eta<\gaptilde(\nbigp_{\ast}\nbigv)$.

For any $0<10\rank(\nbigv)\epsilon\leq\eta$,
let $\psi_{\epsilon,i}$ be a map
$\Partilde(\nbigp_{\ast}\nbigv,i)\lrarr \real$
such that
$(i)$~$|\psi_{\epsilon,i}(b)-b|\allowbreak<2\epsilon$,
$(ii)$~if $e(b_1-b_2)\in\seisuu$
then $\psi_{\epsilon,i}(b_1)-b_1=\psi_{\epsilon,i}(b_2)-b_2$. We~define
$\varphi_{\epsilon,i}\colon
 \seisuu\times\Par(\nbigp_{\ast}\nbigv,i)
\lrarr \real$
by
$\varphi_{\epsilon,i}(k,b):=
 \psi_{\epsilon,i}(b)+\epsilon k$.

We take
$\veca\in\real^{\Lambda}$
for $\Partilde(\nbigp_{\ast}\nbigv,i)$ $(i\in\Lambda)$
as in Lemma~\ref{lem;20.2.12.20}.
The eigenvalues of
the endomorphism $\Res_i\big(\DDlambda\big)$
on $\Gr^{F}_b\big(\nbigp_{\veca}\nbigv_{|H_i}\big)$
are constant on $H_i$
because $H_i$ are compact. We~obtain the well defined
nilpotent part $N_{i,b}$ of
$\Res_i\big(\DDlambda\big)$.
There exists a finite subset $Z_i\subset H_i$
such that the conjugacy classes
of the nilpotent part of $N_{i,b|Q}$
$(Q\in H_i\setminus Z_i)$ are constant. We~obtain the weight filtration $W$
of $\Gr^F_b\big(\nbigp_{\veca}\nbigv_{|H_i\setminus Z_i}\big)$
by algebraic vector subbundles
whose restriction to~$Q\in H_i\setminus Z_i$
are the weight filtration of
$N_{i,b|Q}$.
It~uniquely extends to a filtration of
$\Gr^F_b\big(\nbigp_{\veca}\nbigv_{|H_i}\big)$
by algebraic subbundles,
which is also denoted by $W$.

For any $(k,b)\in\seisuu\times\Par(\nbigp_{\ast}\nbigv,\veca,i)$,
let
$W_kF_b\big(\nbigp_{\veca}\nbigv_{|H_i}\big)$
denote the subbundle
of $\nbigp_{\veca}\nbigv_{|H_i}$
obtai\-ned as the pull back of
$W_k\Gr^F_b\big(\nbigp_{\veca}\nbigv_{|H_i}\big)$
by the projection
$F_b\big(\nbigp_{\veca}\nbigv_{|H_i}\big)\lrarr
 \Gr^F_b\big(\nbigp_{\veca}\nbigv_{|H_i}\big)$. We~define the filtration $F^{(\epsilon)}$
on $\nbigp_{\veca}\nbigv_{|H_i}$
indexed by $\openclosed{a_i-1}{a_i}$
as follows:
\begin{gather*}
 F^{(\epsilon)}_c\nbigp_{\veca}\nbigv_{|H_i}:=
 \sum_{\substack{
 (k,b)\in \seisuu\times\Par(\nbigp_{\ast}\nbigv,\veca,i)
\\
 \varphi_{\epsilon,i}(k,b)\leq c
}}
 W_kF_b\nbigp_{\veca}\nbigv_{|H_i}.
\end{gather*}
We obtain the corresponding filtered bundle
$\nbigp^{(\epsilon)}_{\ast}\nbigv$ over $\nbigv$.
As in the curve case (see Lemma~\ref{lem;20.2.12.30}),
we obtain the following.
\begin{Lemma}
$\big(\nbigp^{(\epsilon)}_{\ast}\nbigv,\DDlambda\big)$
is a good filtered $\lambda$-flat bundle.
\end{Lemma}

We clearly have
$\lim\limits_{\epsilon\to 0}
 c_1\big(\nbigp^{(\epsilon)}_{\ast}\nbigv\big)=
 c_1(\nbigp_{\ast}\nbigv)$
and
$\lim\limits_{\epsilon\to 0}
 \ch_2\big(\nbigp^{(\epsilon)}_{\ast}\nbigv\big)=
 \ch_2(\nbigp_{\ast}\nbigv)$.
The following is standard,
and similar to Lemma~\ref{lem;19.1.28.30}.
(See also~\cite[Proposition~3.28]{mochi4}.)
\begin{Lemma}
\label{lem;21.5.5.1}
 Let $L$ be an ample line bundle on $X$.
Suppose that $\big(\nbigp_{\ast}\nbigv,\DDlambda\big)$ is
$\mu_L$-stable.
Then, if $\epsilon$ is sufficiently small,
$\big(\nbigp^{(\epsilon)}_{\ast}\nbigv,\DDlambda\big)$
is also $\mu_L$-stable.
\end{Lemma}

\subsection{Some families of auxiliary metrics
on a punctured disc}

\subsubsection{Regular model case}
\label{subsection;20.2.6.1}

Let $V$ be a finite dimensional vector space over $\cnum$
with a nilpotent endomorphism $f$.
Let $(a,\alpha)\in\real\times\cnum$.
Let $X$ be a neighbourhood of $0$ in $\cnum$. We~set $H:=\{0\}$. We~set $\nbigv=\nbigo_X(\ast H)\otimes V$.
From $(a,\alpha)$ and $(V,f)$,
we obtain a model filtered $\lambda$-flat bundle
$\big(\nbigp_{\ast}\nbigv,\DDlambda\big)$
by applying the construction
in Section~\ref{subsection;20.2.9.1}
in the case $\gminia=0$ and $e=1$.

Fix $0<\eta<1$.
For any $0\leq 10\rank(\nbigv)\epsilon\leq \eta$,
we take $a(\epsilon)$
such that $|a-a(\epsilon)|<2\epsilon$,
and we obtain the regular filtered $\lambda$-flat bundle
$\big(\nbigp^{(\epsilon)}_{\ast}\nbigv,\DDlambda\big)$
as in Section~\ref{subsection;19.1.28.1}. We~set $E:=\nbigv_{|X\setminus H}$. We~consider the K\"ahler metric $g_{\epsilon}:=
 \big(\eta^2|z|^{2\eta-2}+\epsilon^2|z|^{2\epsilon-2}\big)
 \,{\rm d}z\,{\rm d}\zbar$
of $X\setminus H$.

\begin{Proposition}
\label{prop;20.2.6.2}
There exist Hermitian metrics $h^{(\epsilon)}$
of $E$
for $0\leq 10\rank(\nbigv)\epsilon\leq\eta$
such that the following holds:
\begin{itemize}\itemsep=0pt
\item
$\nbigp^{h^{(\epsilon)}}_{\ast}E=\nbigp^{(\epsilon)}_{\ast}\nbigv$.
\item
$h^{(0)}$ is a Hermitian metric of $\big(E,\DDlambda\big)$,
and $h^{(\epsilon)}\to h^{(0)}$ in the $C^{\infty}$-sense
locally on $X\setminus H$.
\item
There exist $C_i>0$ $(i=0,1,2)$
such that
the following conditions are satisfied for any $\epsilon$:
\begin{gather*}
 \bigl| G\big(h^{(\epsilon)}\big) \bigr|_{g_{\epsilon},h^{(\epsilon)}}\leq C_0,\qquad
 C_1^{-1}|z|^{2\rank(\nbigv)\epsilon} h^{(0)}\leq
 h^{(\epsilon)} \leq C_1|z|^{-2\rank(\nbigv)\epsilon}h^{(0)},
\\
 C_2^{-1}\det(h_0)\leq \det\big(h^{(\epsilon)}\big)\leq C_2\det\big(h^{(0)}\big).
\end{gather*}
\item
Let $B_i^{(\epsilon)}$ $(i=1,2)$
be the $C^{\infty}$-endomorphisms of $E$
determined by the condition
$\DD^{\lambda\star}_{h^{(\epsilon)}}(v)
=
B^{(\epsilon)}_{1}(v)\,{\rm d}z/z+
B^{(\epsilon)}_2(v)\,{\rm d}\zbar/\zbar$
for $v\in V$.
Then,
there exists $C_3>0$ such that
$\big|B^{(\epsilon)}_i\big|_{h^{(\epsilon)}}\leq C_3$
holds for any $\epsilon$.
\end{itemize}
Moreover,
for any $v_1,v_2\in V$,
$h^{(\epsilon)}(v_1,v_2)$
depends only on $|z|$,
where $v_j$ are naturally regarded as holomorphic sections
of $\nbigv$.
\end{Proposition}
In the case $\lambda\neq 0$,
such a family of Hermitian metrics $h^{(\epsilon)}$
is constructed in
\cite[Sections~4.3 and~4.4.1]{mochi5}. We~shall explain the case $\lambda=0$
in Section~\ref{subsection;19.2.11.30}.

\subsubsection{General case}

Let $S_0$ be a finite set.
Let $S_1$ be a finite subset of
$\openclosed{0}{1}$.
Let $S_2$ be a finite subset of $\cnum$.
Let $V_{i,a,\alpha}$
($(i,a,\alpha)\in S_0\times S_1\times S_2$)
be finite dimensional $\cnum$-vector spaces
equipped with a nilpotent endomorphism $f_{i,a,\alpha}$.
Set $r:=\sum \dim V_{i,a,\alpha}$ and $e:=r!$.
Take $\eta>0$ such that
\begin{gather*}
10e\eta< \min\bigl\{ |a-b|\mid a,b\in S_1\cup\{0,1\},\, a\neq b \bigr\}.
\end{gather*}

As in Section~\ref{subsection;20.2.6.1},
we obtain the regular filtered $\lambda$-flat bundles
$\big(\nbigp_{\ast}\nbigv_{i,a,\alpha},\DDlambda_{i,a,\alpha}\big)$
from $(V_{i,a,\alpha},f_{i,a,\alpha})$.
For $0\leq \epsilon<\eta/10r$,
we take $a(\epsilon)$ such that
$|a-a(\epsilon)|<2\epsilon$,
and we obtain
$\big(\nbigp^{(\epsilon)}_{\ast}\nbigv_{i,a,\alpha},\DDlambda_{i,a,\alpha}\big)$. We~set
$E_{i,a,\alpha}:=\nbigv_{i,a,\alpha|X\setminus H}$. We~obtain the metrics
$h^{(\epsilon)}_{i,a,\alpha}$ of $E_{i,a,\alpha}$
as in~Pro\-po\-sition~\ref{prop;20.2.6.2}. We~set
$\nbigp^{(\epsilon)}_{\ast}\nbigv:=
 \bigoplus\nbigp^{(\epsilon)}_{\ast}\nbigv_{i,a,\alpha}$
and
$h^{(\epsilon)}:=\bigoplus h^{(\epsilon)}_{i,a,\alpha}$.

Fix a positive integer $m$ and
a positive number $C$. We~consider the following data:
\begin{itemize}\itemsep=0pt
\item
For each $i\in S_0$,
let $\gminia(i)$ denote a polynomial
$\sum_{j=1}^m \gminia(i)_jz^{-j}
 \in z^{-1}\cnum\big[z^{-1}\big]$
such that \mbox{$|\gminia(i)_j|\leq C$}.
\item
Let $A$ be a holomorphic section of
$\nend(\nbigv)$
with the decomposition
$A=\sum A_{(j,b,\beta),(i,a,\alpha)}$,
where
\begin{gather*}
A_{(j,b,\beta),(i,a,\alpha)}
 \in\nhom(\nbigv_{i,a,\alpha},\nbigv_{j,b,\beta})
\end{gather*}
If $i\neq j$,
 we obtain $|A_{(j,b,\beta),(i,a,\alpha)}|_{h^{(0)}}\leq C|z|^{10m}$,
and if $i=j$,
we obtain $|A_{(i,b,\beta),(i,a,\alpha)}|_{h^{(0)}}\allowbreak\leq C|z|^{4\eta}$.
\end{itemize}

We define the flat $\lambda$-connection
$\DDlambda$ on $\nbigv$
as follows:
\begin{gather*}
 \DDlambda=
 \bigoplus_{i,a,\alpha}
 \bigl( {\rm d}\gminia(i)\id_{\nbigv_{i,a,\alpha}} +\DDlambda_{i,a,\alpha} \bigr)
+A\frac{{\rm d}z}{z}.
\end{gather*}
Let $g_{\epsilon}$
be the K\"ahler metrics of $X\setminus H$
as in Section~\ref{subsection;20.2.6.1}.

\begin{Proposition}
\label{prop;20.2.6.10}
There exists a constant $C'$ depending only on $m$ and $C$
such that
\begin{gather*}
\big|G\big(h^{(\epsilon)}\big)\big|_{g_{\epsilon},h^{(\epsilon)}}
\leq C'.
\end{gather*}
\end{Proposition}

\subsubsection{A consequence}

Let $Y$ be a neighbourhood of $0$ in $\cnum$. We~set $H:=\{0\}$.
Let $\big(\nbigp_{\ast}\nbigv,\DDlambda\big)$
be any good filtered $\lambda$-flat bundle
on $(Y,H)$.
Let $\big(E,\DDlambda\big)$ be the $\lambda$-flat bundle
obtained as the restriction of
$\big(\nbigp_{\ast}\nbigv,\DDlambda\big)$
to~$Y\setminus H$. We~set $e:=\rank(\nbigv)!$.
Let $\varphi$ be as in Section~\ref{subsection;20.2.6.20},
and we set
$\big(\nbigp_{\ast}\nbigvtilde,\DDtilde^{\lambda}\big):=
 \varphi^{\ast}\big(\nbigp_{\ast}\nbigv,\DDlambda\big)$. We~take $\eta>0$ such that
$10e\eta<
 \gaptilde(\nbigp_{\ast}\nbigv)$.
Let $g$ be the K\"ahler metric
$\eta^2|z|^{2\eta-2}\,{\rm d}z\,{\rm d}\zbar$
of $Y\setminus H$. By~using a special case of Proposition~\ref{prop;20.2.6.10},
we obtain the following corollary.

\begin{Corollary}
\label{cor;19.1.20.2}
There exists a Hermitian metric
$h$ of $\big(E,\DDlambda\big)$
such that
$(i)$~$\nbigp^{h}_{\ast}E=
 \nbigp_{\ast}\nbigv$,
$(ii)$~$|G(h)|_{h,g}$
is bounded on $Y\setminus H$.
\end{Corollary}

\begin{proof}
Let
$\big(\nbigp_{\ast}\nbigv_0,\DD_0^{\lambda}\big)$
be a model filtered $\lambda$-flat bundle
with an isomorphism
$\Phi_m\colon \nbigp_{\ast}\nbigv_0\allowbreak\simeq\nbigp_{\ast}\nbigv$
as in Lemma~\ref{lem;20.2.14.20},
where $m$ is a sufficiently large integer. We~recall that
$\big(\nbigp_{\ast}\nbigv_0,\DD_0^{\lambda}\big)$
is obtained as the descent of
the $G$-equivariant model filtered $\lambda$-flat bundle
$\big(\nbigp_{\ast}\nbigvtilde_0,\DDtilde^{\lambda}_0\big)
=\bigoplus_{\gminia,a,\alpha}
 \big(\nbigp_{\ast}\nbigv_{\gminia,a,\alpha}, \DDlambda_{\gminia,a,\alpha}\big)$,
and $\Phi_m$ is induced by
a $G$-equivariant isomorphism
$\nbigp_{\ast}\nbigvtilde_0\simeq \nbigp_{\ast}\nbigvtilde$.
Let~$\htilde_0^{(0)}$ be the Hermitian metric of
$\nbigvtilde_{0|\Ytilde\setminus \Htilde}$
as in Proposition~\ref{prop;20.2.6.2}
with $\epsilon=0$. By~the isomorphism $\Phitilde_m$,
it induces a $G$-equivariant Hermitian metric
$\htilde$
of $\nbigvtilde_{|\Ytilde\setminus\Htilde}$.
Applying Proposition~\ref{prop;20.2.6.10}
to $\big(\nbigp_{\ast}\nbigvtilde,\DDtilde^{\lambda}\big)$
with $\htilde$,
we obtain the boundedness of $G\big(\htilde\big)$
with respect to
$\varphi^{\ast}g$
and $\htilde$.
Because $\htilde$ is $G$-equivariant,
we obtain the Hermitian metric $h$ of $E$
which has the desired property.
\end{proof}

\subsubsection{Proof of Proposition~\ref{prop;20.2.6.2}}
\label{subsection;19.2.11.30}

Let $V_2:=\cnum\,v_1\oplus\cnum v_2$
with the nilpotent map
determined by $f_2(v_1)=v_2$ and $f_2(v_2)=0$. We~obtain
$\nbigv_2=V_2\otimes\nbigo_X(\ast H)$
and
the Higgs field $\theta_2$
determined by
$\theta_2(v)=f_2(v)\,{\rm d}z/z$.
Let $\big(E_2,\delbar_{E_2},\theta_2\big)$ be the Higgs bundle
on $X\setminus H$
obtained as the restriction of
$(\nbigv_2,\theta_2)$.
For any $\epsilon>0$,
we set
$L_{\epsilon}(z):=\epsilon^{-1}(|z|^{-\epsilon}-|z|^{\epsilon})$. We~also set
$L_0:=-\log|z|^{2}$.

\begin{Lemma}
We obtain
$L_0(z)\leq L_{\epsilon}(z)\leq |z|^{-\epsilon}L_{0}(z)$.
There exists $C>0$ such that
\begin{gather*}
 \bigl| \del_z \log L_{\epsilon}(z)\bigr| \leq C|z|^{-1}
\end{gather*}
on $\{|z|\leq 1/2\}$ for any $0\leq \epsilon\leq 1/2$.
\end{Lemma}

\begin{proof}
As proved in~\cite[Section~4.2]{mochi5},
$L_0(z)\leq L_{\epsilon}(z)$ holds. We~set
$g_1(\epsilon):=-\epsilon \log|z|^2-\big(1-|z|^{2\epsilon}\big)$
for any $z\in \Delta^{\ast}$ and for $\epsilon>0$.
It~is easy to check that
$\del_{\epsilon}g_1(\epsilon)\geq 0$
and $\lim\limits_{\epsilon \to 0}g_1(\epsilon)=0$.
Hence, we obtain
$L_{\epsilon}(z)\leq |z|^{-\epsilon}L_0(z)$.
For $0<a<1$ and $0<\epsilon$, we set
$g_2(\epsilon,a):= \frac{\epsilon}{2}(a^{-\epsilon}-a^{\epsilon})^{-1}
 (a^{-\epsilon}+a^{\epsilon})$.
Then,
$\del\log L_{\epsilon}(z)
=-g_2(\epsilon,|z|)\frac{{\rm d}z}{z}$.
Then, we can check that
$\del_ag_2(\epsilon,a)\geq 0$
and $\del_{\epsilon}g_2(\epsilon,a)\geq 0$.
Then, we obtain the second claim of the lemma.
\end{proof}

Let $h_2^{(\epsilon)}$
be the $C^{\infty}$-metric of $E_2$
given by
\begin{gather*}
 h_2^{(\epsilon)}(v_1,v_1)=L_{\epsilon},\qquad
 h_2^{(\epsilon)}(v_2,v_2)=L_{\epsilon}^{-1},\qquad
 h_2^{(\epsilon)}(v_1,v_2)=0.
\end{gather*}

\begin{Lemma}
$\big(E_2,\delbar_{E_2},\theta_2,h_2^{(\epsilon)}\big)$
are harmonic bundles.
Moreover,
the family of metrics $h^{(\epsilon)}_2$
satisfies the condition
in Proposition~$\ref{prop;20.2.6.2}$
for $(\nbigp_{\ast}\nbigv_2,\theta_2)$.
\end{Lemma}

\begin{proof}
Let $H_{\epsilon}$ be the matrix valued function
on $X\setminus D$
determined by
$(H_{\epsilon})_{i,j}:=h_2^{(\epsilon)}(v_i,v_j)$.
Then, the following holds:
\begin{align*}
 \delbar \bigl( H_{\epsilon}^{-1} \del H_{\epsilon} \bigr)
&= \begin{pmatrix}
 \delbar\del\log L_{\epsilon} & 0 \\
 0 & -\delbar\del\log L_{\epsilon}
 \end{pmatrix}
=\begin{pmatrix}
 -1 & 0 \\
 0 & 1
 \end{pmatrix}
\frac{\epsilon^2|z|^{-2}\,{\rm d}\zbar\,{\rm d}z}
 {\bigl(|z|^{-\epsilon}-|z|^{\epsilon}\bigr)^2}
 \\
&=-\begin{pmatrix}
 -1 & 0 \\
 0 & 1
 \end{pmatrix}
 L_{\epsilon}^{-2}\frac{{\rm d}z\,{\rm d}\zbar}{|z|^{2}}.
\end{align*}
Let $\Theta$ be the matrix valued function
representing $\theta_2$ with respect to
the frame $(v_1,v_2)$,
i.e.,
$\theta_2(v_1,v_2)=(v_1,v_2)\Theta$.
Let $\theta^{\dagger}_{2,h^{(\epsilon)}}$
denote the adjoint of $\theta_2$
with respect to $h_2^{(\epsilon)}$.
Let $\Theta^{\dagger}_{\epsilon}$
denote the matrix valued function
representing $\theta^{\dagger}_{2,h^{(\epsilon)}}$.
The following holds:
\begin{gather*}
 \Theta=\begin{pmatrix}
 0 & 0 \\ 1 & 0
 \end{pmatrix}
 \frac{{\rm d}z}{z},\qquad
 \Theta^{\dagger}_{\epsilon}
 =\begin{pmatrix}
 0 & 1 \\ 0 & 0
 \end{pmatrix}
 L_{\epsilon}^{-2}\frac{{\rm d}\zbar}{\zbar}.
\end{gather*}
Hence, we obtain
\begin{gather*}
 \bigl[ \Theta,\Theta^{\dagger}_{\epsilon} \bigr]
=\begin{pmatrix}
 -1 & 0 \\ 0 & 1
 \end{pmatrix}\cdot
 L_{\epsilon}^{-2} \cdot
 \frac{{\rm d}z\,{\rm d}\zbar}{|z|^2}.
\end{gather*}
It~implies that
\begin{gather*}
 \delbar\bigl(
 H_{\epsilon}^{-1}\del H_{\epsilon}
 \bigr)
+\bigl[
 \Theta,\Theta^{\dagger}_{\epsilon}
 \bigr]
=0.
\end{gather*}
It~is exactly the Hitchin equation
for $\big(E_2,\delbar_{E_2},\theta_2,h_{2,\epsilon}\big)$.
The other claim is easy to see.
\end{proof}

For each $\ell\in\seisuu_{>0}$,
we set
$(V_{\ell},f_{\ell}):=\Sym^{\ell-1}(V_2,f_2)$. We~set
$\nbigv_{\ell}:=V_{\ell}\otimes\nbigo_X(\ast H)$
and $\theta_{\ell}:=f_{\ell}\,{\rm d}z/z$. We~obtain the regular filtered Higgs bundles
$\big(\nbigp^{(\epsilon)}_{\ast}\nbigv_{\ell},\theta_{\ell}\big)$.
Note that
$\big(\nbigp^{(\epsilon)}_{\ast}\nbigv_{\ell},\theta_{\ell}\big)$
is naturally isomorphic to
the $(\ell-1)$-th symmetric product
of $\big(\nbigp^{(\epsilon)}_{\ast}\nbigv_{2},\theta_2\big)$.
Hence, $h^{(\epsilon)}_{2}$
induce
harmonic metrics
$h^{(\epsilon)}_{\ell}$
of
$\big(E_{\ell},\delbar_{E_{\ell}},\theta_{\ell}\big):=
 (\nbigv_{\ell},\theta_{\ell})_{|X\setminus H}$
satisfying the conditions in Proposition~\ref{prop;20.2.6.2}
for~$\big(\nbigp^{(\epsilon)}_{\ast}\nbigv_{\ell},\theta_{\ell}\big)$.

Let $(V,f)$, $(a,\alpha)$
and $(\nbigv,\theta)$
be as in Section~\ref{subsection;20.2.6.1}.
There exist integers $\ell_1,\ldots,\ell_n$
such that
\begin{gather*}
 (V,f)\simeq\bigoplus_{j=1}^n \big(V_{\ell_j},f_{\ell_j}\big).
\end{gather*}
We obtain
$\nbigp^{(\epsilon)}_{\ast}\nbigv
\simeq
 \bigoplus\nbigp^{(\epsilon)}_{\ast}\nbigv_{\ell_j}$. We~obtain the harmonic metrics
$\bigoplus_{j=1}^n |z|^{-2a(\epsilon)}h^{(\epsilon)}_{\ell_j}$
for
$\big(E,\delbar_E,\theta\big)$. We~can easily check that
they satisfy the conditions in Proposition~\ref{prop;20.2.6.2}
for $(\nbigp_{\ast}\nbigv,\theta)$.
\hfill\qed

\subsubsection{Proof of Proposition~\ref{prop;20.2.6.10}}

We set
$\Phi:=\bigoplus_{i,a,\alpha} {\rm d}\gminia(i)\id_{\nbigv_{i,a,\alpha}}$
and
$\DD^{\lambda\reg}:= \bigoplus \DD^{\lambda}_{i,a,\alpha}$.
Let $\Phi^{\dagger}_{h^{(\epsilon)}}$
and $A^{\dagger}_{h^{(\epsilon)}}$
denote the adjoint of~$\Phi$ and $A$
with respect to $h^{(\epsilon)}$,
respectively. We~obtain the decompositions
$\DDlambda=\DD^{\lambda\reg}+\Phi+A\,{\rm d}z/z$
and
$\DD_{h^{(\epsilon)}}^{\lambda\star}=
 \DD_{h^{(\epsilon)}}^{\lambda\reg\star}
-\Phi_{h^{(\epsilon)}}^{\dagger}
-A_{h^{(\epsilon)}}^{\dagger}\,{\rm d}\zbar/\zbar$.
Note that
$\bigl[ \DD^{\lambda\reg},\Phi_{h^{(\epsilon)}}^{\dagger}\bigr]
=\bigl[ \DD^{\lambda\reg\star}_{h^{(\epsilon)}}, \Phi \bigr]
=\bigl[ \Phi,\Phi^{\dagger}_{h^{(\epsilon)}} \bigr]=0$. By~the assumption,
we obtain
\begin{gather*}
 \bigl| \bigl[ \Phi,A^{\dagger}_{h^{(\epsilon)}} \bigr] \bigr|_{g_{\epsilon},h^{(\epsilon)}}
= \bigl| \bigl[ \Phi^{\dagger}_{h^{(\epsilon)}},A \bigr] \bigr|_{g_{\epsilon},h^{(\epsilon)}}
\leq 2\eta^{-2}C^2C_1^2|z|^{5m-4\rank(\nbigv)\epsilon-2\eta}.
\end{gather*}
We also obtain
$\bigl|\big[A,A^{\dagger}_{h^{(\epsilon)}}\big] \bigr|_{h^{(\epsilon)}}
\leq 2C^2C_1^2|z|^{4\eta-4\rank(\nbigv)\epsilon}$.
Because
\begin{gather*}
 \bigl[ \DD^{\lambda\reg\star}_{h^{(\epsilon)}},A\,{\rm d}z/z \bigr]
=\bigl[ B^{(\epsilon)}_2, A\bigr]\,|z|^{-2}{\rm d}\zbar\,{\rm d}z,
\end{gather*}
we obtain
\begin{gather*}
 \bigl| \bigl[ \DD^{\lambda\reg\star}_{h^{(\epsilon)}},A\,{\rm d}z/z \bigr]
\bigr|_{g_{\epsilon},h^{(\epsilon)}}
= \Bigl| \bigl[ \DD^{\lambda\reg},A^{\dagger}_{h^{(\epsilon)}}\,{\rm d}\zbar/\zbar \bigr]
\bigr|_{g_{\epsilon},h^{(\epsilon)}}
\leq C_3C|z|^{2\eta}.
\end{gather*}
Hence, we obtain the desired estimate for
$G\big(h^{(\epsilon)}\big)= \bigl[ \DDlambda,\DD^{\lambda\star}_{h^{(\epsilon)}} \bigr]$.
\hfill\qed

\subsection[Estimate of the curvature for Hermitian--Einstein metrics
of a Higgs bundle]
{Estimate of the curvature for Hermitian--Einstein metrics\\
of a Higgs bundle}

Let $X$ be a complex surface.
Let $\big(E,\delbar_E,\theta\big)$ be a Higgs bundle on $X$.
Let $g_i$ be a sequence of~K\"ahler metric on $X$
which is convergent to a K\"ahler metric $g_{\infty}$
in the $C^{\infty}$-sense locally on~$X$.
Let~$h_i$ $(i=1,2,\ldots)$ be Hermitian--Einstein metrics
of the Higgs bundle. We~assume the following:
\begin{itemize}\itemsep=0pt
\item
 $\int_X|G(h_i)|^2_{h_i,g_i}\to 0$.
\end{itemize}

Let $\nabla_{h_i}$ be the Chern connection
of $\big(E,\delbar_E,h_i\big)$.
Let $R(h_i)$ denote the curvature of $\nabla_{h_i}$.
The~following proposition is a refinement
contained in the argument
in~\cite[Section~9.1.1]{mochi4}.

\begin{Proposition}
\label{prop;20.1.30.1}
For any relatively compact open subset $U\subset X$,
and for any $p\geq 1$,
the $L_2^p$-norms of $R(h_i)_{|U}$
with respect to $h_i$ and $g_{\infty}$
are bounded.
\end{Proposition}
\begin{proof}
Let $P$ be any point of $X$.
Let $(X_P,z_1,z_2)$ be a holomorphic coordinate neighbourhood
of $X$ around $P$.
Let us describe $\theta$
as $\theta=\sum_{j=1,2} f_j\,{\rm d}z_j$.
Let $X_{P,0}$ be a relatively compact neighbourhood of $P$ in $X_P$.
According to~\cite[Lemma 2.13]{mochi4},
there exist $C_k>0$ $(k=1,2)$,
which are independent of $h_i$,
such that the following inequalities hold on $X_{P,0}$:
\begin{gather*}
|f_j|_{h_i}^2\leq C_1\cdot
 \exp\bigg( C_2\int_{X}|G(h_i)|_{h_i}^2 \bigg).
\end{gather*}
(Note that $G(h_i)$ is denoted as $F(h_i)$
in~\cite{mochi4}.)
Let $\epsilon$ denote a positive number.
After rescaling the coordinate system,
we may assume the following on $X_{P,0}$:
\begin{gather*}
 \sum|f_j|_{h_i}^2\leq \epsilon/100.
\end{gather*}
There exists $i_0$
such that for $i\geq i_0$
we obtain
$\int|G(h_i)|_{h_i,g_i}^2\leq\epsilon/100$.
Because the $L^2$-norms are scale invariant,
we obtain
\begin{gather*}
 \int_{X_{P,0}}\bigl|R(h_i)\bigr|^2_{h_i,g_{i,P}}
\leq \epsilon/10.
\end{gather*}
Let $X_{P,1}$ be a relatively compact neighbourhood of $P$
in $X_{P,0}$.
If $\epsilon$ is sufficiently small,
by the theorem of Uhlenbeck
\cite[Corollary 2.2]{u1},
there exists an orthonormal frame $\vecv_i$ of $(E,h_i)_{|X_{P,1}}$
for each $i$
such that the connection form $A_i$ of $\nabla_{h_i}$
with respect to $\vecv_i$
satisfies $(i)$~$A_i$ is $L_1^2$,
$(ii)$~$\|A_i\|_{L_1^2}\leq C_3\|R(h_i)\|_{L^2}$ on $X_{P,1}$
for a positive constant $C_3$ independently from $i$,
$(iii)$~$A_i$ satisfies
${\rm d}_{g_{i,P}}^{\ast}A_i=0$,
where ${\rm d}_{g_{i,P}}^{\ast}$ denotes the adjoint of ${\rm d}$
with respect to the metric $g_{i,P}$.
Let~$\Theta_i$ and~$\Theta_i^{\dagger}$
represent $\theta$ and $\theta^{\dagger}_{h_i}$
with respect to the frame $\vecv_i$.
Then, $A_i$ satisfies
\begin{gather*}
 \Lambda_{g_{P,i}}({\rm d}A_i+A_i\wedge A_i)
+\Lambda_{g_{P,i}}\big[\Theta_i,\Theta_i^{\dagger}\big]=0.
\end{gather*}
Let $X_{P,2}$ be a relatively compact neighbourhood of $P$
in $X_{P,1}$. By~the argument of Donaldson
in~the proof of~\cite[Corollary 23]{don3},
we obtain that $A_{i}$ are $L_1^p$ for any $p\geq 2$ on $X_{P,2}$,
and that there exists $C_{4,p}>0$
such that $\|A_i\|_{L_1^p}\leq C_{4,p}$ on $X_{P,2}$,
where $C_{4,p}$ are independent of $i$. In~particular,
there exists $C_{5,p}>0$ independently from $i$
such that $\big\|R(h_i)_{|X_{P,2}}\big\|_{L^p}\leq C_{5,p}$.

Let $A_i=A_i^{0,1}+A_i^{1,0}$
be the decomposition into the $(0,1)$-part and the $(1,0)$-part.
Because $\delbar\theta=0$,
we obtain $\delbar\Theta_i+\big[A_i^{0,1},\Theta_i\big]=0$.
Hence, there exist $C_{6,p}>0$ independently from $i$
such that $\big\|\Theta_{i|X_{P,2}}\big\|_{L_2^p}\leq C_{6,p}$.

Note that
$\delbar_E R(h_i)=0$
and $\del_{E,h_i}R(h_i)=0$.
Let $\delbar^{\ast}_{E,h_i,g_{i,P}}$
denote the formal adjoint of
$\delbar_E$ with respect to $h_i$ and $g_{i,P}$.
Because
$\Lambda_{g_{i,P}}R(h_i)+\Lambda_{g_{i,P}}\big[\theta_i,\theta_i^{\dagger}\big]=0$,
there exists $C_{7,p}>0$
such that
$\bigl\|
 \delbar^{\ast}_{E,h_i,g_{i,P}}R(h_i)_{|X_{P,2}}
 \bigr\|_{L_1^p}<C_{7,p}$.
Let $X_{P,3}$ be a relatively compact neighbourhood of $P$
in~$X_{P,2}$.
There exists $C_{8,p}>0$ independently from $i$
such that
$\bigl\|R(h_i)_{|X_{P,3}}\bigr\|_{L_2^p}<C_{8,p}$.
It~implies the claim of the proposition.
\end{proof}

\section[Existence and continuity of harmonic metrics in the curve case]
{Existence and continuity of harmonic metrics\\ in the curve case}

\subsection{Existence of Hermitian--Einstein metric}

Let $X$ be a compact Riemann surface.
Let $D\subset X$ be a finite subset.
Let $\big(\nbigp_{\ast}\nbigv,\DDlambda\big)$ be
a~stable good filtered $\lambda$-flat bundle on $(X,D)$.
Let $\big(E,\DDlambda\big)$ be the $\lambda$-flat bundle
on $X\setminus D$
obtained as the restriction of $\big(\nbigv,\DDlambda\big)$.
Let $\omega$ be any K\"ahler form of~$X$.
Let $h_{\det(E)}$ be a Hermitian metric of~$\det(E)$
such that
$(i)$~$\Lambda_{\omega} R(h_{\det(E)})=
 2\pi\deg(\nbigp_{\ast}\nbigv)\big(\int_X\omega\big)^{-1}$,
$(ii)$~$h_{\det(E)}$ is adapted to~$\nbigp_{\ast}\det(\nbigv)$,
i.e.,
$\nbigp^{h_{\det(E)}}_{\ast}\det(E)
\simeq
 \nbigp_{\ast}\det(\nbigv)$.
(See Proposition~\ref{prop;19.1.30.31}.)

\begin{Theorem}[Biquard--Boalch]
\label{thm;17.11.1.31}
There exists a unique Hermitian--Einstein metric
$h$ of $\big(E,\DDlambda\big)$
adapted to $\nbigp_{\ast}\nbigv$
such that $\det(h)=h_{\det(E)}$.
If $\deg(\nbigp_{\ast}\nbigv)=0$,
 $h$ is a harmonic metric.
\end{Theorem}

\begin{proof}
It~is enough to prove the case $\deg(\nbigp_{\ast}\nbigv)=0$. We~explain an outline of the proof
based on the fundamental theorem of Simpson
\cite[Theorem~1]{s1}
(and its variant~\cite[Proposition~2.49]{mochi5})
because we obtain a consequence
on the Donaldson functional
from the proof,
which will be useful in the proof of
Proposition~\ref{prop;19.1.29.10}
below.
Set $e:=\rank(\nbigv)!$.
Take $\eta>0$ such that
$10e\eta<\gaptilde(\nbigp_{\ast}\nbigv,\veca)$.
(See Section~\ref{subsection;20.7.9.1} for $\gaptilde$.)

Let $\big(X_P,z_P\big)$ be an admissible coordinate neighbourhood
around $P$.
Set $X_P^{\ast}:=X_P\setminus\{P\}$. We~take a K\"ahler metric $g_{X\setminus D}$
of $X\setminus D$
satisfying the following condition:
\begin{itemize}\itemsep=0pt
\item
$g_{X\setminus D|X_P^{\ast}}$
is mutually bounded with
$|z_P|^{-2+\eta}{\rm d}z_P\,{\rm d}\zbar_P$
on $X_P^{\ast}$
for each $P\in D$.
\end{itemize}
Recall that
the K\"ahler manifold
$(X\setminus D,g_{X\setminus D})$
satisfies the assumptions given in
\cite[Section~2]{s1},
according to~\cite[Proposition~2.4]{s1}.

\begin{Lemma}
There exists a Hermitian metric $h_0$ of $E$
such that the following holds:
\begin{enumerate}\itemsep=0pt
\item[$(a)$]
 $\big(E,{\rm d}''_E,h_0\big)$ is acceptable,
 and $\nbigp^{h_0}_{\ast}E=\nbigp_{\ast}\nbigv$.
\item[$(b)$]
$G(h_0)$ is bounded with respect to
$h_0$ and $g_{X\setminus D}$.
\item[$(c)$]
$\det(h_0)=h_{\det(E)}$.
\end{enumerate}
\end{Lemma}

\begin{proof}
By Corollary~\ref{cor;19.1.20.2},
we obtain a Hermitian metric $h_0'$ of $E$
satisfying $(a)$ and $(b)$. We~define the function $\varphi\colon X\setminus D\lrarr \real$
by $h_{\det(E)}=\det(h_0'){\rm e}^{\varphi}$.
Then, $\varphi$ induces a $C^{\infty}$-function on~$X$. We~set $h_0:=h_0'{\rm e}^{\varphi/\rank(E)}$.
Then, the metric $h_0$ has the desired property.
\end{proof}

For any $\lambda$-flat subbundle $E'\subset E$,
let $h'_0$ denote the Hermitian metric of $E'$ induced by $h_0$.
Let $\DDlambda_{E'}$ denote the Higgs field of $E'$
obtained as the restriction of $\DDlambda$. We~obtain the Chern connection
$\nabla_{h'_0}$
from the $(0,1)$-part of $\DDlambda_{E'}$ and $h'_0$.
Let $R(h'_0)$ denote the curvature of $\nabla_{h_0'}$. We~set
\begin{gather*}
 \deg(E',h_0):=
 \frac{\sqrt{-1}}{2\pi}
 \frac{1}{1+|\lambda|^2}
 \int_{X\setminus D}
 \Tr G\big(E',\DDlambda_{E'},h_0'\big)
=\frac{\sqrt{-1}}{2\pi}\int_{X\setminus D}
 \Tr R(h'_0).
\end{gather*}
Let $\Lambda_{g_{X\setminus D,\eta}}$
denote the adjoint of the multiplication by
the K\"ahler form associated with $g_{X\setminus D,\eta}$.
Because $G(h_0)$ is bounded with respect to
$h_0$ and $g_{X\setminus D,\eta}$,
$\deg(E',h_0)$ is well defined in $\real\cup\{-\infty\}$
by the Chern--Weil formula~\cite[Lemma~3.2]{s1}
(see also~\cite[Lemma~2.34]{mochi5}):
\begin{gather*}
 \deg(E',h_0)
=\frac{\sqrt{-1}}{2\pi}
 \frac{1}{1+|\lambda|^2}
 \int
 \Tr\bigl(\Lambda_{g_{X\setminus D,\eta}} G(h_0)\pi_{E'}\bigr)
-\frac{1}{2\pi}
 \frac{1}{1+|\lambda|^2}
 \int\bigl|\DDlambda\pi_{E'}\bigr|^2.
\end{gather*}
Here, $\pi_{E'}$ denotes the orthogonal projection
$E\lrarr E'$ with respect to $h_0$.

\begin{Lemma}
 $\deg(E',h_0)/\rank(E')<\deg(E,h_0)/\rank(E)$ holds.
Namely,
$(E,\DDlambda,h_0)$
is analytically stable in the sense of
{\rm\cite[Section~6]{s2}}
$($see also {\rm\cite[Section~2.3]{mochi5}}$)$.
\end{Lemma}

\begin{proof}
By~\cite[Lemma~6.1]{s2},
we have $\deg(E,h_0)=\deg\big(\nbigp^{h_0}_{\ast}E\big)=0$.
Let $0\neq E'\subsetneq E$ be a~$\lambda$-flat subbundle on $X\setminus D$. By~\cite[Lemma 6.2]{s2},
if $\deg(E',h_0)\neq-\infty$,
$E'$ extends to a~filtered subbundle
$\nbigp^{h_0'}_{\ast}E'\subset\nbigp^{h_0}_{\ast}E$,
and
$\deg(E',h_0)=\deg\big(\nbigp^{h_0'}_{\ast}E'\big)$ holds.
Because $(\nbigp^{h_0}_{\ast}E,\DDlambda)$ is assumed to be stable,
we obtain
$\deg(E',h_0)/\rank E'<\deg\big(\nbigp^{h_0}_{\ast}E\big)/\rank E=0$.
Hence, $\big(E,\delbar_E,\theta,h_0\big)$ is
analytically stable.
\end{proof}

According to the existence theorem of Simpson~\cite[Theorem~1]{s1}
(see also~\cite[Proposition~2.49]{mochi5}),
there exists a Hermitian--Einstein metric $h$ of $\big(E,\DDlambda\big)$
such that $\det(h)=\det(h_0)$
and that $h$ and $h_0$ are mutually bounded. We~already know the uniqueness
as in Proposition~\ref{prop;19.2.12.200}.
Thus, we obtain Theorem~\ref{thm;17.11.1.31}.
\end{proof}

\subsubsection{Complement on the Donaldson functional}

Let $h_0$ and $g_{X\setminus D}$ be
as in the proof of Theorem~\ref{thm;17.11.1.31}.
Let $\nbigh(h_0)$ be the space of
$C^{\infty}$-Hermitian metrics $h_1$ of $E$
satisfying the following condition:
\begin{itemize}\itemsep=0pt
\item
 Let $u_1$ be the endomorphism of $E$
 such that
 $(i)$ $h_1=h_0{\rm e}^{u_1}$,
 $(ii)$ $u_1$ is self-adjoint with respect to both $h_0$ and $h_1$.
 Then,
 $\sup_{Q\in X\setminus D}|u_1|_{h_0}(Q)
 +\|\DDlambda u_1\|_{L^2}
+\bigl\|\DDlambda\DD^{\lambda\star}_{h_0}u_1\bigr\|_{L^1}<\infty$.
Here, we use the $L^p$-norms
induced by $h_0$ and $g_{X\setminus D}$.
\end{itemize}
The Donaldson functional
$M(h_0,\bullet)\colon \nbigh(h_0)\lrarr\real$ is defined
as in~\cite[Section~5]{s1}
and~\cite[Section~2.4]{mochi5}.

\begin{Proposition}
\label{prop;19.1.28.10}
Let $h$ be the Hermitian--Einstein metric
in Theorem~$\ref{thm;17.11.1.31}$.
Then, $h$ is contained in $\nbigh(h_0)$,
and $M(h_0,h)\leq 0$ holds.
\end{Proposition}

\begin{proof}
Let $b$ be the automorphism of $E$
which is self-adjoint with respect to both $h$ and~$h_0$,
and determined by $h=h_0\cdot b$.
The theorem of Simpson~\cite[Theorem~1]{s1}
(see also~\cite[Pro\-position~2.49]{mochi5})
implies that
$b$ and $b^{-1}$ are bounded,
and that $\DDlambda b$ is $L^2$
with respect to $h_0$
and $g_{X\setminus D}$. By~\cite[Lemma 3.1]{s1}
(see also~\cite[Section~2.2.5]{mochi5}),
we also obtain $\DDlambda\DD^{\lambda\star}_{h_0}b$
is $L^1$.
Hence, $h$ is contained in $\nbigh(h_0)$. In~the proof of~\cite[Theorem~1]{s1}
and~\cite[Proposition~2.39]{mochi5},
the metric $h$ is constructed as the limit
of a subsequence of the heat flow $h_t$ $(t\geq 0)$
for which $\del_tM(h_0,h_t)\leq 0$ holds.
Because $M(h_0,h_0)=0$ by the construction,
we obtain $M(h_0,h_t)\leq 0$,
and hence $M(h_0,h)\leq 0$.
\end{proof}

\subsection{Continuities of some families of Hermitian metrics}

\subsubsection{Setting}

\paragraph{Family of curves.}

Let $\Sigma$ be a compact connected oriented real $2$-dimensional
$C^{\infty}$-manifold
with a finite subset $D\subset \Sigma$.
Let $J_i$ $(i=1,2,\ldots)$ be a sequence of complex structures
on~$\Sigma$
such that
the sequence $J_i$ is convergent to $J$ in the $C^{\infty}$-sense.
Assume that there exists a~neighbourhood $N(D)$ of $D$ in $\Sigma$
such that $J_{i|N(D)}$ are independent of $i$.
Let $X_i$ denote the compact Riemann surfaces
$(\Sigma_i,J_i)$.
Similarly, let $X$ denote the compact Riemann surface
$(\Sigma,J)$.
Let $\kappa_i\colon (T\Sigma,J)\simeq (T\Sigma,J_i)$ be
isomorphisms of complex vector bundles on $\Sigma$
such that
$(i)$~$\kappa_i\to \id$,
$(ii)$ $\kappa_{i|N(D)}=\id$. We~regard $\kappa_i$ as isomorphisms
of complex vector bundles $TX\simeq TX_i$.

For $P\in D$,
let $(X_{P},z_P)$ denote an admissible coordinate neighbourhood
of $P$ in $X$
such that $X_P\subset N(D)$. We~may regard $\big(X_P,z_P\big)$
as a holomorphic coordinate neighbourhood
of $P$ in $X_i$.
Let $r$ be a positive integer, and set $e:=r!$.
As in Section~\ref{subsection;20.2.14.11},
let $\varphi_P\colon \Xtilde_{P}\lrarr X_P$ be
the ramified covering given by
$\varphi_P(\zeta_P)=\zeta_P^e$.
Let $G_P$ denote the Galois group of
the ramified covering $\varphi_P$.

\paragraph{Family of good filtered $\boldsymbol\lambda$-flat bundles.}

Let $\big(\nbigp_{\ast}\nbigv,\DDlambda\big)$ be
a stable good filtered $\lambda$-flat bundle of rank $r$
on $(X,D)$ with $\deg(\nbigp_{\ast}\nbigv)=0$.
Let $\big(\nbigp_{\ast}\nbigv_i,\DDlambda_i\big)$
be stable good filtered $\lambda$-flat bundles
of rank $r$ on $(X_i,D)$
with $\deg(\nbigp_{\ast}\nbigv_i)=0$.
For each $P\in D$,
we set
$\big(\nbigp_{\ast}\nbigv_P,\DDlambda_P\big):=
 \big(\nbigp_{\ast}\nbigv,\DDlambda\big)_{|X_P}$
and
$\big(\nbigp_{\ast}\nbigv_{i,P},\DDlambda_{i,P}\big):=
 \big(\nbigp_{\ast}\nbigv_i,\DDlambda\big)_{|X_P}$.
Set
$\big(\nbigp_{\ast}\nbigvtilde_P,\DDtilde^{\lambda}_P\big):=
\varphi_P^{\ast}\big(\nbigp_{\ast}\nbigv_P,\DDlambda_P\big)$
and
$\big(\nbigp_{\ast}\nbigvtilde_{i,P},\DDtilde^{\lambda}_{i,P}\big):=
\varphi_P^{\ast}\big(\nbigp_{\ast}\nbigv_{i,P},\DDlambda_{i,P}\big)$.
There exist $G_P$-invariant subsets
$\nbigi(P),\nbigi(i,P)\subset
\zeta_P^{-1}\cnum[\zeta_P^{-1}]$
and the formal decompositions
\begin{gather*}
 \big(\nbigp_{\ast}\nbigvtilde_P,\DDtilde^{\lambda}_P\big)\otimes
 \cnum[\![\zeta]\!]
=\bigoplus_{\gminia\in\nbigi(P)}
 \big(\nbigp_{\ast}\nbigvtilde_{P,\gminia},\DDtilde^{\lambda}_{P,\gminia}\big),
\\
 \big(\nbigp_{\ast}\nbigvtilde_{i,P},\DDtilde^{\lambda}_{i,P}\big)\otimes
 \cnum[\![\zeta]\!]
=\bigoplus_{\gminia\in\nbigi(i,P)}
 \big(\nbigp_{\ast}\nbigvtilde_{i,P,\gminia},\DDtilde^{\lambda}_{i,P,\gminia}\big),
\end{gather*}
for each $P\in D$.
Suppose moreover that there exist
$G_P$-invariant bijections
 $\rho_{i,P}\colon \nbigi(P)\simeq \nbigi(i,P)$
 such that the following holds:
\begin{itemize}\itemsep=0pt
\item
 $\rank\nbigvtilde_{P,\gminia}
 =\rank\nbigvtilde_{i,P,\rho_{i,P}(\gminia)}$.
\item
 $\ord \gminia=\ord \rho_{i,P}(\gminia)$
and
 $\ord(\gminia-\gminib)=\ord\bigl(
\rho_{i,P}(\gminia)-\rho_{i,P}(\gminib)\bigr)$.
\item
 $\lim\limits_{i\to\infty}\rho_{i,P}(\gminia)=\gminia$
in $\zeta^{-1}\cnum\big[\zeta^{-1}\big]$.
\end{itemize}

We fix such bijections $\rho_{i,P}$.
Let $\pi_{P,\gminia}$ denote
the projection
$\nbigp_{\ast}\nbigvtilde_P\otimes\cnum[\![\zeta]\!]
\lrarr\nbigp_{\ast}\nbigvtilde_{P,\gminia}$.
Similarly,
let $\pi_{i,P,\gminia}$ denote
the projection
$\nbigp_{\ast}\nbigvtilde_{i,P}\otimes\cnum[\![\zeta]\!]
\lrarr\nbigp_{\ast}\nbigvtilde_{i,P,\gminia}$.

\paragraph{$\boldsymbol{C^{\infty}}$-isomorphisms.}

We set $\big(E,\DDlambda\big):=\big(\nbigv,\DDlambda\big)_{|X\setminus D}$
and $\big(E_i,\DDlambda_i\big):=\big(\nbigv_i,\DDlambda_i\big)_{|X_i\setminus D}$.
Let $h_{0}$ denote $C^{\infty}$-metrics of $E$
adapted to $\nbigp_{\ast}\nbigv$
such that $R(\det(h_0))=0$.
Let $h_{0,i}$ denote $C^{\infty}$-metrics of $E_i$
adapted to $\nbigp_{\ast}\nbigv_i$
such that $R(\det(h_{0,i}))=0$.
Let $d''$ and $d_i''$ denote the $(0,1)$-parts
of~$\DDlambda$ and~$\DDlambda_i$.
Suppose that there exist
$C^{\infty}$-isomorphisms $f_i\colon E\simeq E_i$
satisfying the following conditions:
\begin{itemize}\itemsep=0pt
\item
$f_i^{\ast}(h_{0,i})\to h_{0}$ in the $C^{\infty}$-sense
locally on $\Sigma\setminus D$.
\item
 On $N(D)\setminus D$,
 $f_i$ are holomorphic with respect to
 $d''$ and $d''_i$,
 and $f_i$ extend to isomorphisms
of filtered bundles
 $\nbigp_{\ast}\nbigv_{|N(D)}\simeq
 \nbigp_{\ast}\nbigv_{i|N(D)}$.
\item
For each $P\in D$, we obtain
$\Gr^F_c\big(f_{i|X_P}\big)\circ \Res_P(\DDlambda)
=\Res_P\big(\DDlambda_i\big) \circ \Gr^F_{c}\big(f_{i|X_P}\big)$
on $\Gr^F_c(\nbigv_P)$ for any $c\in\real$.
Moreover, there exists $N(P)\geq 10\rank(\nbigv)|\ord\gminia|$
for any $\gminia\in\nbigi(P)$
such that for the induced isomorphisms
$\varphi_P^{\ast}\big(f_{i|X_P}\big)\colon
 \nbigp_{\ast}\nbigvtilde_P\otimes\cnum[\![\zeta]\!]
\simeq \nbigp_{\ast}\nbigvtilde_{i,P}\otimes\cnum[\![\zeta]\!]$,
we obtain
\begin{gather}
 \bigl(\pi_{P,\gminia}-
\varphi_P^{\ast}\big(f_{i|X_P}\big)^{-1}\circ
\pi_{i,P,\rho_{i,P}(\gminia)}\circ
 \varphi_P^{\ast}\big(f_{i|X_P}\big) \bigr)\nbigp_{\ast}\nbigvtilde_{P}\otimes\cnum[\![\zeta]\!]\nonumber
 \\ \qquad
 {}
\subset \nbigp_{\ast-N(P)}\nbigvtilde_{i,P}\otimes\cnum[\![\zeta]\!],
\label{eq;20.2.17.1}
\end{gather}
and the sequences (\ref{eq;20.2.17.1}) are convergent to $0$
as $i\to\infty$.
\item
$\DDlambda-(f_i\otimes\kappa_i)^{-1}\circ\DDlambda_i\circ f_i\to 0$
in the $C^{\infty}$-sense
with respect to $h_{0}$
locally on $\Sigma\setminus D$.
\end{itemize}

\paragraph{Perturbation.}

We take $\eta_i$ $(i=1,2)$ satisfying
$10e\eta_1<
\gaptilde(\nbigp_{\ast}\nbigv)$
and
$10r\eta_2<\eta_1$. We~take $\veca\in\real^D$
for $\Partilde(\nbigp_{\ast}\nbigv,P)$ $(P\in D)$
as in Lemma~\ref{lem;20.2.12.20}.
For any $0<\epsilon<\eta_2$,
by taking $\psi_{P,\epsilon}$ $(P\in D)$
as in Section~\ref{subsection;19.1.28.1},
we obtain families of
good filtered $\lambda$-flat bundles
$\big(\nbigp_{\ast}^{(\epsilon)}\nbigv,\DDlambda\big)$
and
$\big(\nbigp_{\ast}^{(\epsilon)}\nbigv_i,\DDlambda_i\big)$. We~assume the following for each $P\in D$:
\begin{gather*}
 \sum_{a(P)-1<c\leq a(P)}\psi_{P,\epsilon}(c)\,
 \rank \Gr^F_c\big(\nbigp_{\veca}\nbigv_{|P}\big)
=\sum_{a(P)-1<c\leq a(P)}
 c\rank\Gr^F_c\big(\nbigp_{\veca}\nbigv_{|P}\big).
\end{gather*}
In particular,
$\deg(\nbigp_{\ast}\nbigv)=\deg\big(\nbigp^{(\epsilon)}_{\ast}\nbigv\big)$
and
$\deg(\nbigp_{\ast}\nbigv_i)=\deg\big(\nbigp^{(\epsilon)}_{\ast}\nbigv_i\big)$
hold. By~making $\eta_2$ smaller,
we may assume that
$\big(\nbigp_{\ast}^{(\epsilon)}\nbigv,\DDlambda\big)$
are stable for any $0\leq \epsilon\leq\eta_2$.

\subsubsection{Continuity of the family of harmonic metrics}

According to Theorem~\ref{thm;17.11.1.31},
there exists a harmonic metric
$h^{(\epsilon)}_i$
of $\big(E_i,\DDlambda\big)$
adapted to
$\nbigp^{(\epsilon)}_{\ast}\nbigv_i$
such that
$\det h_i^{(\epsilon)}
=\det h_{0,i}$.
Similarly,
there exists a harmonic metric $h^{(0)}$
of $\big(E,\DDlambda\big)$
adapted to~$\nbigp_{\ast}\nbigv$
such that
$\det h^{(0)}=\det h_0$.
The following proposition is
a variant of~\cite[Propositions~4.1 and~4.2]{mochi5}.

\begin{Proposition}\label{prop;19.1.29.10}
For any sequence
$\epsilon_i\to 0$,
the sequence $h_i^{(\epsilon_i)}$ is convergent to $h^{(0)}$
locally on~$X\setminus D$
in the $C^{\infty}$-sense.
\end{Proposition}

\begin{proof}
For $0\leq \epsilon\leq\eta_2$,
let $g_{X\setminus D,\epsilon}$ be
the K\"ahler metric on $X\setminus D$
such that the following holds on $X_P^{\ast}$ for any $P\in D$:
\begin{gather*}
 g_{X\setminus D,\epsilon|X_P^{\ast}}
=\big(\epsilon^2|z_P|^{2\epsilon}+\eta_1^2|z_P|^{2\eta_1}\big)
 |z_P|^{-2}{\rm d}z_P\,{\rm d}\zbar_P.
\end{gather*}
Let $\Lambda_{\epsilon}$
denote the adjoint of the multiplication by
the K\"ahler form $\omega_{X\setminus D,\epsilon}$
associated with $g_{X\setminus D,\epsilon}$.

By the isomorphisms
$\kappa_i\colon (T\Sigma,J)\simeq (T\Sigma,J_i)$
and the metrics $g_{X\setminus D,\epsilon}$,
we obtain the K\"ahler metrics
$g_{X_i\setminus D,\epsilon}$
of $X_i\setminus D$.
Let $\Lambda_{i,\epsilon}$
denote the adjoint of the multiplication by
the K\"ahler form~$\omega_{X_i\setminus D,\epsilon}$
associated with $g_{X_i\setminus D,\epsilon}$.

There exists an approximation of
$\big(\nbigp_{\ast}\nbigv,\DDlambda\big)_{|X_P}$
by a model filtered $\lambda$-flat bundle as in
Section~\ref{subsection;20.2.6.20}. By~using a family of Hermitian metrics for
the model $\lambda$-flat bundle as in Proposition~\ref{prop;20.2.6.2},
and by using Proposition~\ref{prop;20.2.6.10},
we construct a family of metrics
$h^{(\epsilon)}_{\inn}$ $(0\leq \epsilon\leq\eta_2)$
of $E$
such that the following holds:
\begin{itemize}\itemsep=0pt
\item
$h^{(\epsilon)}_{\inn}$ is adapted to
$\nbigp^{(\epsilon)}_{\ast}\nbigv$.
\item
$\det h_{\inn}^{(\epsilon)}=\det h_{0}$.
\item
$h_{\inn}^{(\epsilon)}\lrarr h_{\inn}^{(0)}$
locally on $X\setminus D$ in the $C^{\infty}$-sense
as $\epsilon \lrarr 0$.
\item
There exists $C_1>0$
such that
$\big|G\big(h^{(\epsilon)}_{\inn}\big)\big|_{
g_{X\setminus D,\epsilon},
h^{(\epsilon)}_{\inn}}<C_1$
for any $\epsilon$.
\end{itemize}

Let $\nu_i$ be the function on $X_i\setminus D$
determined by
$\big(f_i^{-1}\big)^{\ast}(\det h_{0})={\rm e}^{\nu_i}\det h_{i,0}$. We~obtain the Hermitian metrics
$h_{i,\inn}^{(\epsilon)}:=
 {\rm e}^{\nu_i/\rank \nbigv}
 \big(f_i^{-1}\big)^{\ast}\big(h^{(\epsilon)}_{\inn}\big)$
$(0\leq \epsilon\leq \eta_2)$
of $E_i$.
Then, by Proposition~\ref{prop;20.2.6.10},
we obtain the following:
\begin{itemize}\itemsep=0pt
\item
$h^{(\epsilon)}_{i,\inn}$ is adapted to
$\nbigp^{(\epsilon)}_{\ast}\nbigv_i$.
\item
$\det h_{i,\inn}^{(\epsilon)}=\det h_{i,0}$.
\item
$h_{i,\inn}^{(\epsilon)}\lrarr h_{i,\inn}^{(0)}$
locally on $X\setminus D$ in the $C^{\infty}$-sense
as $\epsilon \lrarr 0$.
\item
By replacing $C_1$ with a larger constant,
we may assume
$\big|G\big(h^{(\epsilon)}_{i,\inn}\big)\big|_{
g_{X_i\setminus D,\epsilon},
h^{(\epsilon)}_{i,\inn}}<C_1$
for any $\epsilon$ and any $i$.
\end{itemize}

\begin{Lemma}
\label{lem;19.2.12.1}
Let $u^{(i)}$ $(\epsilon_i\to 0)$
be automorphisms of
$E_i$ which are self-adjoint with respect to~$h^{(\epsilon_i)}_{i,\inn}$
such that the following holds:
\begin{itemize}\itemsep=0pt
\item
$\Tr\big(u^{(i)}\big)=0$.
\item
$h_{i,\inn}^{(\epsilon_i)}{\rm e}^{u^{(i)}}
\in\nbigh\big(h_{in}^{(\epsilon_i)}\big)$,
i.e.,
 $\sup\big\|u^{(i)}\big\|_{h^{(\epsilon_i)}_{i,\inn}}
 +\big\|\DDlambda_i u^{(i)}\big\|_{L^2}
 +\big\|\DDlambda_i \DD^{\lambda\star}_{i,h^{(\epsilon_i)}_{i,\inn}} u^{(i)}\big\|_{L^1}
<\infty$,
where the~$L^p$-norms are taken
with respect to $h^{(\epsilon_i)}_{i,\inn}$ and
$g_{X_i\setminus D,\epsilon_i}$. We~do not assume that the estimate is uniform in $i$.
\item
 There exists $C_2>0$ such that
 $\bigl| \Lambda_{i,\epsilon_i} G\big(h^{(\epsilon_i)}_{i,\inn}{\rm e}^{u^{(i)}}\big)
 \bigr|_{h^{(\epsilon_i)}_{i,\inn}}<C_2$
 for any $i$.
\end{itemize}
Then, there exists $C_3,C_4>0$ such that
the following holds for any $\epsilon_i$
\begin{gather*}
 \sup\big|u^{(i)}\big|_{h^{(\epsilon_i)}_{i,\inn}}
<C_3+C_4M\big(h^{(\epsilon_i)}_{i,\inn},h^{(\epsilon_i)}_{i,\inn}{\rm e}^{u^{(i)}}\big).
\end{gather*}
\end{Lemma}

\begin{proof}
By identifying the vector bundles $E_i$ and $E$ by $f_i$,
we apply the same argument
as in the proof of~\cite[Lemma 2.45]{mochi5}.
\end{proof}

Let $b_{i,1}^{(\epsilon)}$ be the automorphism of
$E_i$ which is self-adjoint with respect to
$h_{i,\inn}^{(\epsilon)}$ and $h^{(\epsilon)}_i$
and determined by
$h^{(\epsilon)}_i=h_{i,\inn}^{(\epsilon)}b_{i,1}^{(\epsilon)}$.
Note that $\det(b_{i,1}^{(\epsilon)})=1$.
Take any sequence $\epsilon_i\to 0$. By~Proposition~\ref{prop;19.1.28.10}
and Lemma~\ref{lem;19.2.12.1},
there exists a constant $C_{10}>0$
such that the following holds for any~$i$:
\begin{gather*}
 \sup_{Q\in X_i\setminus D}
 \bigl|b^{(\epsilon_i)}_{i,1|Q}\bigr|_{h_{i,\inn}^{(\epsilon_i)}}<C_{10},\qquad
 \sup_{Q\in X_i\setminus D} \bigl|\big(b^{(\epsilon_i)}_{i,1|Q}\big)^{-1}\bigr|_{h_{i,\inn}^{(\epsilon_i)}}
<C_{10}.
\end{gather*}

\begin{Lemma}
\label{lem;20.2.11.1}
$\int \Lambda_{i,\epsilon_i}\bigl(
 \delbar_{X_i}\del_{X_i}\Tr\big(b_{i,1}^{(\epsilon_i)}\big)
 \bigr)\omega_{X_i\setminus D,\epsilon_i}=0$ holds.
\end{Lemma}
\begin{proof}
We use Proposition~\ref{prop;19.1.28.10}.
Because
$\DDlambda_i \DD^{\lambda\star}_i
 b_{i,1}^{(\epsilon_i)}$
is $L^1$ with respect to $h_{i,\inn}^{(\epsilon_i)}$
and $g_{X_i\setminus D,\epsilon_i}$,
we obtain that
$\Lambda_{i,\epsilon_i}
 \delbar_{X_i}\del_{X_i}\Tr\big(b_{i,1}^{(\epsilon_i)}\big)$ is $L^1$
with respect to
$g_{X_i\setminus D,\epsilon_i}$.
Because
$\DDlambda_i b_{i,1}^{(\epsilon_i)}$ is $L^2$
with respect to $g_{X_i\setminus D,\epsilon_i}$
and $h^{(\epsilon_i)}_{i,\inn}$,
we obtain that $\del_{X_i}\Tr\big(b_{i,1}^{(\epsilon_i)}\big)$ is $L^2$
with respect to $g_{X_i\setminus D,\epsilon_i}$.
Therefore, we~obtain the claim of the lemma
by using~\cite[Lemma 5.2]{s1}.
\end{proof}

By~\cite[Lemma 3.1]{s1},
the following holds:
\begin{gather*}
 \sqrt{-1}\Lambda_{i,\epsilon_i}\delbar_{X_i}\del_{X_i}
 \Tr\big(b_{i,1}^{(\epsilon_i)}\big)
=-\Tr\bigl( b_{i,1}^{(\epsilon_i)} \Lambda_{i,\epsilon_i} G\big(h^{(\epsilon_i)}_{i,\inn}\big) \bigr)
-\bigl|\DDlambda_i\big(b_{i,1}^{(\epsilon_i)}\big)
 \cdot \big(b_{i,1}^{(\epsilon_i)}\big)^{-1/2}
 \bigr|^2_{h^{(\epsilon_i)}_{i,\inn},g_{X_i\setminus D,\epsilon_i}}.
\end{gather*}
Therefore, there exists $C_{12}>0$ such that
the following holds for any $i$:
\begin{gather*}
 \int\bigl| \DDlambda_i b_{i,1}^{(\epsilon_i)}
 \bigr|^2_{h^{(\epsilon_i)}_{i,\inn},g_{X_i\setminus D,\epsilon_i}}
 \omega_{X_i\setminus D,\epsilon_i}
 <C_{12}.
\end{gather*}
We also obtain
\begin{gather}
 \int\bigl| \DD^{\lambda\star}_i b_{i,1}^{(\epsilon_i)}
 \bigr|^2_{h^{(\epsilon_i)}_{i,\inn},g_{X_i\setminus D,\epsilon_i}}
 \omega_{X_i\setminus D,\epsilon_i} <C_{12}.
\end{gather}
Let $\big(E_i,\delbar^{(\epsilon_i)}_{E_i},\theta_i^{(\epsilon_i)}\big)$
be the Higgs bundles underlying
$\big(E_i,\DDlambda,h_i^{(\epsilon_i)}\big)$.
Then, there exists $C_{13}>0$
such that the following holds for any $i$:
\begin{gather*}
 \int\bigl| \theta^{(\epsilon_i)}_{i} \bigr|^2
_{h^{(\epsilon_i)}_{i,\inn},g_{X_i\setminus D,\epsilon_i}}
 \omega_{X_i\setminus D,\epsilon_i} <C_{13}.
\end{gather*}
Then, by applying the argument
in~\cite[Section~4.5.3]{mochi5},
we obtain the desired convergence
of the sequence $h^{(\epsilon_i)}_i$.
\end{proof}

\subsubsection{Continuity of some families of Hermitian metrics}

For $P\in D$,
we set $X_P^{\ast}:=X_P\setminus\{P\}$. We~may naturally regard $X_P^{\ast}$
as a subset of $X_i\setminus D$.
Fix~$N>10$.
Let $g_{i,\epsilon}$ be a sequence of
K\"ahler metrics of $X_i\setminus D$,
such that
$g_{i,\epsilon}\to g_{\epsilon}$ $(i\to\infty)$
and that
\begin{gather*}
 g_{i,\epsilon|X_P^{\ast}}
=
\bigl(
\epsilon^{N+2}|z_P|^{2\epsilon}
+|z_P|^{2}
\bigr)
\frac{{\rm d}z_P\,{\rm d}\zbar_P}{|z_P|^2}.
\end{gather*}
Let $\epsilon_i$ $(i=1,2,\ldots)$
be a sequence such that
$\epsilon_i\to 0$.
The following proposition
is a variant and a refinement
of~\cite[Proposition~5.1]{mochi5}.

\begin{Proposition}
\label{prop;13.1.7.20}
Let $h_{i,1}^{(\epsilon_i)}$ $(i=1,2,\ldots)$
be Hermitian metrics
of $E_i$ satisfying the following conditions:
\begin{itemize}\itemsep=0pt
\item
$\det h_{i,1}^{(\epsilon_i)}=\det h_{i,0}$.
\item
$\big\|G\big(h_{i,1}^{(\epsilon_i)}\big)\big\|_{L^2,g_{i,\epsilon},h_{i,1}^{(\epsilon_i)}}
 \to 0$ as $i\to\infty$.
\item
Let $s^{(i)}$ be the automorphism of $E_i$
which is self-adjoint
with respect to $h_i^{(\epsilon_i)}$
and determined by
 $h_{i,1}^{(\epsilon_i)}=h_i^{(\epsilon_i)}s^{(i)}$.
Then, $s^{(i)}$ and $\big(s^{(i)}\big)^{-1}$ are bounded
with respect to $h_i^{(\epsilon_i)}$
on $X\setminus D$,
and~$\DDlambda_i s^{(i)}$ are $L^2$
with respect to $h_i^{(\epsilon_i)}$
and $g_{i,\epsilon_i}$.
The estimates may depend on $i$.
\end{itemize}
Then,
the sequence
$\big\{f_i^{\ast}\big(s^{(i)}\big)\big\}$ is weakly convergent to $\id_E$
in $L_1^2$ locally on $X\setminus D$.
Moreover,
there exists $A>0$ such that
$\big|s^{(i)}\big|_{h_i^{(\epsilon_i)}}<A$
and $\big|\big(s^{(i)}\big)^{-1}\big|_{h_i^{(\epsilon_i)}}<A$
for any $i$.
\end{Proposition}

\begin{proof}
This is essentially
the same as~\cite[Proposition~5.1]{mochi5}. We~explain an outline of the proof. We~identify $E_i$ with $E$ by $f_i$. We~set
\begin{gather*}
 c_i:= \sup_{\Sigma\setminus D} \big|s^{(i)}\big|_{h^{(\epsilon_i)}_i}.
\end{gather*}
We set
$\stilde^{(i)}:=
 c_i^{-1}s^{(i)}$. We~set
$\htilde^{(\epsilon_i)}_{i,1}:=c_i^{-1}h^{(\epsilon_i)}_{i,1}
=h^{(\epsilon_i)}_{i}\stilde^{(i)}$.
The following holds.
\begin{gather*}
(1+|\lambda|^2)\Delta_{g_{i,0}} \Tr \stilde^{(i)}
= \Tr\bigl( \stilde^{(i)}\sqrt{-1} \Lambda_{g_{i,0}}
 G\big(\htilde^{(\epsilon_i)}_{i,1}\big) \bigr)
+\sqrt{-1}\Lambda_{g_{i,0}} \Tr\bigl( \DDlambda_i\stilde^{(i)}
 \big(\stilde^{(i)}\big)^{-1}
 \DD^{\lambda\star}_{i,h^{(\epsilon_i)}_i}\stilde^{(i)}
 \bigr).
\end{gather*}
From the boundedness of $\stilde^{(i)}$
and the $L^2$-property of
$\DDlambda \stilde^{(i)}$,
we obtain
$\int\Delta_{g_{i,0}}\Tr\big(\stilde^{(i)}\big)\dvol_{g_{i,0}}=0$
as in Lemma~\ref{lem;20.2.11.1}. We~obtain the following for some $A>0$ and $A'>0$:
\begin{align*}
 \int \bigl| \DDlambda\big(\stilde^{(i)}\big)
 \big(\stilde^{(i)}\big)^{-1/2} \bigr|_{g_{i,0},h^{(\epsilon_i)}_i} \dvol_{g_{i,0}}
&\leq A\cdot\int \bigl|
 \Tr\Lambda_{g_{i,0}}G\big(\htilde^{(\epsilon_i)}_{i,1}\big) \bigr| \cdot \dvol_{g_{i,0}}
 \\
&=A\cdot\int \bigl|
 \Tr\Lambda_{g_{i,\epsilon}}G\big(\htilde^{(\epsilon_i)}_{i,1}\big)
 \bigr| \cdot \dvol_{g_{i,\epsilon}}\\
&\leq A'\bigl\|G\big(\htilde^{(\epsilon_i)}_{i,1}\big)
 \bigr\|_{L^2,h^{(\epsilon_i)}_i,g_{i,\epsilon}}.
\end{align*}
Hence,
the sequence
$\stilde^{(i)}$ is $L_1^2$-bounded
on any compact subset of $X\setminus D$. By~taking an appropriate subsequence,
it is weakly convergent in $L_1^2$
locally on $X\setminus D$.
Let $\stilde^{(\infty)}$ denote the weak limit of the sequence. We~obtain $\DDlambda \stilde^{(\infty)}=0$.
Because $\stilde^{(i)}$ are self-adjoint
and uniformly bounded with respect to $h_i^{(\epsilon_i)}$,
$\stilde^{(\infty)}$ is self-adjoint and bounded
with respect to $h^{(0)}$. We~can prove that
$\stilde^{(\infty)}\neq 0$
by the same argument
as in the proof of~\cite[Lemma 5.2]{mochi5}.
Hence, $\stilde^{(\infty)}$ is a non-zero endomorphism of
$\big(\nbigp_{\ast}\nbigv,\DDlambda\big)$.
It~implies that $\stilde^{(\infty)}$ is a multiplication
by a positive constant $u_{\infty}$.

Note that the sequence $\stilde^{(i)}$
is convergent in $L^p$ for any $p$
locally on $X\setminus D$,
and hence
$\det\big(\stilde^{(i)}\big)$ is convergent to
$\det\big(\stilde^{(\infty)}\big)$ in $L^p$ for any $p$
locally on $X\setminus D$.
Because $\det\big(s^{(i)}\big)=1$, we obtain that
the sequence $c_i^{-\rank(\nbigv)}$ is convergent to
$u_{\infty}^{\rank(\nbigv)}$. In~particular, it implies that
the sequence~$c_i$ is bounded.
Then, we obtain the claim of the proposition.
\end{proof}

\subsection[Tensor product of stable filtered lambda-flat sheaves]
{Tensor product of stable filtered $\boldsymbol\lambda$-flat sheaves}
\label{subsection;20.7.10.1}

Let us state a consequence of Theorem~\ref{thm;17.11.1.31}
on the tensor product of reflexive filtered $\lambda$-flat sheaves
on arbitrary dimensional projective varieties.

Let $X$ be an $n$-dimensional non-singular projective variety
equipped with a very ample line bundle $L$.
Let $H$ be a simple normal crossing hypersurface of $X$
 with the irreducible decomposition
 $H=\bigcup_{i\in\Lambda}H_i$.
Let $\big(\nbigp_{\ast}\nbigv_i,\DDlambda_i\big)$ $(i=1,2)$
be reflexive filtered $\lambda$-flat sheaves on $(X,H)$. We~assume the following condition:
\begin{Condition}
\label{condition;20.7.7.10}
 There exists a Zariski closed subset $Z\subset H$
 with $\dim Z<n-1$
 such that
 $\big(\nbigp_{\ast}\nbigv_i,\DDlambda_i\big)_{|X\setminus Z}$
 $(i=1,2)$
 are good filtered $\lambda$-flat bundles on
 $(X\setminus Z,H\setminus Z)$.
\end{Condition}
 For example,
 if $\DDlambda$ is logarithmic and if $\lambda\neq 0$,
 Condition~\ref{condition;20.7.7.10} is satisfied.

We set
$\nbigvtilde:=\nbigv_1\otimes_{\nbigo_X(\ast H)}\nbigv_2$
which is equipped with
the induced flat $\lambda$-connection $\DDlambdatilde$.
Note that $Z\subset H$
and that $\nbigvtilde_{|X\setminus Z}$ is
a locally free $\nbigo_{X\setminus Z}(\ast (H\setminus Z))$-module.
There exists the natural morphism
\begin{gather*}
\varphi\colon\
\nbigvtilde\lrarr\nbigvtilde^{\lor\lor}:=
\nhom_{\nbigo_X(\ast H)}\bigl(
\nhom_{\nbigo_X(\ast H)}\big(\nbigvtilde,\nbigo_X(\ast H)\big),
\nbigo_X(\ast H)\bigr).
\end{gather*}
The $\nbigo_X(\ast H)$-modules
$\Ker\varphi$
and $\Cok\varphi$
are coherent,
and their supports are contained in $Z\subset H$.
Hence, we obtain that $\Ker\varphi=\Cok\varphi=0$.
It~implies that $\nbigvtilde\simeq\nbigvtilde^{\lor\lor}$,
i.e., $\nbigvtilde$
is a reflexive $\nbigo_X(\ast H)$-module.
For $\veca\in\real^{\Lambda}$,
we set
\begin{gather*}
 \nbigp'_{\veca}\nbigvtilde:= \sum_{\vecb_1+\vecb_2=\veca}
 \Image\bigl( \nbigp_{\vecb_1}\nbigv_1
 \otimes_{\nbigo_X}\nbigp_{\vecb_2}\nbigv_2
\lrarr\nbigvtilde \bigr).
\end{gather*}
Let $\nbigp_{\veca}\nbigvtilde$
denote the coherent reflexive subsheaf
of $\nbigvtilde$
generated by
$\nbigp'_{\veca}\nbigvtilde$.
Thus, we obtain a~reflexive filtered $\lambda$-flat sheaf
$\big(\nbigp_{\ast}\nbigvtilde,\DDlambdatilde\big)$
on $(X,H)$.

 \begin{Proposition}
 If $\big(\nbigp_{\ast}\nbigv_i,\DDlambda_i\big)$
 are $\mu_L$-stable,
 then
 $\big(\nbigp_{\ast}\nbigvtilde,\DDlambdatilde\big)$
 is $\mu_L$-polystable.
 \end{Proposition}

\begin{proof}
 According to Propositions~\ref{prop;19.1.30.100},
\ref{prop;20.7.7.1}
 and Condition~\ref{condition;20.7.7.10},
 there exists a positive integer $m$
 such that the following holds
 for any general complete intersection curve $Y$ of
 $L^{\otimes\,m}$.
 \begin{itemize}\itemsep=0pt
 \item $\big(\nbigp_{\ast}\nbigv_i,\DDlambda_i\big)_{|Y}$
 are stable good filtered $\lambda$-flat bundles.
 \item The natural morphism
	 \begin{gather*}
	 \Hom\bigl(\big(\nbigp_{\ast}\nbigvtilde,\DDlambdatilde\big),
	 \big(\nbigp_{\ast}\nbigvtilde,\DDlambdatilde\big)\bigr)
	 \lrarr
	 \Hom\bigl(\big(\nbigp_{\ast}\nbigvtilde,\DDlambdatilde\big)_{|Y},
	 \big(\nbigp_{\ast}\nbigvtilde,\DDlambdatilde\big)_{|Y}\bigr)
	 \end{gather*}
 is an isomorphism.
 \end{itemize}
Because $\big(\nbigp_{\ast}\nbigv_i,\DDlambda_i\big)_{|Y}$
are stable good filtered $\lambda$-flat bundles,
each $\lambda$-flat bundle
$(\nbigv_i,\DDlambda_i)_{|Y\setminus H}$
are equipped with a Hermitian--Einstein metric $h_i$
adapted to the filtered bundle
$\nbigp_{\ast}\nbigv_i$
by Theorem~\ref{thm;17.11.1.31}.
Because $h_1\otimes h_2$ is a Hermitian--Einstein
metric of the $\lambda$-flat bundle
$\big(\nbigvtilde,\DDlambdatilde\big)_{|Y\setminus H}$
adapted to the filtered bundle $\nbigp_{\ast}\nbigvtilde$,
we obtain that
$\big(\nbigp_{\ast}\nbigvtilde,\DDlambdatilde\big)_{|Y}$
is polystable. By~Corollary~\ref{cor;20.7.7.20},
we obtain that
$\big(\nbigp_{\ast}\nbigvtilde,\DDlambdatilde\big)$ is
$\mu_L$-polystable.
\end{proof}

\section[Preliminary existence theorem for Hermitian--Einstein metrics]
{Preliminary existence theorem for Hermitian--Einstein \\metrics}

\subsection{Statements}

\subsubsection{K\"ahler metrics}
\label{subsection;19.1.30.20}

Let $X$ be a smooth projective surface
with a simple normal crossing hypersurface
$H=\bigcup_{i\in\Lambda}H_i$.
Let $L$ be an ample line bundle on $X$.
Let $g_X$ be the K\"ahler metric of $X$
such that the associated K\"ahler form $\omega_X$
represents $c_1(L)$.

We take Hermitian metrics $g_i$ of $\nbigo(H_i)$.
Let $\sigma_i\colon \nbigo_X\lrarr\nbigo_X(H_i)$
denote the canonical section.
Take $N>10$.
There exists $C>0$ such that
the following form defines a K\"ahler form
on $X\setminus H$ for any $0\leq \epsilon<1/10$:
\begin{gather*}
 \omega_{\epsilon}:=
 \omega_X+
\sum_{i\in\Lambda}
 C\cdot\epsilon^{N+2}\cdot
 \sqrt{-1}
 \del\delbar|\sigma_i|^{2\epsilon}_{g_i}.
\end{gather*}
It~is easy to observe that
$\int_X\omega_{\epsilon}^2=\int_X\omega_X^2$ and
that
$\int_X\omega_{\epsilon}\tau
=\int_X\omega_X\tau$
for any closed $C^{\infty}$-$(1,1)$-form $\tau$ on $X$.

\subsubsection[Condition for good filtered lambda-flat bundles and initial metrics]
{Condition for good filtered $\boldsymbol\lambda$-flat bundles and initial metrics}

Let $\big(\nbigp_{\ast}\nbigv,\DDlambda\big)$
be a good filtered $\lambda$-flat bundle on $(X,H)$
satisfying the following condition: We~set $e:=\rank(\nbigv)!$.
\begin{Condition}
\label{condition;19.1.30.1}
\mbox{{}}
\begin{itemize}\itemsep=0pt
\item
There exists $\vecc\in\real^{\Lambda}$
and $m\in e\seisuu_{>0}$ such that
$\Par(\nbigp_{\ast}\nbigv,i)=\{c_i+n/m\mid n\in\seisuu\}$
for each~$i\in\Lambda$.
\item
The nilpotent part of $\Res_i\big(\DDlambda\big)$
on $\lefttop{i}\Gr^F_b(\nbigp_{\ast}\nbigv)$
are $0$
for any $i\in\Lambda$,
$\veca\in\real^{\Lambda}$
and $b\in\openopen{a_i-1}{a_i}$.
\end{itemize}
\end{Condition}

Let $\big(E,\DDlambda\big)$ denote the $\lambda$-flat bundle
on $X\setminus H$
obtained as the restriction of $\big(\nbigp_{\ast}\nbigv,\DDlambda\big)$.

Let $P$ be any point of $H_i\setminus \bigcup_{j\neq i}H_j$.
Let $(X_P,z_1,z_2)$ be an admissible coordinate
neighbourhood around $P$.
There exists an open subset
$X_P'$ in $\cnum^2=\{(\zeta_1,\zeta_2)\}$
such that the map
$\varphi_P\colon X_P'\lrarr X_P$
given by
$\varphi_P(\zeta_1,\zeta_2)=(\zeta_1^m,\zeta_2)$
is a ramified covering. We~set $H_P':=\{\zeta_1=0\}\allowbreak\cap X_P'$. We~obtain the induced good filtered $\lambda$-flat bundle
$\big(\nbigp_{\ast}\varphi_P^{\ast}\nbigv,\varphi_P^{\ast}\DDlambda\big)$
on $(X_P',H_P')$
such that
$\Par\big(\nbigp_{\ast}\varphi_P^{\ast}\nbigv\big)
=\{m\cdot c_i\}+\seisuu$.
\begin{Definition}
\label{df;19.2.12.2}
A Hermitian metric $h_P$ of $E_{|X_P\setminus H}$
is called strongly adapted to
$\nbigp_{\ast}\nbigv_{|X_P}$
if there exists a $C^{\infty}$ Hermitian metric $h_P'$ of
 $\nbigp_{mc_i}\big(\varphi_P^{\ast}\nbigv\big)$ on $X_P'$
 such that
 $\varphi^{-1}(h_P)=|\zeta_1|^{-2mc_i}h_P'$.
\end{Definition}

Let $P$ be any point of $H_i\cap H_j$ $(i\neq j)$.
Let $(X_P,z_1,z_2)$ be a admissible coordinate neighbourhood
around $P$ such that
$X_P\cap H_i=\{z_1=0\}$
and $X_P\cap H_j=\{z_2=0\}$.
There exists an open subset
$X_P'$ in $\cnum^2=\{(\zeta_1,\zeta_2)\}$
such that the map
$\varphi_P\colon X_P'\lrarr X_P$
given by
$\varphi_P(\zeta_1,\zeta_2)=(\zeta_1^m,\zeta_2^m)$
is a ramified covering. We~set $H_P':=\{\zeta_1\zeta_2=0\}\cap X_P'$. We~obtain the induced good filtered $\lambda$-flat bundle
$\big(\nbigp_{\ast}\varphi_P^{\ast}\nbigv,\varphi^{\ast}\DDlambda\big)$
on $(X_P',H_P')$
such that
$\Par\big(\nbigp_{\ast}\varphi_P^{\ast}\nbigv,1\big)=\{m\cdot c_i\}+\seisuu$
and
$\Par\big(\nbigp_{\ast}\varphi_P^{\ast}\nbigv,2\big)=\{m\cdot c_j\}+\seisuu$.
\begin{Definition}
\label{df;19.2.12.3}
A Hermitian metric $h_P$ of
$E_{|X_P\setminus H}$ is called strongly adapted
to $\nbigp_{\ast}\nbigv_{|X_P}$
if there exists a $C^{\infty}$-Hermitian metric $h_P'$
 of $\nbigp_{(mc_i,mc_j)}\varphi_P^{\ast}(\nbigv)$
 such that
 $\varphi^{\ast}(h_P)
 =|\zeta_1|^{-mc_i}|\zeta_2|^{-mc_j}h_P'$.
\end{Definition}

\begin{Definition}
A Hermitian metric $h$ of $E$
is called strongly adapted to $\nbigp_{\ast}\nbigv$
if the following holds:
\begin{itemize}\itemsep=0pt
\item
 For any $P\in H$,
 there exists a neighbourhood $X_P$ of $P$
 such that
 $h_{|X_P\setminus H}$
 is strongly adapted to
 $\nbigp_{\ast}\nbigv_{|X_P}$
in the sense of
Definitions~\ref{df;19.2.12.2}
and~\ref{df;19.2.12.3}.
\end{itemize}
\end{Definition}

\begin{Lemma}
\label{lem;19.2.12.50}
Let $h$ be a Hermitian metric of $E$
strongly adapted to $\nbigp_{\ast}\nbigv$.
Then, the following holds:
\begin{gather*}
 \left(
 \frac{\sqrt{-1}}{2\pi}
 \right)^2
 \int_{X\setminus H}
 \Tr\bigl( R(h)^2\bigr)
=2\int_X\ch_2(\nbigp_{\ast}\nbigv).
\end{gather*}
\end{Lemma}

\begin{proof}
It~is the equality (36) in the proof of
\cite[Proposition~4.18]{mochi4}.
\end{proof}

For each $i\in\Lambda$,
we choose
$b_i\in \Par(\nbigp_{\ast}\det\nbigv,i)$.
Set $\vecb=(b_i)\in\real^{\Lambda}$. We~take a Hermitian metric $h_{\det(E)}$
of $\det(E)$
such that
$h_{\det(E)}\prod_{i\in\Lambda}|\sigma_i|_{g_i}^{2b_i}$
induces a Hermitian metric of
$\nbigp_{\vecb}\det\nbigv$ of~$C^{\infty}$-class.

\begin{Proposition}
\label{prop;19.1.30.10}
There exists a Hermitian metric $h_{\inn}$ of $E$
such that the following holds:
\begin{itemize}\itemsep=0pt
\item
 $h_{\inn}$ is strongly adapted to $\nbigp_{\ast}\nbigv$.
\item
 $G(h_{\inn})$ is bounded
 with respect to $h_{\inn}$ and $\omega_{\epsilon}$,
where $\epsilon:=m^{-1}$.
\item
The following holds:
\begin{gather}
\label{eq;19.2.12.30}
 \int_{X\setminus H}\Tr\bigl(R(h_{\inn})^2\bigr)
=
 \frac{1}{(1+|\lambda|^2)^2}
 \int_{X\setminus H}\Tr\bigl(G(h_{\inn})^2\bigr).
\end{gather}
\item
$\det(h_{\inn})=h_{\det(E)}$.
\end{itemize}
Such a Hermitian metric $h_{\inn}$ is called
an initial metric of
$\big(\nbigp_{\ast}\nbigv,\DDlambda\big)$.
\end{Proposition}

The case $\lambda=1$ was studied
in~\cite[Sections~14.1, 14.2 and Lemma~14.4.2]{Mochizuki-wild}.
The case $\lambda\neq 1$ can be argued
in the essentially same way. We~shall explain the construction in the case $\lambda=0$
in Section~\ref{subsection;19.1.30.11}
after preliminaries in
Sections~\ref{subsection;19.1.30.12}--\ref{subsection;19.1.30.13}.

\subsubsection{Preliminary existence theorem for Hermitian--Einstein metrics}

Let $\big(\nbigp_{\ast}\nbigv,\DDlambda\big)$
be a good filtered $\lambda$-flat bundle
satisfying Condition~\ref{condition;19.1.30.1}.
Let $h_{\inn}$ be an initial metric for $\big(\nbigp_{\ast}\nbigv,\DDlambda\big)$
as in Proposition~\ref{prop;19.1.30.10}. We~shall prove the following theorem
in Section~\ref{subsection;20.7.9.100}.

\begin{Theorem}
\label{thm;19.1.30.30}
Suppose that $\big(\nbigp_{\ast}\nbigv,\DDlambda\big)$
is $\mu_L$-stable.
Then, there exists a Hermitian--Einstein metric $h_{\HE}$
of $\big(E,\DDlambda\big)$
with respect to the K\"ahler form $\omega_{\epsilon}$
$\big(\epsilon:=m^{-1}\big)$
satisfying the following conditions:
\begin{enumerate}\itemsep=0pt
\item[$(i)$]
 $h_{\HE}$ and $h_{\inn}$ are mutually bounded.
\item[$(ii)$]
 $\DDlambda\big(h_{\HE}\cdot h_{\inn}^{-1}\big)$
 is $L^2$
 with respect to $h_{\inn}$ and $\omega_{\epsilon}$.
\item[$(iii)$]
$\det(h_{\HE})=\det(h_{\inn})$
holds. In~particular, the following holds:
\begin{gather*}
\frac{1}{1+|\lambda|^2}\Tr\big(G(h_{\HE})\big)
=\frac{1}{1+|\lambda|^2}\Tr\big(G(h_{\inn})\big)
=\Tr\big(R(h_{\inn})\big).
\end{gather*}
 \item[$(iv)$]
The following equality holds:
\begin{gather}
\label{eq;19.1.30.111}
 \bigg( \frac{\sqrt{-1}}{2\pi} \bigg)^2
 \frac{1}{(1+|\lambda|^2)^2}
 \int_{X\setminus H} \Tr\bigl(G(h_{\HE})^2\bigr)
=2\int_X \ch_{2}(\nbigp_{\ast}\nbigv).
\end{gather}
\end{enumerate}
\end{Theorem}

\subsection{Around cross points}
\label{subsection;19.1.30.12}

Let $X_0:=\bigl\{(z_1,z_2)\in\cnum^2\mid |z_i|<1\bigr\}$. We~set
$H_i:=X_0\cap\{z_i=0\}$ and $H:=H_1\cup H_2$.
Let $(\nbigp_{\ast}\nbigv,\theta)$
be a good filtered Higgs bundle on $(X_0,H)$. We~choose
$b_i\in\Par(\nbigp_{\ast}\nbigv,i)$ $(i=1,2)$,
and set $\vecb=(b_1,b_2)$. We~also choose any Hermitian metric
$h_{\det(E)}$ of $\det(E)$
such that
$h_{\det(E)}|z_1|^{2b_1}|z_2|^{2b_2}$ is a Hermitian metric
of $\nbigp_{\vecb}(\det\nbigv)$ of $C^{\infty}$-class.

\subsubsection{Unramified case}

Suppose that $(\nbigp_{\ast}\nbigv,\theta)$
satisfies the following condition:

\begin{Condition}
\label{condition;19.2.12.20}
\mbox{{}}
\begin{itemize}\itemsep=0pt
\item
There exists $\vecc=(c_1,c_2)\in\real^2$
such that
$(i)$ $-1<c_i\leq 0$,
$(ii)$ $\Par(\nbigp_{\ast}\nbigv,i)
=\{c_i+n\mid n\in\seisuu\}$.
\item
There exists a decomposition of good filtered Higgs bundles
\begin{gather}
\label{eq;20.2.14.30}
(\nbigp_{\ast}\nbigv,\theta)
=\bigoplus_{\gminia\in\nbigi}
 \bigoplus_{\vecalpha\in\cnum^2}
 (\nbigp_{\ast}\nbigv_{\gminia,\vecalpha},
 \theta_{\gminia,\vecalpha})
\end{gather}
such that
$\theta_{\gminia,\alpha}
 -\big({\rm d}\gminia+\sum\alpha_ i{\rm d}z_i/z_i\big)\id_{\nbigv_{\gminia,\vecalpha}}$
induce holomorphic Higgs fields
of $\nbigp_{\vecc}\nbigv_{\gminia,\vecalpha}$.
\end{itemize}
\end{Condition}

We take any holomorphic frame
$\vecv=(v_j)$
of $\nbigp_{\vecc}\nbigv$
compatible with the decomposition (\ref{eq;20.2.14.30}).
For $j=1,\ldots,r$,
we obtain $(\gminia_j,\vecalpha_j)$
determined by
$v_j\in \nbigp_{\vecc}\nbigv_{\gminia_j,\vecalpha_j}$.
Let $h_0$ be the metric of
$\nbigv_{|X_0\setminus H}$
determined by
$h_0(v_i,v_i)=|z_1|^{-2c_1}|z_2|^{-2c_2}$
and
$h_0(v_i,v_j)=0$ $(i\neq j)$.
Note that
$\del_{h_0}\vecv=\vecv
 \bigl(
 -\sum_{k=1,2} c_k {\rm d}z_k/z_k
 \bigr)\,I$,
where $I$ denotes the identity matrix.
Hence,
$\big[\del_{h_0},\delbar\,\big]=0$ holds. We~obtain the description
$\theta\vecv=\vecv\bigl(
 \Lambda_0+\Lambda_1
 \bigr)$
such that the following holds:
\begin{itemize}\itemsep=0pt
\item
$(\Lambda_0)_{ii}
=({\rm d}\gminia_i+\sum_{k=1,2}\alpha_{ik}{\rm d}z_k/z_k)$
and
$(\Lambda_0)_{ij}=0$ $(i\neq j)$.
\item
$(\Lambda_1)_{ij}$ are holomorphic $1$-forms
for any $i$ and $j$.
Moreover,
$(\Lambda_1)_{ij}=0$ holds
unless $(\gminia_i,\vecalpha_i)=(\gminia_j,\vecalpha_j)$.
\end{itemize}
We obtain
$\theta^{\dagger}_{h_0}\vecv
=\vecv\big(\Lambdabar_0+\lefttop{t}\Lambdabar_1\big)$
and
$\big[\theta,\theta_{h_0}^{\dagger}\big]\vecv
=\vecv\big[\Lambda_1,\lefttop{t}\Lambdabar_1\big]$,
where the entries of
$\big[\Lambda_1,\lefttop{t}\Lambdabar_1\big]$
are $C^{\infty}$ on $X_0$. We~have
$(\del_{h_0}\theta)\vecv
=\vecv(\del\Lambda_1)$,
where
any entries of $\del\Lambda_1$ are holomorphic $2$-forms,
and $(\del\Lambda_1)_{ij}=0$ unless
$(\gminia_i,\vecalpha_i)=(\gminia_j,\vecalpha_j)$.

Note that there exists a $C^{\infty}$-function $u$ on $X_0$
such that
$\det(h_0)={\rm e}^{u}h_{\det(E)}$. We~set
$h_{\inn}:=h_0{\rm e}^{-u/\rank E}$.
\begin{Lemma}
$\big[\theta,\theta_{h_{\inn}}^{\dagger}\big]$,
$\del_{h_{\inn}}\theta$ and $\delbar\theta_{h_{\inn}}^{\dagger}$
are bounded with respect to $h_{\inn}$
and
$\sum_{k=1,2}{\rm d}z_k\,{\rm d}\zbar_k$.
\end{Lemma}

\subsubsection{Ramified case}
\label{subsection;14.1.7.10}

Let $\varphi\colon \cnum^2\lrarr \cnum^2$
be given by $\varphi(\zeta_1,\zeta_2)=(\zeta_1^{m},\zeta_2^{m})$. We~set $X_0':=\varphi^{-1}(X_0)$,
$H_i':=X_0'\cap\varphi^{-1}(H_i)$
and $H':=H_1'\cup H_2'$. We~set
$\Gal(\varphi):=
 \{(\kappa_1,\kappa_2)\in\cnum^2\mid
 \kappa_i^{m}=1\}$,
which acts on $X'_0$ by
$(\kappa_1,\kappa_2)(\zeta_1,\zeta_2)
=(\kappa_1\zeta_1,\kappa_2\zeta_2)$.

Suppose that
$\varphi^{\ast}(\nbigp_{\ast}\nbigv,\theta)$
satisfies Condition~\ref{condition;19.2.12.20}
on $(X',H')$. We~construct a $C^{\infty}$-metric $h_0'$
of $\varphi^{\ast}(E)_{|X_0'\setminus H'_0}$
as in the previous subsection. We~may assume that $h_0'$ is
$\Gal(\varphi)$-invariant.
Note that there exists a $\Gal(\varphi)$-invariant
$C^{\infty}$-function $u$ on $X_0'$
such that
$\det(h_0')={\rm e}^{u}\varphi^{-1}(h_{\det(E)})$. We~set
$h_{\inn}':=h_0'{\rm e}^{-u/\rank(E)}$.
Because it is $\Gal(\varphi)$-invariant,
we obtain the induced metric $h_{\inn}$ of~$E$.

Let $g_{X_0'}$ denote the K\"ahler metric
$\sum_{k=1,2}{\rm d}\zeta_k\,{\rm d}\zetabar_k$ on $X_0'$.
Because $\varphi\colon X_0'\setminus H'\lrarr X_0\setminus H$
is a covering map,
it induces a K\"ahler metric
$\varphi_{\ast}(g_{X_0'})$
of $X_0\setminus H$.

\begin{Lemma}
\label{lem;19.2.12.40}
$\bigl[\theta,\theta_{h_{\inn}}^{\dagger}\bigr]$,
$\del_{h_{\inn}}\theta$ and $\delbar\theta_{h_{\inn}}^{\dagger}$
are bounded with respect to
$\big(h_{\inn},\varphi_{\ast}g_{X_0'}\big)$.
\end{Lemma}

\subsubsection{An estimate}

We set
$Y(\epsilon):=\bigl\{ (z_1,z_2)\in X_0\mid \min(|z_i|)=\epsilon \bigr\}$.
\begin{Lemma}
\label{lem;19.2.12.31}
We obtain
$\lim\limits_{\epsilon\to 0}
 \int_{Y(\epsilon)}\Tr\big(\theta\delbar\theta^{\dagger}\big)=0$
and
$\lim\limits_{\epsilon\to 0}
 \int_{Y(\epsilon)}\Tr\big(\theta^{\dagger} \del\theta\big)=0$.
\end{Lemma}
\begin{proof}
It~is enough to consider the case
where Condition~\ref{condition;19.2.12.20}
is satisfied for $(\nbigp_{\ast}\nbigv,\theta)$.
Let~$f$ be any anti-holomorphic function on $X_0$.
Let us consider
$\int_{Y(\epsilon)}{\rm d}\gminia\,f\,{\rm d}\zbar_1{\rm d}\zbar_2$. We~set
$Y_1(\epsilon):=\{|z_1|=\epsilon,|z_2|\geq \epsilon\}$
and
$Y_2(\epsilon):=\{|z_2|=\epsilon,|z_1|\geq \epsilon\}$. We~have
\begin{gather*}
 \int_{Y_1(\epsilon)}{\rm d}\gminia f\,{\rm d}\zbar_1{\rm d}\zbar_2
=\int_{Y_1(\epsilon)}\del_2\gminia {\rm d}z_2 f {\rm d}\zbar_1{\rm d}\zbar_2.
\end{gather*}
It~is of the form
\begin{gather}
\label{eq;19.1.21.11}
 \int_{Y_1(\epsilon)}
 \frac{\gminib(z_1,z_2)}{z_1^{\ell_1}z_2^{\ell_2}}
 f(\zbar_1,\zbar_2){\rm d}\zbar_1 {\rm d}\zbar_2.
\end{gather}
Here, $\gminib$ is a holomorphic function. We~consider the Taylor expansion of
$\gminib$ and $f$.
Then, the contributions of the terms
\begin{gather*}
 \frac{z_1^{k_1}\zbar_1^{m_1}}{z_1^{\ell_1}}\,{\rm d}\zbar_1 \,
 \frac{z_2^{k_2}}{z_2^{\ell_2}}\,\zbar_2^{m_2}\, {\rm d}\zbar_2\,{\rm d}z_2
\end{gather*}
to (\ref{eq;19.1.21.11}) is $0$
unless $k_1-\ell_1-m_1=1$
and $k_2-\ell_2-m_2=0$.
If the equalities hold,
we have \allowbreak$k_1-\ell_1+m_1=2m_1+1\geq 1$
and
$k_2-\ell_2+m_2=2m_2\geq 0$.
Hence, we obtain
$\lim\limits_{\epsilon\to 0}
 \int_{Y_1(\epsilon)} {\rm d}\gminia f\,{\rm d}\zbar_1{\rm d}\zbar_2\allowbreak=0$.
Similarly, we obtain
$\lim\limits_{\epsilon\to 0}
 \int_{Y_2(\epsilon)} {\rm d}\gminia f\,{\rm d}\zbar_1{\rm d}\zbar_2=0$.
Similarly and more easily,
we obtain
$\lim\limits_{\epsilon\to 0}
 \int_{Y(\epsilon)}(\alpha_i {\rm d}z_i/z_i) f {\rm d}\zbar_1 {\rm d}\zbar_2=0$.
Then, the claim of the lemma follows.
\end{proof}

\subsection{Around smooth points}
\label{subsection;19.1.30.13}

We set $X_0:=\bigl\{(z_1,z_2)\in\cnum^2\mid |z_i|<1\bigr\}$
and $H:=\{z_1=0\}$.
Let $\nu\colon X_0\setminus H\lrarr\real_{>0}$
be a~$C^{\infty}$-function
such that
$\nu|z_1|^{-1}$ induces
a nowhere vanishing $C^{\infty}$-function on $X_0$.
Let $(\nbigp_{\ast}\nbigv,\theta)$
be a good filtered Higgs bundle on $(X_0,H)$.
Let $\big(E,\delbar_E,\theta\big)$ be the Higgs bundle
obtained as the restriction of $(\nbigp_{\ast}\nbigv,\theta)$
to $X_0\setminus H$. We~choose $b\in\Par(\nbigp_{\ast}\det\nbigv)$
and a Hermitian metric $h_{\det(E)}$ of
$\det(E)$ such that
$h_{\det(E)}\nu^{2b}$ induces a $C^{\infty}$ metric of
$\nbigp_b(\det\nbigv)$.

\subsubsection{Unramified case}

Suppose that $(\nbigp_{\ast}\nbigv,\theta)$
satisfies Condition~\ref{condition;19.2.12.11}.
\begin{Condition}
\mbox{{}}\label{condition;19.2.12.11}
\begin{itemize}\itemsep=0pt
\item
There exists $-1<c\leq 0$ such that
$\Par(\nbigp_{\ast}\nbigv)=
 \{c+n\mid n\in\seisuu\}$.
\item
There exists a decomposition of good filtered Higgs bundles
\begin{gather*}
 (\nbigp_{\ast}\nbigv,\theta)
=\bigoplus_{\gminia\in\nbigi}
 \bigoplus_{\alpha\in\cnum}
 (\nbigp_{\ast}\nbigv_{\gminia,\alpha},\theta_{\gminia,\alpha}).
\end{gather*}
\item
$\theta_{\gminia,\alpha}
-({\rm d}\gminia+\alpha {\rm d}z_1/z_1)\id_{\nbigv_{\gminia,\alpha}}$
are holomorphic Higgs fields of
$\nbigp_{c}\nbigv_{\gminia,\alpha}$.
\end{itemize}
\end{Condition}

We take $C^{\infty}$-metrics
$h_{\gminia,\alpha}$
of $\nbigp_c\nbigv_{\gminia,\alpha}$,
and we set
$h_0:=\bigoplus \nu^{-2c}h_{\gminia,\alpha}$. We~may assume that
$\det(h_0)=h_{\det(E)}$.

Let $\vecv=(v_1,\ldots,v_r)$ be any holomorphic frame of
$\nbigp_c\nbigv$ compatible with
the decomposition.
For each $i$,
$\gminia_i$ and $\alpha_i$
are determined by the condition that
$v_i$ is a section of
$\nbigp_{c}\nbigv_{\gminia_i,\alpha_i}$.
There exist matrix valued $C^{\infty}$-$(1,0)$-forms
$A_{\gminia,\alpha}$
such that
\begin{gather*}
 \del_{h_0}\vecv
=\vecv\Bigl( (-c\cdot\del\log\nu^2)I+\sum A_{\gminia,\alpha} \Bigr),
\end{gather*}
where $I$ denotes the identity matrix,
and $(A_{\gminia,\alpha})_{i,j}=0$
unless $(\gminia_i,\alpha_i)=(\gminia_j,\alpha_j)=(\gminia,\alpha)$.
Let~$\Lambda$ denote the matrix valued holomorphic $1$-form
determined by
$\theta\vecv=\vecv\Lambda$.
There exists the decomposition
$\Lambda=\Lambda_0+\Lambda_1$
such that the following holds:
\begin{itemize}\itemsep=0pt
\item
$(\Lambda_0)_{ij}=({\rm d}\gminia_i+\alpha_i {\rm d}z_1/z_1)$
if $i=j$,
and
$(\Lambda_0)_{ij}=0$ if $i\neq j$.
\item
$(\Lambda_1)_{ij}$ are holomorphic $1$-forms,
and $(\Lambda_1)_{ij}=0$
unless $(\gminia_i,\alpha_i)=(\gminia_j,\alpha_j)$.
\end{itemize}
There exists
a matrix valued $C^{\infty}$ $(0,1)$-form $\Lambda_2$
such that
$\theta_{h_0}^{\dagger}\vecv
=\vecv\big(\Lambdabar_0+\Lambda_2\big)$.
Moreover,
$(\Lambda_2)_{ij}=0$ holds unless
$(\gminia_i,\alpha_i)=(\gminia_j,\alpha_j)$.

We have
$R(h_0)=
 \big({-}c\delbar\del\log\sigma^2 \big)I
+\bigoplus R(h_{\gminia,\alpha})$,
where $R(h_{\gminia,\alpha})$
are $C^{\infty}$.
Note that
$d\Lambda_0=0$
and
$[\Lambda_0,\Lambda_i]= \big[\Lambda_0,\Lambdabar_i\big]=0$.
Hence,
$[\theta,\theta_{h_0}^{\dagger}]$,
$\del_{h_0}\theta$ and
$\delbar\theta_{h_0}^{\dagger}$
are $C^{\infty}$. We~also have
\begin{gather*}
\bigl(\del_{h_0}\theta\bigr)\vecv=\vecv\Bigl(\del\Lambda_1
+\Bigl[ \bigoplus A_{\gminia,\alpha},\Lambda_1 \Bigr]\Bigr),
\qquad
 \bigl(\delbar\theta_{h_0}^{\dagger} \bigr)\vecv
=\vecv\big(\delbar\Lambda_2\big).
\end{gather*}

We set
$w_1=z_1|z_1|^{-1}\nu$
and $w_2=z_2$.
Then, it is easy to check that
$(w_1,w_2)$ is a $C^{\infty}$ complex coordinate system.
Clearly,
${\rm d}\zbar_2={\rm d}\wbar_2$.
There exists a $C^{\infty}$-function
$\gamma$ and a $C^{\infty}$ $(0,1)$-form $\kappa$
such that
${\rm d}\zbar_1=\gamma {\rm d}\wbar_1+\wbar_1\kappa_1$. We~set $Y(\epsilon)=\{\nu=\epsilon\}=\{|w_1|=\epsilon\}$.

\begin{Lemma}
$\lim\limits_{\epsilon\to 0}
 \int_{Y(\epsilon)}\Tr\big(\theta\delbar\theta^{\dagger}\big)=0$
and
$\lim\limits_{\epsilon\to 0}
 \int_{Y(\epsilon)}\Tr\big(\theta^{\dagger}\del\theta\big)=0$
hold.
\end{Lemma}
\begin{proof}
It~is enough to prove
$\lim\limits_{\epsilon\to 0}
 \int_{Y(\epsilon)}\Tr\big(\theta\delbar\theta^{\dagger}\big)=0$.
It~is easy to see that
\begin{gather*}
 \lim_{\epsilon\to 0}
 \int_{Y(\epsilon)}\Tr\big(\Lambda_1\delbar\Lambda_2\big)=0.
\end{gather*}
Let us study
$\int_{Y(\epsilon)}\Tr\big(\Lambda_0\delbar\Lambda_2\big)$.
For any $C^{\infty}$-function $g$,
we consider the following integral:
\begin{gather}
\label{eq;19.1.21.20}
 \int_{Y(\epsilon)}
 g ({\rm d}\gminia \cdot {\rm d}\zbar_1\,{\rm d}\zbar_2)
=\int_{Y(\epsilon)}
 (g \gamma) \cdot
 {\rm d}\gminia\,{\rm d}\wbar_1\, {\rm d}\wbar_2
+\int_{Y(\epsilon)}
 g \wbar_1\cdot
 {\rm d}\gminia\,\kappa\, {\rm d}\wbar_2.
\end{gather}
We can rewrite the first term in the right hand side of
(\ref{eq;19.1.21.20}) as follows,
for some non-negative integer $\ell$
and for a $C^{\infty}$-function $\gminib$:
\begin{gather*}
 \int_{Y(\epsilon)}
 (g\gamma)\,
 {\rm d}\gminia\,{\rm d}\wbar_1\,{\rm d}\wbar_2
=\int_{Y(\epsilon)}
 (g\gamma\gminib)
 w_1^{-\ell}
 \, {\rm d}\wbar_1
 \,{\rm d}w_2\,{\rm d}\wbar_2
\end{gather*}
Take $N>\ell$. We~consider the expansion
\begin{gather*}
 g\gamma\gminib=
 \sum_{\substack{k,m\geq 0 \\ k+m\leq N}}
 (g\gamma\gminib)_{k,m}(w_2)w_1^k\wbar_1^m
+O\big(|w_1|^N\big).
\end{gather*}
Here,
$(g\gamma\gminib)_{k,m}(w_2)$ are $C^{\infty}$-functions of $w_2$.
The contributions
\begin{gather*}
 \int_{Y(\epsilon)}
 (g\gamma\gminib)_{k,m}(w_2)
 \frac{w_1^k\wbar_1^{m}}{w_1^{\ell}}
 {\rm d}\wbar_1\,{\rm d}w_2\,{\rm d}\wbar_2
\end{gather*}
are $0$
unless $k-\ell-m=1$.
If $k-\ell-m=1$,
then $k-\ell+m=2m+1\geq 1$ holds.
Hence, we obtain
\begin{gather*}
 \lim_{\epsilon\to 0}
 \int_{Y(\epsilon)}
 (g\gamma)\,
 {\rm d}\gminia\,{\rm d}\wbar_1{\rm d}\wbar_2=0.
\end{gather*}
We rewrite the second term in
the right hand side of (\ref{eq;19.1.21.20})
as follows,
for some $C^{\infty}$-functions $f_i$ $(i=1,2)$
and a non-negative integer $\ell$:
\begin{gather*}
 \int_{Y(\epsilon)}
 g\wbar_1\cdot
 {\rm d}\gminia\,\kappa\,{\rm d}\wbar_2
=\int_{Y(\epsilon)}
 f_1w_1^{-\ell-1}\wbar_1\,
 \,{\rm d}w_1\,{\rm d}w_2\,{\rm d}\wbar_2
+\int_{Y(\epsilon)}
f_2w_1^{-\ell}\,\wbar_1
\,{\rm d}\wbar_1\,{\rm d}w_2\,{\rm d}\wbar_2.
\end{gather*}
Take $N>\ell+1$.
Consider the expansions
$f_i=\sum (f_i)_{k,m}(w_2)\,w_1^{k}\wbar_1^m+O\big(|w_1|^N\big)$.
The contributions
\begin{gather*}
\int_{Y(\epsilon)}
 (f_1)_{k,m}(w_2)
 \,\frac{w_1^{k}\wbar_1^{m+1}}{w_1^{\ell+1}}\,
 {\rm d}w_1\,{\rm d}w_2\,{\rm d}\wbar_2
\end{gather*}
are $0$
unless $k-(\ell+1)-(m+1)=-1$.
If $k-(\ell+1)-(m+1)=-1$ holds,
then we have
$k-(\ell+1)+(m+1)=2m+1\geq 1$.
The contributions
\begin{gather*}
\int_{Y(\epsilon)}
 (f_2)_{k,m}(w_2)
 \,\frac{w_1^{k}\wbar_1^{m+1}}{w_1^{\ell}}\,
 {\rm d}\wbar_1\,{\rm d}w_2\,{\rm d}\wbar_2
\end{gather*}
are $0$ unless
$k-\ell-(m+1)=1$.
If $k-\ell-(m+1)=1$ holds,
then we have
$k-\ell+(m+1)=2(m+1)+1\geq 3$.
Hence, we obtain
\begin{gather*}
 \lim_{\epsilon\to 0} \int_{Y(\epsilon)}
 g\wbar_1\cdot {\rm d}\gminia\,\kappa\,{\rm d}\wbar_2=0.
\end{gather*}

Similarly and more easily,
we obtain
$\lim\limits_{\epsilon\to 0}
 \int_{Y(\epsilon)} g(\alpha {\rm d}z_1/z_1)\,{\rm d}\zbar_1{\rm d}\zbar_2=0$
for any $\alpha\in\cnum$
and for any $C^{\infty}$-function $g$.
Thus, we obtain the claim of the lemma.
\end{proof}

\subsubsection{Ramified case}
\label{subsection;14.1.7.11}

Let $\varphi\colon \cnum^2\lrarr\cnum^2$
be given by
$\varphi(\zeta_1,\zeta_2)=\big(\zeta_1^m,\zeta_2\big)$. We~set
$X_0':=\varphi^{-1}(X_0)$
and $H':=\varphi^{-1}(H)$.
Let $\Gal(\varphi):=\{\mu\in\cnum\mid \mu^m=1\}$,
which acts on $X_0'$
by $\mu\cdot (\zeta_1,\zeta_2)=(\mu\zeta_1,\zeta_2)$.

Suppose that
$\varphi^{\ast}(\nbigp_{\ast}\nbigv,\theta)$
satisfies Condition~\ref{condition;19.2.12.11}. We~construct a Hermitian metric $h_0'$
for $\varphi^{\ast}(\nbigp_{\ast}\nbigv,\theta)$
as in the previous subsection. We~may assume that
$h_0'$ is $\Gal(\varphi)$-invariant.
There exists a $C^{\infty}$-function
$f$ on $X_0'$
determined by
$\det(h_0')=e^f\varphi^{-1}(h_{\det(E)})$. We~set
$h_{\inn}':=h_0'e^{-f/\rank(E)}$.
Because $h_{\inn'}$ is $\Gal(\varphi)$-invariant,
we obtain a Hermitian metric $h_{\inn}$ of $E$
induced by $h_{\inn}'$.
Let $\varphi_{\ast}(g_{X_0'})$
denote the K\"ahler metric of
$X_0\setminus H_0$
induced by $\sum_{k=1,2}{\rm d}\zeta_k\,{\rm d}\zetabar_k$.

\begin{Lemma}
\label{lem;19.2.12.32}
$R(h_{\inn})$,
$\big[\theta,\theta_{h_{\inn}}^{\dagger}\big]$,
$\del_{h_{\inn}}\theta$
and $\delbar\theta_{h_{\inn}}^{\dagger}$
are bounded with respect to
$\varphi_{\ast}g_{X'_0}$ and $h_0$. We~also have
\begin{gather*}
\lim_{\epsilon\to 0}\int_{Y(\epsilon)}
 \Tr\big(\theta\delbar\theta_{h_{\inn}}^{\dagger}\big)=0,
\qquad
 \lim_{\epsilon\to 0}
 \int_{Y(\epsilon)}
 \Tr\big(\theta^{\dagger}\del_{h_{\inn}}\theta\big)=0.
\end{gather*}
\end{Lemma}

\subsection{Proof of Proposition~\ref{prop;19.1.30.10}}
\label{subsection;19.1.30.11}

Let $X$, $H$ and $L$ be as in Section~\ref{subsection;19.1.30.20}.
Let $(\nbigp_{\ast}\nbigv,\theta)$
be a good filtered Higgs bundle on $(X,H)$
satisfying Condition~\ref{condition;19.1.30.1}.
Note that
$(\nbigp_{\ast}\nbigv,\theta)$
is as in Section~\ref{subsection;14.1.7.10}
around any cross point of $H$,
and
$(\nbigp_{\ast}\nbigv,\theta)$
is as in Section~\ref{subsection;14.1.7.11}
around any smooth points of $H$.
There exists a Hermitian metric $h_{\inn}$ of $E$
such that
$(i)$ $\det(h_{\inn})=h_{\det(E)}$,
$(ii)$ the restriction of $h_{\inn}$ around any points of~$H$
are as in Section~\ref{subsection;14.1.7.10}
or Section~\ref{subsection;14.1.7.11}. By~the construction,
$h_{\inn}$ is strongly adapted to
$\nbigp_{\ast}\nbigv$. By~Lemmas~\ref{lem;19.2.12.40}
and~\ref{lem;19.2.12.32},
we obtain that
$R(h_{\inn})$,
$\big[\theta,\theta_{h_{\inn}}^{\dagger}\big]$,
$\del_{h_{\inn}}\theta$
and $\delbar\theta_{h_{\inn}}^{\dagger}$
are bounded with respect to
$h_{\inn}$ and $\omega_{\epsilon}$.
As in the proof of~\cite[Proposition~4.18]{mochi4},
we have
\begin{gather*}
 \Tr \bigl( G(h_{\inn})^2\bigr)
=\Tr\bigl( R(h_{\inn})^2\bigr)
+{\rm d}\bigl( \Tr\big(\theta\delbar\theta_{h_{\inn}}^{\dagger}\big)
+\Tr\big(\theta_{h_{\inn}}^{\dagger}\del_{h_{\inn}}\theta\big) \bigr).
\end{gather*}
Then,
we obtain (\ref{eq;19.2.12.30})
from Lemmas~\ref{lem;19.2.12.31} and~\ref{lem;19.2.12.32}.
Thus, we obtain Proposition~\ref{prop;19.1.30.10}.
\hfill\qed

\subsection{Proof of Theorem~\ref{thm;19.1.30.30}}
\label{subsection;20.7.9.100}

Let $E'\subset E$ be
any coherent $\lambda$-flat $\nbigo_{X\setminus H}$-subsheaf. We~assume that
$E'$ is saturated,
i.e., $E/E'$ is torsion-free.
Let $\big(E',\DDlambda_{E'}\big)$ be the induced $\lambda$-flat sheaf
on $X\setminus H$.
There exists a~discrete subset $Z\subset X\setminus H$
such that
$E'_{|X\setminus (H\cup Z)}$ is a subbundle of
$E_{|X\setminus (H\cup Z)}$.
Let $h'$ denote the metric of~$E'_{|X\setminus (H\cup Z)}$
induced by $h_{\inn}$. We~obtain the Chern connection
$\nabla_{h'}$ of $\big(E',d''_{E'},h'\big)$
and the operator
$\DD^{\lambda\star}_{E',h'}$
from $\DD^{\lambda}_{E'}$ and $h'$.
Let $R(E',h')$ denote the curvature of $\nabla_{h'}$. We~obtain
$G(E',h'):=\bigl[ \DDlambda_{E'},\DD^{\lambda\star}_{E',h'} \bigr]$.
Following~\cite{s1}, we define
\begin{gather*}
 \deg_{\omega_{\epsilon}}(E',h_{\inn}):=
 \frac{\sqrt{-1}}{2\pi}
 \frac{1}{1+|\lambda|^2}
 \int_{X\setminus H}
 \Tr\bigl(\Lambda_{\omega_{\epsilon}} G(E',\theta',h')\bigr)
 \dvol_{\omega_{\epsilon}}.
\end{gather*}
It~is well defined in $\real\cup\{-\infty\}$
by the Chern--Weil formula~\cite[Lemma 3.2]{s1}:
\begin{gather*}
 \deg_{\omega_{\epsilon}}(E',h_{\inn})
=\frac{\sqrt{-1}}{2\pi} \frac{1}{1+|\lambda|^2}
 \int_{X\setminus H}
 \Tr\bigl( \pi_{E'}\Lambda_{\omega_{\epsilon}}G(h_{\inn}) \bigr)
-\frac{1}{2\pi} \frac{1}{1+|\lambda|^2} \int_{X\setminus H}
 \bigl| \DDlambda \pi_{E'} \bigr|^2_{h_{\inn},\omega_{\epsilon}}.
\end{gather*}
Here, $\pi_{E'}$ denotes the orthogonal projection
of $E_{|X\setminus (H\cup Z)}$ onto $E'_{|X\setminus (H\cup Z)}$.

\begin{Lemma}
If $\deg_{\omega_{\epsilon}}(E',\theta)\neq-\infty$,
then
$E'$ extends to a filtered subsheaf
$\nbigp^{h'}_{\ast}E'$ of $\nbigp_{\ast}\nbigv$
and
\begin{gather*}
 \deg_{\omega_{\epsilon}}(E',h_{\inn})
=\int_{X}
 c_1\big(\nbigp^{h'}_{\ast}E'\big)\omega_X
\end{gather*}
holds. As a result,
$\big(E,\delbar_E,\theta,h_{\inn}\big)$
is analytically stable
in the sense of~{\rm \cite{s1}}
$($see also {\rm \cite[Section~2.3]{mochi5}}$)$.
\end{Lemma}

\begin{proof}
If $\deg_{\omega_{\epsilon}}(E',h_{\inn})\neq-\infty$,
we obtain
$\int|d''\pi_{E'}|^2<\infty$.
As studied in~\cite{li2, Li-Narasimhan}
on the basis of~\cite{siu},
we obtain a coherent
$\nbigo_X(\ast H)$-submodule
$\nbigp^{h'}(E')\subset\nbigv$
as an extension of $E'$.
Moreover, as proved in~\cite[Lemma 4.20]{mochi4},
we obtain the equality
$\deg_{\omega_{\epsilon}}(E',h_{\inn})
=\int_{X} c_1\big(\nbigp^{h'}_{\ast}E'\big)\omega_X$.
\end{proof}

According to the fundamental theorem of Simpson~\cite[Theorem~1]{s1}
and its variant~\cite[Proposition~2.49]{mochi5},
there exists a Hermitian--Einstein metric $h_{\HE}$
of $\big(E,\DDlambda\big)$
satisfying the conditions $(i)$, $(ii)$ and $(iii)$. By~\cite[Proposition~3.5]{s1}
and~\cite[Lemma~7.4]{s1}
(see also~\cite[Proposition~2.49]{mochi5}), we obtain
\begin{gather*}
 \bigg(\frac{\sqrt{-1}}{2\pi}\bigg)^2
 \frac{1}{(1+|\lambda|^2)^2}
 \int_{X\setminus H} \Tr\bigl(G(h_{\HE})^2\bigr)
= \bigg(\frac{\sqrt{-1}}{2\pi}\bigg)^2 \frac{1}{(1+|\lambda|^2)^2}
 \int_{X\setminus H} \Tr\bigl(G(h_{\inn})^2\bigr).
\end{gather*}
\looseness=1
It~is equal to
$2\int_X\ch_2(\nbigp_{\ast}\nbigv)$
by Lemma~\ref{lem;19.2.12.50}
and Proposition~\ref{prop;19.1.30.10}.
Thus, Theorem~\ref{thm;19.1.30.30}
is proved.
\hfill\qed

\section{Bogomolov--Gieseker inequality}

Let $X$ be any dimensional smooth connected projective variety
with a simple normal crossing hypersurface
$H=\bigcup_{i\in\Lambda}H_i$.
Let $L$ be any ample line bundle on $X$.

\begin{Theorem}
\label{thm;19.1.30.131}
Let $\big(\nbigp_{\ast}\nbigv,\DDlambda\big)$
be a $\mu_L$-polystable good filtered $\lambda$-flat bundle on $(X,H)$.
Then, the Bogomolov--Gieseker inequality holds:
\begin{gather*}
 \int_X\ch_2(\nbigp_{\ast}\nbigv)
 c_1(L)^{\dim X-2}
\leq
 \frac{\int_Xc_1(\nbigp_{\ast}\nbigv)^2c_1(L)^{\dim X-2}}{2\rank \nbigv}.
\end{gather*}
\end{Theorem}

\begin{proof}
By the Mehta--Ramanathan type theorem
(Proposition~\ref{prop;19.1.30.100}),
it is enough to study the case $\dim X=2$,
which we shall assume in the rest of the proof. We~use the notation in~Section~\ref{subsection;19.1.30.101}.
Let
$\big(\nbigp_{\ast}\nbigv,\DDlambda\big)
=\bigoplus \big(\nbigp_{\ast}\nbigv_j,\DDlambda_j\big)$
be the decomposition
into the stable components.

We set $e:=\rank(\nbigv)!$. We~choose $\eta>0$ such that
$0<10e\eta< \gaptilde(\nbigp_{\ast}\nbigv)$. We~take $\veca\in\real^{\Lambda}$
for~$\Partilde(\nbigp_{\ast}\nbigv,i)$ $(i\in\Lambda)$
as in Lemma~\ref{lem;20.2.12.20}.

Let $m\in e\seisuu_{>0}$
such that $\epsilon:=m^{-1}<\eta/10\rank(\nbigv)$.
For any $b\in\Partilde(\nbigp_{\ast}\nbigv,i)$, we set
$b(\epsilon):=\max\bigl\{ d\in\epsilon\seisuu\mid
 d<b \bigr\}$. We~set
\begin{gather*}
 c_{i,j}:=\frac{1}{\rank \nbigv_j}
 \sum_{b\in\Par(\nbigp_{\ast}\nbigv_j,\veca,i)}
 (b-b(\epsilon))
 \rank \lefttop{i}\Gr^F_{b}(\nbigp_{\veca}\nbigv_j).
\end{gather*}
We have $0\leq c_{i,j}\leq\epsilon$.
For any $b\in\Partilde(\nbigp_{\ast}\nbigv_j,i)$,
we set
$\psi_{\epsilon,i,j}(b):=b(\epsilon)+c_{i,j}$.
Then,
we obtain
$|\psi_{\epsilon,i,j}(b)-b|<2\epsilon$
and the following equalities:
\begin{gather*}
 \sum_{b\in\Par(\nbigp_{\ast}\nbigv_j,\veca,i)}
 \psi_{\epsilon,i,j}(b)\rank\lefttop{i}\Gr^F_b(\nbigp_{\veca}\nbigv_j)
=\sum_{b\in\Par(\nbigp_{\ast}\nbigv_j,\veca,i)}
 b\rank\lefttop{i}\Gr^F_b(\nbigp_{\veca}\nbigv_j).
\end{gather*}
Moreover, we have
$\psi_{\epsilon,i,j}(b)-c_{i,j}\in \epsilon\seisuu$.

Applying the construction in Section~\ref{subsection;19.1.30.101},
we obtain good filtered $\lambda$-flat bundles
$\big(\nbigp^{(\epsilon)}_{\ast}\nbigv_j,\DDlambda_j\big)$
on $(X,H)$. By~the construction, they satisfy
Condition~\ref{condition;19.1.30.1}. By~Lemma~\ref{lem;21.5.5.1},
there exists $m_0$ such that
$\big(\nbigp^{(\epsilon)}_{\ast}\nbigv_j,\DDlambda_j\big)$
are $\mu_L$-stable if $m\geq m_0$.
Let $\big(E_j,\DDlambda_j\big)$ be the $\lambda$-flat bundle
obtained as the restriction of
$\big(\nbigp_{\ast}\nbigv_j,\DDlambda_j\big)$ to $X\setminus H$. We~use the K\"ahler metric $g_{\epsilon}$
of $X\setminus H$
as in Section~\ref{subsection;19.1.30.20}.
There exist Hermitian--Einstein metrics
$h^{(\epsilon)}_{j,\HE}$
of the $\lambda$-flat bundles $\big(E_j,\DDlambda_j\big)$
as in Theorem~\ref{thm;19.1.30.30}
for the good filtered $\lambda$-flat bundles
$\big(\nbigp^{(\epsilon)}_{\ast}\nbigv_j,\DDlambda_j\big)$.
Note that
$\bigoplus h^{(\epsilon)}_{j,\HE}$
is a Hermitian--Einstein metric of
$\bigoplus \big(E_j,\DDlambda_j\big)$.

By Proposition~\ref{prop;19.1.30.110},
the equality (\ref{eq;19.1.30.111}),
and the equality
$\frac{\sqrt{-1}}{2\pi}\frac{1}{1+|\lambda|^2}
 \Tr G\big(h^{(\epsilon)}_{\HE}\big)
=\frac{\sqrt{-1}}{2\pi}R(h_{\det E})$,
we obtain
\begin{gather*}
 \int_X\ch_2\big(\nbigp^{(\epsilon)}_{\ast}\nbigv\big)
\leq
 \frac{\int_Xc_1\big(\nbigp^{(\epsilon)}_{\ast}\nbigv\big)^2}
 {2\rank \nbigv}.
\end{gather*}
By taking the limit as $m\to\infty$, i.e., $\epsilon\to 0$,
we obtain the desired inequality.
\end{proof}

\begin{Corollary}
\label{cor;19.1.30.120}
Let $\big(\nbigp_{\ast}\nbigv,\DDlambda\big)$ be
a $\mu_L$-polystable good filtered
$\lambda$-flat bundle on $(X,H)$.
Suppose that
\begin{gather*}
\int_Xc_1(\nbigp_{\ast}\nbigv)
 c_1(L)^{\dim X-1}=0,
\qquad
\int_X\ch_2(\nbigp_{\ast}\nbigv)
 c_1(L)^{\dim X-2}=0.
\end{gather*}
Then, $c_1(\nbigp_{\ast}\nbigv)=0$ holds.

Moreover,
for any decomposition
$\big(\nbigp_{\ast}\nbigv,\DDlambda\big)
=\bigoplus\big(\nbigp_{\ast}\nbigv_j,\DDlambda_j\big)$
into $\mu_L$-stable good filtered $\lambda$-flat bundles,
we obtain
$c_1(\nbigp_{\ast}\nbigv_j)=0$
and
$\int_X\ch_2(\nbigp_{\ast}\nbigv_j)c_1(L)^{\dim X-2}=0$.
\end{Corollary}

\begin{proof}
On one hand,
because of the Hodge index theorem
and
$\int_Xc_1(\nbigp_{\ast}\nbigv)
 c_1(L)^{\dim X-1}=0$,
we obtain
\begin{gather*}
 \int_Xc_1(\nbigp_{\ast}\nbigv)^2c_1(L)^{\dim X-2}\leq 0,
\end{gather*}
and the equality holds if and only if
$c_1(\nbigp_{\ast}\nbigv)=0$.
On the other hand,
by the Bogomolov--Gieseker inequality
and $\int_X\ch_2(\nbigp_{\ast}\nbigv)c_1(L)^{\dim X-2}=0$,
we obtain
\begin{gather*}
 \int_Xc_1(\nbigp_{\ast}\nbigv)^2c_1(L)^{\dim X-2}\geq 0.
\end{gather*}
Hence, we obtain $c_1(\nbigp_{\ast}\nbigv)=0$.

Let
$\big(\nbigp_{\ast}\nbigv,\DDlambda\big)
=\bigoplus\big(\nbigp_{\ast}\nbigv_j,\DDlambda_j\big)$
be a decomposition
into $\mu_L$-stable good filtered $\lambda$-flat bundles. We~have
$\int_Xc_1(\nbigp_{\ast}\nbigv_j)c_1(L)^{\dim X-1}=0$.
Hence, by the Hodge index theorem,
we obtain
\begin{gather*}
 \int_Xc_1(\nbigp_{\ast}\nbigv_j)^2c_1(L)^{\dim X-2}
\leq 0.
\end{gather*}
By the Bogomolov--Gieseker type inequality,
we obtain
\begin{gather*}
 \int_X\ch_2(\nbigp_{\ast}\nbigv_j)
 c_1(L)^{\dim X-2}\leq 0.
\end{gather*}
Because
$\sum_j\int_X\ch_2(\nbigp_{\ast}\nbigv_j)
 c_1(L)^{\dim X-2}=
 \int_X\ch_2(\nbigp_{\ast}\nbigv)=0$,
we obtain
\begin{gather*}
 \int_X\ch_2(\nbigp_{\ast}\nbigv_j)
 c_1(L)^{\dim X-2}=0.
\end{gather*}
Thus, we obtain the claim of the corollary.
\end{proof}

\begin{Remark}
Although $H$ was assumed to be ample in
\cite[Section~14.4, Corollary~14.5.1]{Mochizuki-wild},
it~is not essential.
Indeed,
for any simple normal crossing hypersurface $H$,
there exists an ample simple normal crossing hypersurface $H'$
such that $H\subset H'$.
Let $\big(\nbigp'_{\ast}\nbigv,\DDlambda\big)$
be the filtered $\lambda$-flat bundle on $(X,H')$
naturally induced by $\big(\nbigp_{\ast}\nbigv,\DDlambda\big)$.
The Chern characters of
$\nbigp_{\ast}\nbigv$
and $\nbigp_{\ast}'\nbigv$ are equal,
and hence the Bogomolov--Gieseker inequalities for
$\nbigp_{\ast}\nbigv$
and $\nbigp_{\ast}'\nbigv$ are equivalent.
\end{Remark}

\section{Existence theorem of pluri-harmonic metrics}\label{section;19.2.11.11}

\subsection{Statement}\label{subsection;19.1.30.130}

Let us prove Theorem~\ref{thm;19.2.11.20}.
According to Corollary~\ref{cor;19.1.30.120},
it is enough to study the case
where $\big(\nbigp_{\ast}\nbigv,\DDlambda\big)$
is a {\em $\mu_L$-stable} good filtered $\lambda$-flat bundle
on $(X,H)$ such that
\begin{gather*}
c_1(\nbigp_{\ast}\nbigv)=0,
\quad
\int_X\ch_2(\nbigp_{\ast}\nbigv)c_1(L)^{\dim X-2}=0.
\end{gather*}
Let $\big(E,\DDlambda\big)$ be the $\lambda$-flat bundle
obtained as the restriction
$\big(\nbigp_{\ast}\nbigv,\DDlambda\big)_{|X\setminus H}$.
Let $h_{\det(E)}$ denote the pluri-harmonic metric
of $\big(\det(E),\DDlambda_{\det(E)}\big)$
strongly adapted to $\nbigp_{\ast}(\det(E))$.
For the proof of~Theorem~\ref{thm;19.2.11.20},
it is enough to prove the following theorem.
\begin{Theorem}
\label{thm;19.1.31.100}
There exists a unique pluri-harmonic metric $h$
of the $\lambda$-flat bundle $\big(E,\DDlambda\big)$
such that
$\nbigp^{h}_{\ast}E=\nbigp_{\ast}\nbigv$
and $\det(h)=h_{\det(E)}$.
\end{Theorem}
The proof is given in the rest of this section.

\subsection{Surface case}
\label{subsection;19.2.1.20}

Let us study the case
$\dim X=2$.
The following argument is essentially
the same as the proof of
\cite[Theorem~5.5]{mochi5}.
Let $\big(\nbigp_{\ast}\nbigv,\DDlambda\big)$
be as in Section~\ref{subsection;19.1.30.130}. We~use the notation in the proof of~Theorem~\ref{thm;19.1.30.131}.
For large $m\in e\seisuu_{>0}$,
we set $\epsilon:=m^{-1}$. We~have the perturbations
$\big(\nbigp^{(\epsilon)}_{\ast}\nbigv,\DDlambda\big)$. We~use the K\"ahler metrics $g_{\epsilon}$
of $X\setminus H$
as in Section~\ref{subsection;19.1.30.20}.
There exist the Hermitian--Einstein metrics
$h^{(\epsilon)}_{\HE}$
of $\big(E,\DDlambda\big)$
adapted to $\big(\nbigp^{(\epsilon)}_{\ast}\nbigv,\DDlambda\big)$
such that $\det\big(h^{(\epsilon)}_{\HE}\big)=h_{\det(E)}$.
\begin{Proposition}
\label{prop;19.1.31.1}
For any sequence $m_i\to\infty$,
we set $\epsilon_i:=m_i^{-1}$.
Then, after going to a~sub\-se\-qu\-ence,
$h^{(\epsilon_i)}_{\HE}$
is convergent almost everywhere on $X\setminus H$,
and the limit $h$ is a pluri-harmonic metric of
the $\lambda$-flat bundle $\big(E,\DDlambda\big)$
adapted to $\nbigp_{\ast}\nbigv$
such that $\det(h)=h_{\det E}$.
\end{Proposition}

\subsubsection{Family of ample hypersurfaces}
\label{subsection;19.2.1.1}

There exists a $0$-dimensional closed subset
$Z\subset H$ such that
$(i)$ $Z$ contains the singular points of~$H$,
$(ii)$ any $P\in H\setminus Z$ has
a neighbourhood $H_P$ in $H$
on which the conjugacy classes of~$\Res(\DDlambda)_{|Q}$ $(Q\in H_P)$ are constant.

Take a sufficiently large integer $M$. We~set
$\gbigz_M:=H^0\big(X,L^{\otimes\,M}\big)\setminus\{0\}$.
It~is equipped with a natural $\cnum^{\ast}$-action.
Let $p_i$ denote the projection of
$X\times\gbigz_M$ onto the $i$-th component.
There exists the universal section
$\gminis$ of $p_1^{\ast}\big(L^{\otimes M}\big)$.
Let $\gbigx_M$ denote the scheme
obtained as $\gminis^{-1}(0)$.
Let~$\ttP_1\colon \gbigx_M\lrarr X$
and $\ttP_2\colon \gbigx_M\lrarr\gbigz_M$
denote the morphism induced by $p_i$.
For each $s\in\gbigz_M$,
let~$X_{s}$ denote the fiber product of
$\ttP_2$
and the inclusion $\{s\}\lrarr\gbigz_M$.

There exists the $\cnum^{\ast}$-invariant
maximal Zariski open subset
$\gbigz^{\circ}_M\subset\gbigz_M$
such that
$(i)$~the induced morphism
 $\ttP_2^{\circ}\colon \gbigx_M^{\circ}:=
\gbigx_M\times_{\gbigz_M}\gbigz_M^{\circ}
\lrarr\gbigz_M^{\circ}$
is smooth,
$(ii)$~$X_s\cup H$ is normal crossing for any~$s\in\gbigz_M^{\circ}$,
$(iii)$~$(X_s\cap H)\cap Z=\varnothing$.
Let $\ttP_1^{\circ}$ denote the restriction of
$\ttP_1$ to $\gbigx^{\circ}_M$.
For any $Q\in \gbigx_M^{\circ}$,
we obtain the subspace
$T_{\ttP_1(Q)}X_{\ttP_2(Q)}
\subset T_{\ttP_1(Q)}X$
of codimension $1$.
It~determines a~point
in $\proj\big(T^{\ast}_{\ttP_1(Q)}X\big)$.
Hence, we obtain the natural morphism
$\ttPtilde^{\circ}_1\colon \gbigx_M^{\circ}\lrarr \proj(T^{\ast}X)$.
If $M$ is sufficiently large,
$\ttP_1^{\circ}$ and $\ttPtilde_1^{\circ}$
are surjective.

By the Mehta--Ramanathan type theorem
(Proposition~\ref{prop;19.1.30.100}),
there exists a non-empty
$\cnum^{\ast}$-invariant Zariski open subset
$\gbigz^{\sankaku}_M$ of $\gbigz_M^{\circ}$
such that the following holds:
\begin{itemize}\itemsep=0pt
\item
For each $s\in \gbigz^{\sankaku}_M$,
$\big(\nbigp_{\ast}\nbigv,\DDlambda\big)_{|X_s}$
is stable.
\end{itemize}

We set
$\gbigx_M^{\sankaku}:=
 \gbigx_M^{\circ}\times_{\gbigz^{\circ}_M}\gbigz_M^{\sankaku}$.
Note that
$W_M:=X\setminus\ttP^{\circ}_1\big(\gbigx_M^{\sankaku}\big)$
is a finite set.
For each $P\in X\setminus (H\cup W_M)$,
the intersection
$\ttPtilde^{\circ}_1\big(\gbigx_M^{\sankaku}\big)
\cap
 \proj(T^{\ast}_PX)$
in $\proj(T^{\ast}X)$
is Zariski dense in $\proj(T^{\ast}_PX)$.

We set $H_s:=X_s\cap H$.
Let $(E_s,\DDlambda_s)$
denote the $\lambda$-flat bundle on $X_s\setminus H_s$
obtained as the restriction of $\big(E,\DDlambda\big)$.
For each $s\in \gbigz_M^{\sankaku}$,
there exists a pluri-harmonic metric
$h_s$ of $\big(E_s,\DDlambda_s\big)$
such that
$(i)$ $h_s$ is adapted to $\nbigp_{\ast}\nbigv_{|X_s}$,
$(ii)$ $\det(h_s)=h_{\det(E)|X_s\setminus H_s}$.

Let $\ttP_1^{\sankaku}\colon \gbigx^{\sankaku}_M\lrarr X$
be the induced map.
Let $\gbigh_M^{\sankaku}:=\big(\ttP_1^{\sankaku}\big)^{-1}(H)$. We~set
$\big(E^{\sankaku},\DDlambda_{E^{\sankaku}}\big)
:=
 \big(\ttP_1^{\sankaku}\big)^{-1}\big(E,\DDlambda\big)$
on $\gbigx_M^{\sankaku}\setminus \gbigh^{\sankaku}_M$. By~Lemma~\ref{lem;20.2.14.20} and Proposition~\ref{prop;19.1.29.10},
the family of pluri harmonic metrics $h_s$
$(s\in\gbigz_M^{\sankaku})$
induces a continuous Hermitian metric $h^{\sankaku}$ of
$E^{\sankaku}$. We~also obtain Hermitian metrics
$h^{\sankaku(\epsilon_i)}:=
 \big(\ttP_1^{\sankaku}\big)^{-1}\big(h^{(\epsilon_i)}_{\HE}\big)$.

\subsubsection{Local holomorphic coordinate systems}
\label{subsection;19.1.31.20}

Let $P\in X\setminus W_M$. We~take $s_{\infty}\in\gbigz_M^{\sankaku}$
such that $P\not\in X_{s_{\infty}}$.
The following is clear
because
$\ttPtilde_1\big(\gbigx^{\sankaku}_M\big)\cap \proj(T^{\ast}_PX)$
is dense in $\proj(T^{\ast}_PX)$.
\begin{Lemma}
There exist $s_i\in \gbigz^{\sankaku}_M$ $(i=1,2)$
and $\delta>0$
such that the following holds:
\begin{itemize}\itemsep=0pt
\item
$P\in X_{s_i}$ $(i=1,2)$.
\item
$X_{s_1}$ and $X_{s_2}$ are transversal at $P$.
\item
$\{s_1+as_{\infty}\mid |a|<\delta\}$,
$\{s_2+as_{\infty}\mid |a|<\delta\}$,
$\bigl\{s_1+s_2+as_{\infty}\mid |a|<\delta\bigr\}$
and
$\bigl\{s_1+\sqrt{-1}s_2+as_{\infty}\mid |a|<\delta\bigr\}$
are contained in
$\gbigz^{\sankaku}_{M}$.
\end{itemize}
\end{Lemma}

We set $x_i:=s_i/s_{\infty}$ $(i=1,2)$.
There exists a neighbourhood $U_P$
of $P$ in $X\setminus H$
such that
$(x_1,x_2)$ is a holomorphic coordinate system
on $U_P$.
Note that
$\bigl\{\sum b_ix_i=c\bigr\}\cap U_P$
is equal to~$U_P\cap X_{b_1s_1+b_2s_2-cs_{\infty}}$.

\subsubsection{Proof of Proposition~\ref{prop;19.1.31.1}}

Take a sequence $m_i\to\infty$ in $\seisuu$. We~set $\epsilon_i:=m_i^{-1}$. By~Proposition~\ref{prop;19.1.30.110},
we obtain the following convergence:
\begin{gather*}
\lim_{i\to\infty} \int_{X\setminus H}
 \bigl|G\big(h^{(\epsilon_i)}_{\HE}\big)\bigr|^2_{h^{(\epsilon_i)}_{\HE},\omega_{\epsilon_i}}=0.
\end{gather*}
Let $h^{(\epsilon_i)}_{\inn}$ be a Hermitian metric
for $\big(\nbigp^{(\epsilon_i)}_{\ast}\nbigv,\DDlambda\big)$
as in Proposition~\ref{prop;19.1.30.10}.
Let $b_i$ be the automorphism of $E$
determined by
$h^{(\epsilon_i)}_{\HE}=h^{(\epsilon_i)}_{\inn}\cdot b_i$.
Then, for each $i$,
$\DDlambda(b_i)$ is $L^2$
with respect to
$h^{(\epsilon_i)}_{\HE}$ and $\omega_{\epsilon_i}$.

Let $\omega_{\epsilon_i,s}$ denote
the K\"ahler form of $X_s\setminus H_s$
induced by $\omega_{\epsilon_i}$.
Let $h_s^{(\epsilon_i)}$ denote the restriction of
$h^{(\epsilon_i)}_{\HE}$ to $X_s\setminus H_s$.

By Fubini's theorem,
after going to a subsequence,
there exists a $\cnum^{\ast}$-invariant subset
\mbox{$\gbigz_M^{\sharp}\!\subset\!\gbigz_M^{\sankaku}$}
with the following property:
\begin{enumerate}\itemsep=0pt
\item[$(a1)$]
$\lim\limits_{i\to\infty}\int_{X_s\setminus H_s}
 \bigl|
 G\big(h_s^{(\epsilon_i)}\big)
 \bigr|^2_{h_s^{(\epsilon_i)},\omega_{\epsilon_i,s}}=0$
holds
for each $s\in\gbigz_M^{\sharp}$.
\item[$(a2)$]
For each $s\in \gbigz_M^{\sharp}$,
$\DDlambda_s\big(b_{i|X_s\setminus H_s}\big)$
is $L^2$ with respect to
$h^{(\epsilon_i)}_s$
and $\omega_{\epsilon_i,s}$.
\item[$(a3)$]
The Lebesgue measure of
$\gbigz_M^{\sankaku}\setminus\gbigz_M^{\sharp}$
is $0$.
\end{enumerate}

Note that the condition $(a2)$ implies the following.
\begin{Lemma}
Let $s\in \gbigz_M^{\sharp}$.
Let $\htilde^{(\epsilon_i)}_s$
be a harmonic metric of $(E_s,\DDlambda_s)$
adapted to
$\big(\nbigp_{\ast}^{(\epsilon_i)}\nbigv,\DDlambda\big)_{|X_s}$
such that
$\det\big(\htilde^{(\epsilon_i)}_s\big)=h_{\det(E)|X_s\setminus H_s}$.
Let $\btilde_{i,s}$ be the automorphism of
$E_{|X_s\setminus H_s}$
determined by
$h^{(\epsilon_i)}_s=\htilde^{(\epsilon_i)}_s\btilde_{i,s}$.
Then,
$\btilde_{i,s}$
and $\btilde_{i,s}^{-1}$
are bounded with respect to
$\htilde^{(\epsilon_i)}_s$,
and $\DDlambda\big(\btilde_{i,s}\big)$
is $L^2$ with respect to
$\htilde^{(\epsilon_i)}_s$
and $\omega_{\epsilon_i,s}$
\end{Lemma}

\begin{proof}
Let $b'_{i,s}$ be the automorphism of $E_s$
determined by
$\htilde^{(\epsilon_i)}_s=h^{(\epsilon_i)}_{\inn|X_s\setminus H_s}$.
Then,
$b'_{i,s}$ and $(b'_{i,s})^{-1}$ are bounded
with respect to $h^{(\epsilon_i)}_{\inn|X_s\setminus H_s}$,
and
$\DDlambda(b'_{i,s})$ is $L^2$
with respect to
$h^{(\epsilon_i)}_{\inn|X_s\setminus H_s}$
and~$\omega_{\epsilon_i,s}$,
according to Proposition~\ref{prop;19.1.28.10}.
Then, we obtain the claim of the lemma.
\end{proof}

\begin{Lemma}
There exists a $\cnum^{\ast}$-invariant subset
$\gbigx_M^{\sharp}\subset
 \gbigx_M^{\sankaku}
 \times_{\gbigz_M^{\sankaku}}
 \gbigz_M^{\sharp}$
such that
the following holds:
\begin{itemize}\itemsep=0pt
\item
 The measure of
 $\gbigx_M^{\sankaku}\setminus\gbigx_M^{\sharp}$
 is $0$.
 \item
 A subsequence of
 $h^{\sankaku(\epsilon_i)}_{|\gbigx_M^{\sharp}}$
 is convergent to $h^{\sankaku}_{|\gbigx_M^{\sharp}}$
 at any point of $\gbigx_M^{\sharp}$.
\end{itemize}
\end{Lemma}
\begin{proof}
By Proposition~\ref{prop;13.1.7.20},
for any $s\in \gbigz_{M}^{\sharp}$,
the sequence
$h_s^{(\epsilon_i)}$
is weakly convergent to $h_s$ in~$L_1^2$
locally on $X_s\setminus H_s$. We~set
$b^{(\epsilon_i)}_s:=h_s^{(\epsilon_i)} h_s^{-1}$. We~obtain
$\det\big(b_s^{(\epsilon_i)}\big)=1$,
and $b^{(\epsilon_i)}_s$ converges
to the identity
locally on $X_s\setminus H_s$
in $L^p$ for any $p\geq 1$. We~set
$g^{(\epsilon_i)}_s:=
 \bigl|
 b^{(\epsilon_i)}_s
 \bigr|_{h_s}$ on $X_s\setminus H_s$. We~obtain the function
$g^{(\epsilon_i)}$
on
$\gbigx_M^{\sankaku}
 \times_{\gbigz_M^{\sankaku}}
 \gbigz_M^{\sharp}$
from
$g^{(\epsilon_i)}_s$
$(s\in\gbigz_M^{\sharp})$. By~Lemma~\ref{lem;20.2.14.20}
and Proposition~\ref{prop;13.1.7.20},
for any compact subset
$K\subset
\gbigx_M^{\sankaku}
\setminus\gbigh_M^{\sankaku}$,
the restriction of
$g^{(\epsilon_i)}$
to
$K\cap\bigl(
 \gbigx_M^{\sankaku}
 \times_{\gbigz_M^{\sankaku}}
 \gbigz_M^{\sharp}
\bigr)$
are uniformly bounded.
Note that the sequence
$\bigl(
h^{\sankaku(\epsilon_i)}\big(h^{\sankaku}\big)^{-1}
\bigr)_{|K\cap X_s}
=b_{s|K\cap X_s}^{(\epsilon_i)}$
is convergent in $L^p$ for any $p\geq 1$. By~Fubini's theorem and Lebesgue theorem,
we obtain the $L^p$-convergence
of~$h^{\sankaku(\epsilon_i)}\big(h^{\sankaku}\big)^{-1}$
to the identity for any $p$
on $K\cap
 \bigl(\gbigx_M^{\sankaku} \times_{\gbigz_M^{\sankaku}} \gbigz_M^{\sharp}\bigr)$.
Then,
after going to a subsequence,
we obtain the desired convergence.
\end{proof}

\begin{Remark}
If $\lambda=0$,
the argument can be simplified.
Indeed, by Proposition~\ref{prop;20.1.30.1},
the curvature $R\big(h_s^{(\epsilon_i)}\big)$ of $h_s^{(\epsilon_i)}$
are bounded locally on $X_s\setminus H$.
Hence, we obtain that
$h^{(\epsilon_i)}_s$ is weakly convergent to $h_s$ in $L_2^2$. In~particular,
$h_s^{(\epsilon_i)}$
is convergent to $h_s$
in the $C^0$-sense locally on $X_s\setminus H_s$.
\end{Remark}

There exists a subset
$X^{\sharp}\subset\ttP^{\sankaku}_1\big(\gbigx_M^{\sharp}\big)$
such that for any $P\in X^{\sharp}$,
the measure of
$\big(\ttP^{\sankaku}_1\big)^{-1}(P)\setminus
 \bigl(
\gbigx_M^{\sharp}
 \bigr)$
is $0$
in $\big(\ttP^{\sankaku}_1\big)^{-1}(P)$,
and that
the measure of $X\setminus X^{\sharp}$ is $0$
in $X$. We~obtain that the sequence~$h^{(\epsilon_i)}_{\HE|X^{\sharp}}$
is convergent
to a Hermitian metric
$h_{\infty}$ of $E_{|X^{\sharp}}$.

\begin{Lemma}
For any $P\in X^{\sharp}$
and
$s\in \ttP_2\bigl(\big(\ttP_1^{\sankaku}\big)^{-1}(P)\bigr)$,
we obtain $(h_s)_{|P}=h_{\infty|P}$.
\end{Lemma}

\begin{proof}
For any $s\in \ttP_2\bigl(\big(\ttP^{\sankaku}_1\big)^{-1}(P) \bigr)\cap \gbigz_M^{\sharp}$,
we obtain
$(h_s)_{|P}=h_{\infty|P}$.
Then, by using the continuity of $h_s$ on $s$, we obtain
$(h_s)_{|P}=h_{\infty|P}$
for any $s\in \ttP_2\bigl(\big(\ttP^{\sankaku}_1\big)^{-1}(P)\bigr)$.
\end{proof}

\begin{Lemma}
Let $P\in X\setminus (W_M\cup H)$.
Then,
for any $s_1,s_2\in
 \ttP_2\bigl(\big(\ttP^{\sankaku}_1\big)^{-1}(P)\bigr)$,
we obtain
$h_{s_1|P}=h_{s_2|P}$.
\end{Lemma}
\begin{proof}
For any $P\in X^{\sharp}$
and for any
 $s_1,s_2\in
 \ttP_2\bigl(
 \big(\ttP^{\sankaku}_1\big)^{-1}(P)\bigr)$,
we obtain
$h_{s_1|P}=h_{\infty|P}
=h_{s_2|P}$. By~the continuity of $h_s$ on $s$,
we obtain the claim of the lemma.
\end{proof}

Then, $h_{\infty}$
extends to a Hermitian metric
of $E_{|X\setminus (H\cup W_M)}$
by setting
$h_{\infty|P}:=(h_s)_{|P}$
for $s\in \ttP_2\big(\big(\ttP^{\sankaku}_1\big)^{-1}(P)\big)$.

\begin{Lemma}\label{lem;19.2.1.10}
$h_{\infty}$ induces a Hermitian metric of
$E_{|X\setminus (H\cup W_M)}$
of $C^1$-class.
The $C^1$-Hermitian metric is also denoted by
$h_{\infty}$.
\end{Lemma}

\begin{proof}
Let $P$ be any point of $X\setminus W_M$.
Let $(U_P,x_1,x_2)$
be a holomorphic coordinate neighbourhood
as in Section~\ref{subsection;19.1.31.20}. By~using Proposition~\ref{prop;19.1.29.10},
we define the continuous Hermitian metric~$h^{(i)}_P$ of $E_{|U_P}$
by the condition that
$h_{P|\{x_i=a\}}$
is equal to the restriction of
$h_{s_1+as_{\infty}}$. By~the construction of $h_{\infty}$,
we obtain
$h_{\infty|U_P}
=h^{(i)}_{P}$.
Hence, we obtain that $h_{\infty|U_P\cap X^{\sharp}}$ induces
a conti\-nuous Hermitian metric $h_{P,\infty}$ of
$E_{|U_P}$,
and $h^{(1)}_P=h_{P,\infty}=h^{(2)}_P$ hold.
Moreover, by Proposition~\ref{prop;19.1.29.10},
any derivative of $h^{(j)}_P$
with respect to $\del_{z_i}$ and $\del_{\zbar_i}$ $(i\neq j)$
are continuous. We~obtain that $h_{P,\infty}$ is $C^1$.
Thus, we obtain the claim of the lemma.
\end{proof}

We obtain the operator $\DD^{\lambda\star}_{h_{\infty}}$
from $\DDlambda$ and $h_{\infty}$. We~define
$G(h_{\infty}):= \bigl[\DDlambda,\DD^{\lambda\star}_{h_{\infty}}\bigr]$
as a current.

\begin{Lemma}
\label{lem;19.1.31.30}
$G(h_{\infty})^{(1,1)}=0$
on $X\setminus (H\cup W_M)$.
\end{Lemma}

\begin{proof}
Let $P\in X\setminus (H\cup W_M)$.
Let $(U_P,x_1,x_2)$ be a holomorphic coordinate neighbourhood
as in Section~\ref{subsection;19.1.31.20}. We~have the expression
\begin{gather*}
 G(h_{\infty})^{(1,1)}
=G(h_{\infty})_{11}{\rm d}x_1\,{\rm d}\xbar_1
+G(h_{\infty})_{12}{\rm d}x_1\,{\rm d}\xbar_2
+G(h_{\infty})_{21}{\rm d}x_2\,{\rm d}\xbar_1
+G(h_{\infty})_{22}{\rm d}x_2\,{\rm d}\xbar_2.
\end{gather*}
Because $h_{\infty|\{x_i=a\}}$ is equal to
$h_{s_i+as_{\infty}}$,
we obtain $G(h_{\infty})_{ii}=0$ for $i=1,2$.

By considering
the holomorphic coordinate system
$(w_1,w_2)=(x_1+x_2,x_1-x_2)$
and the coefficient of
${\rm d}w_1\,{\rm d}\wbar_1$ in $G(h_{\infty})^{(1,1)}$,
we obtain $G(h_{\infty})_{12}+G(h_{\infty})_{21}=0$. By~considering the holomorphic coordinate system
$(z_1,z_2)=\big(x_1+\sqrt{-1}x_2,x_1-\sqrt{-1}x_2\big)$
and the coefficient of
${\rm d}z_1\,{\rm d}\zbar_1$ in $G(h_{\infty})^{(1,1)}$,
we obtain $G(h_{\infty})_{12}-G(h_{\infty})_{21}=0$.
Therefore, we obtain that
$G(h_{\infty})_{ij}=0$.
\end{proof}

\begin{Lemma}
We obtain $\Lambda G(h_{\infty})=0$ on $X\setminus (H\cup W_M)$.
As a result, $h_{\infty}$ is $C^{\infty}$
on $X\setminus (H\cup W_M)$.
If moreover $\lambda\neq 0$,
then $h_{\infty}$ is a pluri-harmonic metric of
$\big(E,\DDlambda\big)_{|X\setminus H\cup W_M}$.
\end{Lemma}

\begin{proof}
The first claim immediately follows from
Lemma~\ref{lem;19.1.31.30}. We~obtain the second claim by
the elliptic regularity
and a standard bootstrapping argument.
The last claim follows from Corollary~\ref{cor;20.2.3.1}.
\end{proof}

\begin{Lemma}
In the case $\lambda=0$,
we obtain $\del_{h_{\infty}}\theta=0$,
i.e.,
$h_{\infty}$ is a pluri-harmonic metric of
the Higgs bundle
$\big(E,\delbar_E,\theta\big)_{|X\setminus(H\cup W_M)}$.
\end{Lemma}

\begin{proof}
Let us observe that the sequence
$\del_{h^{(\epsilon_i)}_{\HE}}-\del_{h_{\infty}}$
is convergent to $0$
almost everywhere on~$X\setminus H$.
It~is enough to prove that
$\del_{h^{(\epsilon_i)}_{\HE|X_s}}-\del_{h_s}$
is convergent to $0$
for $s\in \gbigz_M^{\sharp}$.
Let $b^{(\epsilon_i)}_s$ be the automorphism of
$E_{|X_s\setminus H}$
which is self-adjoint with respect to
$h_{s}$ and $h_{\HE|X_s}^{(\epsilon_i)}$
determined by~$h_{\HE|X_s}^{(\epsilon_i)}
=h_{s}b^{(\epsilon_i)}_s$. By~Proposition~\ref{prop;13.1.7.20},
the sequence
$\big(b_s^{(\epsilon_i)}\big)^{-1}
 \del_{h_s}\big(b_s^{(\epsilon_i)}\big)$
is convergent to $0$ weakly in $L^2$
locally on $X_s\setminus H$. By~Proposition~\ref{prop;20.1.30.1},
the sequence
$\big(b_s^{(\epsilon_i)}\big)^{-1}
 \del_{h_s}\big(b_s^{(\epsilon_i)}\big)$
is bounded in $L_2^p$
locally on $X_s\setminus H$
for any $p\geq 1$.
\begin{Lemma}
\label{lem;21.5.4.3}
 $\big(b_s^{(\epsilon_i)}\big)^{-1}
 \del_{h_s}\big(b_s^{(\epsilon_i)}\big)$
is convergent to $0$ in $L_1^p$
locally on $X_s\setminus H$.
\end{Lemma}
\begin{proof}
Let
$\big(b_s^{\prime(\epsilon_i)}\big)^{-1}
 \del_{h_s}\big(b_s^{\prime(\epsilon_i)}\big)$
 be any subsequence
of $\big(b_s^{(\epsilon_i)}\big)^{-1}
 \del_{h_s}\big(b_s^{(\epsilon_i)}\big)$.
Because it is bounded in $L_2^p$ locally on $X_s\setminus H$,
it contains a subsequence
$\big(b_s^{\prime\prime(\epsilon_i)}\big)^{-1}
\del_{h_s}\big(b_s^{\prime\prime(\epsilon_i)}\big)$
which is weakly convergent
in $L_2^p$ locally on $X_s\setminus H$
for any $p\geq 2$. By~the Sobolev embedding theorem,
the sequence
$\big(b_s^{\prime\prime(\epsilon_i)}\big)^{-1}
\del_{h_s}\big(b_s^{\prime\prime(\epsilon_i)}\big)$
is convergent in $L_1^p$ locally on $X_s\setminus H$.
Because
$\big(b_s^{(\epsilon_i)}\big)^{-1}
\del_{h_s}\big(b_s^{(\epsilon_i)}\big)$
is convergent to $0$ weakly in $L^2$ locally on $X_s\setminus H$,
the limit should be $0$.
Therefore, we obtain the claim of~Lemma~\ref{lem;21.5.4.3}.
\end{proof}

As a result,
$\del_{h_{\HE}^{(\epsilon_i)}}\theta$
is convergent to
$\del_{h_{\infty}}\theta$ almost everywhere.
Note that
\begin{gather*}
 0\leq
 \int_{X\setminus H}
 \bigl|\del_{h_{\HE}^{(\epsilon_i)}} \theta
 \bigr|_{h_{\HE}^{(\epsilon_i)},\omega_{\epsilon_i}}^2
\leq \int_{X\setminus H}
 \bigl| G\big(h_{\HE}^{(\epsilon_i)}\big) \bigr|_{h_{\HE}^{(\epsilon_i)},\omega_{\epsilon_i}}^2
=-8\pi^2\int_X\ch_2\big(\nbigp^{(\epsilon_i)}_{\ast}\nbigv\big).
\end{gather*}
We also have
$\lim\limits_{i\to\infty}
 \int_X\ch_2\big(\nbigp^{(\epsilon_i)}_{\ast}\nbigv\big)
=0$. We~have the following convergence almost everywhere
on $X\setminus H$:
\begin{gather*}
 \lim_{i\to\infty}
\bigl| \del_{h_{\HE}^{(\epsilon_i)}}\theta
 \bigr|_{h_{\HE}^{(\epsilon_i)},\omega_{\epsilon_i}}^2
=
\bigl| \del_{h_{\infty}}\theta \bigr|_{h_{\infty},\omega_{X}}^2.
\end{gather*}
Therefore, we obtain
$\int\bigl|\del_{h_{\infty}}\theta\bigr|^2_{h_{\infty},\omega_X}=0$
by Fatou's lemma.
\end{proof}

\begin{Lemma}
$h_{\infty}$ induces a $C^{\infty}$-metric
of $E$ on $X\setminus H$,
and hence it is a pluri-harmonic metric of
$\big(E,\DDlambda\big)$.
\end{Lemma}
\begin{proof}
It~is enough to prove that
$h_{\infty}$ is a $C^{\infty}$-metric
around any point of $W_M\setminus H$. We~have only to apply the argument
in~\cite[Lemma~5.15]{mochi5}.
\end{proof}

If $\lambda=0$,
we obtain that
$\big(E,\delbar_E,\theta,h_{\infty}\big)$
is a good wild harmonic bundle on $(X,H)$,
because $(\nbigp_{\ast}\nbigv,\theta)$ is a good filtered Higgs bundle.
If $\lambda\neq 0$,
the associated Higgs bundle
$\big(E,\delbar_E,\theta\big)$
with the pluri-harmonic metric $h_{\infty}$
is a good wild harmonic bundle
by~\cite[Proposition~13.5.1]{Mochizuki-wild}. We~obtain a good filtered $\lambda$-flat bundle
$\big(\nbigp_{\ast}^{h_{\infty}}E,\DDlambda\big)$ on $(X,H)$. We~put $H^{[2]}=\bigcup_{i\neq j}(H_i\cap H_j)$.
For any $P\in H\setminus \big(W_M\cup H^{[2]}\big)$,
there exists $s\in\gbigz_M^{\sankaku}$
such that $P\in X_s$. By~the construction,
$h_{\infty|X_s\setminus H_s}=h_{s}$.
Hence, we obtain
$\nbigp_{\ast}^{h_{\infty}}(E)_{|X_s}
=\nbigp_{\ast}(\nbigv)_{|X_s}$.
Let $Y:=(H\cap W_M)\cup H^{[2]}$,
which is a finite subset of $H$. We~obtain that
$\nbigp_{\ast}^{h_{\infty}}(E)_{|X\setminus Y}
\simeq
 \nbigp_{\ast}\nbigv_{|X\setminus Y}$. By~Hartogs theorem,
we obtain that
$\nbigp_{\ast}^{h_{\infty}}(E)
\simeq
 \nbigp_{\ast}\nbigv$.
Thus, the proof of Proposition~\ref{prop;19.1.31.1}
is completed.
\hfill\qed

\subsection{Higher dimensional case}
\label{subsection;19.2.1.21}

Let us prove Theorem~\ref{thm;19.1.31.100}
in the case $\dim X\geq 3$
by an induction on $\dim X$.
Take a sufficiently large integer $M$. We~set
$\gbigz_M:=
 H^0\big(X,L^{\otimes\,M}\big)\setminus\{0\}$,
and let $\gbigx_M\subset X\times \gbigz_M$
be defined as $\gminis^{-1}(0)$
as in Section~\ref{subsection;19.2.1.1}.
For any $s\in \gbigz_M$,
set $X_s:=s^{-1}(0)$.
Let $\proj(T^{\ast}X)$ denote the projectivization of
the cotangent bundle of $X$.
If $M$ is sufficiently large,
there exists a Zariski dense open subset
$\gbigz_M^{\circ}\subset
 \gbigz_M$
such that the following holds:
\begin{itemize}\itemsep=0pt
\item
 $\ttP_2^{\circ}\colon
 \gbigx_M^{\circ}:=\gbigx_M\times_{\gbigz_M}\gbigz_M^{\circ}
\lrarr\gbigz_M^{\circ}$
is smooth.
\item
$\gbigx_M^{\circ}\cup
 \bigl(H\times\gbigz_M^{\circ}\bigr)$
is simply normal crossing.
Moreover
the intersections of any tuple of~irre\-du\-cible components
are smooth over $\gbigz_M^{\circ}$.
\item
The induced map
$\ttP_1\colon \gbigx_M^{\circ}\lrarr X$ is surjective.
Moreover,
 the induced morphism
 $\ttPtilde_1$: $\gbigx_M^{\circ}\lrarr\proj(T^{\ast}X)$
 is surjective.
\end{itemize}

Let $p_{i,j}$ denote the projection of
$X\times\gbigz_M^{\circ}\times\gbigz_M^{\circ}$
onto the product of the $i$-th component
and the $j$-th component.
For $j=2,3$,
let $\big(\gbigx_M^{\circ}\big)^{(j)}$ denote
the pull back of
$\gbigx_M^{\circ}$
by $p_{1,j}$.
There exists a Zariski dense open subset
$\gbigu_M\subset\gbigz_M^{\circ}\times\gbigz_M^{\circ}$
such that the following holds:
\begin{itemize}\itemsep=0pt
\item
Let $\big(\gbigx_M^{\circ}\big)^{(j)}_{\gbigu_M}$
denote the fiber product of
$\big(\gbigx_M^{\circ}\big)^{(j)}$
and $\gbigu_M$
over $\gbigz_M^{\circ}\times\gbigz_M^{\circ}$.
Then,
$\big(\gbigx_M^{\circ}\big)^{(2)}_{\gbigu_M}
 \cup
 \big(\gbigx_M^{\circ}\big)^{(3)}_{\gbigu_M}
 \cup
 (H\times \gbigu_M)$
is simply normal crossing.
Moreover,
the intersection of any tuple of~irre\-ducible components
are smooth over $\gbigu$.
\end{itemize}
By the Mehta--Ramanathan type theorem
(Proposition~\ref{prop;19.1.30.100}),
there exists a Zariski dense open subset
$\gbigu_M^{\sankaku}\subset\gbigu_M$
such that the following holds:
\begin{itemize}\itemsep=0pt
\item
 For $\vecs=(s_1,s_2)\in\gbigu_M^{\sankaku}$,
 we set
 $X_{\vecs}:=X_{s_1}\cap X_{s_2}$.
 Then, the restriction
 $\big(\nbigp_{\ast}\nbigv,\DDlambda\big)_{|X_{\vecs}}$
 is a~$\mu_L$-stable good filtered $\lambda$-flat bundle
 on $(X_{\vecs},H\cap X_{\vecs})$.
\end{itemize}
Hence, there exists a Zariski dense open subset
$\gbigz_M^{\sankaku}\subset\gbigz_M^{\circ}$
such that the following holds:
\begin{itemize}\itemsep=0pt
\item
For any $s\in\gbigz_M^{\sankaku}$,
$\big(\nbigp_{\ast}\nbigv,\DDlambda\big)_{|X_s}$
is a~$\mu_L$-stable good filtered $\lambda$-flat bundle on
$(X_s,H\cap X_s)$.
\item
For any $s_1,s_2\in\gbigz_M^{\sankaku}$,
there exists a Zariski open subset
$\gbigv(s_1,s_2)\subset\gbigz_M^{\sankaku}$
such that the restrictions
$\big(\nbigp_{\ast}\nbigv,\DDlambda\big)_{|X_{(s_i,s_3)}}$ $(i=1,2)$
are $\mu_L$-stable for any $s_3\in\gbigv(s_1,s_2)$.
\end{itemize}

We set $\gbigx_M^{\sankaku}:=
 \gbigx_M\times_{\gbigz_M}\gbigz_M^{\sankaku}$.
Let $\ttP_2^{\sankaku}\colon \gbigx_M^{\sankaku}\lrarr X$
denote the naturally induced morphism.
Then, $W_M:=X\setminus\ttP_2^{\sankaku}\big(\gbigx_M^{\sankaku}\big)$
is a finite subset.

For any $P\in X\setminus (H\cup W_M)$,
there exists $s\in\gbigz_M^{\sankaku}$
such that $P\in X_s$.
Then,
$\big(\nbigp_{\ast}\nbigv_s,\DDlambda_s\big):=
 \big(\nbigp_{\ast}\nbigv,\DDlambda\big)_{|X_s}$
is $\mu_L$-stable,
and the following holds:
\begin{gather*}
\int_{X_s}c_1(\nbigp_{\ast}\nbigv_s)c_1\big(L_{|X_s}\big)^{\dim X_s-1}=0,
\qquad
\int_{X_s}\ch_2(\nbigp_{\ast}\nbigv_s)c_1\big(L_{|X_s}\big)^{\dim X_s-2}=0.
\end{gather*}
There exists a pluri-harmonic metric $h_s$
of $\big(E_s,\delbar_{E_s},\DDlambda_s\big)
 :=\big(E,\delbar_E,\DDlambda\big)_{|X_s\setminus H}$
adapted to
$\nbigp_{\ast}\nbigv_s$
such that $\det(h_s)=h_{\det(E)|X_s\setminus H}$.
Take another $s'\in \gbigz_M^{\sankaku}$
such that $P\in X_{s'}$.
There exists a~pluri-harmonic metric $h_{s'}$
of $\big(E_{s'},\delbar_{E_{s'}},\DDlambda_{s'}\big)$
adapted to $\nbigp_{\ast}\nbigv_{s'}$
such that $\det(h_{s'})=h_{\det(E)|X_{s'}\setminus H}$.
\begin{Lemma}
$h_{s|P}=h_{s'|P}$.
\end{Lemma}

\begin{proof}
Suppose that $X_{s}\cup X_{s'}\cup H$ is
simply normal crossing. We~set $X_{s,s'}:=X_{s}\cap X_{s'}$.
It~is smooth and connected. We~obtain a good filtered $\lambda$-flat bundle
$\big(\nbigp_{\ast}\nbigv,\DDlambda\big)_{|X_{s,s'}}$,
and~$h_{s|X_{s,s'}}$ and~$h_{s'|X_{s,s'}}$
are adapted to
$\nbigp_{\ast}\nbigv_{|X_{s,s'}}$.
Let $b_{s,s'}$ be the automorphism of
$E_{|X_{s,s'}}$ which is self-adjoint
with respect to both $h_{s|X_{s,s'}}$ and $h_{s'|X_{s,s'}}$,
and determined by
$h_{s'|X_{s,s'}}=h_{s|X_{s,s'}}\cdot b_{s,s'}$.
There exists a decomposition
\begin{gather*}
 \big(\nbigp_{\ast}\nbigv,\DDlambda\big)_{|X_{s,s'}}
=\bigoplus
 \big(\nbigp_{\ast}\nbigv_i,\DDlambda_i\big),
\end{gather*}
which is orthogonal with respect to both
$h_{s|X_{s,s'}}$ and $h_{s'|X_{s,s'}}$,
and
$b_{s,s'}=\bigoplus a_i\id_{\nbigv_i}$
for some positive constants $a_i$.

There exists $s_1\in \gbigv(s,s')$.
Then,
$(\nbigp_{\ast}\nbigv,\theta)_{|X_{s_1,s}}$
and
$(\nbigp_{\ast}\nbigv,\theta)_{|X_{s_1,s'}}$
are $\mu_L$-stable.
Therefore,
we have
$h_{s|X_{s_1,s}}=h_{s_1|X_{s_1,s}}$
and
$h_{s'|X_{s_1,s'}}=h_{s_1|X_{s_1,s'}}$. We~obtain that
$h_{s|X_{s_1}\cap X_{s}\cap X_{s'}}
=h_{s'|X_{s_1}\cap X_{s}\cap X_{s'}}$.
It~implies that $a_i$ are $1$,
and hence $h_{s|P}=h_{s'|P}$.

In general,
there exists $s_2\in\gbigz_M^{\sankaku}$
such that
$(i)$ $P\in X_{s_2}$,
$(ii)$ $X_{s}\cup X_{s_2}\cup H$
 and $X_{s'}\cup X_{s_2}\cup H$
 are simply normal crossing. By~the above consideration,
we obtain
$h_{s|P}=h_{s_2|P}=h_{s'|P}$.
\end{proof}

Therefore,
we obtain Hermitian metrics
$h_P$ of $E_{|P}$ $(P\in X\setminus (H\cup W_M))$. By~using the argument in Lemma~\ref{lem;19.2.1.10},
we can prove that they induce a Hermitian metric $h$
of~$E_{|X\setminus(H\cup W_M)}$
of~$C^1$-class. We~obtain $G(h)$ from $\DDlambda$ and $h$
as a current.
Because $h_{|X_s}$ $\big(s\in \gbigu_M^{\sankaku}\big)$
are pluri-harmonic metrics of
$\big(E,\DDlambda\big)_{|X_s\setminus H}$,
we obtain that $G(h)=0$.
It~also implies that $h$ is $C^{\infty}$ on
$X\setminus (H\cup W_M)$. By~using the argument in~\cite[Lemma 5.15]{mochi5},
we obtain that
$h$ induces a pluri-harmonic metric
of $\big(E,\DDlambda\big)$ on $X\setminus H$.
Then, as in the proof of Proposition~\ref{prop;19.1.31.1},
we can conclude that
$\big(E,\DDlambda,h\big)$ is a good wild harmonic bundle,
and that $\nbigp^h_{\ast}(E)=\nbigp_{\ast}\nbigv$.
Thus, we obtain Theorem~\ref{thm;19.1.31.100}.
\hfill\qed

\section{Homogeneity with respect to group actions}
\label{section;19.2.16.30}

\subsection{Preliminary}

\subsubsection{Homogeneous harmonic bundles}

Let $Y$ be a complex manifold.
Let $K$ be a compact Lie group.
Let $\rho\colon K\times Y\lrarr Y$
be a~$K$-action on $Y$
such that $\rho_k\colon Y\lrarr Y$ is holomorphic
for any $k\in K$.
Let $\kappa\colon K\lrarr S^1$ be
a~homomorphism of Lie groups.

Let $\big(E,\delbar_E,\theta,h\big)$ be a harmonic bundle
on $Y$.
It~is called $(K,\rho,\kappa)$-homogeneous
if $\big(E,\delbar_E,h\big)$ is $K$-equivariant
and $k^{\ast}\theta=\kappa(k)\theta$.

\begin{Remark}
According to {s3},
harmonic bundles are equivalent to
polarized variation of pure twistor structure of weight $w$,
 for any given integer $w$.
If $\kappa$ is non-trivial,
as studied in \cite[Section~3]{Mochizuki-TodaII},
by choosing a vector $\gminiv$ in the Lie algebra of $K$
such that $d\kappa(\gminiv)\neq 0$,
we obtain the integrability of the variation
of pure twistor structure
from the homogeneity of harmonic bundles.
\end{Remark}

\subsubsection[Homogeneous filtered Higgs sheaves
and the stability condition with respect to the action]
{Homogeneous filtered Higgs sheaves
and the stability condition \\with respect to the action}\label{subsection;19.2.18.21}

Let $X$ be a connected complex projective manifold
with a simple normal crossing hypersur\-face~$H$.
Let $G$ be a complex reductive algebraic group.
Let $\rho\colon G\times Y\lrarr Y$ be an algebraic $G$-action
on $Y$ which preserves $H$.
Let $\kappa\colon G\lrarr \cnum^{\ast}$
be a homomorphism of complex algebraic groups.

Let $(\nbigp_{\ast}\nbigv,\theta)$
be a filtered Higgs sheaf on $(Y,H)$.
It~is called $(G,\rho,\kappa)$-homogeneous
if $\nbigp_{\ast}\nbigv$ is~$G$-equivariant
and $g^{\ast}\theta=\kappa(g)\theta$
for any $g\in G$.

Let $L$ be a $G$-equivariant ample line bundle on $X$.
A $(G,\rho,\kappa)$-homogeneous
filtered Higgs sheaf $(\nbigp_{\ast}\nbigv,\theta)$ on $(X,H)$
is called $\mu_L$-stable (resp.~$\mu_L$-semistable)
with respect to the $G$-action
if the following holds:
\begin{itemize}\itemsep=0pt
\item
Let $\nbigv'$ be a $G$-invariant
saturated Higgs subsheaf of $\nbigv$
such that $0<\rank\nbigv'<\rank\nbigv$.
Then,
$\mu_L(\nbigp_{\ast}\nbigv')
<\mu_L(\nbigp_{\ast}\nbigv)$
(resp.~$\mu_L(\nbigp_{\ast}\nbigv')
\leq\mu_L(\nbigp_{\ast}\nbigv)$)
holds.
\end{itemize}
A $(G,\rho,\kappa)$-homogeneous
filtered Higgs sheaf $(\nbigp_{\ast}\nbigv,\theta)$ on $(X,H)$
is called $\mu_L$-polystable
with respect to the $G$-action
if it is $\mu_L$-semistable with respect to the $G$-action
and isomorphic to a~direct sum
of $(G,\rho,\kappa)$-homogeneous filtered sheaves
$\bigoplus(\nbigp_{\ast}\nbigv_i,\theta_i)$,
where each $(\nbigp_{\ast}\nbigv_i,\theta_i)$
is $\mu_L$-stable with respect to the $G$-action.

\begin{Lemma}
\label{lem;19.2.1.100}
$(\nbigp_{\ast}\nbigv,\theta)$ is
$\mu_L$-semistable
if and only if
$(\nbigp_{\ast}\nbigv,\theta)$
is $\mu_L$-semistable with respect to
the $G$-action.
\end{Lemma}

\begin{proof}
The ``only if'' part is clear.
Let us prove that the ``if'' part.
Let $\nbigv_0\subset\nbigv$
be the $\beta$-subobject as in Proposition~\ref{prop;19.2.18.110}.
Because $g^{\ast}\nbigv_0$ also has the same property,
we obtain that $\nbigv_0$ is~$G$-invariant.
Then, the claim of the proposition is clear.
\end{proof}

The following lemma is clear.

\begin{Lemma}
If $(\nbigp_{\ast}\nbigv,\theta)$ is
$\mu_L$-stable,
then
$(\nbigp_{\ast}\nbigv,\theta)$ is
$\mu_L$-stable
with respect to the $G$-action.
\end{Lemma}

\begin{Lemma}
\label{lem;19.2.2.1}
If $(\nbigp_{\ast}\nbigv,\theta)$ is
$\mu_L$-stable
with respect to the $G$-action,
then
$(\nbigp_{\ast}\nbigv,\theta)$ is
$\mu_L$-polystable.
\end{Lemma}

\begin{proof}
According to Lemma~\ref{lem;19.2.1.100},
$(\nbigp_{\ast}\nbigv,\theta)$ is $\mu_L$-semistable.
Let $\nbigv_1$ be the socle of
$(\nbigp_{\ast}\nbigv,\theta)$
as in Proposition~\ref{prop;19.2.18.140}.
Because $g^{\ast}\nbigv_1$
also has the same property,
$\nbigv_1$ is $G$-invariant.
Moreover,
$\mu_L(\nbigp_{\ast}\nbigv_1)=\mu_L(\nbigp_{\ast}\nbigv)$
holds.
Hence, we obtain
$\nbigv_1=\nbigv$.
According to Proposition~\ref{prop;19.2.18.140},
$(\nbigp_{\ast}\nbigv,\theta)$
is $\mu_L$-polystable.
\end{proof}

\begin{Remark}
In general,
even if $(\nbigp_{\ast}\nbigv,\theta)$ is $\mu_L$-stable
with respect to the $G$-action,
$(\nbigp_{\ast}\nbigv,\theta)$ is not necessarily $\mu_L$-stable.
\end{Remark}

\subsubsection{Actions of a complex reductive group
and its compact real form}

Let $X$ be a complex projective manifold
equipped with an algebraic action of
a complex reductive group $G$.
Let $L$ be a $G$-equivariant ample line bundle on $X$.
Let $K$ be a compact real form of~$G$.

Let $\big(E,\delbar_E\big)$ be a $G$-equivariant holomorphic vector bundle on $X$.
Then, as the restriction,
we may naturally regard $\big(E,\delbar_E\big)$
as a $K$-equivariant holomorphic vector bundle on $X$.
\begin{Lemma}
\label{lem;19.2.12.120}
The above procedure induces an equivalence
between $G$-equivariant holomorphic vector bundles
and $K$-equivariant holomorphic vector bundles on $X$.
\end{Lemma}

\begin{proof}
Let $\big(E,\delbar_E\big)$ be a $K$-equivariant holomorphic vector bundle
on $X$.
There exists $m_0>0$ such that
$E\otimes L^{\otimes\,m_0}$ is globally generating. We~set
$\nbigg_0:=H^0\big(X,E\otimes L^{\otimes\,m_0}\big)\otimes \big(L^{\otimes\,m_0}\big)^{-1}$.
There exists a naturally induced epimorphism
of $\nbigo_X$-modules $\nbigg_0\lrarr E$.
Let $\nbigk$ denote the kernel.
There exists $m_1>0$ such that
$\nbigk\otimes L^{\otimes\,m_1}$ is globally generating. We~set
$\nbigg_1:=
 H^0\big(X,\nbigk\otimes L^{\otimes\,m_1}\big)
 \otimes \big(L^{\otimes\,m_1}\big)^{-1}$.
There exists a naturally induced epimorphism
$\nbigg_1\lrarr\nbigk$.
Thus, we obtain a resolution
$\nbigg_1\lrarr\nbigg_0$ of $E$.
Because $E$ is $K$-equivariant,
$H^0\big(X,E\otimes L^{\otimes m_0}\big)$
is naturally a~$K$-representation,
$\nbigg_0$ is a $K$-equivariant holomorphic vector bundle
on $X$,
and $\nbigg_0\lrarr E$ is $K$-equivariant.
Hence, $\nbigk$ is a $K$-equivariant holomorphic vector bundle.
Similarly
$H^0\big(X,\nbigk\otimes L^{\otimes m_2}\big)$
is a $K$-representation,
and $\nbigg_1$ is $K$-equivariant holomorphic vector bundle,
and $\nbigg_1\lrarr\nbigk_2$ is $K$-equivariant.

The $K$-representations
on $H^0\big(X,E\otimes L^{\otimes\,m_1}\big)$
and $H^0\big(X,\nbigk\otimes L^{\otimes m_2}\big)$
naturally induce
$G$-representations
on $H^0\big(X,E\otimes L^{\otimes\,m_1}\big)$
and $H^0\big(X,\nbigk\otimes L^{\otimes m_2}\big)$.
Hence,
$\nbigg_i$ are naturally
algebraic $G$-equivariant vector bundles on $X$.
Moreover,
the morphism $\nbigg_1\lrarr\nbigg_0$
is $G$-equivariant and algebraic.
Hence, $E$ is a $G$-equivariant algebraic vector bundle
on $X$.
\end{proof}

\subsection{An equivalence}

\subsubsection{Good filtered Higgs bundles associated with
homogeneous good wild Higgs bundles}

Let $X$ be a connected complex projective manifold
with a simple normal crossing hypersurface~$H$.
Let $G$ be a complex reductive group acting on $(X,H)$.
Let $K$ be a compact real form of $G$.
The actions of $G$ and $K$ on $X$ are denoted by $\rho$.
Let $\kappa\colon G\lrarr \cnum^{\ast}$
be a character.
The induced homomorphism $K\lrarr S^1$ is also denoted by $\kappa$.

Let $\big(E,\delbar_E,\theta,h\big)$ be a
$(K,\rho,\kappa)$-homogeneous harmonic bundle on $X\setminus H$
which is good wild on $(X,H)$. We~obtain a good filtered Higgs bundle
$\big(\nbigp^h_{\ast}E,\theta\big)$ on $(X,H)$.
Because each $\nbigp^h_{\veca}E$ is naturally
a $K$-equivariant holomorphic vector bundle on $X$,
$\nbigp_{\ast}^hE$ is naturally $G$-equivariant
by Lemma~\ref{lem;19.2.12.120}.
Because $k^{\ast}\theta=\kappa(k)\theta$ for any $k\in K$,
we obtain
$g^{\ast}\theta=\kappa(g)\theta$ for any $g\in G$.
Therefore,
$\big(\nbigp_{\ast}^hE,\theta\big)$
is a $(G,\rho,\kappa)$-homogeneous
good filtered Higgs bundle on $(X,H)$.

Let $L$ be a $G$-equivariant ample line bundle on $X$.
\begin{Proposition}
$\big(\nbigp_{\ast}^hE,\theta\big)$
is $\mu_L$-polystable
with respect to the $G$-action,
i.e.,
there exists a~decomposition
$\big(E,\delbar_E,\theta,h\big)
=\bigoplus(E_i,\delbar_{E_i},\theta_i,h_i)$
of $(G,\rho,\kappa)$-homogeneous harmonic bundles
such that each
$\big(\nbigp_{\ast}^{h_i}E_i,\theta_i\big)$ is $\mu_L$-stable
with respect to the $G$-action.
\end{Proposition}

\begin{proof}
Because $\big(\nbigp^h_{\ast}E,\theta\big)$
is $\mu_L$-polystable,
we obtain that
$\big(\nbigp^h_{\ast}E,\theta\big)$
is $\mu_L$-semistable
with res\-pect to the $G$-action.
Let $\nbigv_1\subset\nbigp^hE$
be a $G$-invariant saturated Higgs $\nbigo_X(\ast H)$-submodule
such that
$\mu_L(\nbigp_{\ast}\nbigv_1)=\mu_L\big(\nbigp_{\ast}^hE\big)=0$.
Let $E_1$ be the Higgs subsheaf of $E$
obtained as the restriction of~$\nbigv_1$ to $X\setminus H$.
Then, by the argument in the proof of~\cite[Proposition~13.6.1]{Mochizuki-wild},
we obtain that~$E_1$ is a subbundle,
and the orthogonal complement
$E_2:=E_1^{\bot}$ is also a holomorphic subbundle.
Moreover,
$\theta(E_2)\subset E_2\otimes\Omega^1_{X\setminus H}$,
and $E_2$ is $K$-equivariant.
Hence, we obtain a decomposition
$\big(E,\delbar_E,\theta,h\big)
=\big(E_1,\delbar_{E_1},\theta_1,h_1\big)
\oplus
\big(E_2\delbar_{E_2},\theta_2,h_2\big)$
of $(K,\rho,\kappa)$-homogeneous harmonic bundles.
Then, the claim of the proposition is clear.
\end{proof}

\subsubsection{Uniqueness}

Let $\big(E,\delbar_E,\theta,h\big)$ be a
$(K,\rho,\kappa)$-homogeneous harmonic bundle on $X\setminus H$
which is good wild on $(X,H)$.
Let $h'$ be another pluri-harmonic metric of
$\big(E,\delbar_E,\theta\big)$
such that
$(i)$ $h'$ is $K$-invariant,
$(ii)$ $\nbigp_{\ast}^{h'}E=\nbigp_{\ast}^hE$.
The following is clear from
Proposition~\ref{prop;19.2.12.200}.
\begin{Proposition}
There exists a decomposition
$\big(E,\delbar_E,\theta\big)=\bigoplus_{i=1}^m
 \big(E_{i},\delbar_{E_i},\theta_i\big)$
such that
$(i)$ the decomposition is orthogonal with respect to
both $h$ and $h'$,
$(ii)$ there exist $a_i>0$ $(i=1,\ldots,m)$
 such that
 $h'_{|E_i}=a_ih_{E_i}$,
$(iii)$ the decomposition $E=\bigoplus E_i$
is preserved by the $K$-action.
\end{Proposition}

\subsubsection{Existence theorem}

Let $(\nbigp_{\ast}\nbigv,\theta)$
be a $(G,\rho,\kappa)$-homogeneous good filtered
Higgs bundle on $(X,H)$
such that
\begin{gather*}
 \int_Xc_1(\nbigp_{\ast}\nbigv)c_1(L)^{\dim X-1}
=0,
 \quad
 \int_X\ch_2(\nbigp_{\ast}\nbigv)c_1(L)^{\dim X-2}=0.
\end{gather*}
Let $\big(E,\delbar_E,\theta\big)$
be the Higgs bundle on $X\setminus H$
obtained as the restriction of $(\nbigp_{\ast}\nbigv,\theta)$.

\begin{Theorem}
Suppose that
$(\nbigp_{\ast}\nbigv,\theta)$ is $\mu_L$-stable
with respect to the $G$-action.
Then, there exists a $K$-invariant pluri-harmonic metric $h$
of $\big(E,\delbar_E,\theta\big)$
such that $\nbigp^h_{\ast}E=\nbigp_{\ast}\nbigv$.
If $h'$ is another $K$-invariant pluri-harmonic metric
of $\big(E,\delbar_E,\theta\big)$,
there exists a positive constant $a$
such that $h'=ah$.
\end{Theorem}
\begin{proof}
By Lemma~\ref{lem;19.2.2.1},
$(\nbigp_{\ast}\nbigv,\theta)$ is $\mu_L$-polystable.
There exists the canonical decomposition
\begin{gather*}
 (\nbigp_{\ast}\nbigv,\theta)=
 \bigoplus_{i=1}^m
 (\nbigp_{\ast}\nbigv_i,\theta_i)
 \otimes
 U_i,
\end{gather*}
where $(\nbigp_{\ast}\nbigv_i,\theta_i)$
are $\mu_L$-stable good filtered Higgs bundles
such that
$(\nbigp_{\ast}\nbigv_i,\theta_i)
\not\simeq
(\nbigp_{\ast}\nbigv_j,\theta_j)$ $(i\neq j)$,
and $U_i$ are finite dimensional complex vector spaces.
Let
$\big(E_i,\delbar_{E_i},\theta_i\big)$
denote the Higgs bundle obtained as the restriction of
$(\nbigv_i,\theta_i)$ to $X\setminus H$.
There exist pluri-harmonic metrics~$h_i$
of~$\big(E_i,\delbar_{E_i},\theta_i\big)$
adapted to the filtered bundles $\nbigp_{\ast}\nbigv_i$.
Let $h^{(0)}_{U_i}$ be Hermitian metrics of~$U_i$. We~obtain a pluri-harmonic metric
$h^{(0)}=\bigoplus \big(h_i\otimes h^{(0)}_{U_i}\big)$
of $\big(E,\delbar_E,\theta\big)$
adapted to $\nbigp_{\ast}\nbigv$. By~Proposition~\ref{prop;19.2.12.200}
and the uniqueness of the canonical decomposition,
we obtain the following lemma.
\begin{Lemma}
For any pluri-harmonic metric $h^{(1)}$ of
$\big(E,\delbar_E,\theta\big)$
adapted to $\nbigp_{\ast}\nbigv$,
 there uniquely exist Hermitian metrics $h^{(1)}_{U_i}$ of $U_i$
 such that
$h^{(1)}=\bigoplus\big(h_{i}\otimes h^{(1)}_{U_i}\big)$.
\end{Lemma}

For any $k\in K$,
we obtain a pluri-harmonic metric $k^{\ast}h^{(0)}$
of $\big(E,\delbar_E,\kappa(k)\theta\big)$
adapted to $\nbigp_{\ast}\nbigv$.
Because $|\kappa(k)|=1$,
$k^{\ast}\big(h^{(0)}\big)$ is also a pluri-harmonic metric
of $\big(E,\delbar_E,\theta\big)$
adapted to $\nbigp_{\ast}\nbigv$.
Hence, there uniquely exist
Hermitian metrics $h_{U_i}(k)$
of $U_i$
such that
$k^{\ast}(h)=\bigoplus_{i=1}^m(h_i\otimes h_{U_i}(k))$. By~using the Haar measure ${\rm d}k$ on $K$
with $\int_K{\rm d}k=1$,
we define the Hermitian metric $h$ of $E$ as follows:
\begin{gather*}
 h:=\int_Kk^{\ast}\big(h^{(0)}\big){\rm d}k
 =\bigoplus_{i=1}^m
 \Bigl(h_i\otimes\int_Kh_{U_i}(k){\rm d}k\Bigr).
\end{gather*}
Then, $h$ is also a pluri-harmonic metric. By~the construction, $h$ is $K$-invariant.

Let $h'$ be another $K$-invariant pluri-harmonic metric
of $\big(E,\delbar_E,\theta\big)$ adapted to $\nbigp_{\ast}\nbigv$. We~obtain the decomposition
$\big(E,\delbar_E,\theta\big)=\bigoplus \big(E_i,\delbar_{E_i},\theta_i\big)$
as in Proposition~\ref{prop;19.2.12.200},
which induces a~decom\-position of
good filtered Higgs bundles
$(\nbigp_{\ast}\nbigv)
=\bigoplus(\nbigp_{\ast}\nbigv_i,\theta_i)$.
Because both $h$ and $h'$ are $K$-invariant,
the decompositions are also $K$-invariant.
Hence, the decomposition
$(\nbigp_{\ast}\nbigv)
=\bigoplus(\nbigp_{\ast}\nbigv_i,\theta_i)$
is $G$-invariant. By~the $\mu_L$-stability of
$(\nbigp_{\ast}\nbigv)$,
we obtain $m=1$,
i.e., $h'=ah$ for~$a>0$.
\end{proof}

\begin{Corollary}
\label{cor;19.2.18.20}
We obtain the equivalence between
the isomorphism classes of the following objects:
\begin{itemize}\itemsep=0pt
\item
 $(K,\rho,\kappa)$-homogeneous good wild harmonic bundles
on $(X,H)$.
\item
 $(G,\rho,\kappa)$-homogeneous good filtered Higgs bundles
 $(\nbigp_{\ast}\nbigv,\theta)$
 such that
 $(i)$ it is $\mu_L$-polystable with respect to the $G$-action,
 $(ii)$ $\mu_L(\nbigp_{\ast}\nbigv)=0$,
 $\int_X\ch_2(\nbigp_{\ast}\nbigv)c_1(L)^{\dim X-2}=0$.
\end{itemize}
\end{Corollary}

\subsection*{Acknowledgement}

It~is my great pleasure to dedicate this paper to
Professor Kyoji Saito on the occasion of his 77th birthday.
Not to mention his pioneering and profound achievements,
I have always been impressed with his kindness and generosity
towards the younger generation of mathematicians,
as well as his insatiable quest for mathematics.
I appreciate the editors of this special issue for the invitation.

I thank Carlos Simpson
for his fundamental works on harmonic bundles
which are most important in this study.
I am grateful to the referees
for their careful readings
and helpful comments to improve this manuscript.
I thank Philip Boalch and Andy Neitzke
for their kind comments to a preliminary note
on the proof in the one dimensional case.
I thank Claude Sabbah for discussions
on many occasions
and for his kindness.
I thank Fran\c{c}ois Labourie
for his comment on the definition of wild harmonic bundles.
I thank Pengfei Huang for discussions.
I am grateful to Ya Deng
for asking questions
related to Mehta--Ramanathan type theorems
in~Section~\ref{subsection;20.7.7.21}
and tensor products in Section~\ref{subsection;20.7.10.1}.
I thank Akira Ishii and Yoshifumi Tsuchimoto
for their constant encouragement.
A part of this manuscript was prepared
for lectures in the Oka symposium
and
ICTS program ``Quantum Fields, Geometry and Representation Theory''.
I thank the organizers for the opportunities.

I am partially supported by
the Grant-in-Aid for Scientific Research (S) (No.~17H06127),
the Grant-in-Aid for Scientific Research (S) (No.~16H06335),
the Grant-in-Aid for Scientific Research (A) (No.~21H04429),
the Grant-in-Aid for Scientific Research (C) (No.~15K04843),
and
the Grant-in-Aid for Scientific Research (C) (No.~20K03609),
Japan Society for the Promotion of Science.

\pdfbookmark[1]{References}{ref}
\LastPageEnding

\end{document}